\documentclass[12pt,twoside,leqno]{amsart}
\usepackage{amssymb,amsbsy,amsmath,amsfonts,amssymb,amscd,epsfig,times,graphics,color,xypic}

\usepackage[latin1]{inputenc}
\newcommand{\C}{\mathbb{C}}

\newcommand{\K}{\mathbb{K}}
\newcommand{\R}{\mathbb{R}}
\newcommand{\N}{\mathbb{N}}

\def\dl{{[\![}}
\def\dr{{]\!]}}

\newtheorem{problem}{Problem}
\newtheorem{theorem}{Theorem}
\newtheorem{proposition}{Proposition}
\newtheorem{lemma}{Lemma}
\newtheorem{corollary}{Corollary}
\newtheorem{definition}{Definition}
\newtheorem{example}{Example}
\newtheorem{convention}{Convention}

\newtheorem{openquestion}{Open question}
\newtheorem{openproblem}{Open problem}
\newtheorem{notationalconvention}{Notational Convention}

\begin{document}


\title[Lie symmetries and CR geometry]{
{\Large\bf I:}~Lie symmetries of partial differential equations
\\
and CR geometry}

\author{Jo\"el Merker}

\address{
CNRS, Universit\'e de Provence, LATP, UMR 6632, CMI, 
39 rue Joliot-Curie, F-13453 Marseille Cedex 13, France. 
{\it Internet}:
{\tt http://www.cmi.univ-mrs.fr/$\sim$merker/index.html}}

\email{merker@cmi.univ-mrs.fr} 

\date{\number\year-\number\month-\number\day}

\maketitle

\begin{center}
\begin{minipage}[t]{12cm}
\baselineskip =0.35cm
{\scriptsize

\vskip -0.5cm

\centerline{\bf Table of contents}

\smallskip

{\bf 1.~Completely integrable systems of partial differential 
equations \dotfill 1.}

{\bf 2.~Submanifold of solutions \dotfill 6.}

{\bf 3.~Classification problems \dotfill 10.}

{\bf 4.~Punctual and infinitesimal Lie symmetries \dotfill 12.}

{\bf 5.~Examples \dotfill 15.}

{\bf 6.~Transfer of Lie symmetries to the parameter space \dotfill 18.}

{\bf 7.~Equivalence problems and normal forms \dotfill 25.}

{\bf 8.~Study of two specific examples \dotfill 27.}

{\bf 9.~Dual system of partial differential equations \dotfill 45.}

{\bf 10.~Fundamental pair of foliations and covering property 
\dotfill 46.}

{\bf 11.~Formal and smooth equivalences \dotfill 50.}


}\end{minipage}
\end{center}


{\small

\bigskip

This memoir is divided in three parts\footnote{ Part~II of~\cite{
me2005a} already appeared as \cite{ me2005b}.}. Part~I endeavours a
general, new theory (inspired by modern
CR geometry) of Lie symmetries of completely integrable {\sc
pde} systems, viewed from their associated submanifold of solutions.
Part~II builds general combinatorial formulas for the
prolongations of vector fields to jet spaces.
Part~III characterizes explicitly flatness of some systems of second
order. The results presented here are original and did not
appear in print elsewhere; most formulas of 
Parts~II and~III were checked by means of Maple Release 7.

}

\section*{\S1.~Completely integrable systems of 
partial differential equations}

\subsection*{1.1.~General systems} 
Let $\K= \R$ or $\C$. Let $n \in \N$ with $n \geqslant 1$ and let $x =
(x^1, \dots, x^n) \in \K^n$. Also, let $m \in \N$ with $m \geqslant 1$
and let $y= (y^1, \dots, y^m) \in \K^m$. For $\alpha \in \N^n$, we
denote by a subscript $y_{ x^\alpha}$ the partial derivative
$\partial^{ \vert \alpha \vert} y / \partial x^\alpha$ of a local map
$\K^n \ni x \mapsto y (x) \in \K^m$.

Let $\kappa \in \N$ with $\kappa \geqslant 1$, let $p \in \N$ with $p
\geqslant 1$, choose a collection of $p$ multiindices $\beta (1),
\dots, \beta(p) \in \N^n$ with $\vert \beta(q) \vert \geqslant 1$ for
$q=1,\dots,p$ and $\max_{1\leqslant q \leqslant p}\, \vert \beta(q)
\vert =\kappa$, and choose integers $j(1), \dots, j(p)$ with $1
\leqslant j(q) \leqslant m$ for $q = 1, \dots, p$. In the present
Part~I, we study the Lie symmetries of a general system of
analytic partial differential equations of the form:
\def\theequation{$\mathcal{E}$}\begin{equation}
y_{x^\alpha}^j(x)
=
F_\alpha^j
\Big(
x,y(x),
\big(y_{x^{\beta(q)}}^{j(q)}(x)
\big)_{
1\leqslant q\leqslant p}
\Big),
\end{equation}
where $j$ with $1 \leqslant j \leqslant m$ and $\alpha \in \N^n$
satisfy
\def\theequation{1.2}\begin{equation}
\big(
j,\alpha
\big) 
\neq 
(
j,0
)
\ \ \ \ \
\text{\rm and}
\ \ \ \ \
\big(
j,\alpha
\big) 
\neq
\big(
j(q), \beta(q)
\big).
\end{equation}
In particular,
all $(\kappa +1)$-th 
partial derivatives of the unknown $y = y(x)$ depend
on a certain precise set of derivatives of order $\leqslant \kappa$: 
the system is {\it complete}. In addition, all the other partial
derivatives of order $\leqslant \kappa$ do also depend on the same
precise set of derivatives.

Here, we assume that $u = 0$ is a local solution of the
system~\thetag{ $\mathcal{ E }$} and that the functions $F_\alpha^j$
are $\K$-algebraic (in the sense of Nash) or $\K$-analytic, in a
neighborhood of the origin in $\K^{ n+ m+ p }$. Even if our
concern will be local throughout, we will not introduce any
special notation to speak of open subsets and simply refer to various
$\K^\mu$. We will study five concrete instances, the first three ones
being classical.

\def\theexample{1.3}\begin{example}{\rm
With $n = m = \kappa = 1$, 
a second order ordinary differential equation 
\def\theequation{$\mathcal{E}_1$}\begin{equation}
y_{xx} 
= 
F(x,y,y_x),
\end{equation}
and more generally $y_{ x^{ \kappa +1}} = F ( x, y, y_x, \dots, y_{
x^\kappa } \big)$, where $x, y \in \K$, see~\cite{lie1883, el1890,
tr1896, ca1924, se1931, ca1932a, ol1986, ar1988, bk1989, gtw1989,
hk1989, ib1992, ol1995, n2003}.
}\end{example}

\def\theexample{1.4}\begin{example}{\rm
With $n \geqslant 2$, $m=1$ and $\kappa =1$, 
a complete system of second order
equations
\def\theequation{$\mathcal{E}_2$}\begin{equation}
y_{x^{i_1}x^{i_2}} 
= 
F_{i_1,i_2} 
\big(
x^i,y,y_{x^k}
\big),
\ \ \ \ \ \
1\leqslant i_1,i_2\leqslant n, 
\end{equation}
see~\cite{ha1937, ch1975, su2001} and Part~III below.
}\end{example}

\def\theexample{1.5}\begin{example}{\rm
Dually, with $n = 1$, $m \geqslant 2$ and
$\kappa = 1$, an ordinary system of second
order
\def\theequation{$\mathcal{E}_3$}\begin{equation}
y_{xx}^j 
=
F^j
\big( 
x,y^{j_1},y_x^{j_1}
\big), 
\ \ \ \ \ \ \
j=1,\dots,m, 
\end{equation}
see~\cite{ fe1995, me2004} and the references therein.
}\end{example}

\def\theexample{1.6}\begin{example}{\rm
With $n=1$, $m=2$ and $\kappa = 1$, a system of the form
\def\theequation{$\mathcal{E}_4$}\begin{equation}
\left\{
\aligned
y_x^2
&
=
F
\big(
x,y^1,y^2,y_x^1
\big)
\\
y_{xx}^1
&
=
G
\big(
x,y^1,y^2,y_x^1
\big).
\endaligned\right.
\end{equation}
Differentiating the first equation with respect to $x$ and
substituting, we get the missing equation:
\def\theequation{1.7}\begin{equation}
\aligned
y_{xx}^2
&
=
F_x
+
y_x^1\,F_{y^1}
+
y_x^2\,F_{y^2}
+
y_{xx}^1\,F_{y_x^1}
\\
&
=
F_x
+
y_x^1\,F_{y^1}
+
y_x^2\,F_{y^2}
+
G\,F_{y_x^1}
\\
&
=:
H
\big(
x,y^1,y^2,y_x^1
\big).
\endaligned
\end{equation}
}\end{example}

\def\theexample{1.8}\begin{example}{\rm
With $n=2$, $m=1$ and $\kappa = 2$, a system of the form
\def\theequation{$\mathcal{E}_5$}\begin{equation}
\left\{
\aligned
y_{x^2}
&
=
F
\big(
x^1,x^2,y,y_{x^1},y_{x^1x^1}
\big)
\\
y_{x^1x^1x^1}
&
=
G
\big(
x^1,x^2,y,y_{x^1},y_{x^1x^1}
\big).
\endaligned\right.
\end{equation}
Here, five equations are missing.
Differentiating the first equation with respect to $x^1$ and
substituting:
\def\theequation{1.9}\begin{equation}
\aligned
y_{x^1x^2}
&
=
F_{x^1}
+
y_{x^1}\,F_y
+
y_{x^1x^1}\,F_{y_{x^1}}
+
y_{x^1x^1x^1}\,F_{y_{x^1x^1}}
\\
&
=
F_{x^1}
+
y_{x^1}\,F_y
+
y_{x^1x^1}\,F_{y_{x^1}}
+
G\,F_{y_{x^1x^1}}
\\
&
=:
H
\big(
x^1,x^2,y,y_{x^1},y_{x^1x^1}
\big),
\endaligned
\end{equation}
and then similarly for $y_{ x^2 x^2}$, $y_{ x^1 x^1 x^2}$, $y_{ x^1
x^2 x^2}$, $y_{ x^2 x^2 x^2 }$.
}\end{example}

\subsection*{1.10.~Finitely nondegenerate generic submanifolds 
of $\C^{ n + m}$} Examples~1.3, 1.4, 1.6 and~1.8 (but {\it not}\,
1.5) are intrinsically linked to real submanifolds of complex
submanifolds.

Let $M$ be a real algebraic or analytic local generic CR\footnote{
Fundamentals about Cauchy-Riemann geometry may be found in~\cite{
bo1991, ber1999, me2005a, me2005b, mp2005}.} submanifold of $\C^{ n +
m}$ of codimension $m \geqslant 1$ and of CR dimension $n \geqslant
1$, and let $p \in M$. Classically, there exists local holomorphic
coordinates $t = (z, w) \in \C^n \times \C^m$ centered at $p$ in which
$M$ is represented by
\def\theequation{1.11}\begin{equation}
w^j
=
\overline{\Theta}^j
(z,\bar z,\bar w),
\ \ \ \ \ 
j=1,\dots,m,
\end{equation}
for some local $\C$-analytic map $\Theta = ( \Theta^1, \dots,
\Theta^m)$ satisfying the identity 
\def\theequation{1.12}\begin{equation}
w\equiv 
\overline{\Theta} 
\big(
z,\bar z,\Theta(\bar z,z,w)
\big),
\end{equation}
reflecting the fact that $M$ is real.

\def\thedefinition{1.13}\begin{definition}{\rm
(\cite{ ber1999, me2005a, me2005b, mp2005}) $M$ is {\sl finitely
nondegenerate}\, if there exists an integer $\kappa \geqslant 1$ such
that the local holomorphic map
\def\theequation{1.14}\begin{equation}
(\bar z,\bar w)
\longmapsto
\big(
\overline{\Theta}_{z^\beta}^j
(0,\bar z,\bar w)
\big)_{\vert\beta\vert\leqslant\kappa}^{1\leqslant j\leqslant m}
\end{equation}
is of rank $n+m$ at $(\bar z, \bar w) = (0, 0)$.

}\end{definition}

From~\thetag{ 1.12}, the map $\bar w \mapsto \overline{ \Theta} (0, 0,
\bar w)$ is already of rank $m$ at $\bar w = 0$. One then verifies
(\cite{ ber1999, me2005a, me2005b, mp2005}) that there exist
multiindices $\beta(1), \dots, \beta(n) \in \N^n$ with $\vert \beta(k)
\vert \geqslant 1$ for $k=1, \dots, n$ and $\max_{1 \leqslant k
\leqslant n}\, \vert \beta (k) \vert = \kappa$ together with
integers $j(1), \dots,j(n)$ with $1 \leqslant j(k) \leqslant m$ 
such that the local holomorphic map
\def\theequation{1.15}\begin{equation}
\C^{n+m}\ni(\bar z, \bar w) 
\longmapsto 
\Big(
\big(
\overline{\Theta}^j(0,\bar z, \bar w)\big)^{
1\leqslant j\leqslant m}, 
\Big(
\overline{\Theta}_{z^{\beta(k)}}^{j(k)}
(0,\bar z,\bar w)
\Big)_{
1\leqslant k\leqslant n}
\Big)
\in\C^{m+n}
\end{equation}
is of rank $n+m$ at $(\bar z, \bar w) = (0, 0)$.

\subsection*{ 1.16.~Associated system of partial differential
equations} Generalizing an idea which goes back to B.~Segre
in~\cite{se1931, se1932} ($n = m = 1$), applied by \'E.~Cartan
in~\cite{ca1932a} and studied more recently in~\cite{su2001, gm2003a},
we may associate to $M$ a system of partial differential equations of
the form \thetag{ $\mathcal{ E}$} as follows. Complexifying the
variables $\bar z$ and $\bar w$, we introduce new independent
variables $\zeta \in \C^n$ and $\xi \in \C^m$ together with the
complex algebraic or analytic $m$-codimensional submanifold $\mathcal{
M}$ of $\C^{2 (n+m)}$ defined by
\def\theequation{1.17}\begin{equation}
w^j
=
\overline{\Theta}^j
(z,\zeta,\xi), 
\ \ \ \ \ \ \ \ \
j=1,\dots,m.
\end{equation}
We then consider the ``dependent variables'' $w^j$ as algebraic or
analytic functions of the ``independent variables'' $z^k$, with
additional dependence on the extra ``parameters'' $(\zeta, \xi)$.
Then by applying the differentiation $\partial^{ \vert \alpha \vert}
/\partial z^\alpha$ to~\thetag{ 1.17}, we get $w_{ z^\alpha }^j (z) =
\overline{ \Theta}_{ z^\alpha}^j (z, \zeta, \xi)$. Assuming finite
nondegeneracy and writing these equations for $(j, \alpha) = (j(k),
\beta( k ))$, we obtain a system of $m + n$
equations:
\def\theequation{1.18}\begin{equation}
\left\{
\aligned
w^j(z)
&
=
\overline{\Theta}^j(z,\zeta,\xi),
\ \ \ \ \ j=1,\dots,m,\\
w_{z^{\beta(k)}}^{j(k)}(z)
&
=
\overline{\Theta}_{z^{\beta(k)}}^{j(k)}
(z,\zeta,\xi), \ \ \ \ \ k=1,\dots,n.
\endaligned\right.
\end{equation}
By means of the implicit function theorem we can solve:
\def\theequation{1.19}\begin{equation}
(\zeta,\xi)
= 
R
\big(
z^k,w^j(z),w_{z^{\beta(k)}}^{j(k)}(z) 
\big). 
\end{equation}
Finally, for every pair $(j, \alpha)$ different from $(j,0)$ and from
$(j( k), \beta (k))$, we may replace $( \zeta, \xi)$ by $R$ in the
differentiated expression $w_{ z^\alpha}^j (z) = \overline{ \Theta}_{
z^\alpha}^j (z, \zeta, \xi)$, which yields
\def\theequation{1.20}\begin{equation}
\aligned
w_{z^\alpha}^j(z)
& \
=
\overline{\Theta}_{z^\alpha}^j
\left(z,R
\big(
z^k,w^j(z),w_{z^{\beta(k)}}^{j(k)}(z)
\big)
\right)
\\
& \
=: 
F_\alpha^j
\left(
z^k,w^j(z),w_{z^{\beta(k)}}^{j(k)}(z)
\right).
\endaligned
\end{equation}
This is the {\sl system of partial differential equations associated
to $M$}. 

\def\theexample{1.21}\begin{example}{\rm
(Continued) With $n = m = 1$, {\it i.e.} $M \subset \C^2$ and $\kappa
= 1$, {\it i.e.} $M$ is Levi nondegenerate of equation
\def\theequation{1.22}\begin{equation}
w
=
\bar w 
+ 
i\,z\bar z
+
{\rm O}_3,
\end{equation}
where $z$, $\bar z$ are assigned weight 1 and $w$, $\bar w$ weight 2,
B.~Segre~\cite{ se1931} obtained $w_{ zz} = F (z, w,
w_z)$. J.~Faran~\cite{ fa1980} found some examples of such equations
that {\it cannot}\, come from a $M \subset \C^2$. But
the following was left unsolved.

\def\theopenproblem{1.23}\begin{openproblem}
Characterize equations $y_{ xx} = F (x, y, y_x)$ associated to a real
analytic, Levi nondegenerate ({\it i.e.} $\kappa = 1$) hypersurface $M
\subset \C^2$. Can on read the reality condition~\thetag{ 1.12} on
$F$~? In case of success, generalize to arbitrary $M \subset \C^{ n+
m}$.
\end{openproblem}

}\end{example}

\def\theexample{1.24}\begin{example}{\rm
(Continued) Similarly, the system~\thetag{ $\mathcal{ E}_2$} 
comes from a Levi
nondegenerate hypersurface $M \subset \C^{ n + 1}$ (\cite{ ha1937,
cm1974, ch1975, su2001}. Exercise: why \thetag{ $\mathcal{ E}_3$ }
cannot come from any $M \subset \C^\nu$~?
}\end{example}

\def\theexample{1.25}\begin{example}{\rm
(Continued) With $n=1$, $m=2$ and $\kappa = 1$, the system~\thetag{
$\mathcal{ E}_4$ } comes from a $M \subset \C^3$ which is Levi
nondegenerate and satisfies
\def\theequation{1.26}\begin{equation}
T^cM
+
[T^cM,T^cM]
+
\big[ 
T^c M,[T^c M,T^c M]
\big] 
=
TM 
\end{equation}
at the origin, namely which has equations of the following form,
after some elementary transformations (\cite{ be1997, bes2005}):
\def\theequation{1.27}\begin{equation}
\aligned
w^1
&
=
\bar w^1
+
i\,z\bar z
+
{\rm O}_4,
\\
w^2
&
=
\bar w^2
+
i\,z\bar z(z+\bar z)
+
{\rm O}_4,
\endaligned
\end{equation}
where $z$, $\bar z$ are assigned weight 1 and
$w^1$, $w^2$, $\bar w^1$, $\bar w^2$
weight 2.
}\end{example}

\def\theexample{1.28}\begin{example}{\rm
(Continued) With $n=2$, $m=1$ and $\kappa = 2$, the system~\thetag{
$\mathcal{ E}_5$} comes from a hypersurface 
$M \subset \C^3$ of equation (\cite{ eb1998,
gm2003b, fk2005a, fk2005b, eb2006, gm2006}):
\def\theequation{1.29}\begin{equation}
w
=
\bar w 
+
i\,
\frac{
2\,z^1\bar z^1
+
z^1z^1\bar z^2
+
\bar z^1\bar z^1z^2
}{1-z^2\bar z^2}
+
{\rm O}_4,
\end{equation}
where $z_1$, $\bar z_1$, $z_2$, $\bar z_2$ are assigned weight 1 and
$w$, $\bar w$ weight 2, 
with the assumption that the Levi form has
rank exactly one at every point, 
and with the assumption that
$M$ is 2-nondegenerate at $0$.
}\end{example}

\subsection*{1.30.~Jet spaces, contact forms and Frobenius 
integrability} Throughout the present Part~I, we assume that
the system~\thetag{ $\mathcal{ E}$} is {\sl completely integrable},
namely that the Pfaffian system naturally associated to 
\thetag{ $\mathcal{ E
}$} in the appropriate jet space is involutive in the sense of
Frobenius. This holds automatically in case \thetag{ $\mathcal{ E }$}
comes from a generic submanifold $M \subset \C^{ n+ m}$. In general, we
will construct a {\sl submanifold of solutions}\, associated to
\thetag{ $\mathcal{ E }$}. So, we must explain complete integrability.

We denote by $\mathcal{ J}_{ n, m }^\kappa$ the space of $\kappa$-th
jets of maps $\K^n \ni x \mapsto y(x) \in \K^m$. Let
\def\theequation{1.31}\begin{equation}\left(
x^i, y^j, y_{i_1}^j, y_{i_1,i_2}^j, 
\dots\dots, 
y_{i_1,i_2,\dots,i_\kappa}^j
\right)
\in \K^{n+m+mn+mn^2+\cdots+mn^\kappa},
\end{equation}
denote the natural coordinates on $\mathcal{ J}_{ n, m}^\kappa \simeq
\K^{n+m(1+ n+\cdots + n^\kappa )}$. 
For instance, $(x, y, y_1) \in 
\mathcal{ J}_{ 1, 1}^1$.
We shall sometimes write them
shortly:
\def\theequation{1.32}\begin{equation}
\big(
x^i,y^j, y_\beta^j
\big)
\in
\K^{n+m+m(n+\cdots+n^\kappa)},
\end{equation}
where $\beta \in \N^n$ varies
and satisfies $\vert \beta \vert \leqslant \kappa$.
Sometimes also, we consider these jet coordinates only up to their
symmetries $y_{ i_1, i_2, \dots, i_\lambda }^j = y_{ i_{\sigma (1) },
i_{ \sigma (2)}, \dots, i_{ \sigma( \lambda )}}^j$, where $\sigma$ is
a permutation of $\{ 1, 2, \dots, \lambda
\}$, so that $\mathcal{ J}_{ n, m
}^\kappa \simeq \K^{ n + m \, C_{ n+ \kappa }^\kappa}$, with $C_{
n+\kappa }^\kappa := \frac{ (n+ \kappa)! }{ \kappa ! \ n!}$.

Having these notations at hand, we may develope the canonical system
of contact forms on $\mathcal{ J}_{ n, m}^\kappa$ (\cite{ol1995},
\cite{stk2000}):
\def\theequation{1.33}\begin{equation}
\left\{
\aligned
\theta^j 
:= 
& \ 
dy^j
-
\sum_{k=1}^n\,y_k^j\, dx^k,
\\
\theta_{i_1}^j 
:=
& \
dy_{i_1}^j
-
\sum_{k=1}^n\, 
y_{i_1,k}^j\,dx^k,
\\
\cdots\cdots\cdots\cdots
& \
\cdots\cdots\cdots\cdots\cdots\cdots\cdots\cdots \\
\theta_{i_1,\dots,i_{\kappa-1}}^j 
:=
& \
dy_{i_1,\dots,i_{\kappa-1}}^j
-
\sum_{k=1}^n\,y_{i_1,\dots,i_{\kappa-1},k}^j\,dx^k.
\endaligned\right.
\end{equation}
For instance, with $n=m=1$ and 
$\kappa = 2$, we have
$\theta^1 = dy - y_1\, dx$ and
$\theta_1^1 = 
dy_1 - y_2 \, dx$.
These (linearly independent) one-forms generate a subspace $\mathcal{
CT }_{ n, m }^\kappa$ of the cotangent
$T^* \mathcal{ J}_{ n, m}^\kappa$ whose
dimension equals $m\, C_{ n+ \kappa - 1 }^{ \kappa - 1}$. For the
duality between forms and vectors, the orthogonal $(\mathcal{ CT }_{n,
m }^\kappa)^\bot$ in $T \mathcal{ J}_{ n, m }^\kappa$ is spanned by
the $n + m\, C_{n + \kappa -1 }^\kappa$ vector fields:
\def\theequation{1.34}\begin{equation}
\left\{
\aligned
D_i
:=
& \
\frac{\partial}{\partial x^i}
+
\sum_{j_1=1}^m\,y_i^{j_1}\, 
\frac{\partial}{\partial y^{j_1}}
+\cdots+
\sum_{j_1=1}^m\,
\sum_{k_1,\dots,k_{\kappa-1}=1}^n\, 
y_{i,k_1,\dots,k_{\kappa-1}}^{j_1}\, 
\frac{\partial}{\partial y_{k_1,\dots,k_{\kappa-1}}^{j_1}}, 
\\
T_{i_1,\dots,i_\kappa}^{j_1}
:=
& \
\frac{\partial}{\partial y_{i_1,\dots,i_\kappa}^{j_1}},
\endaligned\right.
\end{equation}
the first $n$ ones being the total differentiation operators, 
considered in Part~II. For $n=m=1$, $\kappa = 2$, we get $\frac{
\partial }{ \partial x} + y_1 \, \frac{ \partial }{ \partial y} + y_2
\, \frac{ \partial }{ \partial y_1}$ and $\frac{ \partial }{ \partial
y_2}$.

Classically (\cite{ ol1986, bk1989, ol1995}), one associates to
\thetag{ $\mathcal{ E }$} its {\sl skeleton} $\Delta_\mathcal{ E}$,
namely the $(n+ m +p)$-dimensional submanifold of $\mathcal{ J}_{ n, m
}^{ \kappa + 1}$ simply defined by the graphed equations:
\def\theequation{1.35}\begin{equation}
y_\alpha^j
=
F_\alpha^j
\Big(x,y,
\big(
y_{\beta(q)}^{j(q)}
\big)_{1\leqslant q \leqslant p}
\Big),
\end{equation}
for $\big( j, \alpha \big) \neq (j, 0)$
and $\neq \big( j(q), \beta(q) \big)$
with $\vert \alpha \vert \leqslant \kappa +1$.
Clearly, the natural coordinates on
$\Delta_\mathcal{ E }$ are:
\def\theequation{1.36}\begin{equation}
\Big(
x,y,
\big(
y_{\beta(q)}^{j(q)}
\big)_{
1\leqslant q \leqslant p}
\Big)
\equiv
\Big(
x,y,
\big(
y_{l_1(q),\dots,l_{\lambda_q}(q)}^{j(q)}
\big)_{
1\leqslant q\leqslant p}\Big)
\in
\K^n\times\K^m\times\K^p,
\end{equation}
where $\lambda_q := \vert \beta (q) \vert$ and $\big( l_1 (q),
\dots,l_{ \lambda_q} (q) \big) := \beta (q )$.

Next, in view of the form~\thetag{1.34} of the generators of
$(\mathcal{ CT }_{ n, m}^{ \kappa +1} 
)^\bot$ and in view of the equations
of $\Delta_\mathcal{ E}$, the intersection
\def\theequation{1.37}\begin{equation}
(\mathcal{CT}_{n,m}^{\kappa+1})^\bot
\cap 
T\Delta_\mathcal{E}
\end{equation} 
is a vector subbundle of $T \Delta_\mathcal{E}$ that is generated by
$n$ linearly independent vector fields
obtained by restricting the $D_i$ to 
$\Delta_{ \mathcal{ E}}$, which yields:
\def\theequation{1.38}\begin{equation}
\left\{
\aligned
{\sf D}_i
=
\frac{\partial }{\partial x^i}+
& \
\sum_{j=1}^m\,{\bf A}_i^j
\Big(
x^{i_1},y^{j_1}, 
y_{\beta(q_1)}^{j(q_1)}
\Big)
\frac{\partial}{\partial y^j}+ \\ 
& \
+
\sum_{q=1}^p\, {\bf B}_i^q
\Big(
x^{i_1},y^{j_1}, 
y_{\beta(q_1)}^{j(q_1)}
\Big)\, 
\frac{\partial }{\partial y_{\beta(q)}^{j(q)}},
\endaligned\right.
\end{equation}
$i = 1, \dots, n$, where the coefficients 
${\bf A}_i^j$ and ${\bf B}_i^q$ are
given by:
\def\theequation{1.39}\begin{equation}
\aligned
{\bf A}_i^j:=
& \
\left\{
\aligned
{}
&
y_i^j \ \text{\rm if the variable} \ y_i^j \ \text{\rm 
appears among the $p$ variables} \ 
y_{\beta(q_1)}^{j(q_1)};
\\
{}
&
F_i^j \ \text{\rm otherwise};\\
\endaligned\right.\\
{\bf B}_i^q:=
& \
\left\{
\aligned
{}
&
y_{i, l_1(q), \dots, l_{\lambda_q} (q)}^{ j(q)} \ \text{\rm if} 
\ y_{l_1 (q), \dots, l_{ \lambda_q (q)}}^{ j(q)} \ \text{\rm
appears among the $p$ variables} \
y_{\beta(q_1)}^{j(q_1)};\\
{}
&
F_{i,l_1(q),\dots,l_{\lambda_q}(q)}^{j(q)} \ 
\text{\rm otherwise}.\\
\endaligned\right.\\
\endaligned
\end{equation}

\def\theexample{1.40}\begin{example}{\rm
For \thetag{ $\mathcal{ E}_1$}, we get ${\sf D} = \frac{ \partial }{
\partial x} +y_1\, \frac{ \partial }{ \partial x} + F (x, y, y_1) \,
\frac{ \partial }{ \partial y_2}$; exercise: treat \thetag{ $\mathcal{
E}_2$} and \thetag{ $\mathcal{ E}_3$}. For \thetag{ $\mathcal{
E}_4$}, we get ${\sf D} = \frac{ \partial }{ \partial x} + y_1^1 \,
\frac{ \partial }{ \partial y^1} + F \, \frac{ \partial }{ \partial
y^2} + G\, \frac{ \partial }{ \partial y_1^1}$. For \thetag{
$\mathcal{ E}_5$}, whose skeleton is written $y_2 = F$, $y_{ 1, 1, 1}
= G$, $y_{ 1, 2} = H$, $y_{ 1, 1, 2} = K$, with $F$, $G$, $H$, $K$
being functions of $\big( x^1, x^2, y, y_1, y_{ 1, 1} \big)$, we get
\def\theequation{1.41}\begin{equation}
\aligned
{\sf D}_1
&
=
\frac{\partial}{\partial x^1}
+
y_1\,\frac{\partial}{\partial y}
+
y_{1,1}\,\frac{\partial}{\partial y_1}
+
G\,\frac{\partial}{\partial y_{1,1}},
\\
{\sf D}_2
&
=
\frac{\partial}{\partial x^2}
+
F\,\frac{\partial}{\partial y}
+
H\,\frac{\partial}{\partial y_1}
+
K\,\frac{\partial}{\partial y_{1,1}}.
\endaligned
\end{equation}
}\end{example}

\def\thedefinition{1.42}\begin{definition}{\rm
The system \thetag{$\mathcal{ E }$} is {\sl completely
integrable} if the $n$ vector fields~\thetag{1.38} satisfy the
Frobenius integrability condition, namely every Lie bracket $[{\sf D}_{
i_1}, {\sf D}_{i_2 }]$, $1 \leqslant i_1, i_2 \leqslant n$, is a linear
combination of the vector fields ${\sf D}_1, \dots, {\sf D}_n$. 
}\end{definition}

Because of their specific form~\thetag{1.38}, we must then have in fact
$[{\sf D}_{i_1}, {\sf D}_{i_2 }] = 0$. For $n = 1$, the condition is
of course void.

\section*{\S2.~Submanifold of solutions}

\subsection*{2.1.~Fundamental foliation of the skeleton} 
As the vector fields ${\sf D}_i$ commute, they equip the skeleton
$\Delta_\mathcal{ E} \simeq \K^{ n+ m+ p}$ with a foliation ${\sf
F}_{\! \Delta_\mathcal{ E} }$ by $n$-dimensional integral manifolds
which are (approximately) directed along the $x$-axis. We draw a
diagram (see only the left side).

\bigskip
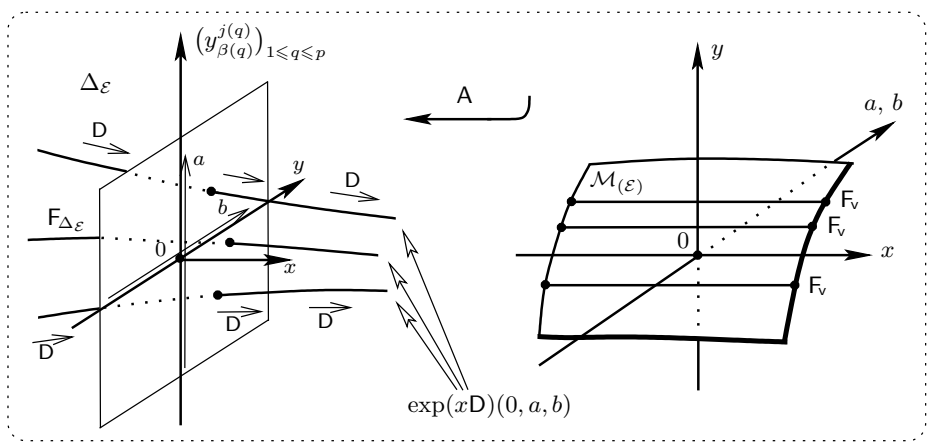
\begin{center}
\begin{picture}(0,0)%
\includegraphics{redressement-A.pstex}%
\end{picture}%
\setlength{\unitlength}{4144sp}%
\begingroup\makeatletter\ifx\SetFigFont\undefined
\def\x#1#2#3#4#5#6#7\relax{\def\x{#1#2#3#4#5#6}}%
\expandafter\x\fmtname xxxxxx\relax \def\y{splain}%
\ifx\x\y   
\gdef\SetFigFont#1#2#3{%
  \ifnum #1<17\tiny\else \ifnum #1<20\small\else
  \ifnum #1<24\normalsize\else \ifnum #1<29\large\else
  \ifnum #1<34\Large\else \ifnum #1<41\LARGE\else
     \huge\fi\fi\fi\fi\fi\fi
  \csname #3\endcsname}%
\else
\gdef\SetFigFont#1#2#3{\begingroup
  \count@#1\relax \ifnum 25<\count@\count@25\fi
  \def\x{\endgroup\@setsize\SetFigFont{#2pt}}%
  \expandafter\x
    \csname \romannumeral\the\count@ pt\expandafter\endcsname
    \csname @\romannumeral\the\count@ pt\endcsname
  \csname #3\endcsname}%
\fi
\fi\endgroup
\begin{picture}(5514,2589)(484,-2233)
\put(5720,-1121){\makebox(0,0)[lb]{\smash{\SetFigFont{9}{10.8}{rm}{\color[rgb]{0,0,0}$x$}%
}}}
\put(4501,-1066){\makebox(0,0)[lb]{\smash{\SetFigFont{9}{10.8}{rm}{\color[rgb]{0,0,0}$0$}%
}}}
\put(1753,-862){\makebox(0,0)[lb]{\smash{\SetFigFont{8}{9.6}{rm}{\color[rgb]{0,0,0}$b$}%
}}}
\put(5619,-239){\makebox(0,0)[lb]{\smash{\SetFigFont{9}{10.8}{rm}{\color[rgb]{0,0,0}$a,\,b$}%
}}}
\put(4705, 86){\makebox(0,0)[lb]{\smash{\SetFigFont{9}{10.8}{rm}{\color[rgb]{0,0,0}$y$}%
}}}
\put(2148,-1217){\makebox(0,0)[lb]{\smash{\SetFigFont{8}{9.6}{rm}{\color[rgb]{0,0,0}$x$}%
}}}
\put(2201,-600){\makebox(0,0)[lb]{\smash{\SetFigFont{8}{9.6}{rm}{\color[rgb]{0,0,0}$y$}%
}}}
\put(1605,-568){\makebox(0,0)[lb]{\smash{\SetFigFont{8}{9.6}{rm}{\color[rgb]{0,0,0}$a$}%
}}}
\put(1378,-1117){\makebox(0,0)[lb]{\smash{\SetFigFont{8}{9.6}{rm}{\color[rgb]{0,0,0}$0$}%
}}}
\put(1606,134){\makebox(0,0)[lb]{\smash{\SetFigFont{9}{10.8}{rm}{\color[rgb]{0,0,0}$\big(y_{\beta(q)}^{j(q)}\big)_{1\leqslant q\leqslant p}$}%
}}}
\put(732,-933){\makebox(0,0)[lb]{\smash{\SetFigFont{9}{10.8}{rm}{\color[rgb]{0,0,0}${\sf F}_{\!\Delta_\mathcal{E}}$}%
}}}
\put(924,-126){\makebox(0,0)[lb]{\smash{\SetFigFont{9}{10.8}{rm}{\color[rgb]{0,0,0}$\Delta_\mathcal{E}$}%
}}}
\put(5476,-823){\makebox(0,0)[lb]{\smash{\SetFigFont{9}{10.8}{rm}{\color[rgb]{0,0,0}${\sf F}_{\!\sf v}$}%
}}}
\put(5398,-980){\makebox(0,0)[lb]{\smash{\SetFigFont{9}{10.8}{rm}{\color[rgb]{0,0,0}${\sf F}_{\!\sf v}$}%
}}}
\put(5284,-1337){\makebox(0,0)[lb]{\smash{\SetFigFont{9}{10.8}{rm}{\color[rgb]{0,0,0}${\sf F}_{\!\sf v}$}%
}}}
\put(997,-402){\makebox(0,0)[lb]{\smash{\SetFigFont{8}{9.6}{rm}{\color[rgb]{0,0,0}${\sf D}$}%
}}}
\put(678,-1733){\makebox(0,0)[lb]{\smash{\SetFigFont{8}{9.6}{rm}{\color[rgb]{0,0,0}${\sf D}$}%
}}}
\put(1798,-1565){\makebox(0,0)[lb]{\smash{\SetFigFont{8}{9.6}{rm}{\color[rgb]{0,0,0}${\sf D}$}%
}}}
\put(2348,-1546){\makebox(0,0)[lb]{\smash{\SetFigFont{8}{9.6}{rm}{\color[rgb]{0,0,0}${\sf D}$}%
}}}
\put(2509,-702){\makebox(0,0)[lb]{\smash{\SetFigFont{8}{9.6}{rm}{\color[rgb]{0,0,0}${\sf D}$}%
}}}
\put(2889,-2050){\makebox(0,0)[lb]{\smash{\SetFigFont{9}{10.8}{rm}{\color[rgb]{0,0,0}$\exp(x{\sf D})(0,a,b)$}%
}}}
\put(3978,-698){\makebox(0,0)[lb]{\smash{\SetFigFont{8}{9.6}{rm}{\color[rgb]{0,0,0}$\mathcal{M}_{(\mathcal{E})}$}%
}}}
\put(3180,-204){\makebox(0,0)[lb]{\smash{\SetFigFont{9}{10.8}{rm}{\color[rgb]{0,0,0}${\sf A}$}%
}}}
\end{picture}

\end{center}

The (abstract, not numerical) integration of \thetag{ $\mathcal{ E}$}
is thus straightforwardly completed: the set of solutions coincides
with the set of leaves of ${\sf F}_{ \Delta_\mathcal{ E }}$. This is
the true geometric content, viewed in the appropriate jet space, of the
assumption of complete integrability.

\subsection*{2.2.~General solution and submanifold of solutions}
To construct the submanifold of solutions $\mathcal{ M}_{( \mathcal{
E})}$ associated to \thetag{ $\mathcal{ E}$}
(sketched in the right hand side), we execute some elementary
analytico-geometric constructions.

At first, we {\it duplicate}\, the coordinates $\big( y_{ \beta (q)}^{
j(q) }, y^j \big) \in \K^p \times \K^m$ by introducing a new subspace
of coordinates $(a, b) \in \K^p \times \K^m$; thus, on the left
diagram, we draw a vertical plane together with $a$- and $b$-axes.
The leaves of the foliation ${\sf
F}_{\! \Delta_\mathcal{ E} }$ are uniquely determined by their
intersections with this plane, consisting of points of coordinates
$(0, a, b) \in \K^n \times \K^p \times \K^m$.

Such points $(0, a, b)$ correspond to the {\sl initial conditions}
$\big( y_{x^{ \beta (q)}}^{ j(q)} (0), 
y (0) \big)$ for the general solution of $(\mathcal{ E
})$. In fact, the (concatenated, multiple) flow
of $\{ {\sf D}_1, \dots, {\sf D}_n \}$ is given by 
\def\theequation{2.3}\begin{equation}
\exp
\big(
x^n{\sf D}_n\,
\big( 
\cdots
(\exp(x^1{\sf D}_1(0,a,b)))
\cdots 
\big)
\big)
=
\big(
x,
\Pi(x,a,b), 
\Omega(x,a,b)
\big)\in\K^n\times\K^m\times\K^p,
\end{equation}
for some two local analytic maps $\Pi = (\Pi^1, \dots, \Pi^m)$ and
$\Omega = (\Omega^1, \dots, \Omega^p)$ and the next lemma is
straightforward.

\def\thelemma{2.4}\begin{lemma}
The general solution of \thetag{ $\mathcal{ E}$}
is 
\def\theequation{2.5}\begin{equation}
y(x)
:=
\Pi(x,a,b),
\end{equation}
where $(a,b)$ varies in $\K^p \times \K^m$.
Furthermore, for $q=1, \dots, p${\rm :}
\def\theequation{2.6}\begin{equation}
\Omega^q(x,a,b)
\equiv
\Pi_{x^{\beta(q)}}^{j(q)}
(x,a,b).
\end{equation}
\end{lemma}

This leads to introducing a fundamental geometric object.

\def\thedefinition{2.7}\begin{definition}{\rm
The {\sl submanifold of solutions} $\mathcal{ V}_{\mathcal{ S}}
(\mathcal{ E})$ associated to \thetag{ $\mathcal{ E}$}
is
the analytic submanifold of $\K_x^n \times \K_a^p \times \K_b^m$
defined by the Cartesian equations{\rm :}
\def\theequation{2.8}\begin{equation}
0
=
-
y^j
+
\Pi^j(x,a,b), 
\ \ \ \ \ \ \ \ 
j=1,\dots,m.
\end{equation}
}\end{definition}

There is a strong interplay between the study of $(\mathcal{ E })$ and
the geometry of $\mathcal{ V}_{ \mathcal{S}}( \mathcal{ E})$. By
construction, the diffeomorphism:
\def\theequation{2.9}\begin{equation}
\left\{
\aligned
{}
&
{\sf A} : 
\K^{n+p+m} \ [{\rm coordinates} \ (x^i,a^q,b^j)] \longrightarrow
\K^{n+m+p}\ \left[{\rm coordinates} \ 
\left(
x^i,y^j,y_{\beta(q)}^{j(q)}
\right)\right]\\
& 
{\sf A}\,(x^i,a^q,b^j)
:=
\left(
x^i,\
\Pi^j(x,a,b),\
\Pi_{x^{\beta(q)}}^{j(q)}(x,a,b),
\right),
\endaligned\right.
\end{equation}
sends the {\sl foliation ${\sf F}_{\sf v}$ by the {\sf v}ariables $x$}
whose leaves are $\{ a = {\rm cst.}, \ b = {\rm cst}. \}$ (see the
diagram), to the previous foliation ${\sf F}_{\! \Delta_\mathcal{ E }
}$.

\subsection*{2.10.~{\sc pde} system associated to a submanifold}
Inversely, let $\mathcal{ M}$ be a submanifold of $\K_x^n \times
\K_y^m \times \K_a^p \times \K_b^m$ of the form
\def\theequation{2.11}\begin{equation}
y^j
=
\Pi^j(x,a,b), 
\ \ \ \ \ \ \ \ \ \ \ \ 
j=1,\dots,m.
\end{equation}
A necessary
condition for it to be the complexification of a generic $M \subset
\C^{ n + m}$ is that $p = n$ (answer to an exercise above).

\def\thedefinition{2.12}\begin{definition}{\rm 
$\mathcal{M}$ is {\sl solvable with respect to the parameters} if
$b \mapsto \Pi (0, 0, b)$ of rank $m$ at $b = 0$ and if
there exist $\kappa \geqslant 1$, multiindices $\beta(1), \dots,
\beta(p) \in \N^n$ with $\vert \beta(q) \vert \geqslant 1$ for $q = 1,
\dots, p$ and $\max_{ 1\leqslant q \leqslant p}\, \vert \beta (q)
\vert = \kappa$, together with integers $j(1), \dots, j(p)$ with $1
\leqslant j(q) \leqslant m$ such that the local
$\K$-analytic map
\def\theequation{2.13}\begin{equation}
\K^{m+p}\ni 
(a,b)
\longmapsto 
\Big(
\big(
\Pi^j(0,a,b)
\big)^{
1\leqslant j\leqslant m},\ 
\Big(
\Pi_{x^{\beta(q)}}^{j(q)}
(0,a,b)
\Big)_{1\leqslant q \leqslant p}
\Big)
\in\K^{m+p}
\end{equation}
is of rank equal to $m+p$ at $(a, b) = (0, 0)$
}\end{definition}

When $\mathcal{ M}$ is the submanifold of solutions of a system
\thetag{ $\mathcal{ E}$}, it is automatically solvable with respect to
the variables, the pairs $(j(q), \beta ( q ))$ being the same as in the
arguments of the right hand sides $F_\alpha^j$ in \thetag{ $\mathcal{
E }$}. Proceeding as in~\S1.16, we may associate to $\mathcal{ M}$ a
system of the form~\thetag{ $\mathcal{ E }$}. Since
we need introduce some new notation, 
let us repeat the argument.

Considering $y = y(x) = \Pi (x, a, b)$ as a function of $x$ with extra
parameters $(a, b)$ and applying $\partial^{ \vert \alpha \vert} \big/
\partial x^\alpha$, we get $y_{ x^\alpha}^j (x) = \Pi_{ x^\alpha}^j
(x, a, b)$. Writing only the relevant $( m + p)$ equations:
\def\theequation{2.14}\begin{equation}
\left\{
\aligned
y^j(x)
&
=
\Pi^j(x,a,b),
\\
y_{x^{\beta(q)}}^{j(q)}
&
=
\Pi_{x^{\beta(q)}}^{j(q)}
(x,a,b),
\endaligned\right.
\end{equation}
the assumption of solvability with respect to 
parameters enables to 
get
\def\theequation{2.15}\begin{equation}
\left\{
\aligned
a^q
&
=
A^q
\big(
x^i,y^j,y_{\beta(q_1)}^{j(q_1)}
\big),
\\
b^j
&
=
B^j
\big(
x^i,y^{j_1},y_{\beta(q)}^{j(q)}
\big).
\endaligned\right.
\end{equation}
For every $(j, \alpha) \neq (j, 0)$ and $\neq (j(q), \beta (q))$, we
then replace $(a, b)$ in $y_{ x^\alpha}^j = \Pi_{ x^\alpha}^j$:
\def\theequation{2.16}\begin{equation}
\aligned
y_{x^\alpha}^j(x)
&
=
\Pi_{x^\alpha}^j
\Big(
x,
A
\big(
x^i,y^{j_1}(x),y_{\beta(q)}^{j(q)}(x)
\big),
B
\big(
x^i,y^{j_1}(x),y_{\beta(q)}^{j(q)}(x)
\big)
\Big)
\\
&
=:
F_\alpha^j
\big(
x^i,
y^{j_1}(x),
y_{x^{\beta(q)}}^{j(q)}(x)
\big).
\endaligned
\end{equation}

\def\theproposition{2.17}\begin{proposition}
There is a one-to-one correspondence 
\def\theequation{2.18}\begin{equation}
(\mathcal{E}_\mathcal{M})
=
(\mathcal{E})
\longleftrightarrow
\mathcal{M}
=
\mathcal{M}_{(\mathcal{E})},
\end{equation}
between completely integrable systems of partial differential
equations of the general form \thetag{ $\mathcal{ E}$} and
submanifolds {\rm (}of solutions{\rm )} $\mathcal{ M}$ of the
form~\thetag{2.11} which are solvable with respect to the parameters.
Of course
\def\theequation{2.19}\begin{equation}
\big(
\mathcal{E}_{\mathcal{M}_{(\mathcal{E})}}
\big)
=
(\mathcal{E})
\ \ \ \ \ \ \ \ \ \ \ 
\text{\rm and}
\ \ \ \ \ \ \ \ \ \ \
\mathcal{M}_{(\mathcal{E}_\mathcal{M})}
=
\mathcal{M}.
\end{equation}
\end{proposition}

\subsection*{2.20.~Transfer of total differentiations}
We notice that the auxiliary functions $A^q$ and $B^j$ enable to
express the inverse of ${\sf A}$:
\def\theequation{2.21}\begin{equation}
{\sf A}^{-1}:\ \
\big(
x^{i_1},
y^{j_1},
y_{\beta(q_1)}^{j(q_1)}
\big)
\longmapsto
\Big(
x^i,
A^q
\big(
x^{i_1},
y^{j_1},
y_{\beta(q_1)}^{j(q_1)}
\big),
B^j
\big(
x^{i_1},
y^{j_1},
y_{\beta(q_1)}^{j(q_1)}
\big)
\Big).
\end{equation}
More importantly, the total differentiation operator considerably
simplifies when viewed on $\mathcal{ M}$. This observation is useful
for translating differential invariants of \thetag{ $\mathcal{ E}$} as
differential invariants of $\mathcal{ M}$.

\def\thelemma{2.22}\begin{lemma}
Through ${\sf A}$, for $i=1, \dots, n$, the pull-back of the total
differentiation operator ${\sf D}_i$ is simply $\frac{ \partial }{
\partial x^i}$, or equivalently{\rm :}
\def\theequation{2.23}\begin{equation}
{\sf A}_*
\Big(
\frac{\partial}{\partial x^i}
\Big)
=
{\sf D}_i.
\end{equation}
\end{lemma}

\proof
Let $\ell = \ell \big( x^i, y^j, 
y_{ \beta (q)}^{ j(q)} \big)$ be any 
function defined on $\Delta_\mathcal{ E}$. 
Composing with ${\sf A}$ yields the function
$\Lambda := \ell \circ {\sf A}$, {\it i.e.}
\def\theequation{2.24}\begin{equation}
\Lambda(x,a,b)
\equiv
\ell
\Big(
x^i,
\Pi^j(x,a,b),
\Pi_{x^{\beta(q)}}^{j(q)}
(x,a,b)
\Big).
\end{equation}
Differentiating with respect to $x^i$, we get, 
dropping the arguments:
\def\theequation{2.25}\begin{equation}
\frac{\partial\Lambda}{\partial x^i}
=
\frac{\partial\ell}{\partial x^i}
+
\sum_{j=1}^m\,\Pi_{x^i}^j\,
\frac{\partial\ell}{\partial y^j}
+
\sum_{q=1}^p\,\Pi_{x^ix^{\beta(q)}}^{j(q)}\,
\frac{\partial\ell}{\partial y_{x^{\beta(q)}}^{j(q)}}.
\end{equation}
Replacing the appearing $\Pi_{ x^\alpha }^j$ for which $(j, \alpha)
\neq (j, 0)$ and $\neq (j(q), \beta (q))$ by $F_\alpha^j$, we recover
${\sf D}_i$ as defined by~\thetag{ 1.38}, whence $\frac{ \partial
\Lambda}{ \partial x^i} = {\sf D}_i \ell$.
\endproof

\subsection*{2.26.~Transfer of algebrico-differential expressions}
The diffeomorphism ${\sf A}$ may be used to translate
algebrico-differential expressions from $\mathcal{ M }$ to \thetag{
$\mathcal{ E }$} and vice-versa:
\def\theequation{2.27}\begin{equation}
{\sf I}_\mathcal{M}
\big(
J_{x,a,b}^{\lambda+\kappa+1}\,
\Pi
\big)
\longleftrightarrow
{\sf I}_{(\mathcal{E})}
\big(
J_{x,y,y_1}^\lambda\,F
\big).
\end{equation}
Here, $\lambda \in \N$, the letter $J$ is used to denote jets, and
${\sf I} = {\sf I}_\mathcal{ M}$ or $= {\sf I}_{ (\mathcal{ E})}$ is a
polynomial or more generally, a quotient of polynomials 
with respect to its jet
arguments. Notice the shift by $\kappa +1$ of the
jet orders. 

\def\theexample{2.28}\begin{example}{\rm
Suppose $n=m=1$ and $\kappa = 1$. Then $F = \Pi_{ xx}$. As an
exercise, let us compute $F_x$, $F_y$, $F_{ y_1}$ in terms of $J_{ x,
a, b}^3 \, \Pi$. We start with the identity
\def\theequation{2.29}\begin{equation}
F(x,y,y_1)
\equiv
\Pi_{xx}
\big(
x,
A(x,y,y_1),
B(x,y,y_1)
\big),
\end{equation}
that we differentiate with respect to 
$x$, to $y$ and to $y_1$:
\def\theequation{2.30}\begin{equation}
\aligned
F_x
&
=
\Pi_{xxx}
+
\Pi_{xxa}\,A_x
+
\Pi_{xxb}\,B_x,
\\
F_y
&
=
\ \ \ \ \ \ \ \ \ \ \ \ \
\Pi_{xxa}\,A_y
+
\Pi_{xxb}\,B_y,
\\
F_{y_1}
&
=
\ \ \ \ \ \ \ \ \ \ \ \ \
\Pi_{xxa}\,A_{y_1}
+
\Pi_{xxb}\,B_{y_1}.
\endaligned
\end{equation}
Thus, we need to compute $A_x$, $A_y$, $A_{ y_1}$, $B_x$, $B_y$, $B_{
y_1}$. This is easy: it suffices to differentiate the two identities
that define $A$ and $B$ as implicit functions, namely:
\def\theequation{2.31}\begin{equation}
\aligned
y
&
\equiv
\Pi
\big(
x,
A(x,y,y_1),
B(x,y,y_1)
\big)
\ \ \ \ \ \ \
\text{\rm and}
\\
y_1
&
\equiv
\Pi_x
\big(
x,
A(x,y,y_1),
B(x,y,y_1)
\big)
\endaligned
\end{equation}
with respect to $x$, to $y$ and to $y_1$, 
which gives six new identities:
\def\theequation{2.32}\begin{equation}
\aligned
0
&
=
\Pi_x
+
\Pi_a\,A_x
+
\Pi_b\,B_x,
\ \ \ \ \ \ \ \ \ \ \ \ \ \
0
=
\Pi_{xx}
+
\Pi_{xa}\,A_x
+
\Pi_{xb}\,B_x,
\\
1
&
=
\ \ \ \ \ \ \ \ \ \ 
\Pi_a\,A_y
+
\Pi_b\,B_y,
\ \ \ \ \ \ \ \ \ \ \ \ \ \ 
0
=
\ \ \ \ \ \ \ \ \ \ \
\Pi_{xa}\,A_y
+
\Pi_{xb}\,B_y
\\
0
&
=
\ \ \ \ \ \ \ \ \ 
\Pi_a\,A_{y_1}
+
\Pi_b\,B_{y_1},
\ \ \ \ \ \ \ \ \ \ \ \
1
=
\ \ \ \ \ \ \ \ \ \ \
\Pi_{xa}\,A_{y_1}
+
\Pi_{xb}\,B_{y_1},
\endaligned
\end{equation}
and to solve each of the three linear systems of two 
equations located in a line,
noticing that their common determinant
$\Pi_b \, \Pi_{ xa} - \Pi_a \, \Pi_{ xb}$
does not vanish at the origin, since $\Pi = b + xa + {\rm O}_3$.
By elementary Cramer formulas, we get:
\def\theequation{2.33}\begin{equation}
\left\{
\aligned
A_x
&
=
\frac{
-\Pi_b\,\Pi_{xx}
+
\Pi_x\,\Pi_{xb}
}{
\Pi_b\,\Pi_{xa}-\Pi_a\,\Pi_{xb}},
\ \ \ \ \ \ \ \ \ \ \ \ \ \ \
B_x
=
\frac{
-\Pi_x\,\Pi_{xa}
+
\Pi_a\,\Pi_{xx}
}{
\Pi_b\,\Pi_{xa}-\Pi_a\,\Pi_{xb}},
\\
A_y
&
=
\frac{
-\Pi_{xb}
}{
\Pi_b\,\Pi_{xa}-\Pi_a\,\Pi_{xb}},
\ \ \ \ \ \ \ \ \ \ \ \ \ \ \ \ \ \
B_y
=
\frac{
\Pi_{xa}
}{
\Pi_b\,\Pi_{xa}-\Pi_a\,\Pi_{xb}},
\\
A_{y_1}
&
=
\frac{
\Pi_b
}{
\Pi_b\,\Pi_{xa}-\Pi_a\,\Pi_{xb}},
\ \ \ \ \ \ \ \ \ \ \ \ \ \ \ \ \
B_{y_1}
=
\frac{
-\Pi_a
}{
\Pi_b\,\Pi_{xa}-\Pi_a\,\Pi_{xb}}.
\endaligned\right.
\end{equation} 
Replacing in~\thetag{ 2.30}, no simplification occurs and we get what
we wanted:
\def\theequation{2.34}\begin{equation}
\left\{
\aligned 
F_x
&
=
\Pi_{xxx}
+
\frac{
\Pi_{xxa}
\big[
-\Pi_b\,\Pi_{xx}
+
\Pi_x\,\Pi_{xb}
\big]
+
\Pi_{xxb}
\big[
-\Pi_x\,\Pi_{xa}
+
\Pi_a\,\Pi_{xx}
\big]
}{
\Pi_b\,\Pi_{xa}-\Pi_a\,\Pi_{xb}},
\\
F_y
&
=
\frac{
-\Pi_{xxa}\,\Pi_{xb}
+
\Pi_{xxb}\,\Pi_{xa}
}{
\Pi_b\,\Pi_{xa}-\Pi_a\,\Pi_{xb}},
\\
F_{y_1}
&
=
\frac{
\Pi_{xxa}\,\Pi_b
-
\Pi_{xxb}\,\Pi_a
}{
\Pi_b\,\Pi_{xa}-\Pi_a\,\Pi_{xb}},
\endaligned\right.
\end{equation}
One sees ${\sf D}\, F = F_x + \Pi_x \, F_y + \Pi_{ xx}\, F_{ y_1} =
\Pi_{ xxx}$ simply, as predicted by Lemma~2.22.

Second order derivatives $F_{ xx}$, $F_{ xy}$, $F_{ xy_1}$, $F_{ yy}$,
$F_{ yy_1}$, $F_{ y_1 y_1}$ have still reasonable complexity, when
expressed in terms of $J_{ x,a,b }^4 \, \Pi$. Beyond, the
computations explode.
}\end{example}

\def\theopenquestion{2.35}\begin{openquestion}
A second order ordinary differential equation $y_{ xx} = F (x, y,
y_x)$ has two fundamental differential invariants, namely {\rm (\cite{
tr1896, ca1924, gtw1989, ol1995}):}
\def\theequation{2.36}\begin{equation}
\aligned
{\sf I}_{(\mathcal{E}_1)}^1
&
:=
\frac{\partial^4F}{\partial y_1^4}
\ \ \ \ \ \ \ \
\text{\rm and}
\\
{\sf I}_{(\mathcal{E}_1)}^2
&
:=
{\sf D}{\sf D}\big(F_{y_1y_1}\big)
-
F_{y_1}\,{\sf D}\big(F_{y_1y_1}\big)
-
4\,{\sf D}\big(F_{yy_1}\big)
+
6\,F_{yy}
-
3\,F_y\,F_{y_1y_1}
+
4\,F_{y_1}\,F_{yy_1}.
\endaligned
\end{equation}
Compute ${\sf I}_{\mathcal{ M}_1}^1$ and 
${\sf I}_{\mathcal{ M}_1}^2$.
\end{openquestion}

Although the notion of diffeomorphism is clear and apparently obvious
from the intuitive, geometric and conceptual viewpoints, in concrete
applications and in explicit computations, it almost never
straightforward to transfer algebrico-differential objects.

\def\theopenproblem{2.37}\begin{openproblem}
For general \thetag{ $\mathcal{ E}$} and $\mathcal{ M}$, build closed
combinatorial formulas executing the double translation~\thetag{ 2.27}.
\end{openproblem}

\subsection*{2.38.~Plan for the sequel}
We will endeavour a general theory showing that the study of systems
\thetag{ $\mathcal{ E}$} and the study of submanifolds of solutions
$\mathcal{ M }$ gives complementary views on the same object. In fact,
Lie symmetries, equivalence problems, Cartan connections, normal forms
and classification lists may be endeavoured on both sides, yielding
essentially equivalent results, though the translation is seldom
straightforward. In Section~3, 4 and~5, we review some features from
the side \thetag{ $\mathcal{ E}$}, before studying some aspects from
the side of $\mathcal{ M}$. A more systematic and complete approach
shall appear as a monography.

\section*{\S3.~Classification problems}

\subsection*{3.1.~Transformations of {\sc pde} systems} 
Through a local $\K$-analytic change of variables close to the
identity $(x, y) \mapsto \varphi (x, y) =: (x', y')$, the system
\thetag{ $\mathcal{ E}$} transforms to a similar system, with primes:
\def\theequation{$\mathcal{E}'$}\begin{equation}
{y'}_{{x'}^\alpha}^j(x')
=
{F'}_\alpha^j
\Big(
x',y'(x'),
\big({y'}_{{x'}^{\beta(q)}}^{j(q)}(x')
\big)_{
1\leqslant q\leqslant p}
\Big).
\end{equation}

\def\theexample{3.2}\begin{example}{\rm
Coming back temporarily to the
notations of \S1.12(II), with $n = m = \kappa = 1$, assume that $y_{
xx} = f (x, y, y_x)$ transforms to $Y_{ XX} = F (X, Y, Y_X)$ through a
local diffeomorphism $(x, y) \mapsto (X, Y) = \big( X (x, y), Y (x, y)
\big)$. How $F$ is related to $f$~? By symmetry, it suffices to
compute $f$ in terms of $F$, $X$, $Y$. The prolongation to $\mathcal{
J}_{ 1, 1}^2$ of the diffeomorphism has components (\cite{ bk1989,
me2004}):
\def\theequation{3.3}\begin{equation}
Y_X
=
\frac{Y_x+y_x\,Y_y}{
X_x+y_xX_y},
\end{equation}
and
\def\theequation{3.4}\begin{equation}
\aligned
Y_{XX}
=
\frac{1}{
\big[
X_x+y_xX_y
\big]^3}
& \
\left(
y_{xx} \cdot 
\left\vert
\begin{array}{cc}
X_x & X_y 
\\
Y_x & Y_y 
\end{array}
\right\vert
+
\left\vert
\begin{array}{cc}
X_x & X_{xx} \\
Y_x & Y_{xx} 
\end{array}
\right\vert
+
\right.
\\
&\ \ \ \ \
\left.
+
y_x\cdot
\left\{
2
\left\vert
\begin{array}{cc}
X_x & X_{xy} \\
Y_x & Y_{xy} 
\end{array}
\right\vert-
\left\vert
\begin{array}{cc}
X_{xx} & X_y \\
Y_{xx} & Y_y 
\end{array}
\right\vert
\right\}+ 
\right.
\\
& \ \ \ \ \
\left.
+ y_xy_x\cdot \left\{
\left\vert
\begin{array}{cc}
X_x & X_{yy} \\
Y_x & Y_{yy} 
\end{array}
\right\vert
-2
\left\vert
\begin{array}{cc}
X_{xy} & X_y \\
Y_{xy} & Y_y
\end{array}
\right\vert
\right\}+
\right.
\\
&\ \ \ \ \ 
\left.
+
y_xy_xy_x\cdot 
\left\{
-
\left\vert
\begin{array}{cc}
X_{yy} & X_y \\
Y_{yy} & Y_y
\end{array}
\right\vert
\right\}
\right).
\endaligned
\end{equation}
It then suffices to replace $Y_{XX}$ above by $F (X, Y, Y_X)$ and to
solve $y_{ xx}$:
\def\theequation{3.5}\begin{equation}
\aligned
y_{xx}
&
=
\frac{1}{
\left\vert
\begin{array}{cc}
X_x & X_y \\
Y_x & Y_y
\end{array}
\right\vert
}
\left(
\big[X_x+y_xX_y\big]^3\,
F
\left(
X,Y,
\frac{Y_x+y_x\,Y_y}{
X_x+y_xX_y}
\right)
-
\left\vert
\begin{array}{cc}
X_x & X_{xx} \\
Y_x & Y_{xx} 
\end{array}
\right\vert
+
\right.
\\
&\ \ \ \ \ \ \ \ \ \ \ \ \ \ \ \ \ \ \ \ \ \ \ \ \ \ \
\left.
+
y_x\cdot
\left\{
-2
\left\vert
\begin{array}{cc}
X_x & X_{xy} \\
Y_x & Y_{xy} 
\end{array}
\right\vert
+
\left\vert
\begin{array}{cc}
X_{xx} & X_y \\
Y_{xx} & Y_y 
\end{array}
\right\vert
\right\}+ 
\right.
\\
& \ \ \ \ \ \ \ \ \ \ \ \ \ \ \ \ \ \ \ \ \ \ \ \ \ \ \
\left.
+ 
y_xy_x\cdot\left\{
-
\left\vert
\begin{array}{cc}
X_x & X_{yy} \\
Y_x & Y_{yy} 
\end{array}
\right\vert
+2
\left\vert
\begin{array}{cc}
X_{xy} & X_y \\
Y_{xy} & Y_y
\end{array}
\right\vert
\right\}+
\right.
\\
&\ \ \ \ \ \ \ \ \ \ \ \ \ \ \ \ \ \ \ \ \ \ \ \ \ \ \
\left.
+
y_xy_xy_x\cdot 
\left\{
\left\vert
\begin{array}{cc}
X_{yy} & X_y \\
Y_{yy} & Y_y
\end{array}
\right\vert
\right\}
\right)
\\
&
=:
f(x,y,y_x).
\endaligned
\end{equation} 
}\end{example}

\def\theopenproblem{3.6}\begin{openproblem}
Find general formulas expressing the 
$F_\alpha^j$ in terms
of ${ F'}_\alpha^j$, ${x'}^i$, ${y'}^j$.
\end{openproblem}

Conversely, given two such systems \thetag{ $\mathcal{ E}$} and
\thetag{ $\mathcal{ E}'$}, {\it when do they transform to each
other}~? Let $\pi_{ \kappa, p} '$ denote the projection from
${\mathcal{ J}'}_{ n,m }^{ \kappa +1}$ 
to $\Delta_{ \mathcal{ E }'}$ defined
by
\def\theequation{3.7}\begin{equation}
\pi_{\kappa,p}'
\big(
{x'}^i,{y'}^j,{y'}_{i_1}^j,\dots,{y'}_{i_1,\dots,i_{\kappa+1}}^j\big)
:=
\Big(
{x'}^i,{y'}^j,{y'}_{\beta(q)}^{j(q)}
\Big).
\end{equation}
Let $\varphi^{ (\kappa +1 )}$ 
be the $(\kappa +1)$-th prolongation of $\varphi$
(Section~1(II)).

\def\thelemma{3.8}\begin{lemma}
{\rm (\cite{ ol1986, bk1989, ol1995})}
The following three conditions are equivalent{\rm :}

\begin{itemize}

\smallskip\item[{\bf (1)}]
$\varphi$ transforms \thetag{ $\mathcal{ E}$} to \thetag{ $\mathcal{
E}' $}{\rm ;}

\smallskip\item[{\bf (2)}]
its $(\kappa +1)$-th 
prolongation $\varphi^{ (\kappa +1 )} : \mathcal{ J}_{ n,
m}^{ \kappa +1} 
\to {\mathcal{ J}'}_{ n,m }^{ \kappa +1}$ 
maps $\Delta_\mathcal{ E}$ to $\Delta_{ \mathcal{ E}'}${\rm ;}

\smallskip\item[{\bf (3)}]
$\varphi^{ (\kappa +1)} : \mathcal{ J}_{ n, m }^{
\kappa +1} \to {\mathcal{
J}'}_{ n,m}^{ \kappa +1}$ 
maps $\Delta_\mathcal{ E}$ to $\Delta_{\mathcal{
E}'}$ and the associated map 
\def\theequation{3.9}\begin{equation}
\Phi_{\mathcal{E},\mathcal{ E}'} 
:=
\pi_{\kappa,p}'
\circ
\big(\varphi^{(\kappa +1)}\big\vert_{
\Delta_\mathcal{ E}} \big)
\end{equation} 
sends every leaf of ${\sf F}_{ \! \Delta_\mathcal{ E}}$ to some leaf
of ${\sf F}_{ \! \Delta_{ \mathcal{ E}' }}$.

\end{itemize}\smallskip
\end{lemma}

\medskip\noindent
{\bf Equivalence problem 3.10.} {\it Find an algorithm to decide whether
two given \thetag{ $\mathcal{ E}$} and \thetag{ $\mathcal{ E}'$} are
equivalent.}

\medskip
\'Elie Cartan's widely applicable method (not reviewed here; \cite{
ca1937, ste1983, g1989, hk1989, fe1995, ol1995}) provides an answer ``in
principle'' to this question by reducing to an $\{ e \}$-structure an
initial G-structure associated to \thetag{ $\mathcal{ E }$}. Due to
the incredible size-length-complexity of the underlying computations,
this approach almost never abutes: it is forced to incompleteness. But
in fact, the main question is to classify.

\medskip\noindent
{\bf Classification problem 3.11.} {\it Classify systems \thetag{
$\mathcal{ E}$}, namely provide complete lists of all possible such
equations written in simplified ``normal'', easily recognizable forms.}

\medskip
Both problems are deeply linked to the classification of Lie algebras
of local vector fields. For $n=1$, $m=1$ and $\kappa = 1$, namely
\thetag{ $\mathcal{ E }_1$}: $y_{ xx} = F (x, y, y_x)$, Lie and Tresse
solved the two problems\footnote{
The author knows no complete confirmation 
of the Lie-Tresse classification by means
of \'E.~Cartan's method of equivalence.}. 
Table~7 of~\cite{ ol1986}, below reproduced,
describes the results.

\medskip
\begin{center}
\begin{tabular}[t]{|| l || l | l | l ||}
\hline
\hline
& \ 
{\rm Symmetry group} \ & \
{\rm Dimension} \ & \ {\rm Invariant equation} \ \\
\hline
\hline
{\bf (1)} 
& 
& 
0
& 
$y_{xx}=F(x,y,y_x)$
\\
\hline
{\bf (2)} 
& 
$\partial_y$
& 
1
& 
$y_{xx}=F(x,y_x)$
\\
\hline
{\bf (3)} 
& 
$\partial_x$, 
$\partial_y$
& 
2
& 
$y_{xx}=F(y_x)$
\\
\hline
{\bf (4)} 
& 
$\partial_x$, $e^x\partial_y$
& 
2
& 
$y_{xx}-y_x=F(y_x-y)$
\\
\hline
{\bf (5)} 
& 
$\partial_x$, $\partial_x-y\partial_y$, 
$x^2\partial_x-2xy\partial_y$ 
& 
3
& 
$y_{xx}=\frac{3y_x^2}{2y}+cy^3$
\\
\hline
{\bf (6)} 
& 
$\partial_x$, $x\partial_x-y\partial_y$, 
& 
3
& 
$y_{xx}=6yy_x-4y^3+$
\\
&
\ \ \ \ \ \ 
$x^2\partial_x-(2xy+1)\partial_y$
& 
& 
\ \ \ \ \ 
$+c(y_x-y^2)^{3/2}$
\\
\hline
{\bf (7)} 
& 
$\partial_x$, $\partial_y$, 
$x\partial_x+\alpha y\partial_y$,
& 
3
& 
$y_{xx}=c(y_x)^{\frac{\alpha-2}{\alpha-1}}$
\\
& 
\ \ \ \ \ 
$\alpha\neq 0,\frac{1}{2}, 1, 2$
& 
&
\\
\hline
{\bf (8)} 
& 
$\partial_x$, $\partial_y$, $x\partial_x+
(x+y)\partial_y$
& 
3
& 
$y_{xx}=ce^{-y_x}$
\\
\hline
{\bf (9)}
& 
$\partial_x$, $\partial_y$, 
$y\partial_x$, $x\partial_y$, $y\partial_y$,
& 
8
& 
$y_{xx}=0$
\\
& 
$x^2\partial_x+xy\partial_y$, 
$xy\partial_x+y^2\partial_y$
& 
& 
\\
\hline
\hline
\end{tabular}

\medskip

{\bf Table~1.}

\end{center}

However, the author knows no modern reference offering a complete
proof of this classification, with precise insight on the assumptions
(some normal forms hold true only at a generic point). In addition,
the above Lie-Tresse list is still slightly incomplete in the sense
that it does not precise which are the conditions satisfied by $F$
(Table~7 in~\cite{ ol1986}) insuring in the first four lines that
$\mathfrak{ SYM} (\mathcal{ E}_1)$ is indeed of small dimension $0$,
$1$ or $2$.

\def\theopenquestion{3.12}\begin{openquestion} 
Specify some precise nondegeneracy conditions upon $F$ in the first
four lines of Table~1.
\end{openquestion}

\section*{\S4.~Punctual and infinitesimal Lie symmetries}

\subsection*{4.1.~Lie symmetries of \thetag{ $\mathcal{ E}$}} 
Let $\varphi = (\phi, \psi)$ be a diffeomorphism of $\K_x^n \times
\K_y^n$ as in~\thetag{ 1.7}(II). 

\def\thedefinition{4.2}\begin{definition}{\rm
(\cite{ ol1986, ol1995, bk1989}) $\varphi$ is a {\sl {\rm (}local{\rm
)} Lie symmetry} of~\thetag{ $\mathcal{E }$} if it transforms the
graph of every solution of~\thetag{ $\mathcal{E }$} into the graph of
another solution.
}\end{definition}

To explain, we must pass to jet spaces. Denote the components of the
$(\kappa +1)$-th prolongation $\varphi^{ (\kappa +1 )} 
:\ \mathcal{ J }_{ 
n,m }^{ \kappa +1} \to \mathcal{ J}_{ n, m }^{ \kappa +1}$ by
\def\theequation{4.3}\begin{equation}
\varphi^{(\kappa+1)}
=
\big(
\phi^{i_1},
\psi^{j_1},
\Phi_{i_1}^j,
\Phi_{i_1,i_2}^j,
\dots\dots,
\Phi_{i_1,i_2,\dots,i_{\kappa+1}}^j
\big).
\end{equation}
The restriction $\varphi^{ (\kappa +1 )} \big \vert_{ \Delta_\mathcal{ E
}}$ is obtained by replacing each jet variable $y_\alpha^j$ by
$F_\alpha^j$, whenever $(j, \alpha) \neq (j,0)$ and $\neq (j (q),
\beta (q))$, and wherever it appears\footnote{ Remind from
Section~1(II) that we have not (open problem) provided a complete
explicit expression of $\Phi_{ i_1, \dots, i_\lambda }^j$ for general
$n \geqslant 1$, $m \geqslant 1$ and $\lambda \geqslant 1$.} in the
$\Phi_{ i_1, \dots, i_\lambda }^j$.

Let $\pi_{ \kappa, p}$ denote the projection from $\mathcal{ J}_{n,
m}^{ \kappa +1}$ 
to $\Delta_\mathcal{ E } \simeq \K^{m+ n +p } $ defined by
\def\theequation{4.4}\begin{equation}
\pi_{\kappa,p}
\big(
x^i,y^j,y_{i_1}^j,\dots,y_{i_1,\dots,i_{\kappa+1}}^j\big)
:=
\Big(
x^i,y^j,y_{\beta(q)}^{j(q)}
\Big),
\end{equation}
and introduce the map
\def\theequation{4.5}\begin{equation}
\varphi_{\Delta_\mathcal{E}}
:=
\pi_{\kappa,p}
\circ
\big(
\varphi^{(\kappa+1)}\big\vert_{\Delta_{\mathcal{E}}}
\big)
\equiv
\Big(
\varphi(x^i,y^j), 
\Phi_{\beta(q)}^{j(q)}
\big(
x^i,y^j,y_{\beta(q_1)}^{j(q_1)}
\big)
\Big).
\end{equation}

\def\thelemma{4.6}\begin{lemma}
{\rm (\cite{ ol1986, ol1995, bk1989}, [$*$])}
The following three conditions are equivalent{\rm :}

\begin{itemize}

\smallskip\item[{\bf (1)}]
the diffeomorphism $\varphi$ is a Lie symmetry of \thetag{
$\mathcal{E}$ }{\rm ;}

\smallskip\item[{\bf (2)}]
$\varphi^{ (\kappa+1 )} \big\vert_{ \Delta_\mathcal{ E}}$ sends
$\Delta_\mathcal{ E}$ to $\Delta_\mathcal{ E}${\rm ;}

\smallskip\item[{\bf (3)}]
$\varphi^{ (\kappa+1 )} \big\vert_{ \Delta_\mathcal{ E }}$ sends
$\Delta_\mathcal{ E }$ to $\Delta_\mathcal{ E }$ {\rm and} $\varphi_{
\Delta_\mathcal{ E}} = \pi_{ \kappa, p} \big( \varphi^{ (\kappa +1 )}
\big\vert_{ \Delta_\mathcal{ E }} \big)$ is a symmetry of the
foliation ${\sf F}_{ \! \Delta_\mathcal{ E }}$, namely it sends every
leaf to some other leaf.
\end{itemize}

Then the set of Lie symmetries of \thetag{
$\mathcal{E}$ } constitutes a local Lie (pseudo)group.
\end{lemma}

\subsection*{ 4.7.~Infinitesimal Lie symmetries of \thetag{ 
$\mathcal{ E}$}} Let 
\def\theequation{4.8}\begin{equation}
\mathcal{ L} 
= 
\sum_{i=1}^n\mathcal{X}^i(x,y)\,
\frac{\partial}{\partial x^i}
+
\sum_{j=1}^m\,\mathcal{Y}^j(x,y)\,
\frac{\partial}{\partial y^j},
\end{equation}
be a (local) vector field on $\K^{ n+m }$ having analytic
coefficients. Denote its flow by $\varphi_t (x,y) :=\exp(t \,
\mathcal{ L}) (x,y)$, $t\in \K$. As in Section~1(II), by
differentiating the prolongation $(\varphi_t)^{ ( \kappa +1 ) }$ with
respect to $t$ at $t=0$, we get the prolonged vector field $\mathcal{
L}^{ (\kappa +1 )}$ on $\mathcal{ J}_{ n, m}^{ \kappa +1}$, having the
general form (Part~II):
\def\theequation{4.9}\begin{equation}
\mathcal{L}^{(\kappa+1)}
=
\mathcal{L}
+
\sum_{j=1}^m\,\sum_{i_1=1}^n\,{\bf Y}_{i_1}^j\,
\frac{\partial}{\partial y_{i_1}^j}
+\cdots+
\sum_{j=1}^m\,\sum_{i_1,\dots,i_{\kappa+1}=1}^n\,
{\bf Y}_{i_1,\dots,i_{\kappa+1}}^j\,
\frac{\partial}{\partial y_{i_1,\dots,i_{\kappa+1}}^j},
\end{equation}
with known explicit expressions for the ${\bf Y}_{ i_1, \dots,
i_\lambda}^j$.

\def\thedefinition{4.10}\begin{definition}{\rm
$\mathcal{ L}$ is an {\sl infinitesimal symmetry} of \thetag{
$\mathcal{ E }$} if for every small $t$, its time-$t$ flow map
$\varphi_t$ is a Lie symmetry of \thetag{ $\mathcal{ E }$ }.
}\end{definition}

The restriction $\mathcal{ L}^{ (\kappa +1)} \big \vert_{
\Delta_\mathcal{ E}}$ is obtained by replacing every $y_\alpha^j$ by
$F_\alpha^j$ in all coefficients ${\bf Y}_{ i_1 }^j, \dots, {\bf Y}_{
i_1, \dots, i_{\kappa +1} 
}^j$. Then the coefficients become functions of
$\big (x^{ i_1}, y^{ j_1}, y_{ \beta (q_1)}^{ j(q_1) } \big)$ only.

\def\thelemma{4.11}\begin{lemma}
{\rm (\cite{ ol1986, ol1995, bk1989}, [$*$])}
The following three conditions are equivalent{\rm :}

\begin{itemize}

\smallskip\item[{\bf (1)}]
the vector field $\mathcal{ L}$ is an infinitesimal Lie symmetry of 
\thetag{ $\mathcal{ E }$}{\rm ;}

\smallskip\item[{\bf (2)}]
its $(\kappa+1)$-th 
prolongation $\mathcal{ L}^{ (\kappa +1 )}$ is tangent to
the skeleton $\Delta_\mathcal{ E}${\rm ;}

\smallskip\item[{\bf (3)}]
$\mathcal{ L}^{ (\kappa +1 )}$ is tangent to $\Delta_\mathcal{ E}$ and
the push-forward
\def\theequation{4.12}\begin{equation}
\mathcal{ L}_{\Delta_\mathcal{E}}
:=
(\pi_{\kappa,p})_* 
\big(
\mathcal{ L}^{(\kappa+1)}
\big\vert_{ \Delta_\mathcal{E}}
\big)
\end{equation}
is an infinitesimal symmetry of the foliation ${\sf F}_{ \!
\Delta_\mathcal{ E}}$, namely for every $i = 1, \dots,
n$, the Lie bracket $\big[ \mathcal{ L }_{ \Delta_\mathcal{ E}}, 
\, {\sf D}_i \big]$ is a linear combination of $\{ {\sf D}_1, \ldots, {\sf D}_n \}$.
\end{itemize}
\end{lemma}

According to~\cite{ ol1986, bk1989, ol1995}, the set of infinitesimal
Lie symmetries constitutes a Lie algebra, with the property $\big[ \mathcal{ L}^{ (\kappa +1)}, \,
{\mathcal{ L}'}^{ (\kappa +1)} \big] = \big[ \mathcal{ L}, \mathcal{
L}' \big]^{ (\kappa +1)}$. We summarize by a diagram.

\bigskip
\begin{center}
\begin{picture}(0,0)%
\includegraphics{lifts.pstex}%
\end{picture}%
\setlength{\unitlength}{4144sp}%
\begingroup\makeatletter\ifx\SetFigFont\undefined
\def\x#1#2#3#4#5#6#7\relax{\def\x{#1#2#3#4#5#6}}%
\expandafter\x\fmtname xxxxxx\relax \def\y{splain}%
\ifx\x\y   
\gdef\SetFigFont#1#2#3{%
  \ifnum #1<17\tiny\else \ifnum #1<20\small\else
  \ifnum #1<24\normalsize\else \ifnum #1<29\large\else
  \ifnum #1<34\Large\else \ifnum #1<41\LARGE\else
     \huge\fi\fi\fi\fi\fi\fi
  \csname #3\endcsname}%
\else
\gdef\SetFigFont#1#2#3{\begingroup
  \count@#1\relax \ifnum 25<\count@\count@25\fi
  \def\x{\endgroup\@setsize\SetFigFont{#2pt}}%
  \expandafter\x
    \csname \romannumeral\the\count@ pt\expandafter\endcsname
    \csname @\romannumeral\the\count@ pt\endcsname
  \csname #3\endcsname}%
\fi
\fi\endgroup
\begin{picture}(3856,1721)(745,-1496)
\put(2894,-657){\makebox(0,0)[lb]{\smash{\SetFigFont{9}{10.8}{rm}{\color[rgb]{0,0,0}$\Delta_\mathcal{E}$}%
}}}
\put(3163,-894){\makebox(0,0)[lb]{\smash{\SetFigFont{9}{10.8}{rm}{\color[rgb]{0,0,0}$\pi_\kappa$}%
}}}
\put(1029,-635){\makebox(0,0)[lb]{\smash{\SetFigFont{9}{10.8}{rm}{\color[rgb]{0,0,0}$\Delta_\mathcal{E}$}%
}}}
\put(1187,-295){\makebox(0,0)[lb]{\smash{\SetFigFont{9}{10.8}{rm}{\color[rgb]{0,0,0}$\pi_{\kappa,p}$}%
}}}
\put(2662,-283){\makebox(0,0)[lb]{\smash{\SetFigFont{9}{10.8}{rm}{\color[rgb]{0,0,0}$\pi_{\kappa,p}$}%
}}}
\put(1178,-1009){\makebox(0,0)[lb]{\smash{\SetFigFont{9}{10.8}{rm}{\color[rgb]{0,0,0}$\pi_p$}%
}}}
\put(2742,-1016){\makebox(0,0)[lb]{\smash{\SetFigFont{9}{10.8}{rm}{\color[rgb]{0,0,0}$\pi_p$}%
}}}
\put(745,-883){\makebox(0,0)[lb]{\smash{\SetFigFont{9}{10.8}{rm}{\color[rgb]{0,0,0}$\pi_\kappa$}%
}}}
\put(4486,-666){\makebox(0,0)[lb]{\smash{\SetFigFont{9}{10.8}{rm}{\color[rgb]{0,0,0}$\mathcal{L}_{\Delta_\mathcal{E}}$}%
}}}
\put(2812,-1454){\makebox(0,0)[lb]{\smash{\SetFigFont{9}{10.8}{rm}{\color[rgb]{0,0,0}$\K_x^n\times\K_y^m$}%
}}}
\put(4492,-1450){\makebox(0,0)[lb]{\smash{\SetFigFont{9}{10.8}{rm}{\color[rgb]{0,0,0}$\mathcal{L}$}%
}}}
\put(766,-1456){\makebox(0,0)[lb]{\smash{\SetFigFont{9}{10.8}{rm}{\color[rgb]{0,0,0}$\K_x^n\times\K_y^m$}%
}}}
\put(1882,-541){\makebox(0,0)[lb]{\smash{\SetFigFont{9}{10.8}{rm}{\color[rgb]{0,0,0}$\varphi_{\Delta_\mathcal{E}}$}%
}}}
\put(1913,-1313){\makebox(0,0)[lb]{\smash{\SetFigFont{9}{10.8}{rm}{\color[rgb]{0,0,0}$\varphi$}%
}}}
\put(1955,117){\makebox(0,0)[lb]{\smash{\SetFigFont{9}{10.8}{rm}{\color[rgb]{0,0,0}$\varphi^{(\kappa+1)}$}%
}}}
\put(908, 15){\makebox(0,0)[lb]{\smash{\SetFigFont{9}{10.8}{rm}{\color[rgb]{0,0,0}$\mathcal{J}_{n,m}^{\kappa+1}$}%
}}}
\put(2889, 24){\makebox(0,0)[lb]{\smash{\SetFigFont{9}{10.8}{rm}{\color[rgb]{0,0,0}$\mathcal{J}_{n,m}^{\kappa+1}$}%
}}}
\put(4509, 28){\makebox(0,0)[lb]{\smash{\SetFigFont{9}{10.8}{rm}{\color[rgb]{0,0,0}$\mathcal{L}^{(\kappa+1)}$}%
}}}
\end{picture}

\end{center}

\subsection*{4.13.~Sophus Lie's algorithm} We describe the general
process. Its complexity will be exemplified in
Section~5 (to be read simultaneously).

The tangency of $\mathcal{ L}^{ ( \kappa +1)}$ 
to $\Delta_\mathcal{ E}$ is
expressed by applying $\mathcal{ L}^{ (\kappa +1 )}$ to the equations $0
= - y_\alpha^j + F_\alpha^j$, which yields:
\def\theequation{4.14}\begin{equation}
0
=
-
{\bf Y}_\alpha^j
+
\sum_{i=1}^n\,\mathcal{X}^i\,
\frac{\partial F_\alpha^j}{\partial x^i}
+
\sum_{l=1}^n\,\mathcal{Y}^l\,
\frac{\partial F_\alpha^j}{\partial y^l}
+
\sum_{q=1}^p\,{\bf Y}_{\beta(q)}^{j(q)}\,
\frac{\partial F_\alpha^j}{\partial y_{\beta(q)}^{j(q)}},
\end{equation}
for $(j, \alpha) \neq (j, 0)$ and $\neq (j(q), \beta (q))$.
Restricting a coefficient ${\bf Y}_{i_1, \dots, i_\lambda}^j$ to
$\Delta_\mathcal{ E}$, namely replacing everywhere in it each
$y_\alpha^j$ by $F_\alpha^j$, provides a specialized coefficient
\def\theequation{4.15}\begin{equation}
\widehat{\bf Y}_{i_1,\dots,i_\lambda}^j
=
\widehat{\bf Y}_{i_1,\dots,i_\lambda}^j
\big(
x^{i_1},y^{j_1},
y_{\beta(q_1)}^{j(q_1)},
J_{x,y}^\lambda\mathcal{X}^{i_1},
J_{x,y}^\lambda\mathcal{Y}^{j_1}
\big),
\end{equation}
that depends {\bf linearly}\, on the $\lambda$-th jet of the
coefficients of $\mathcal{ L}$, as confirmed by an inspection of
Part~II's formulas. Here, we use the jet notation $J_{ x, y}^\lambda
\mathcal{ Z} := \big( \partial_x^{ \alpha_1} \partial_y^{ \beta_1}
\mathcal{ Z} \big)_{ \vert \alpha_1 \vert + \vert \beta_1 \vert
\leqslant \lambda}$. We thus get equations
\def\theequation{4.16}\begin{equation}
0
\equiv
-
\widehat{\bf Y}_\alpha^j
+
\sum_{i=1}^n\,\mathcal{X}^i\,
\frac{\partial F_\alpha^j}{\partial x^i}
+
\sum_{l=1}^n\,
\mathcal{Y}^l\,
\frac{\partial F_\alpha^j}{\partial y^l}
+
\sum_{q=1}^p\,
\widehat{\bf Y}_{\beta(q)}^{j(q)}\,
\frac{\partial F_\alpha^j}{\partial y_{\beta(q)}^{j(q)}},
\end{equation}
involving only the variables $\big( x^{ i_1}, y^{ j_1}, 
y_{ \beta( q_1)}^{ j(q_1)} \big)$.

Next, we develope every such equation with respect to the powers of
$y_{ \beta( q_1)}^{ j(q_1)}$:
\def\theequation{4.17}\begin{equation}
0
\equiv
\sum_{\mu_1,\dots,\mu_p\geqslant 0} \, 
(y_{\beta(1)}^{j(1)})^{\mu_1}\cdots 
(y_{\beta(p)}^{j(p)})^{\mu_p} \, 
\Psi_{\alpha,\mu_1,\dots,\mu_p}^j
\big(
x^{i_1},y^{j_1},
J_{x,y}^{\kappa+1}\mathcal{X}^{i_1},
J_{x,y}^{\kappa+1}\mathcal{Y}^{j_1}
\big).
\end{equation}
The $\Psi_{ \alpha, \mu_1, \dots, \mu_p}^j$ are {\bf linear} with
respect to $\big( J_{ x, y }^{\kappa+1} \mathcal{ X }^{ i_1}, J_{ x, y
}^{\kappa+1} \mathcal{ Y }^{j_1} \big)$, with certain coefficients
analytic with respect to $(x, y)$, which depend intrinsically (but in
a complex manner) on the right hand sides $F_\alpha^j$.

\def\theproposition{4.18}\begin{proposition}
The vector field $\mathcal{ L}$ is an infinitesimal Lie symmetry of
\thetag{ $\mathcal{ E}$ } {\rm if and only if} its coefficients
$\mathcal{ X}^{ i_1}$, $\mathcal{ Y}^{ j_1}$ satisfy the linear {\sc pde}
system{\rm :}
\def\theequation{4.19}\begin{equation}
0
=
\Psi_{\alpha,\mu_1,\dots,\mu_p}^j
\big(
x^{i_1},y^{j_1},
J_{x,y}^{\kappa+1}\mathcal{X}^{i_1},
J_{x,y}^{\kappa+1}\mathcal{Y}^{j_1}
\big)
\end{equation}
for all $(j, \alpha) \neq (j, 0)$ and $\neq (j(q), \beta (q ))$ and for
all $(\mu_1, \dots, \mu_p) \in \N^p$.
\end{proposition}

In all known instances, a finite number of these equations
suffices. 

\def\theexample{4.20}\begin{example}{\rm
With $n = m = \kappa = 1$, a second prolongation $\mathcal{ L}^{ (2)} =
\mathcal{ X}\, \frac{ \partial }{ \partial x} + \mathcal{ Y} \, \frac{
\partial }{ \partial y} + {\bf Y}_1 \, \frac{ \partial }{ \partial
y_1} + {\bf Y}_2 \, \frac{ \partial }{ \partial y_2}$ is tangent to
the skeleton $0 = - y_2 + F( x, y, y_1)$ of \thetag{ $\mathcal{ E}_1$
} if and only if $0 = - {\bf Y}_2 + \mathcal{ X}\, F_x + \mathcal{ Y}
\, F_y + {\bf Y}_1 \, F_{ y_1}$, or, developing:
\def\theequation{4.21}\begin{equation}
\left\{
\aligned
0
&
=
-
\mathcal{Y}_{xx}
+
\big[
-2\,\mathcal{Y}_{xy}
+
\mathcal{X}_{xx}
\big]\,y_1
+
\big[
-\mathcal{Y}_{yy}
+
2\,\mathcal{X}_{xy}
\big]\,(y_1)^2
+
\big[
\mathcal{X}_{yy}
\big]\,(y_1)^3
+
\\
&\ \ \ \ \
+
\big[
-\mathcal{Y}_y
+
2\,\mathcal{X}_x
\big]\,F
+
\big[
3\,\mathcal{X}_y
\big]\,y_1\,F
+
\big[\mathcal{X}\big]\,F_x
+
\big[\mathcal{Y}\big]\,F_y
+
\\
&\ \ \ \ \
+
\big[
\mathcal{Y}_x
\big]\,F_{y_1}
+
\big[
\mathcal{Y}_y
-
\mathcal{X}_x
\big]\,y_1\,F_{y_1}
+
\big[
-\mathcal{X}_y
\big]\,(y_1)^2\,F_{y_1}.
\endaligned\right.
\end{equation}
Developing $F = \sum_{ k\geqslant 0}\, 
(y_1)^k \, F_k (x,y)$, we may obtain equations~\thetag{ 4.19}.
}\end{example}

\section*{\S5.~Examples}

\subsection*{ 5.1.~Second order ordinary differential equation}
Pursuing the study of \thetag{ $\mathcal{ E}_1$}, according to
Section~7 below, we may assume that $F = {\rm O} (y_x)$, or
equivalently $F (x, y, 0) \equiv 0$.

\def\theconvention{5.2}\begin{convention}{\rm
The letters ${\sf R}$ will denote various functions of $(x, y, y_1)$,
changing with the context. Similarly, ${\sf r} = {\sf r} (x, y)$,
excluding the pure jet variable $y_1$. Hence, symbolically:
\def\theequation{5.3}\begin{equation}
{\sf R}
=
{\sf r}
+
y_1\,{\sf r}
+
(y_1)^2\,{\sf r}
+
(y_1)^3\,{\sf r}
+
\cdots.
\end{equation} 
}\end{convention}

So the skeleton is 
\def\theequation{5.4}\begin{equation}
y_2
=
F(x,y,y_1)
= 
y_1\,{\sf R}
=
y_1\,{\sf r}
+
(y_1)^2\,{\sf r}
+
(y_1)^3\,{\sf r}
+
\cdots.
\end{equation}
Applying $\mathcal{ L }^{ (2)}$, see~\thetag{ 2.3}(II) for its
expression, we get:
\def\theequation{5.5}\begin{equation}
0 
= 
-
{\bf Y}_2
+ 
\mathcal{X}\,F_x 
+
\mathcal{Y}\,F_y 
+
{\sf Y}_1\,F_{y_1}.
\end{equation}
Observe that $F_x = (y_1\, {\sf R})_x = {\sf r} \, y_1 + {\sf r} \,
(y_1)^2 + \cdots $ and similarly for $F_y$, but that $(y_1 \, {\sf R}
)_{ y_1} = {\sf r} + {\sf r} \, y_1 + {\sf r} \, (y_1)^2 + \cdots$.
Inserting above ${\sf Y}_1$, ${\sf Y}_2$ given by ~\thetag{ 2.6}(II),
replacing $y_2$ by $y_1 \, {\sf R}$ and computing ${\rm mod}\, (y_1
)^4$, we get:
\def\theequation{5.6}\begin{equation}
\aligned
0
\equiv
&
-
\mathcal{Y}_{xx}
+
\big[
-
2\,\mathcal{Y}_{xy}
+
\mathcal{X}_{xx}
\big]\,y_1
+
\big[
-
\mathcal{Y}_{yy}
+
2\,\mathcal{X}_{xy}
\big]\,(y_1)^2
+
\big[
\mathcal{X}_{yy}
\big]\,(y_1)^3
+
\\
&\
+
\big[
-
\mathcal{Y}_y
+
2\,\mathcal{X}_x
\big]\,
\big(
y_1\,{\sf r}
+
(y_1)^2\,{\sf r}
+
(y_1)^3\,{\sf r}
\big)
+
\big[
3\,\mathcal{X}_y
\big]\,
\big(
(y_1)^2\,{\sf r}
+
(y_1)^3\,{\sf r}
\big)
+
\\
&\
+
\big[
\mathcal{X}
\big]\,
\big(
y_1\,{\sf r}
+
(y_1)^2\,{\sf r}
+
(y_1)^3\,{\sf r}
\big)
+
\big[
\mathcal{Y}
\big]\,
\big(
y_1\,{\sf r}
+
(y_1)^2\,{\sf r}
+
(y_1)^3\,{\sf r}
\big)
+
\\
&\
+
\big[
\mathcal{Y}_x
\big]\,
\big(
{\sf r}
+
y_1\,{\sf r}
+
(y_1)^2\,{\sf r}
+
(y_1)^3\,{\sf r}
\big)
+
\\
&\
+
\big[
\mathcal{Y}_y
-
\mathcal{X}_x
\big]\,
\big(
y_1\,{\sf r}
+
(y_1)^2\,{\sf r}
+
(y_1)^3\,{\sf r}
\big)
+
\big[
-
\mathcal{X}_y
\big]\,
\big(
(y_1)^2\,{\sf r}
+
(y_1)^3\,{\sf r}
\big).
\endaligned
\end{equation}
We gather the powers ${\rm cst.}$, $y_1$, 
$(y_1)^2$ and $(y_1)^3$, equating their coefficients to 
$0$:
\def\theequation{5.7}\begin{equation}
\aligned
0
&
=
-
\mathcal{Y}_{xx}
+
{\sf P}\big(
\mathcal{Y}_x
\big),
\\
0
&
=
-
2\,\mathcal{Y}_{xy}
+
\mathcal{X}_{xx}
+
{\sf P}\big(
\mathcal{Y}_y,
\mathcal{X}_x,
\mathcal{X},
\mathcal{Y},
\mathcal{Y}_x
\big),
\\
0
&
=
-
\mathcal{Y}_{yy}
+
2\,\mathcal{X}_{xy}
+
{\sf P}\big(
\mathcal{Y}_y,
\mathcal{X}_x,
\mathcal{X}_y,
\mathcal{X},
\mathcal{Y},
\mathcal{Y}_x
\big),
\\
0
&
=
\mathcal{X}_{yy}
+
{\sf P}\big(
\mathcal{Y}_y,
\mathcal{X}_x,
\mathcal{X}_y,
\mathcal{X},
\mathcal{Y},
\mathcal{Y}_x
\big)
\endaligned
\end{equation}

\def\theconvention{5.8}\begin{convention}{\rm
The letter ${\sf P}$ will denote various {\it linear} combinations of
some precise partial derivatives of $\mathcal{ X}$, $\mathcal{ Y}$
which have analytic coefficients in $(x, y)$.
}\end{convention}

By cross-differentiations and substitutions in the above system, all
third, fourth, fifth, {\it etc.} order derivatives of $\mathcal{ X},
\mathcal{ Y}$ may be expressed as ${\sf P} \big( \mathcal{ X},
\mathcal{ Y}, \mathcal{ X}_x, \mathcal{ X}_y, \mathcal{ Y}_x,
\mathcal{ Y}_y, \mathcal{ Y}_{ xy}, \mathcal{ Y}_{ yy} \big)$.

\def\theproposition{5.9}\begin{proposition}
An infinitesimal Lie symmetry $\mathcal{ X}\, \frac{ \partial}{
\partial x} + \mathcal{ Y}\, \frac{ \partial}{ \partial y}$ of
\thetag{ $\mathcal{ E }_1$} is uniquely determined by the eight
initial Taylor coefficients{\rm :}
\def\theequation{5.10}\begin{equation}
\mathcal{X}(0),\,
\mathcal{Y}(0),\,
\mathcal{X}_x(0),\,
\mathcal{X}_y(0),\,
\mathcal{Y}_x(0),\,
\mathcal{Y}_y(0),\,
\mathcal{Y}_{xy}(0),\,
\mathcal{Y}_{yy}(0).
\end{equation}
\end{proposition}

The bound $\dim \, \mathfrak{ SYM } (\mathcal{ E}_1) \leqslant 8$ is
attained with $F = 0$, whence all ${\sf P} = 0$ and
\def\theequation{5.11}\begin{equation}
\left\{
\aligned
A
:= 
&\ 
\partial_y, \ \ \ \ \ \ & 
E:= 
&\ 
y\,\partial_y, 
\\
B
:=
&\ 
\partial_x, \ \ \ \ \ \ & 
F
:= 
&\ 
y\,\partial_x, 
\\
C
:=
&\ 
x\,\partial_y, \ \ \ \ \ \ & 
G
:= 
&\ 
xx\,\partial_x
+
xy\,\partial_y, 
\\
D
:=
&\ 
x\,\partial_x, \ \ \ \ \ \ & 
H
:= 
&\ 
xy\,\partial_x
+
yy\,\partial_y. 
\endaligned\right.
\end{equation}
are infinitesimal generators of the group $\text{ \rm PGL}_3( \K) =
\text{\rm Aut} (P_2 (\K ))$ of projective transformations
\def\theequation{5.12}\begin{equation}
(x,y)
\mapsto 
\left(
\frac{\alpha x+\beta y+\gamma}
{\lambda x+\mu y+\nu}, \
\frac{\delta x+\eta y+\epsilon} 
{\lambda x+\mu y+\nu},
\right)
\end{equation}
stabilizing the collections of all affine lines of $\K^2$, namely the
solutions of the {\sl model} equation $y_{ xx} = 0$. The {\sl model
Lie algebra} $\mathfrak{ pgl }_3 ( \K) \simeq \mathfrak{ sl}_3 (\K)$
is simple.

\def\thetheorem{5.13}\begin{theorem}
The bound $\dim \, \mathfrak{ SYM } (\mathcal{ E}_1) \leqslant 8$ is
attained {\rm if and only if} \thetag{ $\mathcal{ E}_1$} is
equivalent, through a diffeomorphism $(x, y) \mapsto (X, Y)$, to $Y_{
XX} = 0$.
\end{theorem}

\proof
The statement is well known (\cite{ lie1883, el1890, tr1896, se1931,
ca1932a, ol1986, hk1989, ib1992, ol1995, sh1997, su2001, n2003,
me2004}). We provide a (new?) proof which has the advantage to enjoy
direct generalizations to all {\sc pde} systems whose model Lie algebras are
semisimple, for instance \thetag{ $\mathcal{ E }_2$}, \thetag{
$\mathcal{ E }_3$} and \thetag{ $\mathcal{ E }_5$}.

The Lie brackets between the eight generators~\thetag{ 5.11}
are:

\medskip
\begin{center}
\begin{tabular}[t]{ l || l | l | l | l | l | l | l | l }
& $A$ & $B$ & $C$ & $D$ & $E$ & $F$ &
$G$ & $H$ 
\\
\hline 
\hline 
$A$ & 
$0$ & $0$ & $0$ & $0$ & $A$ & $B$ & $C$ & $D+2E$ 
\\
\hline
$B$ & 
$0$ & $0$ & $A$ & $B$ & $0$ & $0$ & $E+2D$ & $F$ 
\\
\hline
$C$ & 
$0$ & $-A$ & $0$ & $-C$ & $C$ & $D-E$ & $0$ & $G$ 
\\
\hline
$D$ & 
$0$ & $-B$ & $C$ & $0$ & $0$ & $-F$ & $G$ & $0$ 
\\ 
\hline
$E$ & 
$-A$ & $0$ & $-C$ & $0$ & $0$ & $F$ & $0$ & $H$ 
\\
\hline
$F$ & 
$-B$ & $0$ & $-D+E$ & $F$ & $-F$ & $0$ & $H$ & $0$ 
\\
\hline
$G$ & 
$-C$ & $-E-2D$ & $0$ & $-G$ & $0$ & $H$ & $0$ & $0$ 
\\
\hline
$H$ & 
$-D-2E$ & $-F$ & $-G$ & $0$ & $-H$ & $0$ & $0$ & $0$
\end{tabular}

\medskip
{\bf Table~2.}

\end{center}

Assuming that $\dim \mathfrak{ SYM} (\mathcal{ E}_1) = 8$, taking
account of~\thetag{ 5.7}, after making some linear combinations, there
must exist eight generators of the form
\def\theequation{5.14}\begin{equation}
\left\{
\aligned
A'
:= 
&\ 
\partial_y
+
{\rm O}(1), \ \ \ \ \ \ & 
E':= 
&\ 
y\,\partial_y
+
{\rm O}(2), 
\\
B'
:=
&\ 
\partial_x
+
{\rm O}(1), \ \ \ \ \ \ & 
F'
:= 
&\ 
y\,\partial_x
+
{\rm O}(2), 
\\
C'
:=
&\ 
x\,\partial_y
+
{\rm O}(2), \ \ \ \ \ \ & 
G'
:= 
&\ 
xx\,\partial_x
+
xy\,\partial_y
+
{\rm O}(3), 
\\
D'
:=
&\ 
x\,\partial_x
+
{\rm O}(2), \ \ \ \ \ \ & 
H'
:= 
&\ 
xy\,\partial_x
+
yy\,\partial_y
+
{\rm O}(3).
\endaligned\right.
\end{equation}
To insure that the Lie brackets between these vector fields are small
perturbations of the model ones, we can in advance replace $( x, y)$
by $(\varepsilon x, \varepsilon y)$, so that $y_{ xx} = \varepsilon \,
F \big( \varepsilon x, \varepsilon y, y_x \big)$ is an ${\rm O}
(\varepsilon )$, hence all the remainders ${\rm O} (1)$, ${\rm O} (2)$
and ${\rm O} (3)$ above are also ${\rm O} ( \varepsilon )$. It
follows that the structure constants for $A', \dots, H'$ are
$\varepsilon$-close to those of Table~2.

\def\thetheorem{5.15}\begin{theorem}
{\rm (\cite{ ov1994})}
Every semisimple Lie algebra over $\R$ or $\C$ is rigid{\rm :}
small deformations of the structure constants just give isomorphic
Lie algebras.
\end{theorem}

Consequently, there exists a change of basis close to the identity
leading to new generators $A'', B'', \dots, G'', H''$ having exactly
the same structure constants as in Table~2. Then $A'' (0)$ and $B''
(0)$ are still linearly independent. Since $\big[ A'', B'' \big] = [
A, B] = 0$, there exist local coordinates $(X, Y)$ centered at $0$ in
which $A '' = \partial_X$ and $B'' = \partial_Y$. Since $\big[ A'',
C'' \big] = [ A, C] = 0$ and $\big[ B'', C'' \big] = [ B, C] = A$, it
follows that $C '' = X \partial_Y$. The tangency to $0 = - Y_2 + F (X,
Y, Y_1)$ (with $F (0) = 0$) of $\big( \partial_X \big)^{ (2 )} =
\partial_X$, of $\big( \partial_Y \big)^{ (2 )} = \partial_Y$ and of
$\big( X \partial_Y \big)^{ (2 )} = X \partial_Y + \partial_{ Y_1}$
yields $F = 0$.
\endproof

\def\theopenquestion{5.16}\begin{openquestion}
Does this proof generalize to $y_{ x^{ \kappa +1}} = F \big( x, y,
y_x, \dots, y_{ x^\kappa} \big)$~{\rm ?}
\end{openquestion}

\subsection*{ 5.17.~Complete system of second order}
We now summarize a generalization to \thetag{ $\mathcal{ E}_2$}.
According to Section~7 below, one may assume that the submanifold of
solutions is $y = b + \sum_{ i=1}^n \, a^i \big[ x^i + {\rm O} (\vert
x \vert^2) + {\rm O} (a) + {\rm O} (b) \big]$, whence $y_{x^{ i_1} x^{
i_2}} = F_{ i_1, i_2} \big( x^i, y, y_{ x^k} \big)$ with $F (x, y, 0)
\equiv 0$. Applying to the skeleton $0 = - y_{ i_1, i_2} + F_{ i_1,
i_2} \big( x^i, y, y_k\big)$ a second prolongation $\mathcal{ L}^{
(2)}$ having coefficients ${\bf Y}_{ i_1}$ given by~\thetag{3.9}(II)
and ${\bf Y}_{ i_1, i_2}$ given by~\thetag{ 3.20 }(II), we get
\def\theequation{5.18}\begin{equation}
0
=
-
{\sf Y}_{i_1,i_2}
+
\sum_{k=1}^n\,
\big[\mathcal{X}^k\big]
\frac{\partial F_{i_1,i_2}}{\partial x^k}
+
\big[
\mathcal{Y}
\big]
\frac{\partial F_{i_1,i_2}}{\partial y}
+
\sum_{k=1}^n\,
\big[{\bf Y}_k\big]
\frac{\partial F_{i_1,i_2}}{\partial y_k}.
\end{equation}
Replacing $y_{ i_1, i_2}$ everywhere by $F_{i_1, i_2} = y_1\, {\sf R}+
\cdots + y_n\, {\sf R}$, developping in powers of the pure jet
variables $y_l$ and picking the coefficients of ${\rm cst.}$,
of $y_k$, of
$(y_k)^2$ and of $(y_k)^3$, we get the linear system
\def\theequation{5.19}\begin{equation}
\left\{
\aligned
\mathcal{Y}_{x^{i_1}x^{i_2}}
&
=
{\sf P}\big(
\mathcal{Y}_{x^l}
\big)
\\
\delta_{i_1}^k\,\mathcal{Y}_{x^{i_2}y}
+
\delta_{i_2}^k\,\mathcal{Y}_{x^{i_1}y}
-
\mathcal{X}_{x^{i_1}x^{i_2}}^k
&
=
{\sf P}\big(
\mathcal{Y}_y,
\mathcal{X}_{x^{l_1}}^{l_2},
\mathcal{X}^l,
\mathcal{Y},
\mathcal{Y}_{x^l}
\big)
\\
\delta_{i_1,i_2}^{k,k}\,\mathcal{Y}_{yy}
-
\delta_{i_1}^k\,\mathcal{X}_{x^{i_2}y}^k
-
\delta_{i_2}^k\,\mathcal{X}_{x^{i_1}y}^k
&
=
{\sf P}\big(
\mathcal{Y}_y,
\mathcal{X}_{x^{l_1}}^{l_2},
\mathcal{X}_y^l,
\mathcal{X}^l,
\mathcal{Y},
\mathcal{Y}_{x^l}
\big)
\\
\delta_{i_1,i_2}^{k,k}\,\mathcal{X}_{yy}^k
&
=
{\sf P}\big(
\mathcal{Y}_y,
\mathcal{X}_{x^{l_1}}^{l_2},
\mathcal{X}_y^l,
\mathcal{X}^l,
\mathcal{Y},
\mathcal{Y}_{x^l}
\big),
\endaligned\right.
\end{equation}
upon which obvious linear combinations yield a known
generalization of Proposition~5.9.

\def\theproposition{5.20}\begin{proposition}
{\rm (\cite{ su2001, gm2003a})} An infinitesimal Lie symmetry $\sum_{
k=1}^n \, \mathcal{ X}^k \, \frac{ \partial }{ \partial x^k} +
\mathcal{ Y} \, \frac{ \partial }{ \partial y}$ is uniquely determined
by the $n^2 + 4 n +3$ initial Taylor coefficients{\rm :}
\def\theequation{5.21}\begin{equation}
\mathcal{X}^l(0),\
\mathcal{Y}(0),\
\mathcal{X}_{x^{l_1}}^{l_2}(0),\
\mathcal{X}_y^l(0),\
\mathcal{Y}_{x^l}(0),\
\mathcal{Y}_y(0),\
\mathcal{Y}_{x^ly}(0),\
\mathcal{Y}_{yy}(0).
\end{equation}
\end{proposition}

The bound $\dim \mathfrak{ SYM} (\mathcal{ E}_2) \leqslant n^2 + 4n +
3$ is attained with $F_{ i_1, i_2} = 0$, whence all ${\sf P} = 0$ and
\def\theequation{5.22}\begin{equation}
\left\{
\aligned
A
:= 
&\ 
\partial_y, \ \ \ \ \ \ & 
E:= 
&\ 
y\,\partial_y, 
\\
B_i
:=
&\ 
\partial_{x^i}, \ \ \ \ \ \ & 
F_i
:= 
&\ 
y\,\partial_{x^i}, 
\\
C_i
:=
&\ 
x^i\,\partial_y, \ \ \ \ \ \ & 
G_i
:= 
&\ 
x^i\big(
x^1\,\partial_{x^1}+\cdots+
x^n\,\partial_{x^n}+
y\,\partial_y
\big)
+
xy\,\partial_y, 
\\
D_{i,k}
:=
&\ 
x^i\,\partial_{x^k}, \ \ \ \ \ \ & 
H
:= 
&\ 
y\big(
x^1\,\partial_{x^1}+\cdots+
x^n\,\partial_{x^n}+
y\,\partial_y
\big). 
\endaligned\right.
\end{equation}
are infinitesimal generators of the group $\text{ \rm PGL}_{ n+2} (
\K) = \text{\rm Aut} (P_{ n+1} (\K ))$ of projective transformations
\def\theequation{5.23}\begin{equation}
(x,y)
\mapsto 
\left(
\frac{\alpha_1 x^1+\cdots+\alpha_n x^n+\beta y+\gamma}
{\lambda_1 x^1+\cdots+\lambda_n x^n+\mu y+\nu},\ \
\frac{\delta_1 x^1+\cdots+\delta_n x^n+\eta y+\epsilon} 
{\lambda_1 x^1+\cdots+\lambda_n x^n+\mu y+\nu},
\right)
\end{equation}
stabilizing the collections of all affine planes of $\K^{ n+1}$,
namely the solutions of the {\sl model} equation $y_{ x^{ i_1} x^{
i_2}} = 0$. The {\sl model Lie algebra} $\mathfrak{ pgl }_{ n+2} ( \K)
\simeq \mathfrak{ sl}_{ n+2} (\K)$ is simple, hence rigid.

\def\thetheorem{5.24}\begin{theorem}
The bound $\dim \, \mathfrak{ SYM } (\mathcal{ E}_2) \leqslant n^2 +
4n + 3$ is attained {\rm if and only if} \thetag{ $\mathcal{ E}_2$} is
equivalent, through a diffeomorphism $(x^i, y) \mapsto (X^k, Y)$, to
$Y_{ X^{ k_1} X^{ k_2}} = 0$.
\end{theorem}

The proof, similar to that of Theorem~5.13, is skipped. 

The study of \thetag{ $\mathcal{ E }_3$} also leads to the model
algebra $\mathfrak{ pgl }_{ n+2} ( \K) \simeq \mathfrak{ sl}_{ n+2}
(\K)$ and an analog to Theorem~5.13 holds. Details are similar.

\section*{\S6.~Transfer of Lie symmetries to the parameter space}

\subsection*{6.1.~Stabilization of foliations}
As announced in \S2.38, we now transfer the theory of Lie symmetries to
submanifolds of solutions.

Restarting from \S4.1, let $\varphi$ a Lie symmetry of \thetag{
$\mathcal{ E }$}, namely $\varphi_{ \Delta_\mathcal{ E }}$ stabilizes
${\sf F}_{ \Delta_\mathcal{ E }}$. The diffeomorphism ${\sf A}$
defined by~\thetag{ 2.9} transforms ${\sf F}_{ \sf v}$ to ${\sf F}_{
\Delta_\mathcal{ E }}$. Conjugating, we get the self-transformation
${\sf A}^{-1} \circ \varphi_{ \Delta_\mathcal{ E }} \circ {\sf A}$ of
the $(x, a, b)$-space that must stabilize also the foliation ${\sf
F}_{ \sf v}$. Equivalently, it must have expression:
\def\theequation{6.2}\begin{equation}
\big[
{\sf A}^{-1} 
\circ
\varphi_{\Delta_\mathcal{E}}
\circ 
{\sf A}
\big]
(x,a,b)
=
\big(
\theta(x,a,b),
f(a,b), 
g(a,b)
\big)
\in\K^n\times \K^p\times\K^m,
\end{equation}
where, importantly, the last two components are independent of the
coordinate $x$, because the leaves of ${\sf F}_{ \sf v}$ are just $\{
a = {\rm cst.}, \ b = {\rm cst} \}$.

\def\thelemma{6.3}\begin{lemma}
To every Lie symmetry $\varphi$ of \thetag{ $\mathcal{ E }$ }, there
corresponds a transformation of the parameters
\def\theequation{6.4}\begin{equation}
(a,b)
\longmapsto 
\big(
f(a,b), 
g(a,b)
\big)
=:
h(a,b)
\end{equation} 
meaning that $\varphi$ transforms the local solution $y_{a, b} (x) :=
\Pi (x, a, b)$ to the local solution $y_{ h(a,b)}(x) = 
\Pi (x, h(a, b))$.
\end{lemma}

Unfortunately, the expression of ${\sf A}^{-1} \circ
\varphi_{ \Delta_\mathcal{ E }} \circ {\sf A}$ does not clearly show
that $f$ and $g$ are independent of $x$. Indeed, reminding the
expressions of ${\sf A}$ and of $\Phi$, we have:
\def\theequation{6.5}\begin{equation}
\varphi_{\Delta_\mathcal{E}}
\circ
{\sf A}
(x,a,b)
=
\Big(
\varphi(x,\Pi(x,a,b)),
\Phi_{\beta(q)}^{j(q)}
\big(
x^{i_1},\Pi^{j_1}(x,a,b),
\Pi_{x^{\beta(q_1)}}^{j(q_1)}(x,a,b)
\big)
\Big).
\end{equation}
To compose with ${\sf A}^{ -1}$ whose expression is given by~\thetag{
2.21}, it is useful to split $\varphi = (\phi, \psi) \in \K^n \times
\K^m$, so above we write
\def\theequation{6.6}\begin{equation}
\varphi(x,\Pi(x,a,b))
=
\big(
\phi(x,\Pi(x,a,b)),
\psi(x,\Pi(x,a,b))
\big),
\end{equation}
and finally, droping the arguments:
\def\theequation{6.7}\begin{equation}
\big[
{\sf A}^{-1}
\circ
\varphi_{\Delta_\mathcal{E}}
\circ
{\sf A}
\big]
(x,a,b)
=
\Big(
\phi^i,
A^q
\big(
\phi^{i_1},\psi^{j_1},\Phi_{\beta(q_1)}^{j(q_1)}
\big),
B^j
\big(
\phi^{i_1},\psi^{j_1},\Phi_{\beta(q_1)}^{j(q_1)}
\big)
\Big).
\end{equation}
In case \thetag{ $\mathcal{ E}$} $=$ \thetag{ $ \mathcal{ E}_1$}, is
an exercise to verify by computations that the $A^q(\cdot)$ and $B^j
(\cdot)$ are independent of $x$. In general however, the explicit
expression of $\Phi_{ i_1, \dots, i_\lambda}^j$ is unknown.
Unfortunately also, nothing shows how $\big( f(a,b), g(a,b) \big)$ is
uniquely associated to $\varphi (x,y)$. 
Further explanations are needed.

\subsection*{ 6.8.~Determination of parameter transformations}
At first, we state a geometric reformulation of the preceding lemma.

\def\thelemma{6.9}\begin{lemma}
Every Lie symmetry $(x, y) \mapsto \varphi (x, y)$ of \thetag{
$\mathcal{ E}$} induces a local $\K$-analytic diffeomorphism
\def\theequation{6.10}\begin{equation}
(x,y,a,b)
\longmapsto
\big(
\varphi(x,y),h(a,b)
\big)
\end{equation}
of $\K_x^n \times \K_y^m \times \K_a^p \times \K_b^m$ that maps to
itself the associated submanifold of solutions
\def\theequation{6.11}\begin{equation}
\mathcal{M}_\mathcal{E} 
=
\big\{
(x,y,a,b):\
y
=
\Pi(x,a,b)
\big\}.
\end{equation}
\end{lemma}

\proof
In fact, we know that the $n$-dimensional leaf $\big\{
\big( x, \Pi (x, a, b) \big): \, x\in \K^n \big\}$ is sent $\big\{
\big( x, \Pi (x, h(a,b)) \big): \, x\in \K^n \big\}$.
\endproof

Equivalently, setting $c := (a, b)$ and writing $(\varphi, h) = (\phi,
\psi, h)$, we have $\psi = \Pi (\phi, h)$ when $y = \Pi (x, c)$,
namely
\def\theequation{6.12}\begin{equation}
\psi(x,\Pi(x,c))
\equiv
\Pi\big(
\phi(x,\Pi(x,c)),h(c)
\big)
\end{equation}

\def\theproposition{6.13}\begin{proposition}
There exists a universal rational map $\widehat{ \sf H}$ such that
\def\theequation{6.14}\begin{equation}
h(c)
\equiv
\widehat{\sf H}
\Big(
J_{x,a,b}^{\kappa+1}\,\Pi(x,c),
J_{x,y}^\kappa
\varphi(x,\Pi(x,c))
\Big)
\end{equation}
\end{proposition}

This shows unique determination of $h$ from $\varphi$, given \thetag{
$\mathcal{ E}$} or equivalently, given $\Pi$.

\proof
Differentiating a function $\chi (x, \Pi (x, c))$ with respect to
$x^k$, $k=1, \dots, n$, corresponds to applying to $\chi$ the vector
field
\def\theequation{6.15}\begin{equation}
{\sf L}_k
:=
\frac{\partial}{\partial x_k}+
\sum_{j=1}^m \, 
\frac{\partial \Pi^j}{\partial x_k}(x,c)\, 
\frac{\partial}{\partial y^j}, 
\ \ \ \ \ \ 
k=1,\dots,n.
\end{equation}
Thus, applying ${\sf L}_k$ to the $m$ 
scalar equations~\thetag{ 6.12}, we get
\def\theequation{6.16}\begin{equation}
{\sf L}_k\,\psi^j
=
\sum_{l=1}^n\, 
\frac{\partial\Pi^j}{\partial x^l}\, 
{\sf L}_k\,\phi^l,
\end{equation}
for $1\leqslant k \leqslant n$ and $1\leqslant j \leqslant m$.
It follows from 
the assumption that $\varphi$ is a local diffeomorphism
that ${\rm det}\,
\big(
{\sf L}_k\,\phi^l( 0 )
\big)_{1 \leqslant k \leqslant n }^{1\leqslant l\leqslant n}
\neq 0$ also.
So we may solve the first derivatives
$\Pi_x$ above: there exist universal polynomials
${\sf S}_l^j$ such that
\def\theequation{6.17}\begin{equation}
\frac{\partial\Pi^j}{\partial x^l}
=
\frac{
{\sf S}_l^j
\left(
\big\{{\sf L}_{k'}\,\varphi^{i'}\big\}_{
1\leqslant k'\leqslant n}^{1\leqslant i'\leqslant n+m}
\right)
}{
{\rm det}
\big(
{\sf L}_{k'}\,\phi^{l'}\big)_{
1\leqslant k'\leqslant n}^{1\leqslant l'\leqslant n}
}.
\end{equation}
Again, we apply the ${\sf L}_k$ to these equations, 
getting, thanks to the chain rule:
\def\theequation{6.18}\begin{equation}
\sum_{l_2=1}^n\, 
\frac{\partial^2\Pi^j}{\partial
x^{l_1}x^{l_2}}\, 
{\sf L}_k\,\phi^{l_2}
=
\frac{{\sf R}_{l_1,k}^j
\left(
\big\{{\sf L}_{k_1'}{\sf L}_{k_2'}\varphi^{i'}\big\}_{
1\leqslant k_1',k_2'\leqslant n}^{
1\leqslant i'\leqslant n+m}
\right)
}{ 
\left[
{\rm det}
\big(
{\sf L}_{k'}\,\phi^{l'}\big)_{
1\leqslant k'\leqslant n}^{1\leqslant l'\leqslant n}
\right]^2
}.
\end{equation} 
Here, ${\sf R}_{ l_1, k}^j$ are universal 
polynomials.
Solving the second derivatives 
$\Pi_{ x^{ l_1} x^{ l_2}}^j$, we get
\def\theequation{6.19}\begin{equation}
\frac{\partial^2\Pi^j}{\partial
x^{l_1}x^{l_2}}
=
\frac{{\sf S}_{l_1,l_2}^j
\left(
\big\{{\sf L}_{k_1'}{\sf L}_{k_2'}\varphi^{i'}\big\}_{
1\leqslant k_1',k_2'\leqslant n}^{
1\leqslant i'\leqslant n+m}
\right)
}{ 
\left[
{\rm det}
\big(
{\sf L}_{k'}\,\phi^{l'}\big)_{
1\leqslant k'\leqslant n}^{1\leqslant l'\leqslant n}
\right]^3
}.
\end{equation}
By induction, for every $\beta \in \N^n$:
\def\theequation{6.20}\begin{equation}
\frac{\partial^{\vert\beta \vert} 
\Pi^j}{
\partial x^\beta}
=
\frac{
{\sf S}_\beta^j\left(
\big\{
{\sf L}^{\beta'}\varphi^{i'}
\big\}_{
\vert\beta'\vert\leqslant\vert\beta\vert}^{
1\leqslant i'\leqslant n+m}
\right)
}{
\left[
{\rm det}
\big(
{\sf L}_{k'}\,\phi^{l'}\big)_{
1\leqslant k'\leqslant n}^{1\leqslant l'\leqslant n}
\right]^{2\vert \beta \vert+1}}, 
\end{equation}
where ${\sf S }_\beta^j$ are universal polynomials.
Here, for $\beta' \in \N^n$, we denote by ${\sf L}^{\beta'}$ the
derivation of order $\vert \beta' \vert$ defined by $( {\sf
L}_1)^{ \beta_1 '} \cdots ({\sf L}_n )^{ \beta_n '}$. 

Next, thanks to the assumption that $\mathcal{ M}$ is solvable with respect
to the parameters, there exist integers $j(1), \dots, j(p)$ with
$1\leqslant j(q)\leqslant m$ and multiindices $\beta (1), \dots, \beta
(p) \in \N^n$ with $\vert \beta (q) \vert \geqslant 1$ and $\max_{1
\leqslant q \leqslant p} \, \vert \beta (q) \vert = \kappa$ such that
the local $\K$-analytic map
\def\theequation{6.21}\begin{equation}
\K^{p+m}\ni
c
\longmapsto 
\left(
\big(
\Pi^j(0,c)
\big)^{1\leqslant j\leqslant m},\
\left(
\frac{\partial^{\vert \beta(q)\vert} 
\Pi^{j(q)}}{\partial
x^{\beta(q)}}(0,c)
\right)_{1\leqslant q\leqslant p}
\right)\in\K^{p+m}
\end{equation}
has rank $p+m$ at $c=0$. We then consider in~\thetag{ 6.20} 
only the $(p + m)$ equations written for $(j, 0)$,
$(j(q), \beta (q))$ and we solve $h(c)$ by means of the analytic
implicit function theorem:
\def\theequation{6.22}\begin{equation}
\aligned
h
&
=
\widehat{H}
\left(\phi,
\frac{{\sf S}_{\beta(1)}^{j(1)}\left(
\big\{{\sf L}^{\beta'}\varphi^{i'}\big\}_{
\vert\beta'\vert\leqslant\vert\beta(1)\vert}^{
1\leqslant i'\leqslant n+m}
\right)
}{
{\rm det}
\left[
\big(
{\sf L}_{k'}\,\phi^{l'}\big)_{
1\leqslant k'\leqslant n}^{1\leqslant l'\leqslant n}
\right]^{2\vert\beta(1)\vert+1}
},\dots,
\frac{{\sf S}_{\beta(p)}^{j(p)}\left(
\big\{{\sf L}^{\beta'}\varphi^{i'}\big\}_{
\vert\beta'\vert\leqslant\vert\beta(p)\vert}^{
1\leqslant i'\leqslant n+m}
\right)
}{
{\rm det}
\left[
\big(
{\sf L}_{k'}\,\phi^{l'}\big)_{
1\leqslant k'\leqslant n}^{1\leqslant l'\leqslant n}
\right]^{2\vert\beta(p)\vert+1}
}\right).
\endaligned
\end{equation}
Finally, by developping every derivative ${\sf L}^{ \beta'} \varphi^{
i'}$ (including ${ \sf L}_{ k'} \phi^{ l'}$ as a special case), taking
account of the fact that the coefficients of the ${\sf L}_{ k'}$
depend directly on $\Pi$, we get some universal polynomial ${\sf P}_{
\beta'} \big( J_x^{ \vert \beta' \vert + 1}\, \Pi, J_{ x,y}^{ \vert
\beta ' \vert} \, \varphi^{ i'} \big)$.
Inserting above, we get the map
$\widehat{ \sf H}$.
\endproof

\subsection*{6.23.~Pseudogroup of twin transformations} 
The previous considerations lead to introducing the following. 

\def\thedefinition{6.24}\begin{definition}{\rm
By ${\sf G}_{ \sf v, \sf p}$, we denote the infinite-dimensional
(pseudo)group of local $\K$-analytic diffeomorphisms
\def\theequation{6.25}\begin{equation}
(x,y,a,b)
\longmapsto 
\big( 
\varphi(x,y),h(a,b)
\big) 
\end{equation}
that respect the separation between the variables and the parameters.
}\end{definition}

A converse to Lemma~6.3 holds.

\def\thelemma{6.26}\begin{lemma}
Let $\mathcal{ M}$ be a submanifold $y = \Pi (x, a, b)$ that is
solvable with respect to the parameters $(a, b)$. If a local
$\K$-analytic diffeomorphism $(x,y,a,b)
\longmapsto 
\big( 
\varphi(x,y),h(a,b)
\big)$ of $\K_x^n \times \K_y^m \times \K_a^p
\times \K_b^m$ belonging to ${\sf G}_{ {\sf v}, {\sf p} }$
sends $\mathcal{ M}$ to $\mathcal{ M}$, then $(x, y) \mapsto \varphi
(x, y)$ is a Lie symmetry of the {\sc pde} system $\mathcal{ E}_{
\mathcal{ M} }$ associated to $\mathcal{ M}$.
\end{lemma}

\proof
In fact, since $(\varphi, h)$ respects the separation of
variables and stabilizes $\mathcal{ M}$, it
respects the fundamental pair of foliations
$\big( {\sf F}_{\sf v}, {\sf F}_{\sf p} \big)$, 
namely $\{ (a, b) = (a_0, b_0) \} \cap \mathcal{ M}$
is sent to $\{ (a, b) = h (a_0, b_0) \} \cap \mathcal{ M}$ and
$\{ (x, y) = (x_0, y_0) \} \cap \mathcal{ M}$
is sent to $\{ (x, y) = \varphi (x_0, y_0) \} \cap \mathcal{ M}$.
Hence
$\varphi_{ \Delta_\mathcal{ E_\mathcal{ M} }}$ 
also stabilizes ${\sf F}_{ \Delta_\mathcal{ E}}$.
\endproof

\def\thecorollary{6.27}\begin{corollary}
Through the one-to-one correspondence 
$(\mathcal{ E})
\longleftrightarrow
\mathcal{ M}$
of Proposition~2.17,
Lie symmetries of \thetag{ $\mathcal{ E}$} correspond to elements of
${\sf G}_{ {\sf v}, {\sf p}}$ which stabilize $\mathcal{ M}$.
\end{corollary}

\def\thedefinition{6.28}\begin{definition}{\rm
Let ${\sf Aut}_{{\sf v}, {\sf p}} (\mathcal{ M})$ denote the local
(pseudo)group of $(\varphi, h) \in {\sf G}_{ {\sf v}, {\sf p}}$
stabilizing $\mathcal{ M}$. Let ${\sf Lie} (\mathcal{ E})$ denote
the local (pseudo)group of Lie symmetries
of \thetag{ $\mathcal{ E}$}.
}\end{definition}

In summary:
\def\theequation{6.29}\begin{equation}
\boxed{
{\sf Lie}(\mathcal{E})
\simeq
{\sf Aut}_{{\sf v},{\sf p}}\big(
\mathcal{ M}_{(\mathcal{E})}
\big)
}
\ \ \ \ \ \ \ \ \ \
\text{\rm and}
\ \ \ \ \ \ \ \ \ \
\boxed{
{\sf Aut}_{{\sf v},{\sf p}}(\mathcal{M})
\simeq
{\sf Lie}\big(
\mathcal{E}_\mathcal{M}
\big)
}.
\end{equation}

\subsection*{ 6.30.~Transfer of infinitesimal Lie symmetries}
Let $\mathcal{ L} \in \mathfrak{ SYM} (\mathcal{ E})$, {\it i.e.}
$\mathcal{ L}_{ \Delta_\mathcal{ E}}$ is tangent to $\Delta_\mathcal{
E}$. Through the diffeomorphism ${\sf A}$, the push-forward of
$\mathcal{ L}_{ \Delta_\mathcal{ E} }$ must be of the form
\def\theequation{6.31}\begin{equation}
{\sf A}_*^{-1}(\mathcal{L}_{\Delta_\mathcal{E}})
=
\sum_{k=1}^n\,\Theta^i(x,a,b)\, 
\frac{\partial }{\partial x^i}
+
\sum_{q=1}^p\, 
\mathcal{F}^q(a,b)\, 
\frac{\partial }{\partial a^q}
+
\sum_{j=1}^m\, 
\mathcal{G}^j(a,b)\,
\frac{\partial }{\partial b^j}, 
\end{equation}
where the last two families of $\K$-analytic
coefficients $\mathcal{ F}^q$ and $\mathcal{ G }^j$ depend only on $(a,
b)$.

\def\thelemma{6.32}\begin{lemma}
To every infinitesimal symmetry $\mathcal{ L}$ of \thetag{
$\mathcal{ E}$}, we can associate an infinitesimal symmetry
\def\theequation{6.33}\begin{equation}
\mathcal{L}^*
:=
\sum_{q=1}^p\, 
\mathcal{F}^q(a,b)\, 
\frac{\partial }{\partial a^q}
+
\sum_{j=1}^m\, 
\mathcal{G}^j(a,b)\,
\frac{\partial }{\partial b^j}
\end{equation}
of the space of parameters which tells how the flow of $\mathcal{ L}$
acts infinitesimally on the leaves of ${\sf F}_{ \! \Delta_\mathcal{
E}}$. Furthermore, $\mathcal{ L} + \mathcal{ L }^*$ is tangent to the
submanifold of solutions $\mathcal{ M}_{ (\mathcal{ E } )}$.
\end{lemma}

Considering the flow of $\mathcal{ L} + \mathcal{ L}^*$
reduces these assertions and the next to the arguments
of the preceding paragraphs. So we summarize.

\def\thelemma{6.34}\begin{lemma}
Let $\mathcal{ M}$ be a submanifold $y = \Pi (x, a, b)$ that is
solvable with respect to the parameters $(a, b)$.
If a vector field that respects the separation 
between variables and parameters, namely of the form
\def\theequation{6.35}\begin{equation}
\mathcal{L}
+
\mathcal{L}^*
=
\sum_{i=1}^n\, 
\mathcal{X}^i(x,y)\,
\frac{\partial}{\partial x^i}
+
\sum_{j=1}^m\,
\mathcal{Y}^j(x,y)\,
\frac{\partial}{\partial y^j}
+
\sum_{q=1}^p\,
\mathcal{F}^q(a,b)\,
\frac{\partial}{\partial a^q}
+
\sum_{j=1}^m\,
\mathcal{G}^j(a,b)\,
\frac{\partial}{\partial b^j}
\end{equation}
is tangent to $\mathcal{ M}$, then 
$\mathcal{ L}$ is an infinitesimal Lie
symmetry of $\big( \mathcal{ E}_\mathcal{ M} \big)$
\end{lemma}

\def\thecorollary{6.36}\begin{corollary}
Through the one-to-one correspondence 
$(\mathcal{ E})
\longleftrightarrow
\mathcal{ M}$
of Proposition~2.17,
infinitesimal Lie symmetries of \thetag{ $\mathcal{ E}$} correspond to
vector fields $\mathcal{ L} + \mathcal{ L}^*$ tangent to
$\mathcal{ M }$.
\end{corollary}

\def\thedefinition{6.37}\begin{definition}{\rm
Let $\mathfrak{ SYM } ( \mathcal{ M} )$ denote the Lie algebra of
vector fields $\mathcal{ L} + \mathcal{ L}^*$ tangent to $\mathcal{
M}$. Let $\mathfrak{ SYM} (\mathcal{ E })$ denote the Lie algebra of
infinitesimal Lie symmetries of \thetag{ $\mathcal{ E }$}.
}\end{definition}

In summary:
\def\theequation{6.38}\begin{equation}
\boxed{
\mathfrak{SYM}(\mathcal{E})
\simeq
\mathfrak{SYM}
\big(
\mathcal{M}_{(\mathcal{E})}
\big)
}
\ \ \ \ \ \ \ \ \ \
\text{\rm and}
\ \ \ \ \ \ \ \ \ \
\boxed{
\mathfrak{SYM}
\big(
\mathcal{M}
\big)
\simeq
\mathfrak{SYM}
\big(
\mathcal{E}_\mathcal{M}
\big)
}.
\end{equation}

\subsection*{ 6.39.~Dual defining equations} 
As in \S2.10, let $\mathcal{ M} \subset \K_x^n \times \K_y^m \times
\K_a^p \times \K_b^m$ given by $0 = - y^j + \Pi^j (x, a, b)$ and
assume if to be solvable with respect to the parameters. In
particular, we can solve the $b^j$, obtaining {\sl dual defining
equations}
\def\theequation{6.40}\begin{equation}
b^j
=
{\Pi^*}^j(a,x,y),
\ \ \ \ \ \ \ 
j=1,\dots,m,
\end{equation}
for some local $\K$-analytic map map ${\Pi^*} = ({\Pi^*}^1, \dots,
{\Pi^*}^m )$ satisfying 
\def\theequation{6.41}\begin{equation}
b
\equiv
{\Pi^*}\big(a,x,\Pi(x,a,b)\big)
\ \ \ \ \ \ \
{\rm and}
\ \ \ \ \ \ \
y
\equiv
\Pi\big(x,a,{\Pi^*}(a,x,y)\big).
\end{equation}

\subsection*{6.42.~An algorithm for the computation of $\mathfrak{ 
SYM} (\mathcal{ M })$} The tangency to $\mathcal{ M}$ is expressed 
by applying the vector field~\thetag{ 6.35} to $0 = - y^j + \Pi^j 
(x,a, b)$, which yields:
\def\theequation{6.43}\begin{equation}
\aligned
0
=
-
\mathcal{Y}^j(x,y)
+
\sum_{i=1}^n\,\mathcal{X}^i(x,y)\, 
\Pi_{x^i}^j(x,a,b)
&
+
\sum_{q=1}^p\,
\mathcal{F}^q(a,b)\, 
\Pi_{a^q}^j(x,a,b)
\\
&
+
\sum_{l=1}^m\,
\mathcal{G}^l(a,b)\, 
\Pi_{b^l}^j(x,a,b),
\endaligned
\end{equation}
for $j=1, \dots, m$ and for $(x, y, a, b) \in \mathcal{ M}$. In fact,
after replacing the variable $y$ by $\Pi (x, a, b)$, these equations
should be interpreted as power series identities in $\K\{ x, a, b\}$.

Denote by $\Delta(x, a, b)$ the determinant of the (invertible) matrix
$\big( \Pi_{ b^l}^j (x, a, b ) \big)_{1 \leqslant l, j \leqslant m }$
and by $D(x, a, b)$ its matrix of cofactors, so that $\Pi_b^{ -1} =
[\Delta ]^{-1} \, D$. Hence we can solve $\mathcal{ G}$ from~\thetag{
6.43}:
\def\theequation{6.44}\begin{equation}
\left\{
\aligned
\mathcal{G}(a,b)
\equiv 
\frac{D(x,a,b)}{\Delta(x,a,b)}
&
\left[
\mathcal{Y}
\big(
x,\Pi(x,a,b)
\big)
-
\sum_{i=1}^n\, 
\mathcal{X}^i
\big(
x,\Pi(x,a,b)
\big)\, 
\Pi_{x^i}(x,a,b)
-
\right. 
\\
& 
\left. \ \ \ 
-
\sum_{q=1}^p\, 
\mathcal{F}^q(a,b)\, 
\Pi_{a^q}(x,a,b)
\right].
\endaligned\right.
\end{equation}
Next, we aim to solve the $\mathcal{ F}^q (a,b)$. Consequently, we
gather all the other terms in the brackets as $\Psi_0 \big( J_{ x, a,
b}^1 \Pi, \mathcal{ X}, \mathcal{ Y} \big)$:
\def\theequation{6.45}\begin{equation}
\mathcal{G}(a,b)
\equiv 
\frac{D(x,a,b)}{\Delta(x,a,b)}
\left[
-
\sum_{q=1}^p\,\mathcal{F}^q(a,b)\,
\Pi_{a^q}(x,a,b)
\right]
+
\frac{\Psi_0
\big( J_{ x, a, b}^1 \Pi,
\mathcal{ X}, \mathcal{ Y}
\big)}{
\Delta(x,a,b)}.
\end{equation}
Here, $\Psi_0$ is linear with respect to $(\mathcal{ X}, \mathcal{
Y})$, with polynomial coefficients of degree one in $J_{ x, a, b}^1
\Pi$.

Next, for $k=1, \dots, n$, we differentiate this identity with respect
to $x_k$. Then $\mathcal{ G} (a, b)$ disappears and we chase the
denominator $\Delta^2$:
\def\theequation{6.46}\begin{equation}
\left\{
\aligned
0\equiv 
& \
\left[
\Delta\,D
\right]
\left[
-
\sum_{ q=1}^p\,\mathcal{F}^q(a,b)
\,\Pi_{a^qx^k}(x,a,b)
\right]
+\\
&
\ \ \ \ \ \ \ \ \ \ \ \ \ 
+
\left[
\Delta\,D_{x_k}
-
\Delta_{x_k}\,D
\right]\left[
-
\sum_{q=1}^p\,\mathcal{F}^q(a,b)\, 
\Pi_{a^q}(x,a,b)
\right]
+\\
&
\ \ \ \ \ \ \ \ \ \ \ \ \ 
+
\Psi_k
\big(
J_{x,a,b}^2\Pi,
J_{x,y}^1\mathcal{X},
J_{x,y}^1\mathcal{Y}
\big).
\endaligned\right.
\end{equation}
The $\Psi_k$ are linear with respect to $(J_{ x, y}^1 \mathcal{ X},
J_{ x, y}^1 \mathcal{ Y})$, with polynomial coefficients in $J_{ x, a,
b }^2 \Pi$. Then we further differentiate
with respect to $x$ and by induction, for every $\beta \in
\N^n$, we get:
\def\theequation{6.47}\begin{equation}
\left\{
\aligned
0
\equiv 
&\
\left[
\Delta\,D
\right]
\left[
-
\sum_{q=1}^p\,\mathcal{F}^q(a,b) \, 
\Pi_{a^qx^\beta}(x,a,b)
\right]
+\\
& \
+
\sum_{\vert\beta_1\vert<\vert\beta\vert}\, 
{\sf D}_{\beta,\beta_1}
\big(
J^{\vert\beta_1\vert+1}
\Pi
\big)
\left[
-
\sum_{q=1}^p\,\mathcal{F}^q(a,b)\,
\Pi_{
a^qx^{\beta_1}}
(x,a,b)
\right] 
+\\
&\
+
\Psi_\beta
\big(
J_{x,a,b}^{\vert\beta\vert+1}\Pi,
J_{x,y}^{\vert \beta \vert}\mathcal{X}, 
J_{x,y}^{\vert \beta \vert}\mathcal{Y}),
\endaligned\right.
\end{equation}
where the expressions ${\sf D}_{ \beta, \beta_1 }$ are certain $m
\times m$ matrices with polynomial coefficients in the jet $J_{ x, a,
b}^{ \vert \beta_1 \vert +1} \Pi$, and where the terms $\Psi_\beta
\big( J_{ x, a, b}^{ \vert \beta \vert +1} \Pi, J_{ x, y}^{\vert \beta
\vert} \mathcal{ X}, J_{ x, y}^{\vert \beta \vert} \mathcal{ Y} \big)$
are linear with respect to $ \big( J_{ x, y}^{\vert \beta \vert}
\mathcal{ X}, J_{ x, y}^{\vert \beta \vert} \mathcal{ Y} \big)$, with
polynomial coefficients in $J_{ x, a, b}^{\vert \beta \vert + 1} \Pi$.

Writing these identity for $(j, \beta) = (j(q), \beta (q))$, $q = 1,
\dots, p$, reminding $\max_{ 1 \leqslant q \leqslant p}\, \vert \beta
(q) \vert = \kappa$, it follows from the assumption of solvability
with respect to the parameters (a boring technical check is needed)
that we may solve
\def\theequation{6.48}\begin{equation}
\mathcal{F}^q(a,b)
\equiv 
\Phi^q
\big(J_{x,a,b}^{\kappa+1}
\Pi(x,a,b),
J_{x,y}^\kappa\mathcal{X}
(x,\Pi(x,a,b)), 
J_{x,y}^\kappa
\mathcal{Y}
(x,\Pi(x,a,b))
\big),
\end{equation}
for $q = 1, \dots, p$, where each local $\K$-analytic function
$\Phi_q$ is linear with respect to $(J^\kappa \mathcal{ X },
J^\kappa \mathcal{ Y})$ and rational with respect to $J^{ \kappa
+1} \Pi$, with denominator not vanishing at $(x, a, b) := (0, 0, 0)$.

Pursuing, we differentiate~\thetag{6.48} with respect to $x^l$ for $l=1,
\dots, n$. Then $\mathcal{ F}^q (a, b)$ disappears and we get:
\def\theequation{6.49}\begin{equation}
0
\equiv
\Phi_{q,l}
\big(
J_{x,a,b}^{\kappa+2}
\Pi(x,a,b), 
J_{x,y}^{\kappa+1}
\mathcal{X}(x,\Pi(x,a,b)), 
J_{x,y}^{\kappa+1}
\mathcal{Y}(x,\Pi(x,a,b))
\big),
\end{equation}
for $1 \leqslant q \leqslant p$ and $1 \leqslant l \leqslant
n$. In~\thetag{ 6.46}, we then replace the functions $\mathcal{ F }^q$
by their values $\Phi^q$:
\def\theequation{6.50}\begin{equation}
0
\equiv 
\Psi_{k,j}
\big(
J_{x,a,b}^{\kappa+1}
\Pi(x,a,b), 
J_{x,y}^\kappa
\mathcal{X}
(x,\Pi(x,a,b)),
J_{x,y}^\kappa
\mathcal{Y}
(x,\Pi(x,a,b))
\big),
\end{equation}
for $1 \leqslant k \leqslant n$ and $1 \leqslant j \leqslant m$. Then
we replace the variable $b$ by ${\Pi^*} (a, x, y)$ in the two obtained
systems~\thetag{6.49} and~\thetag{6.50}; taking account of the
functional identity $y \equiv \Pi \big ( x, a, \Pi^* (a, x, y) \big)$
written in~\thetag{ 6.41}, we get
\def\theequation{6.51}\begin{equation}
\left\{
\aligned
0 \equiv
& \ 
\Phi_{q,l}
\big(
J_{x,a,b}^{\kappa+2}
\Pi(x,a,{\Pi^*}(a,x,y)), 
J_{x,y}^{\kappa+1}
\mathcal{X}(x,y), 
J_{x,y}^{\kappa+1}
\mathcal{Y}(x,y)
\big), 
\\
0 
\equiv 
&\
\Psi_{k,j}
\big(
J_{x,a,b}^{\kappa+1}
\Pi(x,a,{\Pi^*}(a,x,y)),
J_{x,y}^\kappa
\mathcal{X}(x,y), 
J_{x,y}^\kappa
\mathcal{Y}(x,y)
\big).
\endaligned\right.
\end{equation}
Finally, we develope these equations in power series
with respect to $a$:
\def\theequation{6.52}\begin{equation}
\left\{
\aligned
0 \equiv 
& \ 
\sum_{\gamma\in\N^p}\,a^\gamma \
\Phi_{q,l,\gamma}
\big(
x,y,
J_{x,y}^{\kappa+1}
\mathcal{X}(x,y), 
J_{x,y}^{\kappa+1}
\mathcal{Y}(x,y)
\big), \\
0 \equiv 
& \
\sum_{\gamma\in\N^p}\,
a^\gamma \
\Psi_{k,j,\gamma}
\big(
x,y,
J_{x,y}^\kappa
\mathcal{X}(x,y), 
J_{x,y}^\kappa
\mathcal{Y}(x,y)
\big),
\endaligned\right.
\end{equation}
where the terms $\Phi_{ q, l, \gamma}$ and $\Psi_{ k, j, \gamma}$ are
linear with respect to the jets of $\mathcal{ X}$, $\mathcal{ Y}$.

\def\theproposition{6.53}\begin{proposition}
A vector field~\thetag{ 6.35} belongs to $\mathfrak{ SYM} ( \mathcal{
M})$ if and only if $\mathcal{ X}^i, \mathcal{ Y}^j$ satisfy the
linear {\sc pde} system
\def\theequation{6.54}\begin{equation}
\left\{
\aligned
0\equiv 
& \
\Phi_{q,l,\gamma}
\big(x,y,
J_{x,y}^{\kappa+1}
\mathcal{X}(x,y), 
J_{x,y}^{\kappa+1}
\mathcal{Y}(x,y)
\big), 
\\
0 
\equiv 
& \
\Pi_{k,j,\gamma}
\big(x,y,
J_{x,y}^\kappa
\mathcal{X}(x,y), 
J_{x,y}^\kappa
\mathcal{Y}(x,y)
\big),
\endaligned\right.
\end{equation}
where $1 \leqslant q\leqslant p$, $1 \leqslant l \leqslant n$, $1
\leqslant k \leqslant n$ and $\gamma \in \N^p$. Then $\mathcal{ F}^q$
defined by~\thetag{ 6.48} and $\mathcal{ G}^j$ defined by~\thetag{ 6.45}
are independent of $x$.
\end{proposition}

This provides a second algorithm, essentially equivalent to Sophus
Lie's.

\def\theexample{6.55}\begin{example}
{\rm 
For $y_{ xx}( x) = F (x, y(x), y_x(x))$, the first line
of~\thetag{ 6.54} is (the second one is redundant):
\def\theequation{6.56}\begin{equation}
\left\{
\aligned
0\equiv &\
\mathcal{X}\left[-
\Pi_{xa}\Pi_{xxx}\Pi_b+
\Pi_a\Pi_{xb}\Pi_{xxx}-
\Pi_x\Pi_{xxa}\Pi_{xb}+
\Pi_{xa}\Pi_x\Pi_{xxb}\right.
+\\
& 
\ \ \ \ \ \ \ \ \ \ \ \ \ \ \ \ \ \ \ \ \ \ \ \ \ \ \ \ \ \ \ \ \ 
\ \ \ \ \ \ \ \ \ \ \ \ \ \ \ \ \ \ \ \ \ \ \ \ \
\left. +\Pi_{xxa}\Pi_{xx}\Pi_b-
\Pi_a\Pi_{xxb}\Pi_{xx}
\right]
+\\
& \
+
\mathcal{Y}\left[
-
\Pi_{xa}\Pi_{xxb}+\Pi_{xxa}
\Pi_{xb}\right]
+\\
& \
+
\mathcal{X}_x\left[
-
2\Pi_{xx}\Pi_{xa}\Pi_b+2\Pi_{xx}
\Pi_a\Pi_{xb}+
\Pi_x\Pi_b\Pi_{xxa}-
\Pi_x\Pi_a\Pi_{xxb}\right]+\\& \
+\mathcal{Y}_x\left[
-\Pi_b\Pi_{xxa}+\Pi_a\Pi_{xxb}
\right]+\\
& \
+\mathcal{X}_y\left[
-3\Pi_x\Pi_{xx}\Pi_{xa}\Pi_b+
3\Pi_x\Pi_a\Pi_{xx}\Pi_{xb}+
(\Pi_x)^2\Pi_b\Pi_{xxa}-(\Pi_x)^2\Pi_a
\Pi_{xxb}\right]
+\\
& \
+\mathcal{Y}_y\left[
\Pi_{xx}\Pi_b\Pi_{xa}-
\Pi_{xx}\Pi_a\Pi_{xb}-
\Pi_x\Pi_b\Pi_{xxa}+\Pi_x
\Pi_a\Pi_{xxb}\right]
+\\
& \
+\mathcal{X}_{xx}\left[
-\Pi_x\Pi_b\Pi_{xa}+\Pi_x\Pi_a
\Pi_{xb}\right]
+\\
& \
+
\mathcal{X}_{xy}\left[
-2(\Pi_x)^2\Pi_b\Pi_{xa}
+2(\Pi_x)^2\Pi_a\Pi_{xb}\right]
+\\
& \
+ \mathcal{X}_{y^2}\left[
-(\Pi_x)^3\Pi_b\Pi_{xa}+
(\Pi_x)^3\Pi_a\Pi_{xb}\right]+
\\
& \
+
\mathcal{Y}_{xx}\left[
\Pi_b\Pi_{xa}-
\Pi_a\Pi_{xb}\right]+\\
& \
+
\mathcal{Y}_{xy}
\left[2\Pi_x\Pi_b\Pi_{xa}-2\Pi_x\Pi_a
\Pi_{xb}\right]
+\\
& \
+
\mathcal{Y}_{y^2}\left[(\Pi_x)^2\Pi_b\Pi_{xa}-
(\Pi_x)^2\Pi_a\Pi_{xb}
\right].
\endaligned\right.
\end{equation}
We observe the similarity with~\thetag{ 4.19}: the expression is linear
in the partial derivatives of $\mathcal{ X}$, $\mathcal{ Y}$ of order
$\leqslant 2$, but the coefficients in the equation above are more
complicated. In fact, after dividing by $-\Pi_b \, \Pi_{ xa} + \Pi_a
\, \Pi_{ xb}$, this equation coincides with~\thetag{ 4.21}, thanks to
$\Pi_x = y_1$ and to the formulas~\thetag{ 2.34} for $F_x$, $F_y$, $F_{
y_1}$.
}\end{example}

\subsection*{6.57.~Infinitesimal CR automorphisms of generic
submanifolds} If the system \thetag{ $\mathcal{ E}$} is associated to
the complexification $\mathcal{ M} = (M)^c$ of a generic $M \subset
\C^{ n + m}$ as in \S1.16, then $a = (\bar z)^c = \zeta$, $b = (\bar
w)^c = \xi$, and the vector field $\mathcal{ L}^*$ associated to an
infinitesimal Lie symmetry 
\def\theequation{6.58}\begin{equation}
\mathcal{L}
=
\sum_{i=1}^n\,\mathcal{X}^i(z,w)\,
\frac{\partial}{\partial z^i}
+
\sum_{j=1}^m\,\mathcal{Y}^j(z,w)\,
\frac{\partial}{\partial w^j}
\end{equation}
of \thetag{ $\mathcal{ E}$} is simply the complexification
$\underline{ \mathcal{ L}}$ of its conjugate
$\overline{ \mathcal{ L}}$, namely
\def\theequation{6.59}\begin{equation}
\mathcal{L}^*
=
\underline{\mathcal{L}}
=
\sum_{i=1}^n\,\overline{\mathcal{X}}^i(\zeta,\xi)\,
\frac{\partial}{\partial \zeta^i}
+
\sum_{j=1}^m\,\overline{\mathcal{Y}}^j(\zeta,\xi)\,
\frac{\partial}{\partial \xi^j}.
\end{equation}
Then the sum $\mathcal{ L} + \underline{ \mathcal{ L}}$ is tangent to
$\mathcal{ M}$ and its flow stabilizes the two invariant foliations,
obtained by intersecting $\mathcal{ M}$ by $\{ (z, w) = {\rm cst.} \}$
or by $\{ (\zeta, \xi) = {\rm cst.} \}$. In~\cite{ me2005a, me2005b},
these two foliations, denoted $\mathcal{ F}$, $\underline{ \mathcal{
F}}$, are called (conjugate) Segre foliations, since its leaves are
the complexifications of the (conjugate) classical Segre varieties
(\cite{ se1931, pi1975, pi1978, we1977, dw1980, bjt1985, df1988,
ber1999, su2001, su2002, su2003, gm2003a}) associated to $M$, viewed in
its ambient space $\C^{ n + m}$. The next definition is also
classical (\cite{ be1979, lo1981, ekv1985, kr1987, kv1987, be1988,
vi1990, st1996, be1997, ber1999, lo2002, fk2005a, fk2005b}):

\def\thedefinition{6.60}\begin{definition}{\rm
By $\mathfrak{ hol} (M)$ is meant the Lie algebra of local holomorphic
vector fields $\mathcal{L} = \sum_{ i=1 }^n\, \mathcal{ X }^i( z, w)\,
\frac{\partial }{ \partial z^i} + \sum_{ j=1 }^m\, \mathcal{Y}^j( z,
w)\, \frac{ \partial }{ \partial w^j}$ whose real flow $\exp \big (t
\mathcal{ L} \big) (z, w)$ induces one-parameter families of local
biholomorphic transformations of $\C^{ n+ m}$ stabilizing $M$.
Equivalently,
\def\theequation{6.61}\begin{equation}
2\,{\rm Re}\,\mathcal{L}
=
\mathcal{L}
+
\overline{\mathcal{L}}
\end{equation}
is tangent to $M$. Again equivalently, $\mathcal{ L} + \underline{
\mathcal{ L}}$ is tangent to $\mathcal{ M} = M^c$.
}\end{definition}

Then obviously $\mathfrak{ hol} (M)$ is a real Lie algebra.

\def\thetheorem{6.62}\begin{theorem}
{\rm (\cite{ ca1932a, ber1999, gm2004})} The complexification
$\mathfrak{ hol} (M) \otimes \C$ identifies with $\mathfrak{ SYM}
\big( \mathcal{ E} (M^c) \big)$. Furthermore, if $M$ is finitely
nondegenerate and minimal at the origin, both are finite-dimensional
and $\mathfrak{ hol} (M)$ is totally real in $\mathfrak{ SYM} \big(
\mathcal{ E} (M^c) \big)$.
\end{theorem}

The minimality assumption is sometimes presented by saying that the
Lie algebra generated by $T^cM$ generates $TM$ at the origin (\cite{
ber1999}). However, it is more natural to proceed with the fundamental
pair of foliations associated to $\mathcal{ M}$ (\cite{ me2001,
gm2004, me2005a, me2005b}). Anticipating Sections~10 and ~11 to which
the reader is referred, we set.

\def\thedefinition{6.63}\begin{definition}{\rm
A real analytic generic submanifold $M \subset \C^{ n+m}$ is {\sl
minimal} at one of its points $p$ if the fundamental pair of
foliations of its complexification $\mathcal{ M}$ is covering at 
$p$ (Definition~10.17).
}\end{definition}

Further informations may be found in Section~10.
We conclude by formulating applications of
Theorems~5.13 and~5.24.

\def\thecorollary{6.64}\begin{corollary}
The bound $\dim \, \mathfrak{ hol} (M) \leqslant 8$ for a Levi
nondegenerate hypersurface $M \subset \C^2$ is attained if and only if
it is locally biholomorphic to the sphere $S^3 \subset \C^2$.
\end{corollary}

\def\thecorollary{6.65}\begin{corollary}
The bound $\dim \, \mathfrak{ hol} (M) \leqslant n^2 + 4n + 3$ for a
Levi nondegenerate hypersurface $M \subset \C^{ n + 1}$ is attained if
and only if it is locally biholomorphic to the sphere $S^{ 2n + 1}
\subset \C^{ n + 1}$.
\end{corollary}

\section*{\S7.~Equivalence problems and normal forms}

\subsection*{7.1.~Equivalences of submanifolds of solutions}
As in~\S3.1, let \thetag{ $\mathcal{ E}$} and \thetag{ $\mathcal{ E
}'$} be two {\sc pde} systems and assume that $\varphi$ transforms
\thetag{ $\mathcal{ E}$} to \thetag{ $\mathcal{ E }'$}. Defining ${\sf
A}'$ similarly as ${\sf A}$, it follows that
\def\theequation{7.2}\begin{equation}
{{\sf A}'}^{-1}
\circ
\Phi_{\mathcal{E},\mathcal{E}'}
\circ
{\sf A}
(x,a,b)
\equiv
\big(
\theta(x,a,b),
f(a,b),
g(a,b)
\big)
=:
(x',a',b')
\end{equation}
transforms ${\sf F}_{\sf v}$ to ${\sf F}_{\sf v}'$, hence induces a
map $(a, b) \mapsto (a', b')$. The arguments of Section~6 apply here
with minor modifications to provide two fundamental lemmas.

\def\thelemma{7.3}\begin{lemma}
Every equivalence $(x, y ) \mapsto (x', y')$ between to {\sc pde}
systems \thetag{ $\mathcal{ E}$} and \thetag{ $\mathcal{ E}'$} comes
with an associated transformation $(a, b) \mapsto (a', b')$ of the
parameter spaces such that
\def\theequation{7.4}\begin{equation}
(x,y,a,b)
\longmapsto
(x',y',a',b')
\end{equation}
is an equivalence between the associated submanifolds of solutions
$\mathcal{ M}_{ (\mathcal{ E})} \to \mathcal{ M}_{ (\mathcal{ E}'
)}'$.
\end{lemma}

Conversely, let $\mathcal{ M }$ and $\mathcal{ M' }$ be two
submanifolds of $\K_x^n \times \K_y^m \times \K_a^p \times \K_b^m$ and
of $\K_{ x'}^n \times \K_{ y'}^m \times \K_{ a'}^p \times \K_{ b'}^m$
represented by $y = \Pi (x, a, b)$ and by $y' = \Pi' ( x', a', b')$,
in the {\it same}\, dimensions. Assume both are solvable with respect
to the parameters.

\def\thelemma{7.5}\begin{lemma}
Every equivalence 
\def\theequation{7.6}\begin{equation}
(x,y,a,b)
\longmapsto
\big(
\varphi(x,y),h(a,b)
\big)
\end{equation}
between $\mathcal{ M}$ and $\mathcal{ M}'$ belonging to ${\sf G}_{
{\sf v}, {\sf p}}$ induces by projection
the equivalence $(x, y) \mapsto \varphi (x,
y)$ between the associated {\sc pde} systems $\big( \mathcal{
E}_{\mathcal{ M}} \big)$ and $\big( \mathcal{ E}_{\mathcal{ M}'} '
\big)$.
\end{lemma}

\subsection*{7.7.~Classification problems}
Consequently, classifying {\sc pde} systems 
under point transformations (Section~3) is equivalent to 
the following.

\medskip\noindent
{\bf Equivalence problem 7.8.}
Find an algorithm to decide whether two given submanifolds (of
solutions) $\mathcal{ M}$ and $\mathcal{ M}'$ are equivalent through
an element of ${\sf G}_{ {\sf v}, {\sf p}}$.

\medskip\noindent
{\bf Classification problem 7.9.}
Classify submanifolds (of solutions) $\mathcal{ M}$, namely provide a
complete list of all possible such equations, including their
automorphism group ${\sf Aut}_{{\sf v}, {\sf p}} (\mathcal{ M})
\subset {\sf G}_{ {\sf v}, {\sf p}}$.

\subsection*{7.10.~Partial normal forms} 
Both problems above are of high complexity. At least as a preliminary
step, it is useful to try to simplify somehow the defining equations
of $\mathcal{ M}$, by appropriate changes of coordinates belonging to
${\sf G}_{ {\sf v}, {\sf p}}$. To begin with, the next lemma holds for
$\mathcal{ M}$ defined by $y = \Pi (x, a, b)$ with the only assumption
that $b \mapsto \Pi (0, 0, b)$ has rank $m$ at $b=0$.

\def\thelemma{7.11}\begin{lemma}
{\rm (\cite{ cm1974, ber1999, me2005a}, [$*$])}
In coordinates $x' = ({ x'}^1, \dots, {\sc x'}^n)$ and $y' = ( {
y'}^1, \dots, {y'}^m)$ an arbitrary submanifold $\mathcal{ M'}$
defined by $y' = \Pi' (x', a', b')$ or dually by $b' = {\Pi '}^* (a',
x', y')$ is equivalent to
\def\theequation{7.12}\begin{equation}
y
=
\Pi(x,a,b)
\ \ \ \ \ \ \ 
\text{\rm or dually to}
\ \ \ \ \ \ \ 
b
=
\Pi^*(a,x,y)
\end{equation}
with 
\def\theequation{7.13}\begin{equation}
\Pi(0,a,b)
\equiv
\Pi(x,0,b)
\equiv
b
\ \ \ \ \ \ \ 
\text{\rm or dually}
\ \ \ \ \ \ \ 
\Pi^*(0,x,y)
\equiv
\Pi^*(a,0,y)
\equiv
y,
\end{equation}
namely $\Pi = b + {\rm O} (xa)$ and
$\Pi^* = y + {\rm O} (ax)$.
\end{lemma}

\proof
We develope
\def\theequation{7.14}\begin{equation}
y'
=
\Pi'(0,a',b')
+
\Lambda'(x')
+
{\rm O}(x'a').
\end{equation}
Since $b' \mapsto \Pi' (0, a', b')$ has rank $m$ at $b' = 0$, the
coordinate change
\def\theequation{7.15}\begin{equation}
b''
:=
\Pi'(0,a',b'),
\ \ \ \ \
a''
:=
a',
\ \ \ \ \ 
x''
:=
x',
\ \ \ \ \ 
y''
:=
y',
\end{equation}
transforms $\mathcal{ M}'$ to $\mathcal{ M}''$ defined by
\def\theequation{7.16}\begin{equation}
y''
=
\Pi''(x'',a'',b'')
:=
b''
+
\Lambda'(x'')
+
{\rm O}(x''a'').
\end{equation}
Solving $b''$ by means of the implicit function theorem, 
we get
\def\theequation{7.17}\begin{equation}
b''
=
{\Pi''}^*(a'',x'',y'')
=
y''
-
\Lambda'(x'')
+
{\rm O}(a''x''),
\end{equation}
and it suffices to set $y := y'' - \Lambda' (x'')$, $x:= x''$ and $a
:= a''$, $b := b''$.
\endproof

Taking account of solvability with respect to the parameters, 
finer normalizations holds.

\def\thelemma{7.18}\begin{lemma}
With $n = m = \kappa = 1$, every submanifold of solutions $y' = b' +
x'a'\big[ 1 + {\rm O}_1 \big]$ of $y_{ x' x'} ' = F '( x', y', y_{ x'}
')$ is equivalent to
\def\theequation{7.19}\begin{equation}
y_{xx}
=
b
+
xa
+
{\rm O}(x^2a^2).
\end{equation}
\end{lemma}

\proof
Writing $y' = b' + x' \big[ a' + a'\, \Lambda ' (a', b') + {\rm O} (x'
a') \big]$, where $\Lambda ' = {\rm O}_1$, we set $a'' := a ' + a' \,
\Lambda' (a', b')$, $b'' := b'$, $x'' := x'$, $y'' := y'$, whence $y''
= b'' + x'' \big[ a'' + {\rm O} (x'' a'' ) \big]$. Dually $b'' = y''
- a'' \big[ x'' + x'' x'' \, \Lambda'' (x'', y'') + {\rm O} (x'' x''
a'') \big]$, so we set $x := x'' + x'' x'' \, \Lambda '' (x'', y'')$,
$y := y''$, $a := a''$, $b := b''$.
\endproof

\def\thecorollary{7.20}\begin{corollary}
Every second order ordinary differential equation 
$y_{ x' x'} ' = F ' (x', y', y_{ x'} ')$ is
equivalent to 
\def\theequation{7.21}\begin{equation}
y_{xx}
=
(y_x)^2\,{\sf R}(x,y,y_x).
\end{equation}
\end{corollary}

\subsection*{7.22.~Complete normal forms} 
The Moser theory of normal forms may be transferred with minor
modifications to submanifolds of solutions associated to
\thetag{ $\mathcal{ E}_1$} and to \thetag{ $\mathcal{ E}_2$}.

\def\thetheorem{7.23}\begin{theorem}
{\rm (\cite{ cm1974, ja1990}, [$*$])} A local $\K$-analytic
submanifold of solutions associated to \thetag{ $\mathcal{ E }_1$
}{\rm :}
\def\theequation{7.24}\begin{equation}
y'
=
b'+x'a'+{\rm O}_3
=
\sum_{k'\geqslant 0}\,\sum_{l'\geqslant 0}\,
\Pi_{k',l'}'(b')\,{x'}^{k'}{a'}^{l'}
\end{equation}
can be mapped, by a transformation $(x', y', a', b') \mapsto (x, y, a,
b)$ belonging to ${\sf G}_{ {\sf v}, {\sf p}}$, to a submanifold of
solutions of the specific form
\def\theequation{7.25}\begin{equation}
y
=
b
+
xa
+
\Pi_{2,4}(b)\,x^2a^4
+
\Pi_{4,2}(b)\,a^2x^4
+
\sum_{k\geqslant 2}\,
\sum_{l\geqslant 2}\,
\sum_{k+l\geqslant 7}\,
\Pi_{k,l}(b)\,x^ka^l.
\end{equation}

\end{theorem}

Solving $(a, b)$ from $y = \Pi$ and $y_x = \Pi_x$ with $\Pi$ as above,
we deduce the following.

\def\thecorollary{7.26}\begin{corollary}
Every $y_{ x' x'}' = F' (x', y', y_{ x'} ')$ is equivalent to
\def\theequation{7.27}\begin{equation}
\aligned
y_{xx}
&
=
(y_x)^2\,
\big[
x^2\,F_{2,2}(y)
+
x^3\,{\sf r}(x,y)
\big]
+
(y_x)^4\,
\big[
F_{0,4}(y)
+
x\,{\sf r}(x,y)
\big]
+
\\
&
\ \ \ \ \
+
\sum_{k\geqslant 0}\,
\sum_{l\geqslant 0}\,
\sum_{k+l\geqslant 5}\,
F_{k,l}(y)\,x^k\,(y_x)^l.
\endaligned
\end{equation}
\end{corollary}

For the completely integrable system \thetag{ $\mathcal{ E}_2$} having
several dependent variables $(x^1, \dots, x^n)$, $n \geqslant 2$, we
have the following.

\def\thetheorem{7.28}\begin{theorem}
{\rm (\cite{ cm1974}, [$*$])} A local $\K$-analytic submanifold of
solutions associated to \thetag{ $\mathcal{ E}_2$}{\rm :}
\def\theequation{7.29}\begin{equation}
y'
=
b'
+
\sum_{1\leqslant k\leqslant n}\,{x'}^k{a'}^k
+
{\rm O}_3
\end{equation}
can be mapped, by a transformation $(x', y', a', b') \mapsto (x, y, a,
b)$ belonging to ${\sf G}_{ { \sf v}, {\sf p }}$, to a submanifold of
solutions of the specific form{\rm :}
\def\theequation{7.30}\begin{equation}
y
=
b
+
\sum_{1\leqslant k\leqslant n}\,
x^ka^k
+
\sum_{k\geqslant 2}\,\sum_{l\geqslant 2}\,
\Pi_{k,l}(x,a,b)
\end{equation}
where
\def\theequation{7.31}\begin{equation}
\Pi_{k,l}(x,a,b)
:=
\sum_{k_1+\cdots+k_n=k}\,\sum_{l_1+\cdots+l_n=l}\,
\Pi_{k_1,\dots,k_n,l_1,\dots,l_n}(b)\,
(x^1)^{k_1}\cdots(x^n)^{k_n}\,
(a^1)^{l_1}\cdots(a^n)^{l_n}
\end{equation}
with the terms $\Pi_{2,2}$, $\Pi_{ 2,3}$ and $\Pi_{ 3, 3}$
satisfying{\rm :}
\def\theequation{7.32}\begin{equation}
0
=
\Delta\,\Pi_{2,2}
=
\Delta\Delta\,\Pi_{2,3}
=
\Delta\Delta\,\Pi_{3,2}
=
\Delta\Delta\Delta\,\Pi_{3,3},
\end{equation}
where
\def\theequation{7.33}\begin{equation}
\Delta
:=
\sum_{1\leqslant k\leqslant n}\,
\frac{\partial^2}{\partial x^k\partial a^k}.
\end{equation}
\end{theorem}

Exercise: solving $(a^k, b)$ from $y = \Pi$ and $y_{ x^l} = \Pi_{ x^l
}$, with $\Pi$ as above, deduce a complete normal form for \thetag{
$\mathcal{ E }_2$}.

\def\theopenproblem{7.34}\begin{openproblem}
Find complete normal forms for submanifolds of solutions associated to
\thetag{ $\mathcal{ E}_4$} and to \thetag{ $\mathcal{ E}_5$}.
\end{openproblem}

\section*{ \S8.~Study of two specific examples}

\subsection*{ 8.1.~Study of the Lie symmetries of
\thetag{ $\mathcal{ E}_4$}}
Its submanifold of solutions possesses two equations:
\def\theequation{8.2}\begin{equation}
y^1
=
\Pi^1(x,a,b^1,b^2)
\ \ \ \ \ \ \ \ \ \ \ \ \ \
y^2
=
\Pi^2(x,a,b^1,b^2).
\end{equation}
For instance, a generic submanifold $M \subset \C^3$ of CR dimension 1
and of codimension 3 has equations of such a form.

Assuming $\mathcal{ V}_\mathcal{ S} (\mathcal{ E}_4)$ to be twin
solvable and having covering submanifold of solutions
({\it see} Definition~10.17), it may be
verified (for $M \subset \C^3$, {\it see} \cite{ be1997}) that at a
Zariski-generic point, its equations are of the form:
\def\theequation{8.3}\begin{equation}
\aligned
y^1
&
=
b^1
+
xa
+
{\rm O}(x^2)
+
{\rm O}(b^1)
+
{\rm O}(b^2),
\\
y^2
&
=
b^2
+
xa(x+a)
+
{\rm O}(x^3)
+
{\rm O}(b^1)
+
{\rm O}(b^2).
\endaligned
\end{equation}
The model has zero remainders with associated system
\def\theequation{8.4}\begin{equation}
y_1^2
=
2x\,y_1^1
+
(y_1^1)^2,
\ \ \ \ \ \ \ \ \ \ \ \ \ \ \
y_2^1
=
0,
\end{equation}
the third equation $y_2^2 = 2\, y_1^1$ being obtained
by differentiating the first.

We may put the submanifold in partial normal form.
Proceeding as in~\cite{ bes2005}, some
partial normalizations belonging to ${\sf G}_{ \sf v, \sf p}$
yield:
\def\theequation{8.5}\begin{equation}
\aligned
y^1
&
=
b^1
+
ax
+
a^2
\big[
\Pi_{3,2}^1(b)\,x^3
+
\Pi_{4,2}^1(b)\,x^4
+\cdots
\big]
+
{\rm O}(a^3\,x^2),
\\
y^2
&
=
b^2
+
a\big[
x^2
+
\Pi_{4,1}^2(b)\,x^4
+\cdots
\big]
+
a^2
\big[
x
+
\Pi_{3,2}^2(b)\,x^3
+\cdots
\big]
+
{\rm O}(a^3\,x^2).
\endaligned
\end{equation}
Redifferentiating, we get an appropriate, partially normalized
system \thetag{ $\mathcal{ E}_4$}:
\def\theequation{8.6}\begin{equation}
\left\{
\aligned
y_1^2
&
=
y_1^1\big(
2x+{\bf g}^1
\big)
+
(y_1^1)^2\big(
1+{\bf g}^2
\big)
+
(y_1^1)^3\,{\sf s}
+
(y_1^1)^4\,{\sf s}
+
(y_1^1)^5\,{\sf s}
+
(y_1^1)^6\,{\sf R},
\\
y_2^1
&
=
(y_1^1)^2\,{\bf h}
+
(y_1^1)^3\,{\sf R},
\\
y_2^2
&
=
y_1^1\big(
2+{\bf g}_x^1
\big)
+
(y_1^1)^2
\big(
{\bf g}_x^2+(2x+{\bf g}^1){\bf h}
\big)
+
(y_1^1)^3\,{\sf r}
+
(y_1^1)^4\,{\sf r}
+
(y_1^1)^5\,{\sf r}
+
(y_1^1)^6\,{\sf R},
\endaligned\right.
\end{equation}
where, precisely:

\begin{itemize}

\smallskip\item[$\bullet$]
${\bf g}^1$, ${\bf g}^2$ and ${\bf h}$ are functions of $(x, y^1,
y^2)$ satisfying ${\bf g}^j = {\rm O} (xx) + {\rm O} (y^1) + {\rm O}
(y^2)$, $j=1, 2$ and ${\bf h} = {\rm O} (x) + {\rm O} (y^1) + {\rm O}
(y^2)$;

\smallskip\item[$\bullet$]
${\sf r}$ and ${\sf s}$ are unspecified functions, varying in the
context, of $(x, y^1, y^2)$ with ${\sf s} = {\rm O} (x) + {\rm O}
(y^1) + {\rm O} (y^2)$, but possibly ${\sf r} (0 ) \neq 0$;

\smallskip\item[$\bullet$]
${\sf R}$ is a remainder function of all the variables $(x, y^1, y^2,
y_1^1)$ parametrizing $\Delta_{ \mathcal{ E}_4}$.

\end{itemize}\smallskip

Letting $\mathcal{ L} = \mathcal{ X} \, \frac{ \partial }{ \partial x}
+ \mathcal{ Y}^1 \, \frac{ \partial }{ \partial y^1} + \mathcal{ Y}^2
\, \frac{ \partial }{ \partial y^2}$ be a candidate infinitesimal Lie
symmetry and applying $\mathcal{ L}^{ (2)} = \mathcal{ L} + {\bf
Y}_1^1 \, \frac{ \partial }{ \partial y_1^1} + {\bf Y}_1^2 \, \frac{
\partial }{ \partial y_1^2} + {\bf Y}_2^1 \, \frac{ \partial }{
\partial y_2^1} + {\bf Y}_2^2 \, \frac{ \partial }{ \partial y_2^2}$
to $\Delta_{ \mathcal{ E}_4}$, we obtain firstly, 
computing ${\rm mod}\, (y_1^1)^5$:
\def\theequation{8.7}\begin{equation}
\aligned
0
&
\equiv
-
{\bf Y}_1^2
+
\big[\mathcal{X}\big]
\big(
y_1^1(2+{\bf g}_x^1)
+
(y_1^1)^2\,{\bf g}_x^2
+
(y_1^1)^3\,{\sf r}
+
(y_1^1)^4\,{\sf r}
\big)
+
\\
&\ \ \ \ \
+
\big[
\mathcal{Y}^1
\big]
\big(
y_1^1\,{\sf r}
+
(y_1^1)^2\,{\sf r}
+
(y_1^1)^3\,{\sf r}
+
(y_1^1)^4\,{\sf r}
\big)
+
\\
&\ \ \ \ \
+
\big[
\mathcal{Y}^2
\big]
\big(
y_1^1\,{\sf r}
+
(y_1^1)^2\,{\sf r}
+
(y_1^1)^3\,{\sf r}
+
(y_1^1)^4\,{\sf r}
\big)
+
\\
&\ \ \ \ \
+
{\bf Y}_1^1
\big(
2x+{\bf g}^1
+
y_1^1(2+2\,{\bf g}^2)
+
(y_1^1)^2\,{\sf s}
+
(y_1^1)^3\,{\sf s}
+
(y_1^1)^4\,{\sf s}
\big),
\endaligned
\end{equation}
and secondly, computing ${\rm mod}\, (y_1^1)^2$:
\def\theequation{8.8}\begin{equation}
0
\equiv
-
{\bf Y}_2^1
+
2\,y_1^1\,{\bf Y}_1^1\,h.
\end{equation}
The third Lie equation involving ${\bf Y}_2^2$ will be
superfluous. Specializing~\thetag{ 4.6}(II) to $m = 2$, we get ${\bf
Y}_1^1$ and ${\bf Y }_1^2$:
\def\theequation{8.9}\begin{equation}
\aligned
{\bf Y}_1^1
&
=
\mathcal{Y}_x^1
+
\big[\mathcal{Y}_{y^1}^1-\mathcal{X}_x\big]\,y_1^1
+
\big[\mathcal{Y}_{y^2}^1\big]\,y_1^2
+
\big[-\mathcal{X}_{y^1}\big]\,(y_1^1)^2
+
\big[-\mathcal{X}_{y^2}\big]\,y_1^1\,y_1^2
\\
{\bf Y}_1^2
&
=
\mathcal{Y}_x^2
+
\big[\mathcal{Y}_{y^1}^2\big]\,y_1^1
+
\big[\mathcal{Y}_{y^2}^2-\mathcal{X}_x\big]\,y_1^2
+
\big[-\mathcal{X}_{y^1}\big]\,y_1^2\,y_1^1
+
\big[-\mathcal{X}_{y^2}\big]\,(y_1^2)^2.
\endaligned
\end{equation}
and also ${\bf Y}_2^1$ and ${\bf Y}_2^2$ (in fact superfluous):
\def\theequation{8.10}\begin{equation}
\aligned
{\bf Y}_2^1
&
=
\mathcal{Y}_{xx}^1
+
\big[2\,\mathcal{Y}_{xy^1}^1-\mathcal{X}_{xx}\big]\,y_1^1
+
\big[2\,\mathcal{Y}_{xy^2}^1\big]\,y_1^2
+
\big[\mathcal{Y}_{y^1y^1}^1-2\,\mathcal{X}_{xy^1}\big]\,(y_1^1)^2
+
\\
&\ \ \
+
\big[2\,\mathcal{Y}_{y^1y^2}^1-2\,\mathcal{X}_{xy^2}\big]\,y_1^1\,y_1^2
+
\big[\mathcal{Y}_{y^2y^2}^1\big]\,(y_1^2)^2
+
\big[-\mathcal{X}_{y^1y^1}\big]\,(y_1^1)^3
+
\\
&\ \ \
+
\big[-2\,\mathcal{X}_{y^1y^2}\big]\,(y_1^1)^2\,y_1^2
+
\big[-\mathcal{X}_{y^2y^2}\big]\,y_1^1\,(y_1^2)^2
+
\big[\mathcal{Y}_{y^1}^1-2\,\mathcal{X}_x\big]\,y_2^1
+
\\
&\ \ \
+
\big[\mathcal{Y}_{y^2}^1\big]\,y_2^2
+
\big[-3\,\mathcal{X}_{y^1}\big]\,y_1^1\,y_2^1
+
\big[-\mathcal{X}_{y^2}\big]\,y_1^1\,y_2^2
+
\big[-2\,\mathcal{X}_{y^2}\big]\,y_1^2\,y_2^1,
\\
{\bf Y}_2^2
&
=
\mathcal{Y}_{xx}^2
+
\big[2\,\mathcal{Y}_{xy^1}^2\big]\,y_1^1
+
\big[2\,\mathcal{Y}_{xy^2}^2-\mathcal{X}_{xx}\big]\,y_1^2
+
\big[\mathcal{Y}_{y^1y^1}^2\big]\,(y_1^1)^2
+
\\
&\ \ \
+
\big[2\,\mathcal{Y}_{y^1y^2}^2-2\,\mathcal{X}_{xy^1}\big]\,y_1^1\,y_1^2
+
\big[\mathcal{Y}_{y^2y^2}^2-2\,\mathcal{X}_{xy^2}\big]\,(y_1^2)^2
+
\big[-\mathcal{X}_{y^1y^1}\big]\,(y_1^1)^2\,y_1^2
+
\\
&\ \ \
+
\big[-2\,\mathcal{X}_{y^1y^2}\big]\,y_1^1\,(y_1^2)^2
+
\big[-\mathcal{X}_{y^2y^2}\big]\,(y_1^2)^3
+
\big[\mathcal{Y}_{y^1}^2\big]\,y_2^1
+
\\
&\ \ \
+
\big[\mathcal{Y}_{y^2}^2-2\,\mathcal{X}_x\big]\,y_2^2
+
\big[-2\,\mathcal{X}_{y^1}\big]\,y_1^1\,y_2^2
+
\big[-\mathcal{X}_{y^1}\big]\,y_1^2\,y_2^1
+
\big[-3\,\mathcal{X}_{y^2}\big]\,y_1^2\,y_2^2.
\endaligned
\end{equation}
Inserting ${\bf Y}_1^2$ and ${\bf Y}_1^1$ in the first Lie
equation~\thetag{ 8.7} in which $y_1^2$ is replaced by the
value $\text{\rm (8.6)}_1$ 
it has on $\Delta_{ \mathcal{ E}_4}$ and still
computing ${\rm mod}\, (y_1^1)^5$, we get, again with ${\sf r}$, 
${\sf s}$ being unspecified functions of $(x, y^1, y^2)$ with 
${\sf s} (0) = 0$:
\def\theequation{8.11}\begin{equation}
\aligned
0
\equiv
&
-
\mathcal{Y}_x^2
+
\big[
-\mathcal{Y}_{y^1}^2
\big]\,y_1^1
+
\\
&
+
\big[
-\mathcal{Y}_{y^2}^2+\mathcal{X}_x
\big]\,
\big(
y_1^1(2x+{\bf g}^1)
+
(y_1^1)^2(1+{\bf g}^2)
+
(y_1^1)^3\,{\sf s}
+
(y_1^1)^4\,{\sf s}
\big)
+
\\
&
+
\big[
\mathcal{X}_{y^1}
\big]\,
\big(
(y_1^1)^2(2x+{\bf g}^1)
+
(y_1^1)^3(1+{\bf g}^2)
+
(y_1^1)^4\,{\sf s}
\big)
+
\\
&
+
\big[
\mathcal{X}_{y^2}
\big]\,
\big(
(y_1^1)^2[2x+{\bf g}^1]^2
+
(y_1^1)^3(4x+2{\bf g}^1)(1+{\bf g}^2)
+
(y_1^1)^4(1+{\sf s})
\big)
+
\\
&
+
\big[
\mathcal{X}
\big]\,
\big(
y_1^1(2+{\bf g}_x^1)
+
(y_1^1)^2\,{\bf g}_x^2
+
(y_1^1)^3\,{\sf r}
+
(y_1^1)^4\,{\sf r}
\big)
+
\\
&
+
\big[
\mathcal{Y}^1
\big]\,
\big(
y_1^1\,{\sf r}
+
(y_1^1)^2\,{\sf r}
+
(y_1^1)^3\,{\sf r}
+
(y_1^1)^4\,{\sf r}
\big)
+
\\
&
+
\big[
\mathcal{Y}^2
\big]\,
\big(
y_1^1\,{\sf r}
+
(y_1^1)^2\,{\sf r}
+
(y_1^1)^3\,{\sf r}
+
(y_1^1)^4\,{\sf r}
\big)
+
\\
&
+
\big[
\mathcal{Y}_x^1
\big]\,
\big(
2x+{\bf g}^1
+
y_1^1(2+2{\bf g}^2)
+
(y_1^1)^2\,{\sf s}
+
(y_1^1)^3\,{\sf s}
+
(y_1^1)^4\,{\sf s}
\big)
+
\\
&
+
\big[
\mathcal{Y}_{y^1}^1-\mathcal{X}_x
\big]\,
\big(
y_1^1(2x+{\bf g}^1)
+
(y_1^1)^2(2+2{\bf g}^2)
+
(y_1^1)^3\,{\sf s}
+
(y_1^1)^4\,{\sf s}
\big)
+
\\
&
+
\big[
\mathcal{Y}_{y^2}^1
\big]\,
\big(
y_1^1[2x+{\bf g}^1]^2
+
(y_1^1)^2(2x+{\bf g}^1)(3+3{\bf g}^2)
+
(y_1^1)^3(2+{\sf s})
+
(y_1^1)^4\,{\sf s}
\big)
+
\\
&
+
\big[
-\mathcal{X}_{y^1}
\big]\,
\big(
(y_1^1)^2(2x+{\bf g}^1)
+
(y_1^1)^3(2+2{\bf g}^2)
+
(y_1^1)^4\,{\sf s}
\big)
+
\\
&
+
\big[
-\mathcal{X}_{y^2}
\big]\,
\big(
(y_1^1)^2[2x+{\bf g}^1]^2
+
(y_1^1)^3(2x+{\bf g}^1)(3+3{\bf g}^2)
+
(y_1^1)^4(2+{\sf s})
\big).
\endaligned
\end{equation}
Collecting the coefficients of the monomials ${\rm cst.}$, $y_1^1$,
$(y_1^1)^2$, $(y_1^1)^3$, $(y_1^1)^4$, we get, after slight
simplification (in the coefficient of $(y_1^1)^2$, the term $(2x +
{\bf g}^1) \mathcal{ X}_x$ annihilates with its opposite; in the
coefficient of $(y_1^1)^3$, two pairs annihilate and then, we divide
by $[1 + {\bf g}^2]$) a system of five linear 
{\sc pde}'s:
\def\theequation{8.12}\begin{equation}
\aligned
0
&
=
-\mathcal{Y}_x^2
+
(2x+{\bf g}^1)\mathcal{Y}_x^1,
\\
0
&
=
-\mathcal{Y}_{y^1}^2
-
(2x+{\bf g}^1)\mathcal{Y}_{y^2}^2
+
(2+{\bf g}_x^1)\mathcal{X}
+
{\sf r}\,\mathcal{Y}^1
+
{\sf r}\,\mathcal{Y}^2
+
\\
&\ \ \ \ \ \
+
(2+2{\bf g}^2)\mathcal{Y}_x^1
+
(2x+{\bf g}^1)\mathcal{Y}_{y^1}^1
+
[2x+{\bf g}^1]^2\mathcal{Y}_{y^2}^1,
\\
0
&
=
-\mathcal{Y}_{y^2}^2
+
\mathcal{X}_x
+
{\bf g}_x^2[1+{\bf g}^2]^{-1}\mathcal{X}
+
{\sf r}\,\mathcal{Y}^1
+
{\sf r}\,\mathcal{Y}^2
+
\\
&\ \ \ \ \ \
+
{\sf s}\,\mathcal{Y}_x^1
+
2\,\mathcal{Y}_{y^1}^1
-
2\,\mathcal{X}_x
+
(6x+3{\bf g}^2)\mathcal{Y}_{y^2}^1,
\\
0
&
=
{\sf s}\,\mathcal{Y}_{y^2}^2
+
{\sf s}\,\mathcal{X}_x
+
(1+{\bf g}^2)\mathcal{X}_{y^1}
+
(2x+{\bf g}^1)(2+2{\bf g}^2)\mathcal{X}_{y^2}
+
\\
&\ \ \ \ \ \
+
{\sf r}\,\mathcal{X}
+
{\sf r}\,\mathcal{Y}^1
+
{\sf r}\,\mathcal{Y}^2
+
{\sf s}\,\mathcal{Y}_x^1
+
{\sf s}\,\mathcal{Y}_{y^1}^1
+
{\sf s}\,\mathcal{X}_x
+
(2+{\sf s})\,\mathcal{Y}_{y^2}^1
-
\\
&\ \ \ \ \ \ 
-
(2+2{\bf g}^2)\mathcal{X}_{y^1}
-
(2x+{\bf g}^1)(3+3{\bf g}^2)\mathcal{X}_{y^2},
\\
0
&
=
{\sf s}\,\mathcal{Y}_{y^2}^2
+
{\sf s}\,\mathcal{X}_x
+
{\sf s}\,\mathcal{X}_{y^1}
+
(1+{\sf s})\mathcal{X}_{y^2}
+
{\sf r}\,\mathcal{X}
+
{\sf r}\,\mathcal{Y}^1
+
{\sf r}\,\mathcal{Y}^2
+
{\sf s}\,\mathcal{Y}_x^1
+
{\sf s}\,\mathcal{Y}_{y^1}^1
+
\\
&\ \ \ \ \ \ 
+
{\sf s}\,\mathcal{X}_x
+
{\sf s}\,\mathcal{Y}_{y^2}^1
+
{\sf s}\,\mathcal{X}_{y^1}
-
(2+{\sf s})\mathcal{X}_{y^2}.
\endaligned
\end{equation}
We then simplify the remainders using ${\sf s} + {\sf s} = {\sf s}$,
${\sf r} + {\sf s} = {\sf r}$ and ${\sf r} + {\sf r} = {\sf r}$; we
divide $\text{\rm (8.12)}_5$ by $(1+ {\sf s})$; we replace $\mathcal{
X}_{ y^2}$ obtained from $\text{\rm (8.12)}_5$ in $\text{\rm
(8.12)}_4$; we divide $\text{\rm (8.12)}_4$ by $(1+ {\bf g}^2)$; we
then solve $\mathcal{ X}_{ y^1}$ from $\text{\rm (8.12)}_4$ and
finally we insert it in $\text{\rm (8.12)}_5$; we get:
\def\theequation{8.13}\begin{equation}
\aligned
0
&
=
-\mathcal{Y}_x^2
+
(2x+{\bf g}^1)\mathcal{Y}_x^1,
\\
0
&
=
-\mathcal{Y}_{y^1}^2
-
(2x+{\bf g}^1)\mathcal{Y}_{y^2}^2
+
(2+{\bf g}_x^1)\mathcal{X}
+
(2+2{\bf g}^2)\mathcal{Y}_x^1
+
(2x+{\bf g}^1)\mathcal{Y}_{y^1}^1
+
\\
&\ \ \ \ \ \
+
[2x+{\bf g}^1]^2\mathcal{Y}_{y^2}^1
+
{\sf r}\,\mathcal{Y}^1
+
{\sf r}\,\mathcal{Y}^2,
\\
0
&
=
-\mathcal{Y}_{y^2}^2
-
\mathcal{X}_x
+
2\,\mathcal{Y}_{y^1}^1
+
(6x+3{\bf g}^2)\mathcal{Y}_{y^2}^1
+
{\sf r}\,\mathcal{Y}^1
+
{\sf r}\,\mathcal{Y}^2
+
{\sf s}\,\mathcal{Y}_x^1
+
\\
&\ \ \ \ \ \ \ \ \ \ \ \ \ \ \ \ \ \ \ \ \ \ \ \
\ \ \ \ \ \ \ \ \ \ \ \ \ \
+
{\bf g}_x^2[1+{\bf g}^2]^{-1}\mathcal{X},
\\
0
&
=
-\mathcal{X}_{y^1}
+
(2+{\sf s})\mathcal{Y}_{y^2}^1
+
{\sf r}\,\mathcal{X}
+
{\sf r}\,\mathcal{Y}^1
+
{\sf r}\,\mathcal{Y}^2
+
{\sf s}\,\mathcal{X}_x
+
{\sf s}\,\mathcal{Y}_x^1
+
{\sf s}\,\mathcal{Y}_{y^1}^1
+
{\sf s}\,\mathcal{Y}_{y^2}^2,
\\
0
&
=
-\mathcal{X}_{y^2}
+
{\sf r}\,\mathcal{X}
+
{\sf r}\,\mathcal{Y}^1
+
{\sf r}\,\mathcal{Y}^2
+
{\sf s}\,\mathcal{X}_x
+
{\sf s}\,\mathcal{Y}_x^1
+
{\sf s}\,\mathcal{Y}_{y^1}^1
+
{\sf s}\,\mathcal{Y}_{y^2}^1
+
{\sf s}\,\mathcal{Y}_{y^2}^2.
\endaligned
\end{equation}
Similarly, developing the second equation~\thetag{ 8.8}
and computing ${\rm mod}\, (y_1^1)^2$, we get:
\def\theequation{8.14}\begin{equation}
\aligned
0
&
\equiv
-
\mathcal{Y}_{xx}^1
+
\big[
-2\,\mathcal{Y}_{xy^1}^1
+
\mathcal{X}_{xx}
\big]\,y_1^1
+
\big[
-(4x+2{\bf g}^1)\mathcal{Y}_{xy^2}
-
(2+{\bf h})\mathcal{Y}_{y^2}^1
\big]
+
\\
&\ \ \ \ \ \ \ \ \ \ \ \ \ \ \
+
\big[
2{\bf h}\,\mathcal{Y}_x^1
\big]\,y_1^1.
\endaligned
\end{equation}
Collecting the coefficients of the monomials ${\rm cst.}$, $y_1^1$, we
get two more linear {\sc pde}'s:
\def\theequation{8.15}\begin{equation}
\aligned
0
&
=
-\mathcal{Y}_{xx}^1,
\\
0
&
=
-2\,\mathcal{Y}_{xy^1}^1
+
\mathcal{X}_{xx}
-
(4x+2{\bf g}^1)\mathcal{Y}_{xy^2}^1
-
(2+{\bf h})\mathcal{Y}_{y^2}^1
+
2{\bf h}\,\mathcal{Y}_x^1.
\endaligned
\end{equation}

\def\theproposition{8.16}\begin{proposition} 
Setting as initial conditions the five specific differential
coefficients
\def\theequation{8.17}\begin{equation}
{\sf P}
:=
{\sf P}
\big(
\mathcal{X},\mathcal{Y}^1,\mathcal{Y}^2,\mathcal{Y}_x^1,
\mathcal{X}_x
\big)
=
{\sf r}\,\mathcal{X}
+
{\sf r}\,\mathcal{Y}^1
+
{\sf r}\,\mathcal{Y}^2
+
{\sf r}\,\mathcal{Y}_x^1
+
{\sf r}\,\mathcal{X}_x,
\end{equation}
it follows by cross differentiations and by linear substitutions from
the seven equations $\text{\rm (8.13)}_i$, $i=1, 2, 3, 4, 5$,
$\text{\rm (8.15)}_j$, $j=1, 2$, that $\mathcal{ X}_{ y^1}$,
$\mathcal{ X}_{ y^2}$, $\mathcal{ Y}_{ y^1}^1$, $\mathcal{ Y}_{
y^2}^1$, $\mathcal{ Y}_x^2$, $\mathcal{ Y}_{ y^1}^2$, $\mathcal{ Y}_{
y^2}^2$ and $\mathcal{ X}_{ xx}$, $\mathcal{ X}_{ xy^1}$, $\mathcal{
X}_{ xy^2}$, $\mathcal{ Y}_{ xx}^1$, $\mathcal{ Y}_{ xy^1}^1$,
$\mathcal{ Y}_{ xy^2}^1$ are uniquely determined as
linear combinations of $(\mathcal{ X}, \mathcal{ Y}^1, 
\mathcal{ Y}^2, \mathcal{ Y}_x^1, \mathcal{ X}_x)$, namely{\rm :}
\def\theequation{8.18}\begin{equation}
\left\{
\aligned
&
\ \ \ \ \ \ \ \ \ \ \ \ \ \ \ \ \ \ \ \ \
\ \ \ \ \ \ \ \ \ \ \ \ \ \ \ \ \ \ \ \ \
\mathcal{Y}_x^2
\overset{1}{=}
{\sf P},\ \ \ \ \ \
\mathcal{X}_{xx}
\overset{2}{=}
{\sf P},\ \ \ \ \ \
\mathcal{Y}_{xx}^1
\overset{3}{=}
{\sf P},
\\
&
\mathcal{X}_{y^1}
\overset{4}{=}
{\sf P},\ \ \ \ \
\mathcal{Y}_{y^1}^1
\overset{5}{=}
{\sf P},\ \ \ \ \
\mathcal{Y}_{y^1}^2
\overset{6}{=}
{\sf P},\ \ \ \ \
\mathcal{X}_{xy^1}
\overset{7}{=}
{\sf P},\ \ \ \ \
\mathcal{Y}_{xy^1}^1
\overset{8}{=}
{\sf P},
\\
&
\mathcal{X}_{y^2}
\overset{9}{=}
{\sf P},\ \ \ \ \
\mathcal{Y}_{y^2}^1
\overset{10}{=}
{\sf P},\ \ \ \ \
\mathcal{Y}_{y^2}^2
\overset{11}{=}
{\sf P},\ \ \ \ \
\mathcal{X}_{xy^2}
\overset{12}{=}
{\sf P},\ \ \ \ \
\mathcal{Y}_{xy^2}^1
\overset{13}{=}
{\sf P}.
\endaligned\right.
\end{equation}
\end{proposition}

Then the expressions ${\sf P}$ are stable under 
differentiation{\rm :}
\def\theequation{8.19}\begin{equation}
\aligned
{\sf P}_x
&
=
{\sf P}
+
{\sf r}\,\mathcal{Y}_x^2
+
{\sf r}\,\mathcal{Y}_{xx}^1
+
{\sf r}\,\mathcal{X}_{xx}
=
{\sf P},
\\
{\sf P}_{y^1}
&
=
{\sf P}
+
{\sf r}\,\mathcal{X}_{y^1}
+
{\sf r}\,\mathcal{Y}_{y^1}^1
+
{\sf r}\,\mathcal{Y}_{y^1}^2
+
{\sf r}\,\mathcal{Y}_{xy^1}^1
+
{\sf r}\,\mathcal{X}_{xy^1}
=
{\sf P},
\\
{\sf P}_{y^2}
&
=
{\sf P}
+
{\sf r}\,\mathcal{X}_{y^2}
+
{\sf r}\,\mathcal{Y}_{y^2}^1
+
{\sf r}\,\mathcal{Y}_{y^2}^2
+
{\sf r}\,\mathcal{Y}_{xy^2}^1
+
{\sf r}\,\mathcal{X}_{xy^2}
=
{\sf P},
\endaligned
\end{equation}
and moreover, all other, higher order partial derivatives of
$\mathcal{ X}$, of $\mathcal{ Y}^1$ and of $\mathcal{ Y}^2$ may be
expressed as ${\sf P} \big( \mathcal{ X}, \mathcal{ Y}^1, \mathcal{
Y}^2, \mathcal{ Y}_x^1, \mathcal{ X}_x \big)$.

\def\thecorollary{8.20}\begin{corollary}
An infinitesimal Lie symmetry of
\thetag{ $\mathcal{ E }_4$} is uniquely determined
by the five initial Taylor coefficients
\def\theequation{8.21}\begin{equation}
\mathcal{X}(0),\
\mathcal{Y}^1(0),\
\mathcal{Y}^2(0),\
\mathcal{Y}_x^1(0),\
\mathcal{X}_x(0).
\end{equation}
\end{corollary}

\proof[Proof of the proposition]
We notice that $\text{\rm (8.18)}_1$ and 
$\text{\rm (8.18)}_3$ are given for free by
$\text{\rm (8.13)}_1$ and by $\text{\rm (8.15)}_1$. 
Differentiating $\text{\rm (8.13)}_3$ with respect to 
$x$, we get:
\def\theequation{8.22}\begin{equation}
\aligned
0
&
=
-
\mathcal{Y}_{xy^2}
-
\mathcal{X}_{xx}
+
2\,\mathcal{Y}_{xy^1}^1
+
(6+3{\bf g}_x^1)\mathcal{Y}_{y^2}^1
+
(6x+3{\bf g}^1)\mathcal{Y}_{xy^2}^1
+
{\sf r}\,\mathcal{Y}^1
+
\\
&\ \ \ \ \ \
+
{\sf r}\,\mathcal{Y}_x^1
+
{\sf r}\,\mathcal{Y}^2
+
{\sf r}\,\mathcal{Y}_x^2
+
{\sf r}\,\mathcal{Y}_x^1
+
{\sf s}\,\mathcal{Y}_{xx}^1
+
{\sf r}\,\mathcal{X}
+
{\bf g}_x^2[1+{\bf g}^2]^{-1}\mathcal{X}_x.
\endaligned
\end{equation}
By $\text{\rm (8.15)}_1$, ${\sf s}\, \mathcal{ Y}_{ xx}^1$
vanishes. We replace $\mathcal{ Y}_x^2$ thanks
to $\text{\rm (8.13)}_1$. 
Differentiating $\text{\rm (8.13)}_1$
with respect to $y^2$, we may substract
$0 = -\mathcal{ Y}_{ xy^2} + 
(2 x + {\bf g}^1) \mathcal{ Y}_{ xy^2}^1 + 
{\sf r}\, \mathcal{ Y}_x^1$. We get:
\def\theequation{8.23}\begin{equation}
\aligned
0
&
=
-\mathcal{X}_{xx}
+
2\,\mathcal{Y}_{xy^1}^1
+
(4x+2{\bf g}^1)\mathcal{Y}_{xy^2}^1
+
\\
&\ \ \ \ \ \
+
(6+3{\bf g}_x^1)\mathcal{Y}_{y^2}^1
+
{\sf r}\,\mathcal{X}
+
{\sf r}\,\mathcal{Y}^1
+
{\sf r}\,\mathcal{Y}^2
+
{\sf r}\,\mathcal{Y}_x^1
+
{\bf g}_x^2[1+{\bf g}^2]^{-1}\mathcal{X}_x.
\endaligned
\end{equation}
By means of $\text{\rm (8.15)}_2$, we replace the 
first three terms and then solve $\mathcal{ Y}_{ y^2}^1$:
\def\theequation{8.24}\begin{equation}
\mathcal{Y}_{y^2}^1
=
{\sf r}\,\mathcal{X}
+
{\sf r}\,\mathcal{Y}^1
+
{\sf r}\,\mathcal{Y}^2
+
{\sf r}\,\mathcal{Y}_x^1
+
{\bf k}^*\,\mathcal{X}_x,
\end{equation}
introducing a notation for a new function
that should be recorded:
\def\theequation{8.25}\begin{equation}
{\bf k}^*
:=
{\bf g}_x^2[1+{\bf g}^2]^{-1}
[4+3{\bf g}_x^1-{\bf h}]^{-1}.
\end{equation}
This is $\text{\rm (8.18)}_{ 10}$.
Next, we differentiate the obtained equation with respect to $x$,
getting:
\def\theequation{8.26}\begin{equation}
\mathcal{Y}_{xy^2}^1
=
{\sf r}\,\mathcal{X}
+
{\sf r}\,\mathcal{Y}^1
+
{\sf r}\,\mathcal{Y}^2
+
{\sf r}\,\mathcal{Y}_x^1
+
{\sf r}\,\mathcal{X}_x
+
{\bf k}^*\,\mathcal{X}_{xx}.
\end{equation}
This is $\text{\rm (8.18)}_{ 13}$.
We replace the obtained value of $\mathcal{ Y}_{ y^2}^1$ in $\text{\rm
(8.13)}_2$, $\text{\rm (8.13)}_3$, $\text{\rm (8.15)}_2$ and the
obtained value of $\mathcal{ Y}_{ x y^2}^1$ in $\text{\rm
(8.15)}_2$. This yields a new, simpler system of seven equations:
\def\theequation{8.27}\begin{equation}
\aligned
0
&
=
-\mathcal{Y}_x^2
+
(2x+{\bf g}^1)\mathcal{Y}_x^1,
\\
0
&
=
-\mathcal{Y}_{y^1}^2
-
(2x+{\bf g}^1)\mathcal{Y}_{y^2}^2
+
(2+{\bf g}_x^1)\mathcal{X}
+
(2+2{\bf g}^2)\mathcal{Y}_x^1
+
(2x+{\bf g}^1)\mathcal{Y}_{y^1}^1
+
\\
&\ \ \ \ \
+
{\sf s}\,\mathcal{X}
+
{\sf r}\,\mathcal{Y}^1
+
{\sf r}\,\mathcal{Y}^2
+
{\sf s}\,\mathcal{Y}_x^1
+
{\sf k}^*[2x+{\bf g}^1]^2\mathcal{X}_x,
\\
0
&
=
-\mathcal{Y}_{y^2}^2
-
\mathcal{X}_x
+
2\,\mathcal{Y}_{y^1}^1
+
{\sf s}\,\mathcal{X}
+
{\sf r}\,\mathcal{Y}^1
+
{\sf r}\,\mathcal{Y}^2
+
{\sf s}\,\mathcal{Y}_x^1
+
{\sf k}^*(6x+3{\bf g}^1)\mathcal{X}_x,
\\
0
&
=
-\mathcal{X}_{y^1}
+
{\sf r}\,\mathcal{X}
+
{\sf r}\,\mathcal{Y}^1
+
{\sf r}\,\mathcal{Y}^2
+
{\sf s}\,\mathcal{Y}_x^1
+
{\sf s}\,\mathcal{X}_x
+
{\sf s}\,\mathcal{Y}_{y^1}^1
+
{\sf s}\,\mathcal{Y}_{y^2}^2,
\\
0
&
=
-\mathcal{X}_{y^2}
+
{\sf r}\,\mathcal{X}
+
{\sf r}\,\mathcal{Y}^1
+
{\sf r}\,\mathcal{Y}^2
+
{\sf s}\,\mathcal{Y}_x^1
+
{\sf s}\,\mathcal{X}_x
+
{\sf s}\,\mathcal{Y}_{y^1}^1
+
{\sf s}\,\mathcal{Y}_{y^2}^2,
\\
0
&
=
-\mathcal{Y}_{xx}^1,
\\
0
&
=
-2\,\mathcal{Y}_{xy^1}^1
+
\mathcal{X}_{xx}
\big(
1
-
{\bf k}^*(4x+2{\bf g}^1)
\big)
+
{\sf r}\,\mathcal{X}
+
{\sf r}\,\mathcal{Y}^1
+
{\sf r}\,\mathcal{Y}^2
+
{\sf r}\,\mathcal{Y}_x^1.
+
{\sf s}\,\mathcal{X}_x.
\endaligned
\end{equation}
Restarting from this system, we differentiate $\text{\rm (8.27)}_3$
with respect to $x$:
\def\theequation{8.28}\begin{equation}
\aligned
0
&
=
-\mathcal{Y}_{xy^2}^2
-
\mathcal{X}_{xx}
+
2\,\mathcal{Y}_{xy^1}^1
+
{\sf r}\,\mathcal{X}
+
{\sf r}\,\mathcal{Y}^1
+
{\sf r}\,\mathcal{Y}^2
+
\\
&\ \ \ \ \ \
+
{\sf r}\,\mathcal{X}_x
+
{\sf r}\,\mathcal{Y}_x^1
+
{\sf r}\,\mathcal{Y}_x^2
+
{\sf s}\,\mathcal{Y}_{xx}^1
+
{\sf k}^*(6x+3{\bf g}^1)\mathcal{X}_{xx}.
\endaligned
\end{equation}
We replace $\mathcal{ Y}_x^2$, we erase $\mathcal{ Y}_{ xx}^1$
and we add $\text{\rm (8.27)}_7$:
\def\theequation{8.29}\begin{equation}
0
=
-\mathcal{Y}_{xy^2}^2
+
{\sf k}^*(2x+{\bf g}^1)\mathcal{X}_{xx}
+
{\sf r}\,\mathcal{X}
+
{\sf r}\,\mathcal{Y}^1
+
{\sf r}\,\mathcal{Y}^2
+
{\sf r}\,\mathcal{Y}_x^1
+
{\sf r}\,\mathcal{X}_x.
\end{equation}
We differentiate $\text{\rm (8.27)}_2$ with respect
to $x$:
\def\theequation{8.30}\begin{equation}
\aligned
0
&
=
-\mathcal{Y}_{xy^1}^2
-
(2+{\bf g}_x^1)\mathcal{Y}_{y^2}^2
-
(2x+{\bf g}^1)\mathcal{Y}_{xy^2}^2
+
{\sf r}\,\mathcal{X}
+
(2+{\bf g}_x^1)\mathcal{X}_x
+
\\
&\ \ \ \ \ \
+
{\sf s}\,\mathcal{Y}_x^1
+
(2+2{\bf g}^2)\mathcal{Y}_{xx}^1
+
(2+{\bf g}_x^1)\mathcal{Y}_{y^1}^1
+
(2x+{\bf g}^1)\mathcal{Y}_{xy^1}^1
+
{\sf r}\,\mathcal{X}_x
+
\\
&\ \ \ \ \ \
+
{\sf s}\,\mathcal{X}_x
+
{\sf r}\,\mathcal{Y}^1
+
{\sf r}\,\mathcal{Y}_x^1
+
{\sf r}\,\mathcal{Y}^2
+
{\sf r}\,\mathcal{Y}_x^2
+
{\sf r}\,\mathcal{Y}_x^1
+
{\sf s}\,\mathcal{Y}_{xx}^1
+
{\bf k}^*[2x+{\bf g}^1]^2\mathcal{X}_{xx}.
\endaligned
\end{equation}
Differentiating $\text{\rm (8.27)}_1$ with respect 
to $y^1$, we may substract 
$0 = - \mathcal{ Y}_{ xy^1}^2 
+ (2x + {\bf g}^1) \mathcal{ Y}_{ x y^1}^1 + 
{\sf r}\, \mathcal{ Y}_x^1$;
we replace $\mathcal{ Y}_x^2$ and
erase $\mathcal{ Y}_{ xx}^1$; we substract
\thetag{ 8.29} multiplied by
$(2x + {\bf g}^1)$; we get:
\def\theequation{8.31}\begin{equation}
0
=
-\mathcal{Y}_{y^2}^2
+
(1+{\sf s})\mathcal{X}_x
+
\mathcal{Y}_{y^1}^1
+
{\sf r}\,\mathcal{X}
+
{\sf r}\,\mathcal{Y}^1
+
{\sf r}\,\mathcal{Y}^2
+
{\sf r}\,\mathcal{Y}_x^1.
\end{equation}
Comparing with $\text{\rm (8.27)}_3$ yields:
\def\theequation{8.32}\begin{equation}
\aligned
\mathcal{Y}_{y^1}^1
&
=
(2+{\sf s})\mathcal{X}_x
+
{\sf r}\,\mathcal{X}
+
{\sf r}\,\mathcal{Y}^1
+
{\sf r}\,\mathcal{Y}^2
+
{\sf r}\,\mathcal{Y}_x^1,
\\
\mathcal{Y}_{y^2}^2
&
=
(3+{\sf s})\mathcal{X}_x
+
{\sf r}\,\mathcal{X}
+
{\sf r}\,\mathcal{Y}^1
+
{\sf r}\,\mathcal{Y}^2
+
{\sf r}\,\mathcal{Y}_x^1.
\endaligned
\end{equation}
These are $\text{\rm (8.18)}_5$
and $\text{\rm (8.18)}_{ 11}$.
Differentiating these two equations
with respect to $x$, replacing
$\mathcal{ Y}_x^2$ and erasing
$\mathcal{ Y}_{ xx}^1$, we get:
\def\theequation{8.33}\begin{equation}
\aligned
\mathcal{Y}_{xy^1}^1
&
=
(2+{\sf s})\mathcal{X}_{xx}
+
{\sf r}\,\mathcal{X}
+
{\sf r}\,\mathcal{Y}^1
+
{\sf r}\,\mathcal{Y}^2
+
{\sf r}\,\mathcal{Y}_x^1
+
{\sf r}\,\mathcal{X}_x,
\\
\mathcal{Y}_{xy^2}^2
&
=
(3+{\sf s})\mathcal{X}_{xx}
+
{\sf r}\,\mathcal{X}
+
{\sf r}\,\mathcal{Y}^1
+
{\sf r}\,\mathcal{Y}^2
+
{\sf r}\,\mathcal{Y}_x^1
+
{\sf r}\,\mathcal{X}_x.
\endaligned
\end{equation}
We then replace this value of $\mathcal{ Y}_{ xy^2}^2$ in~\thetag{
8.29} and solve $\mathcal{ X}_{ xx}$: this yields $\text{\rm
(8.18)}_2$.

To conclude, we replace $\mathcal{ X}_{ xx}$ so obtained in $\text{\rm
(8.27)}_7$: this yields $\text{\rm (8.18)}_8$. We replace $\mathcal{
Y}_{ y^1}^1$ and $\mathcal{ Y}_{ y^2}^2$ from~\thetag{ 8.32} in
$\text{\rm (8.27)}_4$ and in $\text{\rm (8.27)}_5$: this yields
$\text{\rm (8.18)}_4$ and this yields $\text{\rm (8.18)}_9$.
Thanks to $\text{\rm (8.18)}_2$ (got)
we observe that
\def\theequation{8.34}\begin{equation}
{\sf P}_x
=
{\sf P}
+
{\sf r}\,\mathcal{Y}_{xx}^1
+
{\sf r}\,\mathcal{Y}_x^2
+
{\sf r}\,\mathcal{X}_{xx}
=
{\sf P}.
\end{equation}
Differentiating $\text{\rm (8.18)}_4$ (got) and $\text{\rm (8.18)}_9$
(got) with respect to $x$ then yields $\text{\rm (8.18)}_7$ and
$\text{\rm (8.18)}_{ 12 }$. We replace $\mathcal{ Y}_{ y^1}^1$ and
$\mathcal{ Y}_{ y^2}^2$ from $\text{\rm (8.18)}_5$ (got) and
$\text{\rm (8.18)}_{ 11 }$ (got) in $\text{\rm (8.27)}_2$: this yields
$\text{\rm (8.18)}_6$. Finally, to obtain the very last $\text{\rm
(8.18)}_{ 13}$, we differentiate $\text{\rm (8.18)}_{ 10}$ (got) with
respect to $x$.


The proof of Proposition~8.16 is complete.
\endproof

We claim that the bound $\dim \mathfrak{ SYM} (\mathcal{ E}_4)
\leqslant 5$ is attained for the model~\thetag{ 8.4}. Indeed, with $0
= {\sf r} = {\sf s}$ and $0 = {\bf g}^1 = {\bf g}^2 = {\bf h}$ (whence
${\bf k}^* = 0$) \thetag{ 8.24} is $\mathcal{ Y}_{ y^2}^1 = 0$ and
then the seven equations \thetag{ 8.27} are:
\def\theequation{8.35}\begin{equation}
\left\{
\aligned
0
&
=
-\mathcal{Y}_x^2
+
2x\,\mathcal{Y}_x^1,
\\
0
&
=
-\mathcal{Y}_{y^1}^2
-
2x\,\mathcal{Y}_{y^2}^2
+
2\,\mathcal{X}
+
2\,\mathcal{Y}_x^1
+
2x\,\mathcal{Y}_{y^1}^1,
\\
0
&
=
-\mathcal{Y}_{y^2}^2
-
\mathcal{X}_x
+
2\,\mathcal{Y}_{y^1}^1,
\\
0
&
=
-\mathcal{X}_{y^1},
\\
0
&
=
-\mathcal{X}_{y^2},
\\
0
&
=
-\mathcal{Y}_{xx}^1,
\\
0
&
=
-2\,\mathcal{Y}_{xy^1}^1
+
\mathcal{X}_{xx},
\endaligned\right.
\end{equation}
having the general solution
\def\theequation{8.36}\begin{equation}
\left\{
\aligned
\mathcal{X}
&
=
a-d+e\,x,
\\
\mathcal{Y}^1
&
=
b+d\,x+2e\,y^1,
\\
\mathcal{Y}^2
&
=
c+2a\,y^1+3e\,y^2+d\,xx.
\endaligned\right.
\end{equation}
depending on five parameters $a, b, c, d, e \in \K$.
Five generators of $\mathfrak{ SYM} (\mathcal{ E}_4)$ 
are:
\def\theequation{8.37}\begin{equation}
\left\{
\aligned
\mathcal{D}
&
:=
x\,\partial_x
+
2y^1\,\partial_{y^1}
+
3y^2\,\partial_{y^2},
\\
\mathcal{L}_1
&
:=
-\partial_x
+
x\,\partial_{y^1}
+
xx\,\partial_{y^2},
\\
\mathcal{L}_1'
&
:=
\partial_x
+
2y^1\,\partial_{y^2},
\\
\mathcal{L}_2
&
:=
\partial_{y^1},
\\
\mathcal{L}_3
&
:=
\partial_{y^2}.
\endaligned\right.
\end{equation}

\noindent
The commutator table

\medskip
\begin{center}
\begin{tabular}[t]{ l || l | l | l | l | l }
& $\mathcal{D}$ & $\mathcal{L}_1$ & $\mathcal{L}_1'$ 
& $\mathcal{L}_2$ & $\mathcal{L}_3$
\\
\hline 
\hline
$\mathcal{D}$ & 0 & $-\mathcal{L}_1$ & $-\mathcal{L}_1'$ & 
$-2\,\mathcal{L}_2$ & $-3\,\mathcal{L}_3$
\\
\hline
$\mathcal{L}_1$ & $\mathcal{L}_1$ & 0 &
$-\mathcal{L}_2$ & 0 & 0
\\
\hline
$\mathcal{L}_1'$ & $\mathcal{L}_1'$ & $\mathcal{L}_2$ & 0 & 
$-2\,\mathcal{L}_3$ & 0
\\
\hline
$\mathcal{L}_2$ & $2\,\mathcal{L}_2$ & 0 & 
$2\,\mathcal{L}_3$ & 0 & 0
\\
\hline
$\mathcal{L}_3$ & $3\,\mathcal{L}_3$ & 0 & 0 & 0 & 0
\end{tabular}

\medskip

{\bf Table~3.}

\end{center}

\noindent
shows that the subalgebra spanned by $\mathcal{ L}_1$, $\mathcal{
L}_1'$, $\mathcal{ L}_2$, $\mathcal{ L}_3$ is isomorphic to the unique
irreducible 4-dimensional nilpotent Lie algebra $\mathfrak{ n}_4^1$
(\cite{ ov1994, bes2005}). Then $\mathfrak{ SYM} (\mathcal{ E}_4)$ is
a semidirect product of $\K$ with $\mathfrak{ n}_4^1$. The author
ignores whether it is rigid. The following accessible research will
be pursued in a subsequent publication.

\def\theopenproblem{8.38}\begin{openproblem}
Classify systems \thetag{ $\mathcal{ E}_4$} up to point
transformations. Deduce a complete classification, up to local
biholomorphisms, of all real analytic generic submanifolds of
codimension 2 in $\C^3$, valid at a Zariski-generic point.
\end{openproblem}

\subsection*{8.39.~Almost everywhere rigid hypersurfaces}
When studying and classifying differential objects, it is essentially
no restriction to assume their Lie symmetry groups to be of dimension
$\geqslant 1$, the study of objects having no infinitesimal symmetries
being an independent field of research. In particular, if $M \subset
\C^{ n +1}$ ($n \geqslant 1$) is a connected real analytic
hypersurface, we may suppose that $\dim \mathfrak{ hol} (M) \geqslant
1$, at least. So let $\mathcal{ L}$ be a nonzero holomorphic vector
field with $\mathcal{ L} + \overline{ \mathcal{ L}}$ tangent to $M$.

\def\thelemma{8.40}\begin{lemma}
{\rm (\cite{ ca1932a, st1996, ber1999})} If in addition $M$ is finitely
nondegenerate, then
\def\theequation{8.41}\begin{equation}
\Sigma
:=
\big\{
p\in M:\
\mathcal{L}(p)\in T_p^cM
\big\}
\end{equation}
is a proper real analytic subset of $M$.
\end{lemma}

In other words, at every point $p$ belonging to the Zariski-dense
subset $M \backslash \Sigma$, the real nonzero vector $\mathcal{ L}
(p) + \overline{ \mathcal{ L} (p)} \in T_p M$ supplements $T_p^c
M$. Straightening $\mathcal{ L}$ in a neighborhood of $p$, there exist
local coordinates $t = (z_1, \dots, z_n, w)$ with $T_0^c M = \{ w = 0
\}$, $T_0 M = \{ {\rm Im} \, w = 0 \}$, whence $M$ is given by ${\rm
Im}\, w = h (z, \bar z, {\rm Re}\, w)$, and with $\mathcal{ L} =
\frac{ \partial }{ \partial w}$. The tangency of $\frac{ \partial }{
\partial w} + \frac{ \partial }{ \partial \bar w} = \frac{ \partial }{
\partial u}$ to $M$ entails that $h$ is indendepent of $u$. Then the
complex equation of $M$ is of the precise
form
\def\theequation{8.42}\begin{equation}
w
=
\bar w
+
i\,\overline{\Theta}(z,\bar z),
\end{equation}
with $\overline{ \Theta} = 2 h$ simply. The reality of $h$ reads
$\overline{ \Theta} (z, \bar z) \equiv \Theta (\bar z, z)$.

\def\thedefinition{8.43}\begin{definition}{\rm
A real analytic hypersurface $M \subset \C^{ n + 1}$ is called {\sl
rigid} at one of its points $p$ if there exists $\mathcal{ L} \in
\mathfrak{ hol} (M)$ with
\def\theequation{8.44}\begin{equation}
T_pM
=
T_p^cM
\oplus
\R
\big(
\mathcal{L}(p)
+
\overline{\mathcal{L}}(p)
\big).
\end{equation}
}\end{definition}

Similar elementary facts hold for general submanifolds of solutions.

\def\thelemma{8.45}\begin{lemma} 
With $n\geqslant 1$ and $m = 1$, let
$\mathcal{ M}$ be a {\rm (}connected{\rm )} submanifold of solutions
that is solvable with respect to the parameters. If there exists a
nonzero $\mathcal{ L} + \mathcal{ L}^* \in \mathfrak{ SYM} (\mathcal{
M})$, then at Zariski-generic points $p \in \mathcal{ M}$, we have
$\mathcal{ L} (p) \not \in {\sf F}_{ \sf v} (p)$ and there exist local
coordinates centered at $p$ in which $\mathcal{ L} = \frac{ \partial
}{ \partial y}$, $\mathcal{ L}^* = \frac{ \partial }{ \partial b}$,
whence $\mathcal{ M}$ has equation of the form
\def\theequation{8.46}\begin{equation}
y
=
b
+
\Pi(x,a),
\end{equation}
with $\Pi$ independent of $b$.
\end{lemma}

The associated system \thetag{ $\mathcal{ E}_\mathcal{ M}$} has then
equations $F_\alpha$ that are all independent of $y$.

\subsection*{8.47.~Study of the Lie symmetries of \thetag{ 
$\mathcal{ E}_5$}} In Example~1.28, it is thus essentially no
restriction to assume the hypersurface $M \subset \C^3$ to be rigid.

\def\thetheorem{8.48}\begin{theorem}
{\rm (\cite{ gm2003b, fk2005a, fk2005b})} The model hypersurface $M_0$
of equation
\def\theequation{8.49}\begin{equation}
w
=
\bar w
+ 
i\,
\frac{
2\,z^1\bar z^1
+
z^1z^1\bar z^2
+
\bar z^1\bar z^1z^2
}{1-z^2\bar z^2}
\end{equation}
has transitive Lie symmetry algebra $\mathfrak{ hol} (M_0)$ isomorphic
to $\mathfrak{ so} (3, 2)$ and is locally biholomorphic to a
neighborhood of every geometrically smooth point of the tube
\def\theequation{8.50}\begin{equation}
({\rm Re}\,w')^2
=
({\rm Re}\,z_1')^2
+
({\rm Re}\,z_1')^3
\end{equation}
over the standard cone of $\R^3$. Both are Levi-degenerate with Levi
form of rank 1 at every point and are 2-nondegenerate. The associated
{\sc pde} system \thetag{ $\mathcal{ E}_{ \mathcal{ M}_0}$}
\def\theequation{8.51}\begin{equation}
y_{x^2}
=
\frac{1}{4}\,(y_{x^1})^2,
\ \ \ \ \ \ \ 
y_{x^1x^1x^1}
=
0
\end{equation}
{\rm (}plus other equations obtained by cross differentiation{\rm )}
has infinitesimal Lie symmetry algebra isomorphic to $\mathfrak{ so}
(5, \C)$, the complexification $\mathfrak{ so} (3, 2) \otimes \C$.
\end{theorem}

Through tentative issues (\cite{ eb2006, gm2006}), it has been
suspected that $M_0$ is the right model in the category of real
analytic hypersurfaces $M \subset \C^3$ having Levi form of rank 1
that are 2-nondegenerate everywhere. Based on the rigidity of the
simple Lie algebra $\mathfrak{ so} (5, \C)$ (Theorem~5.15),
Theorem~8.105 below will confirm this expectation.

\subsection*{ 8.52.~Preparation}
Thus, translating the considerations to the {\sc pde} language, with
$n=2$ and $m=1$, consider a submanifold of solutions of the form
\def\theequation{8.53}\begin{equation}
\aligned
y
&
=
b
+
\Pi(x,a)
\\
&
=
b
+
\frac{2\,x^1a^1+x^1x^1a^2+a^1a^1x^2}{
1-x^2a^2}
+
{\rm O}_4,
\endaligned
\end{equation}
where ${\rm O}_4$ is a function of $(x, a)$ only. The term $2 \, x^1
a^1$ corresponds to a Levi form of rank $\geqslant 1$ at every
point. The term $x^1 x^1 a^2$ guarantees solvability with respect to the
parameters (compare Definition~2.12). Let us
develope
\def\theequation{8.54}\begin{equation}
\Pi(x,a)
=
\sum_{k_1,k_2\geqslant 0}\,
\sum_{l_1,l_2\geqslant 0}\,
\Pi_{k_1,k_2,l_1,l_2}\,
(x^1)^{k_1}\,
(x^2)^{k_2}\,
(a^1)^{l_1}\,
(a^2)^{l_2},
\end{equation}
with $\Pi_{ k_1, k_2, l_1, l_2} \in \K$. Of course, $\Pi_{ 1, 0, 1,
0} = 2$, $\Pi_{ 2, 0, 0, 1} = 1$ and $\Pi_{ 0, 1, 2, 0} = 1$.

\def\thelemma{8.55}\begin{lemma}
A transformation belonging to ${\sf G}_{ {\sf v}, {\sf p}}$ insures
\def\theequation{8.56}\begin{equation}
\aligned
\Pi_{k_1,k_2,0,0}
&
=
0,
\ \ \ \ \ \ \ 
k_1+k_2\geqslant 0,
\ \ \ \ \ \ \ \ \ \
\Pi_{0,0,l_1,l_2}
&
=
0,
\ \ \ \ \ \ \ 
l_1+l_2\geqslant 0,
\\
\Pi_{k_1,k_2,1,0}
&
=
0,
\ \ \ \ \ \ \
k_1+k_2\geqslant 2,
\ \ \ \ \ \ \ \ \ \
\Pi_{1,0,l_1,l_2}
&
=
0,
\ \ \ \ \ \ \ 
l_1+l_2\geqslant 2,
\\
\Pi_{k_1,k_2,2,0}
&
=
0,
\ \ \ \ \ \ \
k_1+k_2\geqslant 2,
\ \ \ \ \ \ \ \ \ \
\Pi_{2,0,l_1,l_2}
&
=
0,
\ \ \ \ \ \ \ 
l_1+l_2\geqslant 2.
\endaligned
\end{equation} 
\end{lemma}

\proof
Lemma~7.11 achieves the first line. The monomial $x^1$ being factored
by $[a^1 + {\rm O}_2 (a)]$, we set $a^1 := a^1 + {\rm O}_2 (a)$ to
achieve $\Pi_{ 1, 0, l_1, l_2} = 0$, $l_1 + l_2 \geqslant 2$. As in
the proof of Lemma~7.18, we pass to the dual equation $b = y - \Pi (x,
a)$ to complete $\Pi_{ k_1, k_2, 1, 0} = 0$, $k_1 + k_2 \geqslant 2$.
Finally, $x^1 x^1$ is factored by $[ a^2 + {\rm O}_2 (a)]$, so we
proceed similarly to achieve the third line.
\endproof

Since $\Pi (x, a)$ is assumed to be independent of $b$, the assumption
that the Levi form of $M \subset \C^3$ has exactly rank 1 at every
point translates to:
\def\theequation{8.57}\begin{equation}
0
\equiv
\left\vert
\begin{array}{cc}
\Pi_{x^1a^1}
&
\Pi_{x^1a^2}
\\
\Pi_{x^2a^1}
&
\Pi_{x^2a^2}
\end{array}
\right\vert.
\end{equation}
For later use, it is convenient to 
develope somehow $\Pi$ with respect to the powers of
$(a^1, a^2)$:
\def\theequation{8.58}\begin{equation}
\aligned
y
&
=
b
+
\frac{2\,x^1a^1+x^1x^1a^2+a^1a^1x^2}{
1-x^2a^2}
+
\\
&
\ \ \ \ \ \ \ \ \ \ \ \ \
+
a^2\,{\bf b}(x)
+
a^1a^2\,{\bf d}(x)
+
a^2a^2\,{\bf e}(x)
+
a^1a^1a^1\,{\bf f}(x)
+
a^1a^1a^2\,{\bf g}(x)
+
\\
&
\ \ \ \ \ \ \ \ \ \ \ \ \
+
(a^1)^4\,{\sf R}
+
(a^1)^3a^2\,{\sf R}
+
a^1(a^2)^2\,{\sf R}+
(a^2)^3\,{\sf R},
\endaligned
\end{equation}
with ${\sf R} = {\sf R} (x, a)$ being an unspecified remainder.
Thanks to the previous lemma, the coefficients ${\bf a}$ of $a^1$ and
${\bf c}$ of $a^1 a^1$ must vanish. The function ${\bf b}$ is an
${\rm O}_3$.

\def\thelemma{8.59}\begin{lemma}
The function ${\bf b}$ depends only on $x^1$, is an ${\rm O}_3 (x^1)$
and the function ${\bf g}$ satisfies ${\bf g}_{ x^2 x^2} ( 0) = 0$.
\end{lemma}

\proof
Developing $[ 1 - x^2 a^2 ]^{ -1} = 1 + x^2 a^2 + (x^2a^2)^2 + {\rm O}_3
(x^2 a^2)$, inserting the right hand side of
\def\theequation{8.60}\begin{equation}
\aligned
y
-
b
&
=
a^1
\big[
2x^1
\big]
+
a^2
\big[
x^1x^1
+
{\bf b}(x)
\big]
+
a^1a^1
\big[
x^2
\big]
+
a^1a^2
\big[
2x^1x^2
+
{\bf d}(x)
\big]
+
\\
&
\ \ \ \ \ 
+
a^2a^2
\big[
x^1x^1x^2
+
{\bf e}(x)
\big]
+
a^1a^1a^1
\big[
{\bf f}(x)
\big]
+
a^1a^1a^2
\big[
x^2x^2
+
{\bf g}(x)
\big]
+
\\
&
\ \ \ \ \
+
(a^1)^4\,{\sf R}
+
(a^1)^3a^2\,{\sf R}
+
a^1(a^2)^2\,{\sf R}+
(a^2)^3\,{\sf R}
\endaligned
\end{equation}
in the determinant~\thetag{ 8.57} and selecting the
coefficients of ${\rm cst.}$, of $a^1$, of $a^2$ and of $a^1 a^1$,
we get four {\sc pde}s:
\def\theequation{8.61}\begin{equation}
\aligned
0
&
=
2\,{\bf b}_{x^2},
\\
0
&
=
2\,{\bf d}_{x^2}
-
2\,{\bf b}_{x^1},
\\
0
&
=
4\,{\bf e}_{x^2}
-
2x^1\,{\bf d}_{x^2}
-
2x^1\,{\bf b}_{x^1}
-
{\bf d}_{x^2}\,{\bf b}_{x^1},
\\
0
&
=
2\,{\bf g}_{x^2}
-
2\,{\bf d}_{x^1}
-
\big[
6x^1
+
3\,{\bf b}_{x^1}
\big]\,{\bf f}_{x^2}.
\endaligned
\end{equation}
The first one yields ${\bf b} = {\bf b} (x^1)$, which must be an ${\rm
O}_3 (x^1)$, because the whole remainder is an ${\rm O}_4$.
Differentiating the fourth with respect to $x^2$, it then follows that
${\bf g}_{ x^2 x^2} (0) = 0$.

\subsection*{ 8.62.~Associated {\sc pde} system \thetag{ $\mathcal{ E
}_5$}} Next, differentiating~\thetag{ 8.60} 
with respect to $x^1$,
to $x^1 x^1$ and to 
$x^1 x^1 x^1$, we compute $y_{
x^1}$ and $y_{ x^1 x^1}$, we substitute $y_1$ and $y_{ 1, 1}$ and we
push the monomials $a^2 a^2$, $a^1 a^1 a^1$ and $a^1 a^1 a^2$ in the
remainder:
\def\theequation{8.63}\begin{equation}
\aligned
y_1
&
=
2a^1
+
a^2[2x^1+{\bf b}_{x^1}]
+
a^1a^2[2x^2+{\bf d}_{x^1}]
+
(a^2)^2\,{\sf R}
+
(a^1)^3\,{\sf R}
+
(a^1)^2\,a^2\,{\sf R},
\\
y_{1,1}
&
=
a^2[2+{\bf b}_{x^1x^1}]
+
a^1a^2[{\bf d}_{x^1x^1}]
+
(a^2)^2\,{\sf R}
+
(a^1)^3\,{\sf R}
+
(a^1)^2\,a^2\,{\sf R},
\\
y_{1,1,1}
&
=
a^2[{\bf b}_{x^1x^1x^1}]
+
a^1a^2[{\bf d}_{x^1x^1x^1}]
+
(a^2)^2\,{\sf R}
+
(a^1)^3\,{\sf R}
+
(a^1)^2a^2\,{\sf R}.
\endaligned
\end{equation}
Here, the written remainder {\it cannot}\, incorporate $a^1 a^1$, so
it is said that 
the coefficient of $a^1 a^1$ does vanish in each
equation above. Solving for $a^1$ and $a^2$ from the first two
equations, we get
\begin{small}
\def\theequation{8.64}\begin{equation}
\left\{
\aligned
a^1
&
=
\frac{1}{2}\,y_1
-
y_{1,1}\,
\left[
\frac{2x^1+{\bf b}_{x^1}}{4+2\,{\bf b}_{x^1x^1}}
\right]
-
y_1y_{1,1}\,
\left[
\frac{2x^2+{\bf d}_{x^1}}{8+4\,{\bf b}_{x^1x^1}}
\right]
+
\\
&
\ \ \ \ \
+
(y_{1,1})^2\,{\sf R}
+
(y_1)^3\,{\sf R}
+
(y_1)^2y_{1,1}\,{\sf R},
\\
a^2
&
=
y_{1,1}
\left[
\frac{1}{2+{\bf b}_{x^1x^1}}
\right]
-
y_1y_{1,1}
\left[
\frac{{\bf d}_{x^1x^1}}{2(2+{\bf b}_{x^1x^1})^2}
\right]
+
(y_{1,1})^2\,{\sf R}
+
(y_1)^3\,{\sf R}
+
(y_1)^2y_{1,1}\,{\sf R}.
\endaligned\right.
\end{equation}
\end{small}

\noindent
We then get (notice the change of remainder):
\begin{small}
\def\theequation{8.65}\begin{equation}
\aligned
a^1a^1
&
=
\frac{1}{4}\,(y_1)^2
-
y_1y_{1,1}
\left[
\frac{2x^1+{\bf b}_{x^1}}{4+2\,{\bf b}_{x^1x^1}}
\right]
-
(y_1)^2y_{1,1}
\left[
\frac{2x^2+{\bf d}_{x^1}}{8+4\,{\bf b}_{x^1x^1}}
\right]
+
(y_1)^3\,{\sf R}
+
(y_{1,1})^2\,{\sf R},
\\
a^1a^2
&
=
y_1y_{1,1}
\left[
\frac{1}{4+2\,{\bf b}_{x^1x^1}}
\right]
-
(y_1)^2y_{1,1}
\left[
\frac{{\bf d}_{x^1x^1}}{(4+2\,{\bf b}_{x^1x^1})^2}
\right]
+
(y_1)^3\,{\sf R}
+
(y_{1,1})^2\,{\sf R},
\\
a^2a^2
&
=
(y_1)^3\,{\sf R}
+
(y_{1,1})^2\,{\sf R},
\\
a^1a^1a^1
&
=
-
(y_1)^2y_{1,1}
\left[
\frac{6x^1+3\,{\bf b}_{x^1}}{16+8\,{\bf b}_{x^1x^1}}
\right]
(y_1)^3\,{\sf R}
+
(y_{1,1})^2\,{\sf R},
\\
a^1a^1a^2
&
=
(y_1)^2\,y_{1,1}
\left[
\frac{1}{8+4\,{\bf b}_{x^1x^1}}
\right]
+
(y_1)^3\,{\sf R}
+
(y_{1,1})^2\,{\sf R}.
\endaligned
\end{equation}
\end{small}

\noindent
Differentiating~\thetag{ 8.60} with respect to $x^2$, 
substituting $y_2$ for $y_{ x^2}$ and
replacing ${\bf d}_{ x^2}$ by
${\bf b}_{ x^1}$ thanks to $\text{\rm (8.61)}_2$,
we get
\begin{small}
\def\theequation{8.66}\begin{equation}
\aligned
y_2
&
=
a^1a^1
+
a^1a^2[2x^1+{\bf b}_{x^1}]
+
a^2a^2[x^1x^1+{\bf e}_{x^2}]
+
a^1a^1a^1[{\bf f}_{x^2}]
+
a^1a^1a^2[2x^2+{\bf g}_{x^2}]
+
\\
&
\ \ \ \ \
+
(a^1)^4\,{\sf R}
+
(a^1)^3a^2\,{\sf R}
+
a^1(a^2)^2\,{\sf R}
+
(a^2)^3\,{\sf R}.
\endaligned
\end{equation}
\end{small}
Replacing the monomials~\thetag{ 8.65}, we finally obtain:
\def\theequation{8.67}\begin{equation}
\aligned
y_2
&
=
\frac{1}{4}\,(y_1)^2
+
(y_1)^2y_{1,1}
\left[
\frac{2\,{\bf g}_{x^2}-2\,{\bf d}_{x^1}
-
(6x^1+3\,{\bf b}_{x^1}){\bf f}_{x^2}}{
16+8\,{\bf b}_{x^1x^1}}
-
\frac{(2x^1+{\bf b}_{x^1})\,{\bf d}_{x^1x^1}}{
(4+2\,{\bf b}_{x^1x^1})^2}
\right]
+
\\
&
\ \ \ \ \ \ \ \ \ \ \ \ \ \ \ \ \ \ \ \ \ \ 
+
(y_1)^3\,{\sf R}
+
(y_{1,1})^2\,{\sf R}.
\endaligned
\end{equation}
Thanks to $\text{\rm (8.61)}_4$, the first
(big) coefficient of $(y_1)^2 y_{ 1, 1}$ is
zero; then the remainder coefficient is an ${\rm O} (x^1)$, hence
vanishes at $x = 0$, together with its partial first 
derivative with respect
to $x^2$. Accordingly, by ${\sf s}^* = {\sf s}^* (x^1, x^2)$, we will
denote an unspecified function satisfying
\def\theequation{8.68}\begin{equation}
\boxed{
{\sf s}^*(0)
=
0
\ \ \ \ \ \ \
\text{\rm and}
\ \ \ \ \ \ \
{\sf s}_{x^2}^*(0)
=
0
}.
\end{equation}

\def\thelemma{8.69}\begin{lemma}
The skeleton of the {\sc pde} system 
\thetag{ $\mathcal{ E}_5$} associated to
the submanifold~\thetag{ 8.58} possesses three main equations of the
form
\def\theequation{$\Delta_{\mathcal{E}_5}$}\begin{equation}
\left\{
\aligned
y_2
&
=
\frac{1}{4}\,(y_1)^2
+
(y_1)^3\,{\sf r}
+
(y_1)^4\,{\sf r}
+
(y_1)^5\,{\sf r}
+
(y_1)^6\,{\sf R}
+
\\
&
\ \ \ \ \ \ \ \ \ \ \ \ \ \ \ \ \
+
y_{1,1}
\big[
(y_1)^2\,{\sf s}^*
+
(y_1)^3\,{\sf r}
+
(y_1)^4\,{\sf r}
+
(y_1)^5\,{\sf r}
\big]
+
(y_{1,1})^2\,{\sf R},
\\
y_{1,2}
&
=
\frac{1}{2}\,y_1y_{1,1}
+
(y_1)^3\,{\sf r}
+
(y_1)^4\,{\sf r}
+
(y_1)^5\,{\sf r}
+
(y_1)^6\,{\sf R}
+
\\
&
\ \ \ \ \ \ \ \ \ \ \ \ \ \ \ \ \ \
+
y_{1,1}
\big[
(y_1)^2\,{\sf r}
+
(y_1)^3\,{\sf r}
+
(y_1)^4\,{\sf r}
+
(y_1)^5\,{\sf r}
\big]
+
(y_1)^6\,{\sf R},
\\
y_{1,1,1}
&
=
(y_1)^3\,{\sf r}
+
(y_1)^4\,{\sf r}
+
y_{1,1}
\big[
{\sf r}
+
y_1\,{\sf r}
+
(y_1)^2\,{\sf r}
+
(y_1)^3\,{\sf r}
\big]
+
\\
&
\ \ \ \ \ \ \ \ \ \ \ \ \ \ \ \ \ \
+
(y_{1,1})^2
\big[
{\sf r}
+
y_1\,{\sf r}
+
(y_1)^2\,{\sf r}
+
(y_1)^3\,{\sf r}
\big]
+
(y_{1,1})^3\,{\sf R},
\endaligned\right.
\end{equation}
where the letter ${\sf r}$ denotes an unspecified function of $(x^1,
x^2)$, and where the coefficient ${\sf s}^*$ of $(y_1)^2 y_{ 1, 1}$
in the first equation satisfies~\thetag{ 8.68}.
\end{lemma}

\proof
To get the second equation, we compute:
\def\theequation{8.70}\begin{equation}
\aligned
y_{1,2}
&
=
a^1a^2[2+{\bf b}_{x^1x^1}]
+
a^1a^1a^2[{\bf g}_{x^1x^2}]
+
(a^1)^3\,{\sf R}
+
(a^2)^2\,{\sf R}
\\
&
=
\frac{1}{2}\,y_1y_{1,1}
+
(y_1)^2y_{1,1}\,{\sf r}
+
(y_1)^3\,{\sf R}
+
(y_{1,1})^2\,{\sf R}.
\endaligned
\end{equation}
The third equation is got similarly from $\text{\rm (8.63)}_3$.  To
conclude, we develope the first two equations ${\rm mod}\, \big[
(y_1)^6, (y_{ 1, 1})^2 \big]$ and the third one ${\rm mod}\, \big[
(y_1)^4, (y_{ 1, 1})^3 \big]$.
\endproof

This precise skeleton will be referred to as $\Delta_{ \mathcal{
E}_5}$ in the sequel. With the letter ${\sf r}$, the computation rules
are ${\rm cst.} {\sf r} = {\sf r} + {\sf r} = {\sf r} + {\sf s}^* =
{\sf r} \cdot {\sf r} = {\sf r}$; sometimes, ${\sf s}^*$ may be
replaced plainly by ${\sf r}$.

\subsection*{ 8.71.~Infinitesimal Lie symmetries of \thetag{ $\mathcal{ 
E}_5$}} Letting $\mathcal{ L} = \mathcal{ X}^1 \, \frac{ \partial }{
\partial x^1} + \mathcal{ X}^2 \, \frac{ \partial }{ \partial x^2} +
\mathcal{ Y} \, \frac{ \partial }{ \partial y}$ be a candidate
infinitesimal Lie symmetry and applying
\def\theequation{8.72}\begin{equation}
\aligned
\mathcal{L}^{(3)}
&
=
\mathcal{X}^1\,\frac{\partial}{\partial x^1}
+
\mathcal{X}^2\,\frac{\partial}{\partial x^2}
+
\mathcal{Y}\,\frac{\partial}{\partial y}
+
{\bf Y}_1\,\frac{\partial}{\partial y_1}
+
{\bf Y}_2\,\frac{\partial}{\partial y_2}
+
\\
&
\ \ \ \ \
+
{\bf Y}_{1,1}\,\frac{\partial}{\partial y_{1,1}}
+
{\bf Y}_{1,2}\,\frac{\partial}{\partial y_{1,2}}
+
{\bf Y}_{2,1}\,\frac{\partial}{\partial y_{2,1}}
+
{\bf Y}_{2,2}\,\frac{\partial}{\partial y_{2,2}}
+
\\
&
\ \ \ \ \
+
{\bf Y}_{1,1,1}\,\frac{\partial}{\partial y_{1,1,1}}
+
\cdots
+
{\bf Y}_{2,2,2}\,\frac{\partial}{\partial y_{2,2,2}}
\endaligned
\end{equation}
to the skeleton $\Delta_{ \mathcal{ E}_5}$, we obtain firstly,
computing ${\rm mod}\, \big[ (y_1)^5, y_{1, 1} \big]$:
\def\theequation{8.73}\begin{equation}
\aligned
0
&
\equiv
-{\bf Y}_2
+
\frac{1}{2}\,y_1\,{\bf Y}_1
+
\\
&
\ \ \ \ \ 
+
(y_1)^3\,{\sf r}\,\mathcal{X}^1
+
(y_1)^4\,{\sf r}\,\mathcal{X}^1
+
(y_1)^3\,{\sf r}\,\mathcal{X}^2
+
(y_1)^4\,{\sf r}\,\mathcal{X}^2
+
\\
&
\ \ \ \ \
+
{\bf Y}_1
\big[
(y_1)^2\,{\sf r}
+
(y_1)^3\,{\sf r}
+
(y_1)^4\,{\sf r}
\big]
+
\\
&
\ \ \ \ \
+
{\bf Y}_{1,1}
\big[
(y_1)^2\,{\sf s}^*
+
(y_1)^3\,{\sf r}
+
(y_1)^4\,{\sf r}
\big],
\endaligned
\end{equation}
secondly, computing ${\rm mod}\, \big[ (y_1)^5, y_{ 1, 1} \big]$:
\def\theequation{8.74}\begin{equation}
\aligned
0
&
\equiv
-
{\bf Y}_{1,2}
+
\frac{1}{2}\,y_1\,{\bf Y}_{1,1}
+
\\
&
\ \ \ \ \ 
+
(y_1)^3\,{\sf r}\,\mathcal{X}^1
+
(y_1)^4\,{\sf r}\,\mathcal{X}^1
+
(y_1)^3\,{\sf r}\,\mathcal{X}^2
+
(y_1)^4\,{\sf r}\,\mathcal{X}^2
+
\\
&
\ \ \ \ \ \
+
{\bf Y}_1\big[
(y_1)^2\,{\sf r}
+
(y_1)^3\,{\sf r}
+
(y_1)^4\,{\sf r}
\big]
+
\\
&
\ \ \ \ \
+
{\bf Y}_{1,1}
\big[
(y_1)^2\,{\sf r}
+
(y_1)^3\,{\sf r}
+
(y_1)^4\,{\sf r}
\big],
\endaligned
\end{equation}
and thirdly, computing ${\rm mod}\, \big[ (y_1)^3,
(y_{1,1})^2 \big]$:
\def\theequation{8.75}\begin{equation}
\aligned
0
&
\equiv
-{\bf Y}_{1,1,1}
+
y_{1,1}\,\mathcal{X}^1
+
y_{1,1}\,\mathcal{X}^2
+
\\
&
\ \ \ \ \
+
y_{1,1}\,y_1\,\mathcal{X}^1
+
y_{1,1}\,y_1\,\mathcal{X}^2
+
y_{1,1}\,(y_1)^2\,\mathcal{X}^1
+
y_{1,1}\,(y_1)^2\,\mathcal{X}^2
+
\\
&
\ \ \ \ \ 
+
{\bf Y}_1
\big[
(y_1)^2\,{\sf r}
\big]
+
{\bf Y}_{1,1}
\big[
{\sf r}
+
y_1\,{\sf r}
+
(y_1)^2\,{\sf r}
\big]
+
y_{1,1}\,{\bf Y}_1
\big[
{\sf r}
+
y_1\,{\sf r}
+
(y_1)^2\,{\sf r}
\big]
+
\\
&
\ \ \ \ \ \ \ \ \ \ \ \ \ \ \ \ \ \ \
+
y_{1,1}\,{\bf Y}_{1,1}
\big[
{\sf r}
+
y_1\,{\sf r}
+
(y_1)^2\,{\sf r}
\big].
\endaligned
\end{equation}
Specializing to $n=2$ the formulas~\thetag{ 3.9}(II), 
\thetag{ 3.20}(II) and
\thetag{ 3.24}(II), we get ${\bf Y}_1$, ${\bf Y}_2$, ${\bf Y}_{ 1,
1}$, ${\bf Y}_{1, 2}$ and ${\bf Y}_{ 1, 1, 1}$:
\def\theequation{8.76}\begin{equation}
{\bf Y}_1
=
\mathcal{Y}_{x^1}
+
\big[
\mathcal{Y}_y
-
\mathcal{X}_{x^1}^1
\big]\,
y_1
+
\big[
-\mathcal{X}_{x^1}^2
\big]\,y_2
+
\big[
-\mathcal{X}_y^1
\big]\,(y_1)^2
+
\big[
-\mathcal{X}_y^2
\big]\,y_1y_2.
\end{equation}\def\theequation{8.77}\begin{equation}
{\bf Y}_2
=
\mathcal{Y}_{x^2}
+
\big[
-\mathcal{X}_{x^2}^1
\big]\,
y_1
+
\big[
\mathcal{Y}_y
-
\mathcal{X}_{x^2}^2
\big]\,y_2
+
\big[
-\mathcal{X}_y^1
\big]\,y_1y_2
+
\big[
-\mathcal{X}_y^2
\big]\,y_2y_2.
\end{equation}
\def\theequation{8.78}\begin{equation}
\left\{
\aligned
{\bf Y}_{1,1}
&
=
\mathcal{Y}_{x^1x^1}
+
\big[
2\,\mathcal{Y}_{x^1y}
-
\mathcal{X}_{x^1x^1}^1
\big]\,y_1
+
\big[
-\mathcal{X}_{x^1x^1}^2
\big]\,y_2
+
\big[
\mathcal{Y}_{yy}
-
2\,\mathcal{X}_{x^1y}^1
\big]\,(y_1)^2
+
\\
&
\ \ \ \ \
+
\big[
-2\,\mathcal{X}_{x^1y}^2
\big]\,y_1y_2
+
\big[
-\mathcal{X}_{yy}^1
\big]\,(y_1)^3
+
\big[
-\mathcal{X}_{yy}^2
\big]\,(y_1)^2y_2
+
\\
&
\ \ \ \ \
+
\big[
\mathcal{Y}_y-2\,\mathcal{X}_{x^1}^1
\big]\,y_{1,1}
+
\big[
-
2\,\mathcal{X}_{x^1}^2
\big]\,y_{1,2}
+
\big[
-3\,\mathcal{X}_y^1
\big]\,y_1y_{1,1}
+
\\
&
\ \ \ \ \
+
\big[
-\mathcal{X}_y^2
\big]\,y_2y_{1,1}
+
\big[
-2\,\mathcal{X}_y^2
\big]\,y_1y_{1,2}.
\endaligned\right.
\end{equation}
\def\theequation{8.79}\begin{equation}
\left\{
\aligned
{\bf Y}_{1,2}
&
=
\mathcal{Y}_{x^1x^2}
+
\big[
\mathcal{Y}_{x^2y}
-
\mathcal{X}_{x^1x^2}^1
\big]\,y_1
+
\big[
\mathcal{Y}_{x^1y}
-
\mathcal{X}_{x^1x^2}^2
\big]\,y_2
+
\\
&
\ \ \ \ \
+
\big[
-
\mathcal{X}_{x^2y}^1
\big]\,(y_1)^2
+
\big[
\mathcal{Y}_{yy}
-
\mathcal{X}_{x^2y}^2
-
\mathcal{X}_{x^1y}^1
\big]\,y_1y_2
+
\big[
-\mathcal{X}_{x^1y}^2
\big]\,y_2y_2
+
\\
&
\ \ \ \ \
+
\big[
-\mathcal{X}_{yy}^1
\big]\,(y_1)^2y_2
+
\big[
-\mathcal{X}_{yy}^2
\big]\,y_1(y_2)^2
+
\\
&
\ \ \ \ \
+
\big[
-\mathcal{X}_{x^2}^1
\big]\,y_{1,1}
+
\big[
\mathcal{Y}_y
-
\mathcal{X}_{x^2}^2
-
\mathcal{X}_{x^1}^1
\big]\,y_{1,2}
+
\big[
-\mathcal{X}_{x^1}^2
\big]\,y_{2,2}
+
\\
&
\ \ \ \ \
+
\big[
-2\,\mathcal{X}_y^1
\big]\,y_1y_{1,2}
+
\big[
-2\,\mathcal{X}_y^2
\big]\,y_2y_{1,2}.
\endaligned\right.
\end{equation}
\def\theequation{8.80}\begin{equation}
\left\{
\aligned
{\bf Y}_{1,1,1}
&
=
\mathcal{Y}_{x^1x^1x^1}
+
\big[
3\,\mathcal{Y}_{x^1x^1y}
-
\mathcal{X}_{x^1x^1x^1}^1
\big]\,y_1
+
\big[
-\mathcal{X}_{x^1x^1x^1}^2
\big]\,y_2
+
\\
&
\ \ \ \ \
+
\big[
3\,\mathcal{Y}_{x^1yy}
-
3\,\mathcal{X}_{x^1x^1y}^1
\big]\,(y_1)^2
+
\big[
-3\,\mathcal{X}_{x^1x^1y}^2
\big]\,y_1y_2
+
\\
&
\ \ \ \ \
+
\big[
\mathcal{Y}_{yyy}
-
3\,\mathcal{X}_{x^1yy}^1
\big]\,(y_1)^3
+
\big[
-3\,\mathcal{X}_{x^1yy}^2
\big]\,(y_1)^2y_2
+
\\
&
\ \ \ \ \
+
\big[
-\mathcal{X}_{yyy}^1
\big]\,(y_1)^4
+
\big[
-\mathcal{X}_{yyy}^2
\big]\,(y_1)^3y_2
+
\\
&
\ \ \ \ \
+
\big[
3\,\mathcal{Y}_{x^1y}
-
3\,\mathcal{X}_{x^1x^1}^1
\big]\,y_{1,1}
+
\big[
-3\,\mathcal{X}_{x^1x^1}^2
\big]\,y_{1,2}
+
\\
&
\ \ \ \ \
+
\big[
3\,\mathcal{Y}_{yy}
-
9\,\mathcal{X}_{x^1y}^1
\big]\,y_1y_{1,1}
+
\big[
-3\,\mathcal{X}_{x^1y}^2
\big]\,y_2y_{1,1}
+
\\
&
\ \ \ \ \
+
\big[
-6\,\mathcal{X}_{x^1y}^2
\big]\,y_1y_{1,2}
+
\big[
-6\,\mathcal{X}_{yy}^1
\big]\,(y_1)^2y_{1,1}
+
\big[
-3\,\mathcal{X}_{yy}^2
\big]\,y_1y_2y_{1,1}
+
\\
&
\ \ \ \ \
+
\big[
-3\,\mathcal{X}_{yy}^2
\big]\,(y_1)^2\,y_{1,2}
+
\big[
-3\,\mathcal{X}_y^1
\big]\,(y_{1,1})^2
+
\big[
-3\,\mathcal{X}_y^2
\big]\,y_{1,1}y_{1,2}
+
\\
&
\ \ \ \ \
+
\big[
\mathcal{Y}_y
-
3\,\mathcal{X}_{x^1}^1
\big]\,y_{1,1,1}
+
\big[
-3\,\mathcal{X}_{x^1}^2
\big]\,y_{1,1,2}
+
\big[
-4\,\mathcal{X}_y^1
\big]\,y_1y_{1,1,1}
+
\\
&
\ \ \ \ \
+
\big[
-\mathcal{X}_y^2
\big]\,y_2y_{1,1,1}
+
\big[
-3\,\mathcal{X}_y^2
\big]\,y_1y_{1,1,2}.
\endaligned\right.
\end{equation}
Inserting ${\bf Y}_2$, ${\bf Y}_{ 1,2}$, ${\bf Y}_{ 1, 1, 1}$, ${\bf
Y}_1$, ${\bf Y}_{ 1, 1}$ in the three Lie equations~\thetag{ 8.73},
\thetag{ 8.74}, \thetag{ 8.75}, replacing $y_2$, $y_{ 1, 2}$, $y_{ 1,
1, 1}$ by the values they have on $\Delta_{ \mathcal{ E}_5}$, we get
firstly five linear {\sc pde}s by picking the coefficients of ${\rm
cst.}$, of $y_1$, of $(y_1)^2$, of $(y_1)^3$, of $(y_1)^4$ in~\thetag{
8.73}:
\def\theequation{8.81}\begin{equation}
\left\{
\aligned
0
&
=
\mathcal{Y}_{x^2},
\\
0
&
=
\mathcal{X}_{x^2}^1
+
\frac{1}{2}\,\mathcal{Y}_{x^1},
\\
0
&
=
\mathcal{Y}_y
+
\mathcal{X}_{x^2}^2
-
2\,\mathcal{X}_{x^1}^1
+
{\sf r}\,\mathcal{Y}_{x^1}
+
{\sf s}^*\,\mathcal{Y}_{x^1x^1},
\\
0
&
=
2\,\mathcal{X}_y^1
+
\mathcal{X}_{x^1}^2
+
{\sf r}\,\mathcal{X}^1
+
{\sf r}\,\mathcal{X}^2
+
{\sf r}\,\mathcal{Y}_{x^1}
+
{\sf r}\,\mathcal{Y}_y
+
{\sf r}\,\mathcal{X}_{x^1}^1
+
\\
&
\ \ \ \ \ \ \ \ \ \ \ \ \ \ \ \ \ \ \ \ \ \ \ \ \
+
{\sf r}\,\mathcal{Y}_{x^1x^1}
+
{\sf s}^*\,\mathcal{Y}_{x^1y}
+
{\sf s}^*\,\mathcal{X}_{x^1x^1}^1,
\\
0
&
=
\mathcal{X}_y^2
+
{\sf r}\,\mathcal{X}^1
+
{\sf r}\,\mathcal{X}^2
+
{\sf r}\,\mathcal{Y}_{x^1}
+
{\sf r}\,\mathcal{Y}_y
+
{\sf r}\,\mathcal{X}_{x^1}^1
+
{\sf r}\,\mathcal{X}_{x^1}^2
+
{\sf r}\,\mathcal{X}_y^1
+
\\
&
\ \ \ \ \ \ \ \ \ \ \ \ \ \ \ \ \ \ \ \ \ \ \ \ \
+
{\sf r}\,\mathcal{Y}_{x^1x^1}
+
{\sf r}\,\mathcal{Y}_{x^1y}
+
{\sf r}\,\mathcal{X}_{x^1x^1}^1
+
{\sf s}^*\,\mathcal{X}_{x^1x^1}^2
+
{\sf s}^*\,\mathcal{Y}_{yy}
+
{\sf s}^*\,\mathcal{X}_{x^1y}^1,
\endaligned\right.
\end{equation}
secondly, we get three more linear {\sc pde}s by picking the
coefficients of $(y_1)^2$, of $(y_1)^3$, of $(y_1)^4$ in~\thetag{
8.74}:
\def\theequation{8.82}\begin{equation}
\left\{
\aligned
0
&
=
3\,\mathcal{Y}_{x^1y}
+
\mathcal{X}_{x^1x^2}^2
+
4\,\mathcal{X}_{x^2y}^1
-
2\,\mathcal{X}_{x^1x^1}^1
+
{\sf r}\,\mathcal{X}^1
+
{\sf r}\,\mathcal{X}^2
+
{\sf r}\,\mathcal{Y}_{x^1}
+
{\sf r}\,\mathcal{Y}_{x^1x^1},
\\
0
&
=
2\,\mathcal{Y}_{yy}
+
2\,\mathcal{X}_{x^2y}^2
-
6\,\mathcal{X}_{x^1y}^1
-
\mathcal{X}_{x^1x^1}^2
+
{\sf r}\,\mathcal{X}^1
+
{\sf r}\,\mathcal{X}^2
+
{\sf r}\,\mathcal{Y}_{x^1}
+
{\sf r}\,\mathcal{Y}_y
+
{\sf r}\,\mathcal{X}_{x^1}^1
+
\\
&
\ \ \ \ \ \ \ \ \ \ \ \ \ \ \ \ \ \ \ \ \ \ \ \ \ \ \ \ 
\ \ \ \ \ \ \ \ \ \ \ \ \ \ \ \ \ \ \ \ \ \ \ \ \ \ \ \ \ \ \
+
{\sf r}\,\mathcal{Y}_{x^1x^1}
+
{\sf r}\,\mathcal{Y}_{x^1y}
+
{\sf r}\,\mathcal{X}_{x^1x^1}^1,
\\
0
&
=
4\,\mathcal{X}_{yy}^1
+
3\,\mathcal{X}_{x^1y}^2
+
{\sf r}\,\mathcal{X}^1
+
{\sf r}\,\mathcal{X}^2
+
{\sf r}\,\mathcal{Y}_{x^1}
+
{\sf r}\,\mathcal{Y}_y
+
{\sf r}\,\mathcal{X}_{x^1}^1
+
{\sf r}\,\mathcal{X}_{x^1}^2
+
{\sf r}\,\mathcal{X}_y^1
+
\\
&
\ \ \ \ \ \ \ \ \ \ \ \ \ \ \ \ \ \ \ \ \ \ \ \ \ \ \ \ \ 
+
{\sf r}\,\mathcal{Y}_{x^1x^1}
+
{\sf r}\,\mathcal{Y}_{x^1y}
+
{\sf r}\,\mathcal{X}_{x^1x^1}^1
+
{\sf r}\,\mathcal{X}_{x^1x^1}^2
+
{\sf r}\,\mathcal{Y}_{yy}
+
{\sf r}\,\mathcal{X}_{x^1y}^1.
\endaligned\right.
\end{equation}
and thirdly, we get five more linear {\sc pde}s by picking the
coefficients of ${\rm cst.}$, of $y_1$,
of $y_{ 1, 1}$, of $y_1 y_{ 1, 1}$, of $(y_1)^2 y_{ 1, 1}$
in~\thetag{ 8.75}
\def\theequation{8.83}\begin{equation}
\left\{
\aligned
0
&
=
\mathcal{Y}_{x^1x^1x^1}
+
{\sf r}\,\mathcal{Y}_{x^1x^1},
\\
0
&
=
-3\,\mathcal{Y}_{x^1x^1y}
+
\mathcal{X}_{x^1x^1x^1}^1
+
{\sf r}\,\mathcal{Y}_{x^1x^1}
+
{\sf r}\,\mathcal{Y}_{x^1y}
+
{\sf r}\,\mathcal{X}_{x^1x^1}^1,
\\
0
&
=
\mathcal{Y}_{x^1y}
-
\mathcal{X}_{x^1x^1}^1
+
{\sf r}\,\mathcal{X}^1
+
{\sf r}\,\mathcal{X}^2
+
{\sf r}\,\mathcal{Y}_{x^1}
+
{\sf r}\,\mathcal{Y}_y
+
{\sf r}\,\mathcal{X}_{x^1}^1
+
{\sf r}\,\mathcal{Y}_{x^1x^1}
\\
0
&
=
-\frac{3}{2}\,\mathcal{X}_{x^1x^1}^2
+
3\,\mathcal{Y}_{yy}
-
9\,\mathcal{X}_{x^1y}^1
+
{\sf r}\,\mathcal{X}^1
+
{\sf r}\,\mathcal{X}^2
+
{\sf r}\,\mathcal{X}_{x^1}^1
+
{\sf r}\,\mathcal{Y}_{x^1}
+
{\sf r}\,\mathcal{Y}_y
+
\\
&
\ \ \ \ \ \ \ \ \ \ \ \ \ \ \ \ \ \ \ \ \ \ \ \ \ \ \ \ \ \ \ 
\ \ \ \ \ \ \ \ \ \ \ \ \ \ \ \ \ \ 
+
{\sf r}\,\mathcal{X}_{x^1}^2
+
{\sf r}\,\mathcal{X}_y^1
+
{\sf r}\,\mathcal{Y}_{x^1x^1}
+
{\sf r}\,\mathcal{Y}_{x^1y}
+
{\sf r}\,\mathcal{X}_{x^1x^1}^1,
\\
0
&
=
6\,\mathcal{X}_{yy}^1
+
\frac{15}{4}\,\mathcal{X}_{x^1y}^2
+
{\sf r}\,\mathcal{X}^1
+
{\sf r}\,\mathcal{X}^2
+
{\sf r}\,\mathcal{X}_{x^1}^1
+
{\sf r}\,\mathcal{Y}_{x^1}
+
{\sf r}\,\mathcal{Y}_y
+
{\sf r}\,\mathcal{X}_{x^1}^2
+
{\sf r}\,\mathcal{X}_y^1
+
{\sf r}\,\mathcal{X}_y^2
+
\\
&
\ \ \ \ \ \ \ \ \ \ \ \ \ \ \ \ \ \ \ \ \ \ \ \ \ \ \ \ \ \
+
{\sf r}\,\mathcal{Y}_{x^1x^1}
+
{\sf r}\,\mathcal{Y}_{x^1y}
+
{\sf r}\,\mathcal{X}_{x^1x^1}^1
+
{\sf r}\,\mathcal{X}_{x^1x^1}^2
+
{\sf r}\,\mathcal{Y}_{yy}
+
{\sf r}\,\mathcal{X}_{x^1y}^1.
\endaligned\right.
\end{equation}

\def\theproposition{8.84}\begin{proposition}
Setting as initial conditions
the ten specific differential coefficients
\def\theequation{8.85}\begin{equation}
\aligned
{\sf P}
:=
&\
{\sf P}
\big(
\mathcal{X}^1,
\mathcal{X}^2,
\mathcal{Y},
\mathcal{X}_y^1,
\mathcal{X}_{x^2}^2,
\mathcal{Y}_{x^1},
\mathcal{Y}_y,
\mathcal{X}_{x^1x^2}^2,
\mathcal{Y}_{x^1x^1},
\mathcal{Y}_{yy}
\big)
\\
=
&\
{\sf r}\,\mathcal{X}^1
+
{\sf r}\,\mathcal{X}^2
+
{\sf r}\,\mathcal{Y}
+
{\sf r}\,\mathcal{X}_y^1
+
{\sf r}\,\mathcal{X}_{x^2}^2
+
{\sf r}\,\mathcal{Y}_{x^1}
+
{\sf r}\,\mathcal{Y}_y
+
{\sf r}\,\mathcal{X}_{x^1x^2}^2
+
{\sf r}\,\mathcal{Y}_{x^1x^1}
+
{\sf r}\,\mathcal{Y}_{yy},
\endaligned
\end{equation}
it follows by cross differentiations and by linear
substitutions from the 
Lie equations $\text{\rm (8.81)}_i$, $i=1, 2, 3, 4, 5$, 
$\text{\rm (8.82)}_j$, $j=1, 2, 3$, 
$\text{\rm (8.83)}_i$, $i=1, 2, 3, 4, 5$, 
that 
$\mathcal{ X}_{ x^1}^1$,
$\mathcal{ X}_{ x^1}^2$,
$\mathcal{ Y}_{ x^2}$,   
$\mathcal{ X}_{ x^2}^1$, 
$\mathcal{ X}_y^2$, 
$\mathcal{ X}_{ x^1y}^1$, 
$\mathcal{ X}_{ x^2 x^2}^2$,
$\mathcal{ Y}_{ x^1 x^2}$, 
$\mathcal{ X}_{ x^2 y}^1$,
$\mathcal{ X}_{ x^2 y}^2$, 
$\mathcal{ Y}_{ x^1 y}$,  
$\mathcal{ X}_{ yy}^1$,  
$\mathcal{ Y}_{ x^2 y}$, 
$\mathcal{ X}_{ x^1x^1 x^2}^2$,
$\mathcal{ Y}_{ x^1 x^1 x^1}$,  
$\mathcal{ X}_{ x^1 x^2 x^2}^2$ 
$\mathcal{ Y}_{ x^1 x^1 x^2}$, 
$\mathcal{ X}_{ x^1 x^2 y}^2$,
$\mathcal{ Y}_{ x^1 x^1 y}$, 
$\mathcal{ Y}_{ x^1 yy}$, 
$\mathcal{ Y}_{ x^2yy}$, 
$\mathcal{ Y}_{ yyy}$
are uniquely determined as linear combinations
of $\big(
\mathcal{ X}^1,
\mathcal{ X}^2,
\mathcal{ Y},
\mathcal{ X}_y^1,
\mathcal{ X}_{x^2}^2,
\mathcal{ Y}_{x^1},
\mathcal{ Y}_y,
\mathcal{ X}_{x^1x^1}^2,
\mathcal{ Y}_{x^1x^1},
\mathcal{ Y}_{yy}
\big)$, namely{\rm :}
\def\theequation{8.86}\begin{equation}
\left\{
\aligned
&
\mathcal{X}_{x^1}^1
\overset{1}{=}
{\sf P},\ \ \ \ \ \ \ \ \ \ \ \ \ \
\mathcal{X}_{x^1}^2
\overset{2}{=}
{\sf P},\ \ \ \ \ \ \ \ \ \ \ \ \ \ \ \ \
\mathcal{Y}_{x^2}
\overset{3}{=}
{\sf P},
\\
&
\mathcal{X}_{x^2}^1
\overset{4}{=}
{\sf P},\ \ \ \ \ \ \ \ \ \ \ \ \ \ \ \
\mathcal{X}_y^2
\overset{5}{=}
{\sf P},\ \ \ \ \ \ \ \ \ \ \ \ \ \
\\
&
\mathcal{X}_{x^1y}^1
\overset{6}{=}
{\sf P},\ \ \ \ \ \ \ \ \ \ \
\mathcal{X}_{x^2x^2}^2
\overset{7}{=}
{\sf P},\ \ \ \ \ \ \ \ \ \ \ \ \
\mathcal{Y}_{x^1x^2}
\overset{8}{=}
{\sf P},
\\
&
\mathcal{X}_{x^2y}^1
\overset{9}{=}
{\sf P},\ \ \ \ \ \ \ \ \ \ \ \ \,
\mathcal{X}_{x^2y}^2
\overset{10}{=}
{\sf P},\ \ \ \ \ \ \ \ \ \ \ \ \ \ \,
\mathcal{Y}_{x^1y}
\overset{11}{=}
{\sf P},
\\
&
\mathcal{X}_{yy}^1
\overset{12}{=}
{\sf P},
\ \ \ \ \ \ \ \ \ \ \ \ \ \ \ \ \ \ \ \ \ \ \ \ \ \ \ \ \ \ \ \ 
\ \ \ \ \ \ \ \ \ \ \ \ \ \ 
\mathcal{Y}_{x^2y}
\overset{13}{=}
{\sf P},
\\
&
\ \ \ \ \ \ \ \ \ \ \ \ \ \ \ \ \ \ \ \ \ \ \ \ \
\mathcal{X}_{x^1x^1x^2}^2
\overset{14}{=}
{\sf P}, \ \ \ \ \ \ \ \ \ \ \
\mathcal{Y}_{x^1x^1x^1}
\overset{15}{=}
{\sf P},
\\
&
\ \ \ \ \ \ \ \ \ \ \ \ \ \ \ \ \ \ \ \ \ \ \ \ \
\mathcal{X}_{x^1x^2x^2}^2
\overset{16}{=}
{\sf P}, \ \ \ \ \ \ \ \ \ \ \
\mathcal{Y}_{x^1x^1x^2}
\overset{17}{=}
{\sf P},
\\
&
\ \ \ \ \ \ \ \ \ \ \ \ \ \ \ \ \ \ \ \ \ \ \ \ \ \
\mathcal{X}_{x^1x^2y}^2
\overset{18}{=}
{\sf P}, \ \ \ \ \ \ \ \ \ \ \ \ \
\mathcal{Y}_{x^1x^1y} 
\overset{19}{=}
{\sf P},
\\
&
\ \ \ \ \ \ \ \ \ \ \ \ \ \ \ \ \ \ \ \ \ \ \ \ \ \ \ \ \ \ 
\ \ \ \ \ \ \ \ \ \ \ \ \ \ \ \ \ \ \ \ \ \ \ \ \ \ \ \ \ \ \,
\mathcal{Y}_{x^1yy}
\overset{20}{=}
{\sf P},
\\
&
\ \ \ \ \ \ \ \ \ \ \ \ \ \ \ \ \ \ \ \ \ \ \ \ \ \ \ \ \ \ 
\ \ \ \ \ \ \ \ \ \ \ \ \ \ \ \ \ \ \ \ \ \ \ \ \ \ \ \ \ \ \,
\mathcal{Y}_{x^2yy}
\overset{21}{=}
{\sf P},
\\
&
\ \ \ \ \ \ \ \ \ \ \ \ \ \ \ \ \ \ \ \ \ \ \ \ \ \ \ \ \ \ 
\ \ \ \ \ \ \ \ \ \ \ \ \ \ \ \ \ \ \ \ \ \ \ \ \ \ \ \ \ \ \ \
\mathcal{Y}_{yyy}
\overset{22}{=}
{\sf P}.
\endaligned\right.
\end{equation}
\end{proposition}

Then the expressions ${\sf P}$ are stable under differentiation with
respect to $x^1$, to $x^2$, to $y$ and moreover, all other, higher
order partial derivatives of $\mathcal{ X}^1$, of $\mathcal{ X}^2$, of
$\mathcal{ Y}$ may be expressed as ${\sf P}\big( \mathcal{ X}^1,
\mathcal{ X}^2, \mathcal{ Y}, \mathcal{ X}_y^1, \mathcal{ X}_{x^2}^2,
\mathcal{ Y}_{x^1}, \mathcal{ Y}_y, \mathcal{ X}_{x^1x^2}^2, \mathcal{
Y}_{x^1x^1}, \mathcal{ Y}_{yy} \big)$.

\def\thecorollary{8.87}\begin{corollary}
Every infinitesimal Lie symmetry of the {\sc pde} system 
\thetag{ $\mathcal{
E}_5$} is uniquely determined by the ten initial Taylor coefficients
\def\theequation{8.88}\begin{equation}
\mathcal{X}^1(0),
\mathcal{X}^2(0),
\mathcal{Y}(0),
\mathcal{X}_y^1(0),
\mathcal{X}_{x^2}^2(0),
\mathcal{Y}_{x^1}(0),
\mathcal{Y}_y(0),
\mathcal{X}_{x^1x^2}^2(0),
\mathcal{Y}_{x^1x^1}(0),
\mathcal{Y}_{yy}(0).
\end{equation}
\end{corollary}

\proof[Proof of the proposition]
At first, 
$\text{\rm (8.83)}_1$ yields
$\text{\rm (8.86)}_{15}$;
$\text{\rm (8.81)}_1$ yields
$\text{\rm (8.86)}_3$; differentiating
$\text{\rm (8.81)}_1$
with respect to $x^1$ yields
$\text{\rm (8.86)}_8$; differentiating
$\text{\rm (8.81)}_1$
with respect to $y$ yields
$\text{\rm (8.86)}_{13}$; differentiating
$\text{\rm (8.81)}_1$
with respect to $x^1x^1$ yields
$\text{\rm (8.86)}_{17}$;
and differentiating
$\text{\rm (8.81)}_1$
with respect to $yy$ yields
$\text{\rm (8.86)}_{21}$.
Also, rewriting 
$\text{\rm (8.81)}_2$ as
\def\theequation{8.89}\begin{equation}
\mathcal{X}_{x^2}^1
=
-\frac{1}{2}\,\mathcal{Y}_{x^1},
\end{equation}
we get $\text{\rm (8.86)}_4$; and rewriting $\text{\rm (8.81)}_3$ as
\def\theequation{8.90}\begin{equation}
\mathcal{X}_{x^1}^1
=
\frac{1}{2}\,\mathcal{X}_{x^2}^2
+
\frac{1}{2}\,\mathcal{Y}_y
+
{\sf r}\,\mathcal{Y}_{x^1}
+
{\sf s}^*\,\mathcal{Y}_{x^1x^1},
\end{equation}
we get $\text{\rm (8.86)}_1$.

Next, 
differentiating 
$\text{\rm (8.81)}_2$ with respect to $x^1$ and
$\text{\rm (8.81)}_3$ with respect to $x^2$, we get, 
taking account of $0 = \mathcal{ Y}_{ x^2 y} = 
\mathcal{ Y}_{ x^1 x^2} = \mathcal{ Y}_{ x^1 x^1 x^2}$, 
replacing $\mathcal{ X}_{ x^1x^2}$ by
$-\frac{ 1}{ 2}\, \mathcal{ Y}_{ x^1 x^1}$ and
solving for $\mathcal{ X}_{ x^2 x^2}^2$:
\def\theequation{8.91}\begin{equation}
\aligned
0
&
=
\mathcal{X}_{x^1x^2}^1
+
\frac{1}{2}\,\mathcal{Y}_{x^1x^1},
\\
\mathcal{X}_{x^2x^2}^2
&
=
-(1+{\sf s}_{x^2}^*)\,\mathcal{Y}_{x^1x^1}
+
{\sf r}\,\mathcal{Y}_{x^1}.
\endaligned
\end{equation}
This is $\text{\rm (8.86)}_7$.
Differentiating $\text{\rm (8.91)}_2$
with respect to 
$x^1$, taking account of 
$\text{\rm (8.83)}_1$, we get
$\text{\rm (8.86)}_{ 16}$:
\def\theequation{8.92}\begin{equation}
\mathcal{X}_{x^1x^2x^2}^2
=
{\sf r}\,\mathcal{Y}_{x^1}
+
{\sf r}\,\mathcal{Y}_{x^1x^1}.
\end{equation}
We then replace $\mathcal{ X}_{ x^1}^1$ from~\thetag{ 8.90}
in  
$\text{\rm (8.81)}_4$:
\def\theequation{8.93}\begin{equation}
\aligned
0
&
=
\mathcal{X}_{x^1}^2
+
2\,\mathcal{X}_y^1
+
{\sf r}\,\mathcal{X}^1
+
{\sf r}\,\mathcal{X}^2
+
{\sf r}\,\mathcal{X}_{x^2}^2
+
{\sf r}\,\mathcal{Y}_{x^1}
+
{\sf r}\,\mathcal{Y}_y
+
\\
&
\ \ \ \ \ \ \ \ \ \ \ \ \ \ \ \ \ \ \ \ \ \ \ \
+
{\sf r}\,\mathcal{Y}_{x^1x^1}
+
{\sf s}^*\,\mathcal{Y}_{x^1y}
+
{\sf s}^*\,\mathcal{X}_{x^1x^1}^1.
\endaligned
\end{equation}
We differentiate this equation with respect to $x^2$, 
knowing $\mathcal{ Y}_{ x^2} = 0$:
\def\theequation{8.94}\begin{equation}
\aligned
0
&
=
\mathcal{X}_{x^1x^2}^2
+
2\,\mathcal{X}_{x^2y}^1
+
{\sf r}\,\mathcal{X}^1
+
{\sf r}\,\mathcal{X}_{x^2}^1
+
{\sf r}\,\mathcal{X}^2
+
{\sf r}\,\mathcal{X}_{x^2}^2
+
{\sf r}\,\mathcal{X}_{x^2x^2}^2
+
{\sf r}\,\mathcal{Y}_{x^1}
+
{\sf r}\,\mathcal{Y}_y
+
\\
&
\ \ \ \ \ \ \ \ \ \ \ \ \ \ \ \ \ \ \ \ \ \ \ \ \ \ \ \ \ \ 
+
{\sf r}\,\mathcal{Y}_{x^1x^1}
+
{\sf s}_{x^2}^*\,\mathcal{Y}_{x^1y}
+
{\sf s}_{x^2}^*\,\mathcal{X}_{x^1x^1}^1
+
{\sf s}^*\,\mathcal{X}_{x^1x^1x^2}.
\endaligned
\end{equation}
We replace: $\mathcal{ X}_{ x^2}^1$ from~\thetag{ 8.89};
$\mathcal{ X}_{ x^2 x^2 }^2$ from
$\text{\rm (8.91)}_2$; we differentiate
$\text{\rm (8.81)}_2$ with respect to $x^1x^1$ to 
replace $\mathcal{ X}_{ x^1 x^1 x^2}^1$ by ${\sf r}\,
\mathcal{ Y}_{ x^1x^1}$, thanks to $\text{\rm (8.83)}_1$; 
and we reorganize:
\def\theequation{8.95}\begin{equation}
2\,\mathcal{X}_{x^2y}^1
+
{\sf s}_{x^2}^*\,\mathcal{Y}_{x^1y}
+
{\sf s}_{x^2}^*\,\mathcal{X}_{x^1x^1}^1
=
-
\mathcal{X}_{x^1x^2}^2
+
{\sf r}\,\mathcal{X}^1
+
{\sf r}\,\mathcal{X}^2
+
{\sf r}\,\mathcal{X}_{x^2}^2
+
{\sf r}\,\mathcal{Y}_{x^1}
+
{\sf r}\,\mathcal{Y}_y
+
{\sf r}\,\mathcal{Y}_{x^1x^1}.
\end{equation}
We differentiate
$\text{\rm (8.81)}_2$ with respect to 
$y$ and 
$\text{\rm (8.81)}_3$ with respect to $x^1$:
\def\theequation{8.96}\begin{equation}
\aligned
\mathcal{X}_{x^2y}^1
+
\frac{1}{2}\,\mathcal{Y}_{x^1y}
&
=
0,
\\
\mathcal{Y}_{x^1y}
-
2\,\mathcal{X}_{x^1x^1}^1
&
=
-
\mathcal{X}_{x^1x^2}^2
+
{\sf r}\,\mathcal{Y}_{x^1}
+
{\sf r}\,\mathcal{Y}_{x^1x^1}.
\endaligned
\end{equation}
For the three unknowns $\mathcal{ X}_{ x^1 x^1}^1$, 
$\mathcal{ Y}_{ x^1 y}$, $\mathcal{ X}_{ x^ 2 y}^1$,
we solve the three equations
\thetag{ 8.95}, 
$\text{\rm (8.96)}_1$,
$\text{\rm (8.96)}_2$, 
reminding ${\sf s}_{ x^2}^* (0) = 0$:
\def\theequation{8.97}\begin{equation}
\aligned
\mathcal{X}_{x^1x^1}^1
&
=
{\sf r}\,\mathcal{X}^1
+
{\sf r}\,\mathcal{X}^2
+
{\sf r}\,\mathcal{X}_{x^2}^2
+
{\sf r}\,\mathcal{Y}_{x^1}
+
{\sf r}\,\mathcal{Y}_y
+
{\sf r}\,\mathcal{Y}_{x^1x^1}
+
{\sf r}\,\mathcal{X}_{x^1x^2}^2,
\\
\mathcal{Y}_{x^1y}
&
=
{\sf r}\,\mathcal{X}^1
+
{\sf r}\,\mathcal{X}^2
+
{\sf r}\,\mathcal{X}_{x^2}^2
+
{\sf r}\,\mathcal{Y}_{x^1}
+
{\sf r}\,\mathcal{Y}_y
+
{\sf r}\,\mathcal{Y}_{x^1x^1}
+
{\sf r}\,\mathcal{X}_{x^1x^2}^2,
\\
\mathcal{X}_{x^2y}^1
&
=
{\sf r}\,\mathcal{X}^1
+
{\sf r}\,\mathcal{X}^2
+
{\sf r}\,\mathcal{X}_{x^2}^2
+
{\sf r}\,\mathcal{Y}_{x^1}
+
{\sf r}\,\mathcal{Y}_y
+
{\sf r}\,\mathcal{Y}_{x^1x^1}
+
{\sf r}\,\mathcal{X}_{x^1x^2}^2.
\endaligned
\end{equation}
We get $\text{\rm (8.86)}_{11}$ and $\text{\rm (8.86)}_9$.

Thus, we may replace $\mathcal{ X}_{ x^1 x^1}^1$ and $\mathcal{ Y}_{
x^1 y}$ in $\text{\rm (8.81)}_4$ to get $\text{\rm (8.86)}_2$:
\def\theequation{8.98}\begin{equation}
\mathcal{X}_{x^1}^2
=
-2\,\mathcal{X}_y^1
+
{\sf r}\,\mathcal{X}^1
+
{\sf r}\,\mathcal{X}^2
+
{\sf r}\,\mathcal{X}_{x^2}^2
+
{\sf r}\,\mathcal{Y}_{x^1}
+
{\sf r}\,\mathcal{Y}_y
+
{\sf r}\,\mathcal{Y}_{x^1x^1}
+
{\sf r}\,\mathcal{X}_{x^1x^2}^2.
\end{equation}

Next, we differentiate $\text{\rm (8.83)}_3$ with respect to $x^1$ and
we replace: $\mathcal{ X}_{ x^1}^1$ from \thetag{ 8.90}; $\mathcal{
X}_{ x^1}^2$ from \thetag{ 8.98}; $\mathcal{ Y}_{ x^1y}$ from
$\text{\rm (8.97)}_2$; $\mathcal{ X}_{ x^1x^1}^1$ from $\text{\rm
(8.97)}_1$; $\mathcal{ Y}_{ x^1 x^1 x^1}$ from $\text{\rm (8.83)}_1$;
and we compare with $\text{\rm (8.83)}_2$; we differentiate $\text{\rm
(8.96)}_1$ with respect to $x^1$ and $\text{\rm (8.96)}_2$ with
respect to $x^1$; solving, we obtain four new relations:
\def\theequation{8.99}\begin{equation}
\aligned
\mathcal{X}_{x^1x^1x^1}^1
&
=
{\sf r}\,\mathcal{X}^1
+
{\sf r}\,\mathcal{X}^2
+
{\sf r}\,\mathcal{X}_{x^2}^2
+
{\sf r}\,\mathcal{Y}_{x^1}
+
{\sf r}\,\mathcal{Y}_y
+
{\sf r}\,\mathcal{Y}_{x^1x^1}
+
{\sf r}\,\mathcal{X}_{x^1x^2}^2,
\\
\mathcal{Y}_{x^1x^1y}
&
=
{\sf r}\,\mathcal{X}^1
+
{\sf r}\,\mathcal{X}^2
+
{\sf r}\,\mathcal{X}_{x^2}^2
+
{\sf r}\,\mathcal{Y}_{x^1}
+
{\sf r}\,\mathcal{Y}_y
+
{\sf r}\,\mathcal{Y}_{x^1x^1}
+
{\sf r}\,\mathcal{X}_{x^1x^2}^2,
\\
\mathcal{X}_{x^1x^2y}^2
&
=
{\sf r}\,\mathcal{X}^1
+
{\sf r}\,\mathcal{X}^2
+
{\sf r}\,\mathcal{X}_{x^2}^2
+
{\sf r}\,\mathcal{Y}_{x^1}
+
{\sf r}\,\mathcal{Y}_y
+
{\sf r}\,\mathcal{Y}_{x^1x^1}
+
{\sf r}\,\mathcal{X}_{x^1x^2}^2,
\\
\mathcal{X}_{x^1x^1x^2}^2
&
=
{\sf r}\,\mathcal{X}^1
+
{\sf r}\,\mathcal{X}^2
+
{\sf r}\,\mathcal{X}_{x^2}^2
+
{\sf r}\,\mathcal{Y}_{x^1}
+
{\sf r}\,\mathcal{Y}_y
+
{\sf r}\,\mathcal{Y}_{x^1x^1}
+
{\sf r}\,\mathcal{X}_{x^1x^2}^2.
\endaligned
\end{equation}
We get
$\text{\rm (8.86)}_{ 19}$ and $\text{\rm (8.86)}_{ 14}$.

Next, in $\text{\rm (8.81)}_5$, we replace: $\mathcal{ X}_{ x^1}^1$
from~\thetag{ 8.90}; $\mathcal{ X}_{ x^1}^2$ from~\thetag{
8.98}; $\mathcal{ Y}_{ x^1 y}$ from $\text{\rm (8.97)}_2$;
we get:
\def\theequation{8.100}\begin{equation}
\aligned
\mathcal{X}_y^2
&
=
{\sf r}\,\mathcal{X}^1
+
{\sf r}\,\mathcal{X}^2
+
{\sf r}\,\mathcal{X}_y^1
+
{\sf r}\,\mathcal{X}_{x^2}^2
+
{\sf r}\,\mathcal{Y}_{x^1}
+
{\sf r}\,\mathcal{Y}_y
+
{\sf r}\,\mathcal{Y}_{x^1x^1}
+
\\
&
\ \ \ \ \
+
{\sf s}^*\,\mathcal{X}_{x^1x^1}^2
+
{\sf s}^*\,\mathcal{Y}_{yy}
+
{\sf s}^*\,\mathcal{X}_{x^1y}^1.
\endaligned
\end{equation}
We differentiate~\thetag{ 8.98} with respect to $x^1$ and we replace:
$\mathcal{ X}_{ x^1}^1$ from \thetag{ 8.90}; $\mathcal{ X}_{ x^1}^2$
from \thetag{ 8.98}; $\mathcal{ Y}_{ x^1y}$ from $\text{ \rm
(8.97)}_2$; $\mathcal{ Y}_{ x^1 x^1 x^1}$ from $\text{ \rm (8.83)}_1$;
$\mathcal{ X}_{ x^1 x^1 x^2}^2$ from $\text{ \rm (8.99)}_4$; we get:
\def\theequation{8.101}\begin{equation}
\mathcal{X}_{x^1x^1}^2
+
2\,\mathcal{X}_{x^1y}^1
=
{\sf r}\,\mathcal{X}^1
+
{\sf r}\,\mathcal{X}^2
+
{\sf r}\,\mathcal{X}_y^1
+
{\sf r}\,\mathcal{X}_{x^2}^2
+
{\sf r}\,\mathcal{Y}_{x^1}
+
{\sf r}\,\mathcal{Y}_y
+
{\sf r}\,\mathcal{Y}_{x^1x^1}
+
{\sf r}\,\mathcal{X}_{x^1x^2}^2.
\end{equation} 
In $\text{\rm (8.82)}_2$, we replace: $\mathcal{ X}_{ x^1}^1$
from~\thetag{ 8.90}; $\mathcal{ X}_{ x^1 x^1}^1$ from $\text{\rm
(8.97)}_1$; $\mathcal{ Y}_{ x^1 y}$ from $\text{\rm (8.97)}_2$;
and we reorganize:
\def\theequation{8.102}\begin{equation}
2\,\mathcal{X}_{x^2y}^2
-
6\,\mathcal{X}_{x^1y}^1
-
\mathcal{X}_{x^1x^1}^2
=
-
2\,\mathcal{Y}_{yy}
+
{\sf r}\,\mathcal{X}^1
+
{\sf r}\,\mathcal{X}^2
+
{\sf r}\,\mathcal{X}_{x^2}^2
+
{\sf r}\,\mathcal{Y}_{x^1}
+
{\sf r}\,\mathcal{Y}_y
+
{\sf r}\,\mathcal{Y}_{x^1x^1}
+
{\sf r}\,\mathcal{X}_{x^1x^2}^2.
\end{equation}
Differentiating $\text{\rm (8.81)}_3$ with respect to $y$, we replace:
$\mathcal{ Y}_{ x^1 y}$ from $\text{\rm (8.97)}_2$;
$\mathcal{ Y}_{ x^1x^1 y}$ from $\text{\rm (8.99)}_2$; 
and we reorganize:
\def\theequation{8.103}\begin{equation}
\mathcal{X}_{x^2y}^2
-
2\,\mathcal{X}_{x^1y}^1
=
-
\mathcal{Y}_{yy}
+
{\sf r}\,\mathcal{X}^1
+
{\sf r}\,\mathcal{X}^2
+
{\sf r}\,\mathcal{X}_{x^2}^2
+
{\sf r}\,\mathcal{Y}_{x^1}
+
{\sf r}\,\mathcal{Y}_y
+
{\sf r}\,\mathcal{Y}_{x^1x^1}
+
{\sf r}\,\mathcal{X}_{x^1x^2}^2.
\end{equation}
For the three unknowns $\mathcal{ X}_{ x^1 x^1}^2$, 
$\mathcal{ X}_{ x^1y}^1$, $\mathcal{ X}_{ x^2 y}^2$, 
we then solve the four equations
\thetag{ 8.101},
\thetag{ 8.102},
\thetag{ 8.103},
$\text{\rm (8.83)}_4$ (in which we replace: $\mathcal{ X}_{ x^1}^1$
from \thetag{ 8.90}; $\mathcal{ X}_{ x^1}^2$ from \thetag{ 8.98};
$\mathcal{ Y}_{ x^1y}$ from $\text{ \rm (8.97)}_2$; $\mathcal{ X}_{
x^1 x^1}^1$ from$\text{ \rm (8.97)}_1$):
\def\theequation{8.104}\begin{equation}
\aligned
\mathcal{X}_{x^1x^1}^2
&
=
{\sf r}\,\mathcal{X}^1
+
{\sf r}\,\mathcal{X}^2
+
{\sf r}\,\mathcal{X}_y^1
+
{\sf r}\,\mathcal{X}_{x^2}^2
+
{\sf r}\,\mathcal{Y}_{x^1}
+
{\sf r}\,\mathcal{Y}_y
+
{\sf r}\,\mathcal{Y}_{x^1x^1}
+
{\sf r}\,\mathcal{X}_{x^1x^2}^2
+
{\sf r}\,\mathcal{Y}_{yy},
\\
\mathcal{X}_{x^1y}^1
&
=
{\sf r}\,\mathcal{X}^1
+
{\sf r}\,\mathcal{X}^2
+
{\sf r}\,\mathcal{X}_y^1
+
{\sf r}\,\mathcal{X}_{x^2}^2
+
{\sf r}\,\mathcal{Y}_{x^1}
+
{\sf r}\,\mathcal{Y}_y
+
{\sf r}\,\mathcal{Y}_{x^1x^1}
+
{\sf r}\,\mathcal{X}_{x^1x^2}^2
+
{\sf r}\,\mathcal{Y}_{yy},
\\
\mathcal{X}_{x^2y}^2
&
=
{\sf r}\,\mathcal{X}^1
+
{\sf r}\,\mathcal{X}^2
+
{\sf r}\,\mathcal{X}_y^1
+
{\sf r}\,\mathcal{X}_{x^2}^2
+
{\sf r}\,\mathcal{Y}_{x^1}
+
{\sf r}\,\mathcal{Y}_y
+
{\sf r}\,\mathcal{Y}_{x^1x^1}
+
{\sf r}\,\mathcal{X}_{x^1x^2}^2
+
{\sf r}\,\mathcal{Y}_{yy}.
\endaligned
\end{equation}
We get 
$\text{\rm (8.86)}_6$ and
$\text{\rm (8.86)}_{ 10}$.
Replacing then $\mathcal{ X}_{ x^1 x^1}^2$, 
$\mathcal{ X}_{ x^1 y}^1$ in~\thetag{ 
8.100} gives
\def\theequation{8.105}\begin{equation}
\mathcal{X}_y^2
=
{\sf r}\,\mathcal{X}^1
+
{\sf r}\,\mathcal{X}^2
+
{\sf r}\,\mathcal{X}_y^1
+
{\sf r}\,\mathcal{X}_{x^2}^2
+
{\sf r}\,\mathcal{Y}_{x^1}
+
{\sf r}\,\mathcal{Y}_y
+
{\sf r}\,\mathcal{Y}_{x^1x^1}
+
{\sf r}\,\mathcal{X}_{x^1x^2}^2
+
{\sf r}\,\mathcal{Y}_{yy}.
\end{equation}
This is $\text{\rm (8.86)}_5$.

Next, we differentiate \thetag{ 8.103} with respect to 
$x^1$ and we replace: 
$\mathcal{ X}_{ x^1}^1$ from \thetag{ 8.90}; 
$\mathcal{ X}_{ x^1}^2$ from \thetag{ 8.98}; 
$\mathcal{ Y}_{ x^1y}$ from $\text{\rm (8.97)}_2$;
$\mathcal{ Y}_{x^1x^1x^1}$ from $\text{\rm (8.83)}_1$;
$\mathcal{ X}_{x^1x^1x^2}^2$ from $\text{\rm (8.99)}_4$; 
we get:
\def\theequation{8.106}\begin{equation}
\mathcal{Y}_{x^1yy}
+
\mathcal{X}_{x^1x^2y}^2
-
2\,\mathcal{X}_{x^1x^1y}^1
=
{\sf r}\,\mathcal{X}^1
+
{\sf r}\,\mathcal{X}^2
+
{\sf r}\,\mathcal{X}_{x^2}^2
+
{\sf r}\,\mathcal{Y}_{x^1}
+
{\sf r}\,\mathcal{Y}_y
+
{\sf r}\,\mathcal{Y}_{x^1x^1}
+
{\sf r}\,\mathcal{X}_{x^1x^2}^2.
\end{equation}
Also, we differentiate 
$\text{\rm (8.83)}_3$ with respect to $y$ and
we replace: 
$\mathcal{ X}_y^2$ from \thetag{ 8.105}; 
$\mathcal{ Y}_{ x^1y}$ from $\text{\rm (8.97)}_2$;
$\mathcal{ X}_{ x^1y}^1$ from
$\text{\rm (8.104)}_2$; $\mathcal{ Y}_{ x^1 x^1 y}$ from
$\text{\rm (8.99)}_2$; we get:
\def\theequation{8.107}\begin{equation}
\mathcal{Y}_{x^1yy}
-
\mathcal{X}_{x^1x^1y}^1
=
{\sf r}\,\mathcal{X}^1
+
{\sf r}\,\mathcal{X}^2
+
{\sf r}\,\mathcal{X}_y^1
+
{\sf r}\,\mathcal{X}_{x^2}^2
+
{\sf r}\,\mathcal{Y}_{x^1}
+
{\sf r}\,\mathcal{Y}_y
+
{\sf r}\,\mathcal{Y}_{x^1x^1}
+
{\sf r}\,\mathcal{X}_{x^1x^2}^2
+
{\sf r}\,\mathcal{Y}_{yy}.
\end{equation}
Also, we replace in $\text{\rm (8.82)}_3$:
$\mathcal{ X}_{ x^1}^1$
from \thetag{ 8.90};  
$\mathcal{ X}_{ x^1}^2$
from \thetag{ 8.98}; 
$\mathcal{ Y}_{ x^1y}$ from
$\text{\rm (8.97)}_2$;
$\mathcal{ X}_{ x^1x^1}^1$ from 
$\text{\rm (8.97)}_1$;
$\mathcal{ X}_{ x^1x^1}^2$ from 
$\text{\rm (8.104)}_1$;
$\mathcal{ X}_{x^1y}^1$ from 
$\text{\rm (8.104)}_2$; we get:
\def\theequation{8.108}\begin{equation}
4\,\mathcal{X}_{yy}^1
+
3\,\mathcal{X}_{x^1y}^2
=
{\sf r}\,\mathcal{X}^1
+
{\sf r}\,\mathcal{X}^2
+
{\sf r}\,\mathcal{X}_y^1
+
{\sf r}\,\mathcal{X}_{x^2}^2
+
{\sf r}\,\mathcal{Y}_{x^1}
+
{\sf r}\,\mathcal{Y}_y
+
{\sf r}\,\mathcal{Y}_{x^1x^1}
+
{\sf r}\,\mathcal{X}_{x^1x^2}^2
+
{\sf r}\,\mathcal{Y}_{yy}.
\end{equation}
We differentiate this equation with respect to 
$x^2$ and we replace: $4\, \mathcal{ X}_{ x^2 yy}^1$
by $-2\, \mathcal{ Y}_{ x^1 yy}$ from
\thetag{ 8.89};
$\mathcal{ X}_{ x^2}^1$ from \thetag{ 8.98};
$\mathcal{ X}_{ x^2 y}^1$ from 
$\text{\rm (8.97)}_3$;
(notice $0 = \mathcal{ Y}_{ x^1 x^2} = 
\mathcal{ Y}_{ x^2 y}$); 
$\mathcal{ X}_{ x^2 x^2}^2$ from 
$\text{\rm (8.91)}_2$; 
$\mathcal{ X}_{ x^1 x^2 x^2}^2$ from 
\thetag{ 8.92};
we get:
\def\theequation{8.109}\begin{equation}
-2\,\mathcal{Y}_{x^1yy}
+
3\,\mathcal{X}_{x^1x^2y}^2
=
{\sf r}\,\mathcal{X}^1
+
{\sf r}\,\mathcal{X}^2
+
{\sf r}\,\mathcal{X}_y^1
+
{\sf r}\,\mathcal{X}_{x^2}^2
+
{\sf r}\,\mathcal{Y}_{x^1}
+
{\sf r}\,\mathcal{Y}_y
+
{\sf r}\,\mathcal{Y}_{x^1x^1}
+
{\sf r}\,\mathcal{X}_{x^1x^2}^2
+
{\sf r}\,\mathcal{Y}_{yy}.
\end{equation}
For the three unknowns $\mathcal{ X}_{ x^1 x^1 y}^1$, 
$\mathcal{ Y}_{ x^1 yy}$, $\mathcal{ X}_{ x^1 x^2 y}^2$, 
we solve the three equations 
\thetag{ 8.106}, \thetag{ 8.107}, \thetag{ 8.108}; 
we get:
\def\theequation{8.110}\begin{equation}
\aligned
\mathcal{X}_{x^1x^1y}^1
&
=
{\sf r}\,\mathcal{X}^1
+
{\sf r}\,\mathcal{X}^2
+
{\sf r}\,\mathcal{X}_y^1
+
{\sf r}\,\mathcal{X}_{x^2}^2
+
{\sf r}\,\mathcal{Y}_{x^1}
+
{\sf r}\,\mathcal{Y}_y
+
{\sf r}\,\mathcal{Y}_{x^1x^1}
+
{\sf r}\,\mathcal{X}_{x^1x^2}^2
+
{\sf r}\,\mathcal{Y}_{yy}, 
\\
\mathcal{Y}_{x^1yy}
&
=
{\sf r}\,\mathcal{X}^1
+
{\sf r}\,\mathcal{X}^2
+
{\sf r}\,\mathcal{X}_y^1
+
{\sf r}\,\mathcal{X}_{x^2}^2
+
{\sf r}\,\mathcal{Y}_{x^1}
+
{\sf r}\,\mathcal{Y}_y
+
{\sf r}\,\mathcal{Y}_{x^1x^1}
+
{\sf r}\,\mathcal{X}_{x^1x^2}^2
+
{\sf r}\,\mathcal{Y}_{yy},
\\
\mathcal{X}_{x^1x^2y}^2
&
=
{\sf r}\,\mathcal{X}^1
+
{\sf r}\,\mathcal{X}^2
+
{\sf r}\,\mathcal{X}_y^1
+
{\sf r}\,\mathcal{X}_{x^2}^2
+
{\sf r}\,\mathcal{Y}_{x^1}
+
{\sf r}\,\mathcal{Y}_y
+
{\sf r}\,\mathcal{Y}_{x^1x^1}
+
{\sf r}\,\mathcal{X}_{x^1x^2}^2
+
{\sf r}\,\mathcal{Y}_{yy}.
\endaligned
\end{equation}
We get $\text{\rm (8.86)}_{ 20}$ and $\text{\rm (8.86)}_{ 18}$.

Next, in \thetag{ 8.93}, we replace: 
$\mathcal{ Y}_{ x^1y}$ from 
$\text{\rm (8.97)}_2$; 
$\mathcal{ X}_{ x^1 x^1}^1$ from 
$\text{\rm (8.97)}_1$; we get:
\def\theequation{8.111}\begin{equation}
\mathcal{X}_{x^1}^2
+
2\,\mathcal{X}_y^1
=
{\sf r}\,\mathcal{X}^1
+
{\sf r}\,\mathcal{X}^2
+
{\sf r}\,\mathcal{X}_{x^2}^2
+
{\sf r}\,\mathcal{Y}_{x^1}
+
{\sf r}\,\mathcal{Y}_y
+
{\sf r}\,\mathcal{Y}_{x^1x^1}
+
{\sf r}\,\mathcal{X}_{x^1x^2}^2.
\end{equation}
We differentiate this equation with respect to $y$ and
we replace: 
$\mathcal{ X}_y^2$ from 
\thetag{ 8.105};
$\mathcal{ X}_{ x^2y}^2$ from 
$\text{\rm (8.104)}_3$; 
$\mathcal{Y}_{x^1y}$ from  
$\text{\rm (8.97)}_2$;
$\mathcal{Y}_{x^1x^1y}$ from 
$\text{\rm (8.99)}_2$;
$\mathcal{ X}_{ x^1 x^2 y}^2$ from
$\text{\rm (8.99)}_3$; 
we get:
\def\theequation{8.112}\begin{equation}
\mathcal{X}_{x^1y}^2
+
2\,\mathcal{X}_{yy}^1
=
{\sf r}\,\mathcal{X}^1
+
{\sf r}\,\mathcal{X}^2
+
{\sf r}\,\mathcal{X}_y^1
+
{\sf r}\,\mathcal{X}_{x^2}^2
+
{\sf r}\,\mathcal{Y}_{x^1}
+
{\sf r}\,\mathcal{Y}_y
+
{\sf r}\,\mathcal{Y}_{x^1x^1}
+
{\sf r}\,\mathcal{X}_{x^1x^2}^2
+
{\sf r}\,\mathcal{Y}_{yy}.
\end{equation}
For the two unknowns $\mathcal{ X}_{ yy}^1$ and
$\mathcal{ X}_{ x^1 y}^2$, we 
solve the two equations 
\thetag{ 8.108} and \thetag{ 8.112}; 
we get:
\def\theequation{8.113}\begin{equation}
\aligned
\mathcal{X}_{yy}^1
&
=
{\sf r}\,\mathcal{X}^1
+
{\sf r}\,\mathcal{X}^2
+
{\sf r}\,\mathcal{X}_y^1
+
{\sf r}\,\mathcal{X}_{x^2}^2
+
{\sf r}\,\mathcal{Y}_{x^1}
+
{\sf r}\,\mathcal{Y}_y
+
{\sf r}\,\mathcal{Y}_{x^1x^1}
+
{\sf r}\,\mathcal{X}_{x^1x^2}^2
+
{\sf r}\,\mathcal{Y}_{yy},
\\
\mathcal{X}_{x^1y}^2
&
=
{\sf r}\,\mathcal{X}^1
+
{\sf r}\,\mathcal{X}^2
+
{\sf r}\,\mathcal{X}_y^1
+
{\sf r}\,\mathcal{X}_{x^2}^2
+
{\sf r}\,\mathcal{Y}_{x^1}
+
{\sf r}\,\mathcal{Y}_y
+
{\sf r}\,\mathcal{Y}_{x^1x^1}
+
{\sf r}\,\mathcal{X}_{x^1x^2}^2
+
{\sf r}\,\mathcal{Y}_{yy}.
\endaligned
\end{equation}
We get $\text{\rm (8.86)}_{12}$.

Next, we differentiate $\text{\rm (8.113)}_1$
with respect to $x^1$ and we replace:
$\mathcal{ X}_{ x^1}^1$, 
$\mathcal{ X}_{ x^1}^2$, 
$\mathcal{ X}_{ x^1y}^1$, 
$\mathcal{ Y}_{ x^1y}$, 
$\mathcal{ Y}_{ x^1x^1x^1}$, 
$\mathcal{ X}_{ x^1x^1x^2}^2$, 
$\mathcal{ Y}_{ x^1yy}$; we get:
\def\theequation{8.114}\begin{equation}
\mathcal{X}_{x^1yy}^1
=
{\sf r}\,\mathcal{X}^1
+
{\sf r}\,\mathcal{X}^2
+
{\sf r}\,\mathcal{X}_y^1
+
{\sf r}\,\mathcal{X}_{x^2}^2
+
{\sf r}\,\mathcal{Y}_{x^1}
+
{\sf r}\,\mathcal{Y}_y
+
{\sf r}\,\mathcal{Y}_{x^1x^1}
+
{\sf r}\,\mathcal{X}_{x^1x^2}^2
+
{\sf r}\,\mathcal{Y}_{yy}.
\end{equation}
Also, we differentiate $\text{\rm (8.113)}_2$ with respect to $x^1$
and we replace:
\def\theequation{8.115}\begin{equation}
\mathcal{X}_{x^1x^1y}^2
=
{\sf r}\,\mathcal{X}^1
+
{\sf r}\,\mathcal{X}^2
+
{\sf r}\,\mathcal{X}_y^1
+
{\sf r}\,\mathcal{X}_{x^2}^2
+
{\sf r}\,\mathcal{Y}_{x^1}
+
{\sf r}\,\mathcal{Y}_y
+
{\sf r}\,\mathcal{Y}_{x^1x^1}
+
{\sf r}\,\mathcal{X}_{x^1x^2}^2
+
{\sf r}\,\mathcal{Y}_{yy}.
\end{equation}
Also, we differentiate 
$\text{\rm (8.83)}_4$ with respect to $y$; we
replace $\mathcal{X}_{ x^1 yy}^1$ from \thetag{ 8.114}, 
we replace $\mathcal{X}_{x^1x^1y}^2$ from \thetag{ 8.115}; 
and
we achieve other evident replacements; we get:
\def\theequation{8.116}\begin{equation}
\mathcal{Y}_{yyy}
=
{\sf r}\,\mathcal{X}^1
+
{\sf r}\,\mathcal{X}^2
+
{\sf r}\,\mathcal{X}_y^1
+
{\sf r}\,\mathcal{X}_{x^2}^2
+
{\sf r}\,\mathcal{Y}_{x^1}
+
{\sf r}\,\mathcal{Y}_y
+
{\sf r}\,\mathcal{Y}_{x^1x^1}
+
{\sf r}\,\mathcal{X}_{x^1x^2}^2
+
{\sf r}\,\mathcal{Y}_{yy}.
\end{equation}
This is $\text{\rm (8.96)}_{22}$, which completes the proof.
\endproof

\def\thetheorem{8.117}\begin{theorem}
The bound $\dim \mathfrak{ SYM} (\mathcal{ E}_5) \leqslant 10$ is
attained if and only if \thetag{ $\mathcal{ E}_5$} is equivalent, 
through a diffeomorphism $(x^1, x^2, y) \longmapsto 
(X^1, X^2, Y)$, to
the model system 
\def\theequation{8.118}\begin{equation}
Y_{X^2}
=
0, 
\ \ \ \ \ \ \ \ \ \
Y_{X^1X^1X^1}
=
0.
\end{equation}
\end{theorem}

\proof
Firstly, setting ${\sf r} = {\sf s}^*= 0$ everywhere, the
solution to \thetag{ 8.81}, \thetag{ 8.82}, \thetag{ 8.83} is
\def\theequation{8.119}\begin{equation}
\aligned
\mathcal{X}^1
&
=
k
+
(c+j)\,x^1
-
b\,x^2
-
h\,y
+
e\,x^1x^1
-
d\,x^1x^2
+
f\,x^1y
-
e\,x^2y,
\\
\mathcal{X}^2
&
=
g
+
2h\,x^1
+
2j\,x^2
-
d\,x^2x^2
+
2e\,x^1x^2
-
f\,x^1x^1,
\\
\mathcal{Y}
&
=
a
+
2b\,x^1
+
2c\,y
+
d\,x^1x^1
+
2e\,x^1y
+
f\,yy,
\endaligned
\end{equation}
where $a, b, c, d, e, f, g, h, j, k \in \K$ are arbitrary. Computing
the third prolongations of the ten vector fields
\begin{small}
\def\theequation{8.120}\begin{equation}
\aligned
&
\frac{\partial}{\partial x^1},
\ \ \ \ \ \ \ \ \ \ \
\frac{\partial}{\partial x^2},
\ \ \ \ \ \ \ \ \ \ \
\frac{\partial}{\partial y},
\\
&
-x^2\,\frac{\partial}{\partial x^1}
+
2x^1\,\frac{\partial}{\partial y},
\ \ \ \ \ \ \ \ 
x^1\,\frac{\partial}{\partial x^1}
+
2y\,\frac{\partial}{\partial y},
\ \ \ \ \ \ \ \ 
x^1\,\frac{\partial}{\partial x^1}
+
2x^2\,\frac{\partial}{\partial x^2},
\ \ \ \ \ \ \ \ 
-y\,\frac{\partial}{\partial x^1}
+
2x^1\,\frac{\partial}{\partial x^2},
\\
&
-x^1x^2\,\frac{\partial}{\partial x^1}
-
x^2x^2\,\frac{\partial}{\partial x^2}
+
x^1x^1\,\frac{\partial}{\partial y},
\ \ \ \ \ \ \
(x^1x^1-x^2y)\,\frac{\partial}{\partial x^1}
+
2x^1x^2\,\frac{\partial}{\partial x^2}
+
2x^1y\,\frac{\partial}{\partial y},
\\
&
x^1y\,\frac{\partial}{\partial x^1}
-
x^1x^1\,\frac{\partial}{\partial x^2}
+
yy\,\frac{\partial}{\partial y}
\endaligned
\end{equation}
\end{small}

\noindent
one verifies that they all are tangent to the skeleton $y_2 = \frac{
1}{ 4} \, (y_1)^2$, $y_{ 1, 1, 1} = 0$. Thus the bound is attained.
One then verifies (\cite{
fk2005a}) that the spanned Lie algebra is isomorphic to
$\mathfrak{ so} (5, \C)$.

\def\thelemma{8.121}\begin{lemma}
Assuming the normalizations of Lemma~8.54, the remainder
${\sf O}_4$ in~\thetag{ 8.53} is an ${\rm O}_3 (x^1, a^1)${\rm :}
\def\theequation{8.122}\begin{equation}
y
=
b
+
\frac{2\,x^1a^1+x^1x^1a^2+a^1a^1x^2}{1-x^2a^2}
+
(x^1)^3\,{\sf R}
+
(x^1)^2a^1\,{\sf R}
+
x^1(a^1)^2\,{\sf R}
+
(a^1)^3\,{\sf R}.
\end{equation}
\end{lemma}

\proof
Indeed, writing
\def\theequation{8.123}\begin{equation}
y
=
b
+
x^1\,\Lambda^{1,0}
+
a^1\,\Lambda^{0,1}
+
x^1x^1\,\Lambda^{2,0}
+
x^1a^1\,\Lambda^{1,1}
+
a^1a^1\,\Lambda^{0,2}
+
{\rm O}_3(x^1,a^1),
\end{equation}
with $\Lambda^{ i,j} = \Lambda^{ i, j} (x^2, a^2)$, 
and developing the determinant~\thetag{ 8.54} 
with respect to the powers of 
$(x^1, a^1)$, the vanishing of the 
coefficients of ${\rm cst.}$, of $x^1$, of $a^1$ 
yields the system
\def\theequation{8.124}\begin{equation}
\left\{
\aligned
0
&
\equiv
\Lambda_{a^2}^{1,0}\,\Lambda_{x^2}^{0,1},
\\
0
&
\equiv
\Lambda^{1,1}\,\Lambda_{x^2a^2}^{1,0}
-
2\,\Lambda_{a^2}^{2,0}\,\Lambda_{x^2}^{0,1}
-
\Lambda_{x^2}^{1,1}\,\Lambda_{a^2}^{1,0},
\\
0
&
\equiv
\Lambda^{1,1}\,\Lambda_{x^2a^2}^{0,1}
-
\Lambda_{a^2}^{1,1}\,\Lambda_{x^2}^{0,1}
-
2\,\Lambda_{x^2}^{0,2}\,\Lambda_{a^2}^{1,0}.
\endaligned\right.
\end{equation}
If the first equation yields $\Lambda_{ a^2}^{ 1, 0} \equiv 0$,
replacing in the second, using
$\Lambda^{ 2, 0} = a^2 + {\rm O}_2$, we deduce that $\Lambda_{ x^2}^{
0, 1} \equiv 0$ also. Similarly, $\Lambda_{ x^2}^{ 0, 1} \equiv 0$
implies $\Lambda_{ a^2}^{ 1, 0} \equiv 0$. Since the coordinate system
satisfies the normalization $\Pi (0, a ) \equiv \Pi (x, 0) \equiv 0$,
necessarily $\Lambda^{ 1, 0} = {\rm O} (a^2)$ and $\Lambda^{ 0, 1} =
{\rm O} (x^2)$. We deduce:
\def\theequation{8.125}\begin{equation}
0
\equiv
\Lambda^{1,0}
\equiv
\Lambda^{0,1}.
\end{equation}
Redeveloping the determinant, the vanishing of the coefficients of
$x^1 x^1$, of $x^1 a^1$, of $a^1 a^1$ yields the system
\def\theequation{8.126}\begin{equation}
\left\{
\aligned
0
&
\equiv
\Lambda^{1,1}\,\Lambda_{x^2a^2}^{2,0}
-
2\,\Lambda_{a^2}^{2,0}\,\Lambda_{x^2}^{1,1},
\\
0
&
\equiv
\Lambda^{1,1}\,\Lambda_{x^2a^2}^{1,1}
-
\Lambda_{a^2}^{1,1}\,\Lambda_{x^2}^{1,1}
-
4\,\Lambda_{a^2}^{2,0}\,\Lambda_{x^2}^{0,2},
\\
0
&
\equiv
\Lambda^{1,1}\,\Lambda_{x^2a^2}^{0,2}
-
2\,\Lambda_{x^2}^{0,2}\,\Lambda_{a^2}^{1,1}.
\endaligned\right.
\end{equation}
Since $\Lambda^{ 1, 1} (0) = 2 \neq 0$, we may divide by $\Lambda^{ 1,
1}$, obtaining a {\sc pde} system with the three functions $\Lambda_{
x^2 a^2}^{ 2, 0}$, $\Lambda_{ x^2 a^2}^{ 1, 1}$, $\Lambda_{ x^2 a^2}^{
0, 2}$ in the left hand side. We observe that the normalizations of
Lemma~8.55 entail
\def\theequation{8.127}\begin{equation}
\Lambda^{2,0}
=
a^2
+
{\rm O}(x^2a^2),
\ \ \ \ \ \ \ \ \ \ \
\Lambda^{1,1}
=
2
+
{\rm O}(x^2a^2),
\ \ \ \ \ \ \ \ \ \ \
\Lambda^{0,2}
=
x^2
+
{\rm O}(x^2a^2).
\end{equation}
By cross differentiations in the mentioned {\sc pde} system, it
follows that all the Taylor coefficients of $\Lambda^{ 2, 0}$,
$\Lambda^{ 1, 1}$, $\Lambda^{ 0, 2}$ are uniquely determined. As
already discovered in~\cite{ gm2003b}, the unique solution
\def\theequation{8.128}\begin{equation}
\Lambda^{2,0}
=
\frac{a^2}{1-x^2a^2},
\ \ \ \ \ \ \ \ \ \ \
\Lambda^{1,1}
=
\frac{2}{1-x^2a^2},
\ \ \ \ \ \ \ \ \ \ \
\Lambda^{0,2}
=
\frac{x^2}{1-x^2a^2},
\end{equation}
guarantees, when the remainder ${\rm O}_3 (x^1, a^1)$ vanishes, that
the determinant~\thetag{ 8.45} indeed vanishes identically.
\endproof

Conversely, suppose that $\dim \mathfrak{ SYM} (\mathcal{ E}_5) = 
10$ is maximal.

With $\varepsilon \neq 0$ small, replacing $(x^1, x^2, y, a^1, a^2,
b)$ by $(\varepsilon x^1, x^2, \varepsilon^2 y, \varepsilon a^1, a^2,
\varepsilon \varepsilon b)$ in~\thetag{ 8.122} and dividing by
$\varepsilon \varepsilon$, the remainder terms become small:
\def\theequation{8.129}\begin{equation}
y
=
b
+
\frac{2\,x^1a^1+x^1x^1a^2+a^1a^1x^2}{1-x^2a^2}
+
{\rm O}(\varepsilon).
\end{equation}
Then all the remainders in the equations $\Delta_{ \mathcal{ E}_5}$ of
the skeleton are ${\rm O}( \varepsilon )$. We get ten generators
similar to~\thetag{ 8.120}, plus an ${\rm O}( \varepsilon )$
perturbation. Thanks to the rigidity of $\mathfrak{ so} (5, \C)$,
Theorem~5.15 provides a change of generators, close to the $10 \times
10$ identity matrix, leading to the same structure constants as those
of the ten vector fields~\thetag{ 8.120}. As in the end of the proof
of Theorem~5.13, we may then straighten some relevant vector fields
(exercise) and finally check that their tangency to the skeleton
implies that it is the model one. Theorem~8.117 is proved.  
\endproof

\def\thecorollary{8.130}\begin{corollary}
Let $M \subset \C^3$ be a connected real analytic hypersurface
whose Levi form has uniform rank 1 that is 2-nondegenerate at every
point. Then
\def\theequation{8.131}\begin{equation}
\dim\mathfrak{hol}(M)
\leqslant 
10,
\end{equation}
and the bound is attained if and only if $M$ is locally, in a
neighborhood of Zariski-generic points, biholomorphic to the model
$M_0$.
\end{corollary}

\section*{ \S9.~Dual system of partial differential equations}

\subsection*{9.1.~Solvability with respect to the variables}
Let $\mathcal{ M}$ be as in \S2.10 defined by $y = \Pi (x, a, b)$ or
dually by $b = {\Pi^*} (a, x, y)$.

\def\thedefinition{9.2}\begin{definition}{\rm
$\mathcal{M}$ is {\sl solvable with respect to the variables} if there
exist an integer $\kappa^* \geqslant 1$ and multiindices $\delta(1),
\dots, \delta(n) \in \N^p$ with $\vert \delta(l)\vert \geqslant 1$ for
$l = 1, \dots, n$ and $\max_{ 1\leqslant l \leqslant n} \, \vert
\delta (l) \vert = \kappa^*$, together with integers $j(1), \dots,
j(n)$ with $1 \leqslant j(l) \leqslant m$ such that the local
$\K$-analytic map
\def\theequation{9.3}\begin{equation}
\K^{n+m}
\ni
(x,y) 
\longmapsto 
\Big(
\big(
{\Pi^*}^j(0,x,y)
\big)^{1\leqslant j\leqslant m},\
\Big(
{\Pi^*}_{a^{\delta(l)}}^{j(l)}(0,x,y)
\Big)_{1\leqslant l\leqslant n}
\Big)
\in
\K^{m+n}
\end{equation}
is of rank equal to $n+m$ at $(x,y) = (0,0)$
}\end{definition} 

If $\mathcal{ M }$ is a complexified generic submanifold, solvability
with respect to the parameters is equivalent to solvability with
respect to the variables, because ${\Pi^*} = \overline{ \Pi }$. This
is untrue in general: with $n = 2$, $m = 1$, consider the system $y_{
x^2} = 0$, $y_{ x^1 x^1} = 0$, whose general solutions is $y(x) = b +
x_1 a$ with $x^2$ absent.

To characterize generally such a degeneration property, we
develope both
\def\theequation{9.4}\begin{equation}
\aligned
y^j
&
=
\Pi^j(x,a,b)
=
\sum_{\beta\in\N^n}\,
x^\beta\,\Pi_\beta^j(a,b)
\ \ \ \ \ \ \ \
\text{\rm and}
\\
b^j
&
=
{\Pi^*}^j(a,x,y)
=
\sum_{\delta\in\N^p}\,
a^\delta\,{\Pi^*}_\delta^j(x,y),
\endaligned
\end{equation}
with analytic functions $\Pi_\beta^j (a, b)$, ${\Pi^* }_\delta^j
(x,y)$ and we introduce two $\K^\infty$-valued maps
\def\theequation{9.5}\begin{equation}
\aligned
\mathcal{Q}_\infty: \ \ \ 
(a,b)
&
\longmapsto 
\left(
\Pi_\beta^j(a,b)
\right)_{\beta\in\N^n}^{
1\leqslant j\leqslant m}
\ \ \ \ \ \ \ \
\text{\rm and}
\\
\mathcal{Q}_\infty^*: \ \ \
(x,y)
&
\longmapsto
\left(
{\Pi^*}_\delta^j(x,y)
\right)_{
\delta\in\N^p}^{
1\leqslant j\leqslant m}.
\endaligned
\end{equation}
Since $b \mapsto \big( \Pi_0^j (0, b) \big)^{1 \leqslant j \leqslant
m}$ and $y \mapsto \big( {\Pi^*}_0^j (0, y) \big)^{1 \leqslant j
\leqslant m}$ are already both of rank $m$ at the origin, the generic
ranks of these two maps, defined by testing the nonvanishing of minors
of their infinite Jacobian matrices, satisfy
\def\theequation{9.6}\begin{equation}
\aligned
{\rm genrk}\,\mathcal{Q}_\infty
&
=
m+p_\mathcal{M}
\ \ \ \ \ \ \ \
\text{\rm and}
\\
{\rm genrk}\,\mathcal{Q}_\infty^*
&
=
m+n_\mathcal{M}
\endaligned
\end{equation}
for some two integers $0 \leqslant p_\mathcal{ M } \leqslant p$ and $0
\leqslant n_\mathcal{ M } \leqslant n$. So at a Zariski-generic
point, the ranks are {\it equal}\, to $m + p_\mathcal{ M }$ and to $m
+ n_\mathcal{ M }$. 

\def\theproposition{9.7}\begin{proposition}
There exists a local proper $\K$-analytic subset $\Sigma_\mathcal{ M}$
of $\K_x^n \times \K_y^m \times \K_a^p \times \K_b^m$
whose equations, of the specific form
\def\theequation{9.8}\begin{equation}
\Sigma_\mathcal{M}
=
\Big\{
r_\nu(a,b)
=
0,\
\nu\in\N,
\ \ \ \ \
r_\mu^*(x,y)
=
0
\ \ \ \ \ 
\mu\in\N
\Big\},
\end{equation}
are obtained by equating to zero all $(m + p_\mathcal{ M }) \times (m
+ p_\mathcal{ M })$ minors of ${\rm Jac}\, \mathcal{ Q}_\infty$ and
all $(m + n_\mathcal{ M }) \times (m + n_\mathcal{ M })$ minors of
${\rm Jac}\, \mathcal{ Q }_\infty^*$, such that for every point $p = 
(x_p,
y_p, a_p, b_p) \not\in \Sigma_\mathcal{ M}$, there exists a local
change of coordinates respecting the separation of the variables $(x,
y)$ and $(a,b)$
\def\theequation{9.9}\begin{equation}
(x,y,a,b)
\mapsto
\big(
\varphi(x,y),\,
h(a,b)
\big)
=:
(x',y',a',b')
\end{equation}
by which $\mathcal{ M}$ is transformed to a submanifold $\mathcal{
M}'$ centered and localized at $p' = p$ having equations
\def\theequation{9.10}\begin{equation}
y'
=
\Pi'(x',a',b')
\ \ \ \ \
\text{\it and dually}
\ \ \ \ \
b'
=
{\Pi'}^*(a',x',y')
\end{equation}
with $\Pi'$ and ${\Pi'}^*$ independent of 
\def\theequation{9.11}\begin{equation}
\big(
x_{n_\mathcal{M}+1}',\dots,x_n'\big)
\ \ \ \ \
\text{\rm and of}
\ \ \ \ \
\big(
a_{p_\mathcal{M}+1}',\dots,a_p'
\big).
\end{equation}
So $\mathcal{ M }'$, may be considered to be living in $\K_{ x'}^{
n_\mathcal{M }} \times \K_{y '}^m \times \K_{ a'}^{ p_\mathcal{ M }}
\times \K_{ b'}^m$ and in such a smaller space, at $p ' = p$, it is
solvable both with respect to the parameters and to the variables.
\end{proposition}

Interpretation: by forgetting some innocuous variables, at a
Zariski-generic point, any $\mathcal{ M}$ is both solvable with
respect to the parameters and to the variables. These two assumptions
will be held up to the end of this Part~I.

\subsection*{9.12.~Dual system \thetag{ $\mathcal{ E }^*$} and 
isomorphisms $\mathfrak{ SYM}( \mathcal{ E }) \simeq \mathfrak{ SYM}
\big( \mathcal{ \mathcal{ V}_\mathcal{ S}(\mathcal{ E} )} \big) =
\mathfrak{SYM} \big( \mathcal{ \mathcal{ V}_\mathcal{ S} (\mathcal{
E}^* )} \big) \simeq \mathfrak{ SYM}( \mathcal{ E}^*)$} To a system
\thetag{ $\mathcal{ E}$}, we associate its submanifold of solutions
$\mathcal{ M} := \mathcal{V}_\mathcal{ S} ( \mathcal{ E})$. Assuming
it to be solvable with respect to the variables and proceeding as in
\S2.10, we can derive a {\sl dual system of completely integrable
partial differential equations} of the form
\def\theequation{$\mathcal{E}^*$}\begin{equation}
b_{a^\gamma}^j(a)
=
G_\gamma^j
\Big(
a,b(a),
\big(
b_{a^{\delta(l)}}^{j(l)}(a)
\big)_{
1\leqslant l\leqslant n}
\Big),
\end{equation}
where $(j,\gamma) \neq (j, 0)$ and
$\neq (j(l), \delta (l))$.
Its submanifold of solutions $\mathcal{ \mathcal{ V}_\mathcal{ S}
(\mathcal{ E}^* )} \equiv \mathcal{ V }_\mathcal{ S}(\mathcal{ E})$
has equations dual to those of $\mathcal{ \mathcal{ V}_\mathcal{ S}
(\mathcal{ E} )}$.

\def\thetheorem{9.13}\begin{theorem}
Under the assumption of twin
solvability, we have{\rm :}
\def\theequation{9.14}\begin{equation}
\mathfrak{SYM}(\mathcal{E})
\simeq 
\mathfrak{SYM}
\big(\mathcal{V}_\mathcal{S}(\mathcal{ E})\big)
=\,
\mathfrak{SYM}
\big(\mathcal{V}_\mathcal{S}(\mathcal{ E}^*)\big)
\simeq
\mathfrak{SYM}(\mathcal{E}^*),
\end{equation}
through $\mathcal{ L} \longleftrightarrow \mathcal{ L} + \mathcal{ L}^*
= \mathcal{ L}^* + \mathcal{ L} \longleftrightarrow \mathcal{ L}^*$.
\end{theorem}

\section*{\S10.~Fundamental pair of foliations and covering property}

\subsection*{10.1.~Fundamental pair of foliations on $\mathcal{ M}$ }
As in \S2, let \thetag{ $\mathcal{ E}$} and $\mathcal{ M} = \mathcal{
V}_\mathcal{ S} (\mathcal{ E})$ be defined by $y = \Pi (x, a, b)$ or
dually by $b = \Pi^* (a, x, y)$. Abbreviate
\def\theequation{10.2}\begin{equation}
z 
:= 
(x,y)
\ \ \ \ \
\text{\rm and}
\ \ \ \ \
c
:= 
(a,b). 
\end{equation}
Every transformation $(z, c) \mapsto \big( \varphi (z), h (c) \big)$
belonging to ${\sf G}_{\sf v, \sf p}$ 
stabilizes both $\{ z = {\rm cst.} \}$ and $\{ c = {\rm cst.}
\}$. Accordingly, the two foliations of $\mathcal{ M }$
\def\theequation{10.3}\begin{equation}
{\sf F}_{\sf v}
:=
\bigcup_{c_0}\,
\mathcal{M}
\cap
\big\{
c
=
c_0
\big\}
\ \ \ \ \
\text{\rm and}
\ \ \ \ \
{\sf F}_{\sf p}
:=
\bigcup_{z_0}\,
\mathcal{M}
\cap
\big\{
z
=
z_0
\big\}
\end{equation}
are invariant under changes of coordinates. We call $({\sf F}_{\sf
v}, {\sf F}_{\sf p})$ the {\sl fundamental pair of foliations} on
$\mathcal{ M}$. The leaves of the {\sl foliation by {\sf v}ariables}
${\sf F}_{\sf v}$ are $n$-dimensional: 
\def\theequation{10.4}\begin{equation}
{\sf F}_{\sf v}(c_0)
=
\big\{
(z,c_0):\
y
=
\Pi(x,c_0)
\big\}
\end{equation}
The leaves of the {\sl
foliation by {\sf p}arameters} ${\sf F}_{\sf p}$ are $p$-dimensional:
\def\theequation{10.5}\begin{equation}
{\sf F}_{\sf p}(c_0)
=
\big\{
(z_0,c):\
b
=
\Pi^*(a,z_0)
\big\}
\end{equation}
We draw a diagram. In it, the positive codimension is invisible:
\def\theequation{10.6}\begin{equation}
m
=
\dim\mathcal{M}
-
\dim{\sf F}_{\sf v}
-
\dim{\sf F}_{\sf p}
\geqslant 
1
\end{equation}

\bigskip
\begin{center}
\begin{picture}(0,0)%
\includegraphics{double-foliation.pstex}%
\end{picture}%
\setlength{\unitlength}{3947sp}%
\begingroup\makeatletter\ifx\SetFigFont\undefined
\def\x#1#2#3#4#5#6#7\relax{\def\x{#1#2#3#4#5#6}}%
\expandafter\x\fmtname xxxxxx\relax \def\y{splain}%
\ifx\x\y   
\gdef\SetFigFont#1#2#3{%
  \ifnum #1<17\tiny\else \ifnum #1<20\small\else
  \ifnum #1<24\normalsize\else \ifnum #1<29\large\else
  \ifnum #1<34\Large\else \ifnum #1<41\LARGE\else
     \huge\fi\fi\fi\fi\fi\fi
  \csname #3\endcsname}%
\else
\gdef\SetFigFont#1#2#3{\begingroup
  \count@#1\relax \ifnum 25<\count@\count@25\fi
  \def\x{\endgroup\@setsize\SetFigFont{#2pt}}%
  \expandafter\x
    \csname \romannumeral\the\count@ pt\expandafter\endcsname
    \csname @\romannumeral\the\count@ pt\endcsname
  \csname #3\endcsname}%
\fi
\fi\endgroup
\begin{picture}(5724,2351)(56,-1625)
\put(4311,-167){\makebox(0,0)[lb]{\smash{\SetFigFont{8}{9.6}{rm}{\color[rgb]{0,0,0}${\sf F}_{\sf v}$}%
}}}
\put(1140,-1376){\makebox(0,0)[lb]{\smash{\SetFigFont{8}{9.6}{rm}{\color[rgb]{0,0,0}$\mathcal{M}$}%
}}}
\put(4640,-1376){\makebox(0,0)[lb]{\smash{\SetFigFont{8}{9.6}{rm}{\color[rgb]{0,0,0}$\mathcal{M}$}%
}}}
\put(3020,429){\makebox(0,0)[lb]{\smash{\SetFigFont{8}{9.6}{rm}{\color[rgb]{0,0,0}$c$}%
}}}
\put(5275,-526){\makebox(0,0)[lb]{\smash{\SetFigFont{8}{9.6}{rm}{\color[rgb]{0,0,0}$z$}%
}}}
\put(3621,-526){\makebox(0,0)[lb]{\smash{\SetFigFont{8}{9.6}{rm}{\color[rgb]{0,0,0}${\sf F}_{\sf p}$}%
}}}
\put(2085,-861){\makebox(0,0)[lb]{\smash{\SetFigFont{8}{9.6}{rm}{\color[rgb]{0,0,0}${\sf L}$}%
}}}
\put(1680,-116){\makebox(0,0)[lb]{\smash{\SetFigFont{8}{9.6}{rm}{\color[rgb]{0,0,0}${\sf L}^*$}%
}}}
\put(3030,-531){\makebox(0,0)[lb]{\smash{\SetFigFont{8}{9.6}{rm}{\color[rgb]{0,0,0}$0$}%
}}}
\end{picture}

\end{center}

\subsection*{10.7.~Chains $\Gamma_k$ and dual chains $\Gamma_k^*$}
Similarly as in~\cite{ gm2004, me2005a, me2005b, mp2005} (in a CR
context), we introduce two collections $( {\sf L}_k )_{1 \leqslant k
\leqslant n}$ and $({\sf L}_q^*)_{1 \leqslant q \leqslant p}$ of
vector fields whose integral manifolds coincide with the leaves of
${\sf F}_{ \sf v}$ and of ${\sf F}_{ \sf p}$:
\def\theequation{10.8}\begin{equation}
\left\{
\aligned
{\sf L}_k
:=
\frac{\partial}{\partial x_k}+
\sum_{j=1}^m \, 
\frac{\partial \Pi^j}{\partial x_k} (x,a,b)\, 
\frac{\partial}{\partial y^j}, 
\ \ \ \ \ \ 
k=1,\dots,n,
\\
{\sf L}_q^*
:=
\frac{\partial}{\partial a^q}+
\sum_{j=1}^m\, 
\frac{\partial{\Pi^*}^j}{\partial a^q}(a,x,y)\,
\frac{\partial}{\partial b^j},
\ \ \ \ \ \ 
q=1,\dots,p.
\endaligned\right.
\end{equation}
Let $(z_0, c_0) = (x_0, y_0, a_0, b_0) \in \mathcal{ M}$ be a fixed
point, let $x_1 : =(x_1^1, \dots, x_1^n) \in \K^n$ and define the
multiple flow map
\def\theequation{10.9}\begin{equation}
\left\{
\aligned
{\sf L}_{x_1}(x_0,y_0,a_0,b_0)
:=
& \
\exp(x_1{\sf L})(p_0)
:=
\exp\big(x_1^n{\sf L}_n\big(\cdots(
\exp(x_1^1{\sf L}_1(z_0,c_0)))
\cdots)\big)
\\
:=
& \
\big(
x_0+x_1,
\Pi(x_0+x_1,a_0,b_0),
a_0,b_0
\big).
\endaligned\right.
\end{equation}
Similarly, for $a_1 = (a_1^1, \dots, a_1^p ) \in \K^p$, define
the multiple flow map
\def\theequation{10.10}\begin{equation}
{\sf L}_{a_1}^*
(x_0,y_0,a_0,b_0)
:=
\big(
x_0,y_0,a_0+a_1,
\Pi^*(a_0+a_1,x_0,y_0)
\big).
\end{equation}
Starting from the $(z_0, c_0) = (0, 0)$ and moving alternately along
${\sf F}_{\sf v}$, ${\sf F}_{\sf p}$, ${\sf F}_{\sf v}$, 
{\it etc.}, we obtain
\def\theequation{10.11}\begin{equation}
\left\{
\aligned
\Gamma_1(x_1)
:=
& \
{\sf L}_{x_1}(0), \\
\Gamma_2(x_1,a_1)
:=
& \
{\sf L}_{a_1}^*({\sf L}_{x_1}(0)), \\
\Gamma_3(x_1,a_1,x_2)
:=
& \
{\sf L}_{x_2} ({\sf L}_{a_1}^*({\sf L}_{x_1}(0))), \\
\Gamma_4(x_1,a_1,x_2,a_2)
:=
& \
{\sf L}_{a_2}^*(
{\sf L}_{x_2} ({\sf L}_{a_1}^*({\sf L}_{x_1}(0)))),
\endaligned\right.
\end{equation}
and so on. Generally, we get {\sl chains}\, $\Gamma_k := \Gamma_k ([x
a ]_k )$, where $[x a]_k := (x_1, a_1, x_2, a_2, \dots)$ with exactly
$k$ terms, where each $x_l \in \K^n$ and each $a_l \in \K^p$.

If, instead, the first movement consists in moving along ${\sf F}_{\sf
p}$, we start with $\Gamma_1^*( a_1):= {\sf L}_{ a_1}^* (0)$,
$\Gamma_2^*( a_1, x_1) := {\sf L}_{x_1} ( {\sf L}_{ a_1}^* (0))$, {\it
etc.}, and generally we get {\sl dual chains}\, $\Gamma_k^* ([a
x]_k)$, where $[a x]_k := ( a_1, x_1, a_2, x_2, \dots)$, with exactly
$k$ terms. Both $\Gamma_k$ and $\Gamma_k^*$ have range in $\mathcal{ M
}$.

For $k=1, 2,3, \cdots$, integers $e_k$ and $e_k^*$ are defined
inductively by
\def\theequation{10.12}\begin{equation}
\left\{
\aligned
e_1+e_2+e_3+\cdots+e_k={\rm genrk}_\K(\Gamma_k),\\
e_1^*+e_2^*+e_3^*+\cdots+e_k^*={\rm genrk}_\K(\Gamma_k^*).
\endaligned\right.
\end{equation}
By~\thetag{10.9} and~\thetag{ 10.10}, it is clear that $e_1 =n$, $e_2
= p$, $e_1^* =p$, and $e_2^* = n$. 

\def\theexample{10.13}\begin{example}
{\rm 
For $y_{xx} = 0$, the submanifold of solutions $\mathcal{ M}$ is
simply $y = b + xa$, whence
\def\theequation{10.14}\begin{equation}
\left\{
\aligned
\Gamma_1(x_1) 
= 
&\
(x_1,0,0,0), 
\\
\Gamma_2(x_1,a_1) 
= 
&\
(x_1,0,a_1,-x_1a_1), 
\\
\Gamma_3(x_1,a_1,x_2) 
= 
&\
(x_1+x_2,x_2a_1,a_1,-x_1a_1).
\endaligned\right.
\end{equation}
The rank at $(0, 0, 0)$ of $\Gamma_3$ is equal to two, not
more. However, its generic rank is equal to three. Similar
observations hold for the two submanifolds of solutions $y = b+ xx a +
xa a$ and $y = b+ x a^1+ x x a^2$ (in $\K^5$).
}\end{example}

\def\thelemma{10.15}\begin{lemma}
If ${\rm genrk }_\K ( \Gamma_{ k+1}) = {\rm genrk}_\K (\Gamma_{ k})$,
then for each positive integer $l \geqslant 1$, we have ${\rm
genrk}_\K( \Gamma_{ k+l})={\rm genrk }_\K ( \Gamma_k )$. The same
stabilization property holds for $ \Gamma_k^*$.
\end{lemma}

\subsection*{10.16.~Covering property} We now formulate a central 
concept.

\def\thedefinition{10.17}\begin{definition}
{\rm 
The pair of foliations $({\sf F}_{\sf v},{\sf F}_{\sf p})$ is
{\sl covering at the origin} if there exists an integer $k$ such that
the generic rank of $\Gamma_k$ is (maximal possible) equal to $\dim_\K
\, \mathcal{M}$. Since for $a_1=0$, the dual $(k+1)$-th chain
$\Gamma_{k+1}^*$ identifies with the $k$-th chain $\Gamma_k$, 
the same property holds for the dual chains.
}
\end{definition}

\def\theexample{10.18}\begin{example}
{\rm
With $n=1$, $m=2$ and $p=1$ the submanifold defined by $y^1 = b^1$ and
$y^2 = b^2 + x a$ is twin solvable, but its pair of foliations is not
covering at the origin. Then $\mathfrak{ SYM} (\mathcal{ M})$ is
infinite-dimensional, since for $a = a (y^1)$ an arbitrary function,
it contains $a (y^1)\, \frac{ \partial}{ \partial y^1 } + a( b^1)\,
\frac{ \partial }{ \partial b^1 }$.
}
\end{example}

Because we aim only to study finite-dimensional Lie symmetry groups of
partial differential equations, in the remainder of this Part~I, we
will constantly assume the covering property to hold.

By Lemma~10.15, there exist two well defined integers $\mu$ and
$\mu^*$ such that $e_3, e_4, \dots, e_{ \mu+1}>0$, but $e_{
\mu+l}= 0$ for all $l \geqslant 2$ and similarly, $e_3^*, e_4^*,
\dots, e_{ \mu^*+1}^* >0$, but $e_{ \mu^*+l}^* =0$ for all $l
\geqslant 2$. Since the pair of foliations is covering, we have the
two dimension equalities
\def\theequation{10.19}\begin{equation}
\left\{
\aligned
n+p+e_3+\cdots+e_{\mu+1}
=
\dim_\K\, \mathcal{M}=n+m+p, 
\\
p+n+e_3^*+\cdots+e_{\mu^*+1}^*=\dim_\K \, \mathcal{M}=n+m+p.
\endaligned\right.
\end{equation}
By definition, the ranges of $\Gamma_{ \mu + 1}$ and of $\Gamma_{
\mu^* + 1 }^*$ cover (at least; more is true, see: Theorem~10.28) an
open subset of $\mathcal{ M }$. Also, it is elementary to verify the
four inequalities
\def\theequation{10.20}\begin{equation}
\aligned
{}
&
\mu\leqslant 1+m, \ \ \ \ \ \ \ \ \ \ \ \ 
\mu^*\leqslant \ 1+m,\\
&
\mu\leqslant \mu^*+1, \ \ \ \ \ \ \ \ \ \ \ 
\mu^*\leqslant \mu+1.
\endaligned
\end{equation}
In fact, since $\Gamma_{ k+1}$ with $x_1 = 0$ identifies with 
$\Gamma_k^*$, the second line follows.

\def\thedefinition{10.21}\begin{definition}{\rm
The {\sl type of the covering pair of foliations $({\sf F}_{\sf v},
{\sf F}_{\sf p})$} is the pair of integers
\def\theequation{10.22}\begin{equation}
(\mu,\mu^*), 
\ \ \ \ \ \ \ \ \ \ \
\text{\rm with} 
\ \ \ \ \ \ 
\max(\mu,\mu^*)\leqslant 1+m.
\end{equation}
}\end{definition}

\def\theexample{10.23}\begin{example}{\rm
(Continued) We write down the explicit expressions of $\Gamma_4$ and
of $\Gamma_5$:
\def\theequation{10.24}\begin{equation}
\left\{
\aligned
\Gamma_4(x_1,a_1,x_2,a_2; 0) = 
&\
\big(x_1+x_2, \, x_2a_1,
\,a_1+a_2,\, 
-x_1a_1-x_1a_2-x_2a_2,
\big), 
\\
\Gamma_5(x_1,a_1,x_2,a_2,x_3; 0) =
&\
\big(
x_1+x_2+x_3,\, 
x_2a_1+x_3a_1+x_3a_2,\,
a_1+a_2 
\\
& \ \ \ \ \ \ \ \ \ \ \ \ \ \ \ \ \ \ \ \ \ \ \ \ 
\ \ \ \ \ \ \ \ \ \ \ \ \ \ \ \ \ 
-x_1a_1-x_1a_2-x_2a_2
\big).
\endaligned\right.
\end{equation}
Here, $\dim \mathcal{ M} = 3$. By computing its Jacobian matrix,
$\Gamma_5$ is of rank 3 at every point $(x_1, a_1, 0, -a_1,
- x_1) \in \K^5$ with $a_1 \neq 0$. Since (obviously)
\def\theequation{10.25}\begin{equation}
\Gamma_5
\big(
x_1,a_1,0,-a_1,-x_1
\big)
=0\in\mathcal{M},
\end{equation}
we deduce that $\Gamma_5 $ is submersive (``covering'') from a small
neighborhood of $\big( x_1, a_1, 0, -a_1, - x_1 \big)$ in $
\K^5$ onto a neighborhood of the origin in $\mathcal{ M }$.
}\end{example}

\subsection*{10.26.~Covering a neighborhood of the origin in 
$\mathcal{ M}$} For $(z_0, c_0) \in \mathcal{ M}$ fixed and close to
the origin, we denote by $\Gamma_k \big( [ xa]_k ; (z_0, c_0) \big)$
and by $\Gamma_k^* \big( [ ax]_k ; (z_0, c_0) \big)$ the (dual) chains
issued from $(z_0, c_0)$. For given parameters $[xa ]_k = (x_1, a_1,
x_2, \dots)$, we denote by $[ -xa ]_k$ the collection $( \cdots, -x_2,
-a_1, -x_2)$ with minus signs and reverse order; similarly, we
introduce $[- a x]_k$. Notably, we have ${\sf L}_{- x_1} ( {\sf L}_{
x_1}( 0)) = 0$ (because ${\sf L}_{- x_1 + x_1}(\cdot)={\sf
L}_0(\cdot)={\rm Id}$), and also ${\sf L}_{ -x_1}( {\sf L}_{ -a_1 }^*(
{\sf L}_{ a_1}^*( {\sf L}_{ x_1}( 0)))) = 0$ and generally:
\def\theequation{10.27}\begin{equation}
\Gamma_k
\big(
[-xa]_k;\Gamma_k([xa]_k;0)
\big)
\equiv 
0.
\end{equation}
Geometrically speaking, by following backward the $k$-th chain
$\Gamma_k$, we come back to $0$.

\def\thetheorem{10.28}\begin{theorem}
{\rm (\cite{ me2005a, me2005b}, [$*$])} The two maps $\Gamma_{2 \mu
+ 1 }$ and $\Gamma_{2 \mu^* + 1}^*$ are submersive onto a
neighborhood of the origin in $\mathcal{ M}$. Precisely, there exist
two points $[xa]_{2 \mu +1}^0 \in \K^{ (\mu+1)n+\mu p}$ and $[a
x]_{2 \mu^* +1}^{0} \in\K^{ \mu^* n+ (\mu^*+1)p}$ arbitrarily
close to the origin with $\Gamma_{ 2\mu +1}( [xa]_{ 2\mu +1}^0)=0$
and $\Gamma_{ 2\mu^* +1}^*([ a x]_{2 \mu^* +1 }^{0} ) =0$ such
that the two maps
\def\theequation{10.29}\begin{equation}
\left\{
\aligned
\K^{(\mu+1)n+\mu p}
\ni
[xa]_{2\mu+1}
\longmapsto 
& \
\Gamma_{2\mu+1}
\big(
[xa]_{2\mu+1}
\big)
\in
\mathcal{M} 
\ \ \ \ \ \text{\rm and} 
\\
\K^{\mu^*n+(\mu^*+1)p}
\ni
[a x]_{2\mu^*+1}
\longmapsto 
& \
\Gamma_{2\mu^*+1}^*\big(
[a x]_{2\mu^*+1}
\big)
\in
\mathcal{M}
\endaligned\right.
\end{equation} 
are of rank $n+ m+p = \dim_\K \, \mathcal{ M}$ at the points $[x a]_{
2 \mu}^0$ and $[ a x]_{ 2\mu^* }^{0}$ respectively. In particular,
the ranges of the two maps $\Gamma_{2 \mu + 1}$ and $\Gamma_{ 2
\mu^* +1 }^*$ cover a neighborhood of the origin in $\mathcal{ M}$.
\end{theorem}

Let $\pi_z (z, c) := z$ and $\pi_c (z, c) := c$ be the two canonical
projections. The next corollary will be useful in
Section~12. In the example above, it also follows that the map
\def\theequation{10.30}\begin{equation}
[xa]_4\mapsto
\pi_c 
\big(
\Gamma_4([xa_4])
\big)
=
\big(
a_1+a_2,\,
-x_1a_1-x_1a_2-x_2a_2
\big) 
\in\K^2
\end{equation}
is of rank two at all points $[ xa ]_4^0 := \big( x_1^0, a_1^0, 0,
-a_1^0 \big)$ with $a_1^0 \neq 0$.

\def\thecorollary{10.31}\begin{corollary}
{\rm (\cite{ me2005a, me2005b}, [$*$])} There exist two points $[x a
]_{2 \mu }^0 \in \K^{ \mu (n + p)}$ and $[a x ]_{ 2 \mu^* }^{ 0}
\in \K^{ \mu^* (n+p) }$ arbitrarily close to the origin with $\pi_c(
\Gamma_{2 \mu } ([xa]_{ 2 \mu }^0)) = 0$ and $\pi_z \big( \Gamma_{
2\mu^* }^*( [a x]_{ 2 \mu^* }^{0} )) = 0$ such that the two maps
\def\theequation{10.32}\begin{equation}
\left\{
\aligned
\K^{\mu(n+p)}\ni [xa]_{2\mu}
\longmapsto 
& \
\pi_c 
\big(
\Gamma_{2\mu}([xa]_{2\mu})
\big)
\in
\K^{m+p} 
\ \ \ \ \ \text{\rm and} 
\\
\K^{\mu^*(n+p)}\ni [a x]_{2\mu^*}
\longmapsto 
& \
\pi_z
\big(
\Gamma_{2\mu^*}^*([a x]_{2\mu^*})
\big)\in\K^{n+m}
\endaligned\right.
\end{equation}
are of rank $m+p$ at the point $[x a ]_{2 \mu }^0 \in \K^{ \mu
(n+p)}$ and of rank $n+m$ at the point $[a x ]_{ 2 \mu^* }^{ 0}
\in \K^{\mu^* (n+p) }$. 
\end{corollary}

In the case $m=1$ (single dependent variable $y \in \K$), the covering
property always hold with $\mu = \mu^* = 2$.

\section*{\S11.~Formal and smooth equivalences between
submanifolds of solutions}

\subsection*{11.1.~Transformations of submanifolds of solutions}
Lemma~7.3 shows that every equivalence $\varphi$ between two {\sc pde}
systems~\thetag{ $\mathcal{ E}$} and~\thetag{ $\mathcal{ E}'$} lifts
as a transformation which respects the separation between variables
and parameters of the form
\def\theequation{11.2}\begin{equation}
(x,y,a,b)
\longmapsto
\big(
\phi(x,y),\psi(x,y),f(a,b),g(a,b)
\big)
=
\big(
\varphi(x,y),h(a,b)
\big)
=:
(x',y',a',b')
\end{equation}
from the source submanifolds of solutions $\mathcal{ M} := \mathcal{
V}_{ \mathcal{ S}} (\mathcal{ E})$ to the target $\mathcal{ M}' :=
\mathcal{ V}_{ \mathcal{ S}} (\mathcal{ E}')$, whose equations are
\def\theequation{11.3}\begin{equation}
\aligned
y
&
=
\Pi(x,c)
\ \ \ \ \ \ \ \
\text{\rm or dually}
\ \ \ \ \
b
=
{\Pi^*}(a,z)
\ \ \ \ \
\text{\rm and}
\\
{y'}
&
=
{\Pi'}(x',c')
\ \ \ \ \
\text{\rm or dually}
\ \ \ \ \
{b'}
=
{{\Pi'}^*}(a',z').
\endaligned
\end{equation}
The study of transformations between submanifolds of solutions
possesses strong similarities with the study of CR mappings between CR
manifolds (\cite{ pi1975, we1977, dw1980, bjt1985, df1988, ber1999,
me2005a, me2005b}). In fact, one may transfer the whole theory of the
analytic reflection principle to this more general context. In the
present \S10 and in the next \S11, we select and establish some of the
results that are useful to the Lie theory. Some accessible open
questions will also be formulated.

\medskip

Maps of the form~\thetag{ 11.2} send leaves of ${\sf F}_{ \sf v}$ and
of ${\sf F}_{ \sf p}$ to leaves of ${\sf F}_{ \sf v}'$, and of ${\sf
F}_{ \sf p}'$, respectively.

\bigskip
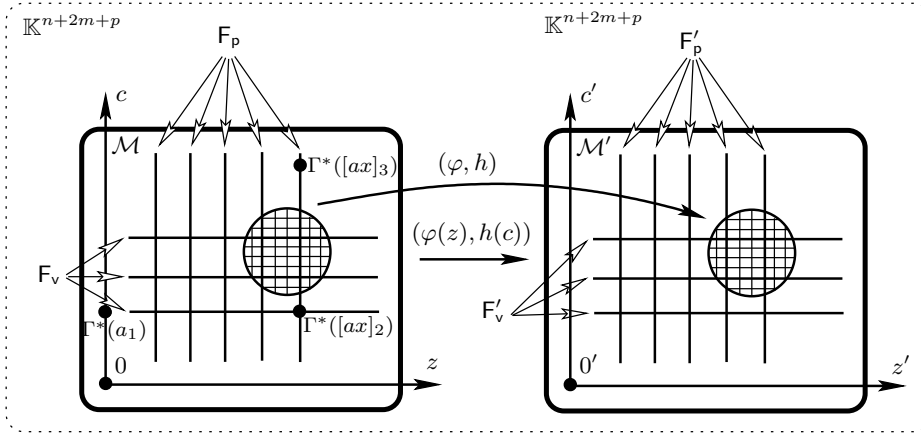
\begin{center}
\begin{picture}(0,0)%
\includegraphics{new-stabilization-foliations.pstex}%
\end{picture}%
\setlength{\unitlength}{4144sp}%
\begingroup\makeatletter\ifx\SetFigFont\undefined
\def\x#1#2#3#4#5#6#7\relax{\def\x{#1#2#3#4#5#6}}%
\expandafter\x\fmtname xxxxxx\relax \def\y{splain}%
\ifx\x\y   
\gdef\SetFigFont#1#2#3{%
  \ifnum #1<17\tiny\else \ifnum #1<20\small\else
  \ifnum #1<24\normalsize\else \ifnum #1<29\large\else
  \ifnum #1<34\Large\else \ifnum #1<41\LARGE\else
     \huge\fi\fi\fi\fi\fi\fi
  \csname #3\endcsname}%
\else
\gdef\SetFigFont#1#2#3{\begingroup
  \count@#1\relax \ifnum 25<\count@\count@25\fi
  \def\x{\endgroup\@setsize\SetFigFont{#2pt}}%
  \expandafter\x
    \csname \romannumeral\the\count@ pt\expandafter\endcsname
    \csname @\romannumeral\the\count@ pt\endcsname
  \csname #3\endcsname}%
\fi
\fi\endgroup
\begin{picture}(5559,2589)(484,-2188)
\put(1127,-510){\makebox(0,0)[lb]{\smash{\SetFigFont{9}{10.8}{rm}{\color[rgb]{0,0,0}$\mathcal{M}$}%
}}}
\put(597,196){\makebox(0,0)[lb]{\smash{\SetFigFont{9}{10.8}{rm}{\color[rgb]{0,0,0}$\K^{n+2m+p}$}%
}}}
\put(1147,-1827){\makebox(0,0)[lb]{\smash{\SetFigFont{9}{10.8}{rm}{\color[rgb]{0,0,0}$0$}%
}}}
\put(3718,214){\makebox(0,0)[lb]{\smash{\SetFigFont{9}{10.8}{rm}{\color[rgb]{0,0,0}$\K^{n+2m+p}$}%
}}}
\put(1758,153){\makebox(0,0)[lb]{\smash{\SetFigFont{9}{10.8}{rm}{\color[rgb]{0,0,0}${\sf F}_{\sf p}$}%
}}}
\put(684,-1278){\makebox(0,0)[lb]{\smash{\SetFigFont{9}{10.8}{rm}{\color[rgb]{0,0,0}${\sf F}_{\sf v}$}%
}}}
\put(1153,-195){\makebox(0,0)[lb]{\smash{\SetFigFont{9}{10.8}{rm}{\color[rgb]{0,0,0}$c$}%
}}}
\put(3009,-1798){\makebox(0,0)[lb]{\smash{\SetFigFont{9}{10.8}{rm}{\color[rgb]{0,0,0}$z$}%
}}}
\put(4523,117){\makebox(0,0)[lb]{\smash{\SetFigFont{9}{10.8}{rm}{\color[rgb]{0,0,0}${\sf F}_{\sf p}'$}%
}}}
\put(3933,-203){\makebox(0,0)[lb]{\smash{\SetFigFont{9}{10.8}{rm}{\color[rgb]{0,0,0}$c'$}%
}}}
\put(5793,-1830){\makebox(0,0)[lb]{\smash{\SetFigFont{9}{10.8}{rm}{\color[rgb]{0,0,0}$z'$}%
}}}
\put(3926,-1825){\makebox(0,0)[lb]{\smash{\SetFigFont{9}{10.8}{rm}{\color[rgb]{0,0,0}$0'$}%
}}}
\put(3340,-1496){\makebox(0,0)[lb]{\smash{\SetFigFont{9}{10.8}{rm}{\color[rgb]{0,0,0}${\sf F}_{\sf v}'$}%
}}}
\put(3078,-610){\makebox(0,0)[lb]{\smash{\SetFigFont{9}{10.8}{rm}{\color[rgb]{0,0,0}$(\varphi,h)$}%
}}}
\put(2920,-1029){\makebox(0,0)[lb]{\smash{\SetFigFont{9}{10.8}{rm}{\color[rgb]{0,0,0}$(\varphi(z),h(c))$}%
}}}
\put(2274,-1578){\makebox(0,0)[lb]{\smash{\SetFigFont{8}{9.6}{rm}{\color[rgb]{0,0,0}$\Gamma^*([a x]_2)$}%
}}}
\put(2299,-623){\makebox(0,0)[lb]{\smash{\SetFigFont{8}{9.6}{rm}{\color[rgb]{0,0,0}$\Gamma^*([a x]_3)$}%
}}}
\put(947,-1605){\makebox(0,0)[lb]{\smash{\SetFigFont{8}{9.6}{rm}{\color[rgb]{0,0,0}$\Gamma^*(a_1)$}%
}}}
\put(3906,-519){\makebox(0,0)[lb]{\smash{\SetFigFont{9}{10.8}{rm}{\color[rgb]{0,0,0}$\mathcal{M}'$}%
}}}
\end{picture}

\end{center}

\subsection*{ 11.4.~Regularity and jet parametrization}
Some strong rigidity properties underly the above diagram.
Especially, the smoothness of the two pairs $\big( {\sf F}_{ \sf v},
\, {\sf F}_{ \sf p} \big)$ and $\big( {\sf F}_{ \sf v}', \, {\sf F}_{
\sf p}' \big)$ governs the smoothness of $(\varphi, h)$.

We shall study the regularity of a {\it purely formal}\, map $(z', c')
= \big( \varphi (z), h(c) \big)$, namely $\varphi (z) \in \K \dl z
\dr^{ n+m}$ and $h (c) \in \K \dl c \dr^{ p + m}$, assuming \thetag{
$\mathcal{ E}$} and \thetag{ $\mathcal{ E}'$} to be analytic.
Concretely, the assumption that
$(\varphi, h)$ maps $\mathcal{ M}$ to $\mathcal{ M}'$
reads as one of the four equivalent identities:
\def\theequation{11.5}\begin{equation}
\left\{
\aligned
\psi
\big(
x,\Pi(x,c)
\big)
&
\equiv
\Pi'
\big(
\phi(x,\Pi(x,c)),h(c)
\big), 
\\
\psi
(z)
&
\equiv
\Pi'
\big(
\phi(z),h(a,\Pi^*(a,z)
\big),
\\
g
\big(
a,\Pi^*(a,z)
\big)
&
\equiv
{\Pi'}^*
\big(
f(a,\Pi^*(a,z)),\varphi(z)
\big),
\\
g(c)
&
\equiv
{\Pi'}^*
\big(
f(c),\varphi(x,\Pi(x,c))
\big),
\endaligned\right.
\end{equation}
in $\K \dl x, c \dr^m$ and in $\K \dl a, z \dr^m$.

\def\thetheorem{11.6}\begin{theorem}
Let $(\varphi, h) := \mathcal{ M} \to \mathcal{ M}'$ be a purely
formal equivalence between two local $\K$-analytic submanifolds of
solutions. Assume that the fundamental pair of foliations $( {\sf
F}_{\sf v}, {\sf F}_{ \sf p})$ is covering at the origin, with type
$(\mu, \mu^* )$ at the origin. Assume that $\mathcal{ M}'$ is both
$\kappa$-solvable with respect to the parameters and
$\kappa^*$-solvable with respect to the variables. Set $\ell := \mu^*
(\kappa + \kappa^*)$ and $\ell^* := \mu (\kappa^* + \kappa)$. Then
there exist two $\K^{ n + m}$-valued and $\K^{ p + m}$-valued local
$\K$-analytic maps $\Phi_{ \ell}$ and $H_{ \ell^*}$, constructible
only by means of $\Pi$, $\Pi^*$, $\Pi'$, ${ \Pi' }^*$, such that the
following two formal power series identities hold{\rm :}
\def\theequation{11.7}\begin{equation}
\left\{
\aligned
\varphi(z)
&
\equiv
\Phi_{\ell}
\big(
z,J_z^{\ell}\varphi(0)
\big), 
\\
h(c) 
&
\equiv
H_{\ell^*}
\big(
c,J_c^{\ell^*}h(0)
\big),
\endaligned\right.
\end{equation}
in $\K \dl z\dr^{ n + m}$ and in $\K \dl c\dr^{ p+m}$, where $J_z^{
\ell} \varphi (0)$ denotes the $\ell$-th jet of $h$ at the origin and
similarly for $J_c^{ \ell^*} h(0)$. In particular, as a corollary, we
have the following two automatic regularity properties:

\begin{itemize}

\smallskip\item[$\bullet$]
$\varphi (z) \in \K \{ z \}^{ n+m}$ and $h (c) \in \K \{ c \}^{ p +
m}$ are in fact convergent{\rm ;}

\smallskip\item[$\bullet$]
if in addition $\mathcal{ M}$ and $\mathcal{ M}'$ are $\K$-algebraic
in the sense of Nash, then $\Phi_{ \ell}$ and $H_{ \ell^*}$ are also
$\K$-algebraic, whence $\varphi (z) \in \mathcal{ A}_\K \{ z \}^{
n+m}$ and $h (c) \in \mathcal{ A}_\K \{ c \}^{ p + m}$ are in fact
$\K$-algebraic.
\end{itemize}
\end{theorem}

\proof
We remind the explicit expressions of the two collections of vector
fields spanning the leaves of the two foliations ${ \sf F}_{\sf v}$
and ${\sf F}_{\sf p}$:
\def\theequation{11.8}\begin{equation}
\left\{
\aligned
{\sf L}_k
:=
\frac{\partial}{\partial x_k}+
\sum_{j=1}^m \, 
\frac{\partial \Pi^j}{\partial x_k}(x,c)\, 
\frac{\partial}{\partial y^j}, 
\ \ \ \ \ \ 
k=1,\dots,n,
\\
{\sf L}_q^*
:=
\frac{\partial}{\partial a^q}+
\sum_{j=1}^m\, 
\frac{\partial{\Pi^*}^j}{\partial a^q}(a,z)\,
\frac{\partial}{\partial b^j},
\ \ \ \ \ \ 
q=1,\dots,p.
\endaligned\right.
\end{equation}
Observe that differentiating the first line of~\thetag{11.5} with
respect to $x^k$ amounts to applying the derivation ${\sf
L}_k$. Similarly, differentiating the third line of~\thetag{11.5} with
respect to $a^q$ amounts to applying ${\sf L}_q^*$. We thus get for
$(z, c) \in \mathcal{ M }$
\def\theequation{11.9}\begin{equation}
\left\{
\aligned
{\sf L}_k\,\psi(z)
&
=
\sum_{l=1}^n\, 
\frac{\partial\Pi'}{\partial{x'}^l}
\big(
\phi(z),h(c)
\big)\, 
{\sf L}_k\,\phi^l(z)
\ \ \ \ \ \text{\rm and}
\\
{\sf L}_q^*\,g(c)
&
=
\sum_{r=1}^p\,
\frac{\partial{\Pi'}^*}{\partial {a'}^r}\, 
\big(
f(c),\varphi(z)
\big)\, 
{\sf L}_q^*\,f^r(c),
\endaligned\right.
\end{equation}
It follows from $\det \big( \frac{ \partial \varphi^k }{ \partial z^l}
\big) (0) \neq 0$ and $\det \big( \frac{ \partial h^k }{ \partial c^l}
\big) (0) \neq 0$ that the two formal determinants
\def\theequation{11.10}\begin{equation}
{\rm det}\,
\big(
{\sf L}_k\,\phi^l(z)
\big)_{1\leqslant k\leqslant n}^{1\leqslant l\leqslant n}
\ \ \ \ \ \ \ 
\text{\rm and}
\ \ \ \ \ \ \
{\rm det}\,
\big(
{\sf L}_q^*\,f^r(c)
\big)_{1\leqslant q\leqslant p}^{1\leqslant r\leqslant p}
\end{equation}
have nonvanishing constant term.
Consequently, these two matrices are invertible in 
$\K\dl z \dr$ and in $\K \dl c \dr$. So there exist universal
polynomials ${\sf S}_l^j$ and ${{\sf S}^*}_r^j$ such that
\def\theequation{11.11}\begin{equation}
\left\{
\aligned
\frac{\partial{\Pi'}^j}{\partial{x'}^l}
\big(\varphi(z),h(c)\big)
&
=
\frac{
{\sf S}_l^j
\left(
\big\{{\sf L}_{k'}\,\varphi^{i'}(z)\big\}_{
1\leqslant k'\leqslant n}^{1\leqslant i'\leqslant n+m}
\right)
}{
{\rm det}
\big(
{\sf L}_{k'}\,\phi^{l'}(z)\big)_{
1\leqslant k'\leqslant n}^{1\leqslant l'\leqslant n}
} 
\ \ \ \ \ \text{\rm and}
\\
\frac{\partial{{\Pi'}^*}^j}{\partial{a'}^r}
\big(f(c),\varphi(z)\big)
&
=
\frac{ 
{{\sf S}^*}_r^j\left(
\big\{{\sf L}_{q'}^*\,h^{i'}(c)\big\}_{
1\leqslant q'\leqslant p}^{1\leqslant i'\leqslant p+m}
\right)
}{
{\rm det}\big(
{\sf L}_{q'}^*\,f^{r'}(c)
\big)_{1\leqslant q'\leqslant p}^{1\leqslant r'\leqslant p}
},
\endaligned\right.
\end{equation}
for $1 \leqslant j \leqslant m$, for $1 \leqslant l \leqslant n$, for
$1 \leqslant r \leqslant p$ and for $(z, c) \in \mathcal{ M}$.

Again, we apply the vector fields ${\sf L}_k$ to the obtained first
line and the vector fields ${\sf L}_q^*$ to the obtained second line,
getting, thanks to the chain rule:
\def\theequation{11.12}\begin{equation}
\left\{
\aligned
\sum_{l_2=1}^n\, 
\frac{\partial^2{\Pi'}^j}{\partial
{x'}^{l_1}{x'}^{l_2}}
\big(\phi(z),h(c)\big)\, 
{\sf L}_k\,\phi^{l_2}(z) 
&
=
\frac{{\sf R}_{l_1,k}^j
\left(
\big\{{\sf L}_{k_1'}{\sf L}_{k_2'}\varphi^{i'}(z)\big\}_{
1\leqslant k_1',k_2'\leqslant n}^{
1\leqslant i'\leqslant n+m}
\right)
}{ 
\left[
{\rm det}
\big(
{\sf L}_{k'}\,\phi^{l'}(z)\big)_{
1\leqslant k'\leqslant n}^{1\leqslant l'\leqslant n}
\right]^2
} 
\ \ \ \ \ \text{\rm and} 
\\
\sum_{r_2=1}^p\, 
\frac{\partial^2{{\Pi'}^*}^j}{
\partial{a'}^{r_1}{a'}^{r_2}}
\big(f(c),\varphi(z)\big)\, 
{\sf L}_q^*\,f^{r_2}(c)
&
=
\frac{{{\sf R}^*}_{r_1,q}^j
\left(
\big\{
{\sf L}_{q_1'}^*{\sf L}_{q_2'}^*h^{i'}(c)
\big\}_{
1\leqslant q_1',q_2'\leqslant p}^{1\leqslant i'\leqslant p+m}
\right)}
{\left[
{\rm det}\big(
{\sf L}_{q'}^*\,f^{r'}(c)
\big)_{1\leqslant q'\leqslant p}^{1\leqslant r'\leqslant p}
\right]^2},
\endaligned\right.
\end{equation}
for $1 \leqslant j \leqslant m$, for $1 \leqslant l_1, l_2 \leqslant
n$, for $1 \leqslant r_1, r_2 \leqslant p$ and for $(z, c) \in
\mathcal{ M}$. Here, ${\sf R}_{ l_1, k}^j$ and ${{\sf R}^*}_{r_1,
q}^j$ are universal polynomials.
Then applying once more Cramer's rule, we get
\def\theequation{11.13}\begin{equation}
\left\{
\aligned 
\frac{\partial^2{\Pi'}^j}{\partial
{x'}^{l_1}{x'}^{l_2}}
\big(\phi(z),h(c)\big) 
&
=
\frac{{\sf S}_{l_1,l_2}^j
\left(
\big\{{\sf L}_{k_1'}{\sf L}_{k_2'}\varphi^{i'}(z)\big\}_{
1\leqslant k_1',k_2'\leqslant n}^{
1\leqslant i'\leqslant n+m}
\right)
}{ 
\left[
{\rm det}
\big(
{\sf L}_{k'}\,\phi^{l'}(z)\big)_{
1\leqslant k'\leqslant n}^{1\leqslant l'\leqslant n}
\right]^3
} 
\ \ \ \ \ \text{\rm and} 
\\
\frac{\partial^2{{\Pi'}^*}^j}{
\partial{a'}^{r_1}{a'}^{r_2}}
\big(f(c),\varphi(z)\big)
= 
& \
\frac{{{\sf S}^*}_{r_1,r_2}^j
\left(
\{
{\sf L}_{q_1'}^*{\sf L}_{q_2'}^*
h^{i'}(c)
\}_{1\leqslant q_1',q_2'\leqslant p}^{1\leqslant i'\leqslant p+m}
\right)}{\left[
{\rm det}\big(
{\sf L}_{q'}^*\,f^{r'}(c)
\big)_{1\leqslant q'\leqslant p}^{1\leqslant r'\leqslant p}
\right]^3
}.
\endaligned\right.
\end{equation}
By induction, for every $j$ with $1\leqslant j \leqslant m$ and every
two multiindices $\beta \in \N^n$ and $\delta \in \N^p$, there exists
two universal polynomials ${\sf S}_\beta^j$ and ${{\sf S}^*
}_\delta^j$ such that
\def\theequation{11.14}\begin{equation}
\left\{
\aligned
\frac{\partial^{\vert\beta \vert} 
{\Pi'}^j}{
\partial{x'}^\beta}
\big(\phi(z),h(c)\big)
&
=
\frac{
{\sf S}_\beta^j\left(
\big\{
{\sf L}^{\beta'}\varphi^{i'}(z)
\big\}_{
\vert\beta'\vert\leqslant\vert\beta\vert}^{
1\leqslant i'\leqslant n+m}
\right)
}{
\left[
{\rm det}
\big(
{\sf L}_{k'}\,\phi^{l'}(z)\big)_{
1\leqslant k'\leqslant n}^{1\leqslant l'\leqslant n}
\right]^{2\vert \beta \vert+1}} 
\ \ \ \ \ \text{\rm and} 
\\
\frac{\partial^{\vert\gamma\vert}{{\Pi'}^*}^j}{
\partial{a'}^\delta}
\big(f(c),\varphi(z)\big)
&
=
\frac{ 
{{\sf S}^*}_\delta^j
\left(
\big\{
{\sf L}^{*\delta'}h^{i'}(c)
\big\}_{\vert \delta' \vert 
\leqslant \vert \delta \vert}^{1\leqslant i'\leqslant p+m}
\right)
}{
\left[
{\rm det}\big(
{\sf L}_{q'}^*\,f^{r'}(c)
\big)_{1\leqslant q'\leqslant p}^{1\leqslant r'\leqslant p}
\right]^{2\vert\delta\vert+1}}.
\endaligned\right.
\end{equation}
Here, for $\beta' \in \N^n$, we denote by ${\sf L}^{\beta'}$ the
derivation of order $\vert \beta' \vert$ defined by $( {\sf
L}_1)^{ \beta_1 '} \cdots ({\sf L}_n )^{ \beta_n '}$. Similarly,
for $\delta'\in \N^p$, ${\sf L}^{* \delta'}$ denotes the derivation
of order $\vert \delta'\vert$ defined by $({\sf L}_1^* )^{
\delta_1' } \cdots ( {\sf L}_p^*)^{ \delta_p' }$.

Next, by the assumption that $\mathcal{ M}'$ is solvable with respect
to the parameters, there exist integers $j(1), \dots, j(p)$ with
$1\leqslant j(q)\leqslant m$ and multiindices $\beta (1), \dots, \beta
(p) \in \N^n$ with $\vert \beta (q) \vert \geqslant 1$ and $\max_{1
\leqslant q \leqslant p} \, \vert \beta (q) \vert = \kappa$ such that
the local $\K$-analytic map
\def\theequation{11.15}\begin{equation}
\K^{p+m}\ni
c' 
\longmapsto 
\left(
\big(
{\Pi'}^j(0,c')
\big)^{1\leqslant j\leqslant m},\
\left(
\frac{\partial^{\vert \beta(q)\vert} 
{\Pi'}^{j(q)}}{\partial
{x'}^{\beta(q)}}(0,c')
\right)_{1\leqslant q\leqslant p}
\right)\in\K^{p+m}
\end{equation}
is of rank $p+m$ at $c'=0$. Similarly, by the assumption that
$\mathcal{ M}'$ is solvable with respect to the variables, there exist
integers $j^\sim (1), \dots, j^\sim(n)$ with $1\leqslant j^\sim (l)
\leqslant m$ and multiindices $\delta (1), \dots, \delta (p) \in \N^n$
with $\vert \delta (q) \vert \geqslant 1$ and $\max_{1 \leqslant q
\leqslant p} \, \vert \delta (q) \vert = \kappa^*$ such that the local
$\K$-analytic map
\def\theequation{11.16}\begin{equation}
\K^{n+m}\ni
z'
\longmapsto 
\left(
\big(
{{\Pi'}^*}^j(0,z')
\big)^{1\leqslant j\leqslant m},\
\left(
\frac{\partial^{\vert\delta(l)\vert}
{{\Pi'}^*}^{j^\sim(l)}}{
\partial{a'}^{\delta(l)}}\big(0,z'\big)
\right)_{1\leqslant l\leqslant n}
\right)\in\K^{n+m}
\end{equation}
is of rank $n+m$ at $z' = 0$. We then consider from the first line
of~\thetag{ 11.14} only the $(p + m)$ equations written for $(j, 0)$,
$(j(q), \beta (q))$ and we solve $h(c)$ by means of the analytic
implicit function theorem; also, in the second line of~\thetag{ 11.14},
we consider the $(n+ m)$ equations written for $(j, 0)$, $(j^\sim (l),
\delta (l ))$ and we solve $\varphi (z)$. 
We get:
\def\theequation{11.17}\begin{equation}
\left\{
\aligned
h(c)
&
=
\widehat{H}
\left(\phi(z),
\frac{{\sf S}_{\beta(1)}^{j(1)}\left(
\big\{{\sf L}^{\beta'}\varphi^{i'}(z)\big\}_{
\vert\beta'\vert\leqslant\vert\beta(1)\vert}^{
1\leqslant i'\leqslant n+m}
\right)
}{
{\rm det}
\left[
\big(
{\sf L}_{k'}\,\phi^{l'}(z)\big)_{
1\leqslant k'\leqslant n}^{1\leqslant l'\leqslant n}
\right]^{2\vert\beta(1)\vert+1}
},\dots
\right. 
\\
&
\ \ \ \ \ \ \ \ \ \ \ \ \ \ \ \ \ \ \ \ \ 
\left.
\dots, \,
\frac{{\sf S}_{\beta(p)}^{j(p)}\left(
\big\{{\sf L}^{\beta'}\varphi^{i'}(z)\big\}_{
\vert\beta'\vert\leqslant\vert\beta(p)\vert}^{
1\leqslant i'\leqslant n+m}
\right)
}{
{\rm det}
\left[
\big(
{\sf L}_{k'}\,\phi^{l'}(z)\big)_{
1\leqslant k'\leqslant n}^{1\leqslant l'\leqslant n}
\right]^{2\vert\beta(p)\vert+1}
}\right), 
\\
\varphi(z) 
&
= 
\widehat{\Phi}
\left(
f(c),\,
\frac{{{\sf S}^*}_{\delta(1)}^{j^\sim(1)}\left(
\big\{
{\sf L}^{*\delta'}h^{i'}(c)
\big\}_{\vert \delta' \vert 
\leqslant \vert \delta(1) \vert}^{1\leqslant i'\leqslant p+m}
\right)
}{
\left[
{\rm det}\big(
{\sf L}_{q'}^*\,f^{r'}(c)
\big)_{1\leqslant q'\leqslant p}^{1\leqslant r'\leqslant p}
\right]^{2\vert\delta(1)\vert+1}},\dots 
\ \ \ \ \ \ \ \right. 
\\
&
\ \ \ \ \ \ \ \ \ \ \ \ \ \ \ \ \ \ \ \ \ 
\left.
\dots,\,
\frac{{{\sf S}^*}_{\delta(n)}^{j^\sim(n)}\left(
\big\{
{\sf L}^{*\delta'}h^{i'}(c)
\big\}_{\vert\delta'\vert 
\leqslant \vert \delta(n)\vert}^{1\leqslant i'\leqslant p+m}
\right)
}{
\left[
{\rm det}\big(
{\sf L}_{q'}^*\,f^{r'}(c)
\big)_{1\leqslant q'\leqslant p}^{1\leqslant r'\leqslant p}
\right]^{2\vert\delta(n)\vert+1}}
\right),
\endaligned\right.
\end{equation}
for $(z, c) \in \mathcal{ M }$. The maps $\widehat{ H}$ and $\widehat{
\Phi }$ depend only on ${\Pi' }$, ${\Pi' }^*$.

\def\thelemma{11.18}\begin{lemma}
For every $\beta' \in \N^n$, there exists a universal polynomial ${\sf
P}_{ \beta'}$ in the jet variables $J_z^{ \vert \beta' \vert}$ having
$\K$-analytic coefficients in $(z, c)$ which depends only on $\Pi$,
$\Pi^*$ such that, for $i' = 1, \dots, n+ m${\rm :}
\def\theequation{11.19}\begin{equation}
{\sf L}^{\beta'} 
\varphi^{i'}(z) 
\equiv
{\sf P}_{\beta'}
\left(
z,\,c,\,J_z^{\vert\beta'\vert}\varphi^{i'}(z)
\right).
\end{equation}
A similar property holds for ${\sf L}^{ * \delta'} h^{ i'} (c)$.
\end{lemma}

We deduce that there exist two local $\K$-analytic mappings $\Phi_0^0$ 
and $H_0^0$ such that we can write
\def\theequation{11.20}\begin{equation}
\left\{
\aligned
\varphi(z)
&
=
\Phi_0^0
\big(
z,\, c, \, J_c^{\kappa^*}h(c)
\big), 
\\
h(c)
&
=
H_0^0
\big(
z,\,c,\,J_z^\kappa\varphi(z)
\big), 
\endaligned\right.
\end{equation}
for $(z, c) \in \mathcal{ M}$. Concretely, this
means that we have two equivalent
pairs of formal identities
\def\theequation{11.21}\begin{equation}
\aligned
\varphi(z)
&
\equiv
\Phi_0^0
\big(z,\,a,\,\Pi^*(a,z),\, 
J_c^{\kappa^*}h(a,\Pi^*(a,z))
\big)
\\
\varphi\big(x,\Pi(x,c)\big)
&
\equiv
\Phi_0^0
\big(
x,\,\Pi(x,c),\,c,\,J_c^{\kappa^*}h(c)
\big)
\\
h(c) 
&
\equiv
H_0^0
\big(
x,\,\Pi(x,c),\,c,\, 
J_z^{\kappa}\varphi(x,\Pi(x,c))
\big)
\\
h\big(a,\Pi^*(a,z)\big)
&
\equiv
H_0^0
\big(
z,\,a,\,\Pi^*(a,z),\,
J_z^\kappa\varphi(z)
\big)
\endaligned
\end{equation}
in $\K\dl a, z\dr^{ n + m }$ and in $\K \dl x, c \dr^{p + m}$. We
notice that, whereas $\varphi$ and $h$ are {\it a priori}\, only
purely formal, by construction, $\Phi_0^0$ and $H_0^0$ are
$\K$-analytic near $\big( 0,0,J_c^{ \kappa^* } h (0) \big)$ and
near $\big( 0, 0, J_z^\kappa \varphi (0) \big)$.

Next, we introduce the following vector fields with $\K$-analytic
coefficients tangent to $\mathcal{ M}$:
\def\theequation{11.22}\begin{equation}
\left\{
\aligned
{\sf V}_j
:= 
& \
\frac{\partial }{\partial y^j}+
\sum_{l=1}^m\, 
\frac{\partial{\Pi^*}^l}{
\partial y^j}(a,z)\, 
\frac{\partial }{\partial b^l}, 
\ \ \ \ \ j=1,\dots,m
\ \ \ \ \ \text{\rm and}
\\
{\sf V}_j^* 
:=
& \
\frac{\partial}{\partial b^j}
+
\sum_{l=1}^m\, 
\frac{\partial\Pi^l}{
\partial b^j}(x,c) \, 
\frac{\partial }{\partial y^l}, 
\ \ \ \ \ j=1,\dots,m. 
\endaligned\right.
\end{equation}
Indeed, we check that ${\sf V}_{ j_1 }[ b^{ j_2} - { \Pi^* }^{ j_2 }(
a, z)] \equiv 0$ and that ${\sf V}_{ j_1 }^*[ y^{ j_2 } - \Pi^{ j_2 }(
x, c)] \equiv 0$.

For $\delta ' \in \N^m$, we observe that ${\sf V}^{ \delta '} \varphi
= \frac{ \partial^{ \vert \delta '\vert} \varphi }{ \partial y^{
\delta'}}$. Applying then ${\sf L}^{ \beta '}$ with $\beta' \in
\N^n$, we get for $ i =1, \dots, n+m$:
\def\theequation{11.23}\begin{equation}
{\sf L}^{\beta'}{\sf V}^{\delta'}
\varphi^i(z)
=
{\sf Q}_{\beta',\delta'}
\big(
z,c,J_z^{\vert\beta'\vert+\vert\delta'\vert}\varphi^i(z)
\big),
\end{equation}
with ${\sf Q}_{ \beta', \delta'}$ universal. Since the $n + m$ vector
fields ${\sf L }_k$ and ${\sf V }_j$, having coefficients depending on
$(z, c)$, span the tangent space to $\K_x^n \times \K_y^m$, the change
of basis of derivations yields, by induction, the following.

\def\thelemma{11.24}\begin{lemma}
For every $\alpha \in \N^{ n+m }$, there exists a universal polynomial
${\sf P}_\alpha$ in its last variables with coefficients being
$\K$-analytic in $(z, c)$ and depending only on $\Pi$, $\Pi^*$ such
that, for $i = 1, \dots, n + m${\rm :}
\def\theequation{11.25}\begin{equation}
\partial_z^\alpha
\varphi^i(z)
\equiv
{\sf P}_\alpha
\Big(z,c, 
\big(
{\sf L}^{\beta'}{\sf V}^{\delta'}
\varphi^i(z)
\big)_{\vert\beta'\vert +
\vert\delta'\vert\leqslant 
\vert\alpha\vert}
\Big).
\end{equation}
\end{lemma}

We are now in position to state and to prove the first fundamental
technical lemma which generalizes the two formulas~\thetag{11.20} to
arbitrary jets.

\def\thelemma{11.26}\begin{lemma}
For every $\lambda \in \N$, there exist two local
$\K$-analytic maps, $\Phi_0^\lambda$ valued in $\K^{ (n+ m) C_{ n+ m+
\lambda }^\lambda }$, and $H_0^\lambda$ valued in $\K^{ (p+m ) C_{ p+m+
\lambda }^\lambda }$, such that{\rm :}
\def\theequation{11.27}\begin{equation}
\left\{
\aligned
J_z^\lambda\varphi(z)
&
\equiv
\Phi_0^\lambda
\big(z,\,c,\,
J_c^{\kappa^*+\lambda}
h(c)\big), 
\\
J_c^\lambda h(c) 
&
\equiv
H_0^\lambda
\big(z,\,c,\, 
J_z^{\kappa+\lambda}\varphi(z)
\big).
\endaligned\right.
\end{equation}
\end{lemma}

\proof
Consider for instance the first line. To obtain it, it suffices to
apply the derivations ${\sf L}^{ \beta'} {\sf V}^{ \delta'}$ with
$\vert \beta'\vert + \vert \delta' \vert \leqslant \lambda$ to the
first line of~\thetag{11.20}, to use the chain rule and to apply
Lemma~11.24.
\endproof

Let $\theta \in \K^l$, $l\in \N$, 
let $Q( \theta)= \big( Q_1( \theta),
\dots, Q_{ n+ 2m+p} ( \theta) \big) \in \K \dl \theta \dr^{ n+ 2m+p}$
and let $a_1 \in \K^p$. As the multiple flow of $\mathcal{\sf L}^*$
given by~\thetag{10.10} does not act on the variables $(x, y)$, we have
the trivial but crucial property:
\def\theequation{11.28}\begin{equation}
\varphi\left(
{\sf L}_{a_1}^*(Q(\theta))
\right)
\equiv
\varphi\left(\pi_z(
{\sf L}_{a_1}^*
(Q(\theta)))\right)
\equiv
\varphi\left(
\pi_z(
Q(\theta))
\right)
\equiv
\varphi\left(
Q(\theta)
\right).
\end{equation}
At the end, we allow to suppress the projection $\pi_z$: this slight
abuse of notation will lighten slightly the writting of further
formulas. More generally, for $\lambda \in \N$, $a_1 \in \K^p$,
$x_1 \in \K^n$:
\def\theequation{11.29}\begin{equation}
\aligned
J_z^\lambda\varphi 
\big(
{\sf L}_{a_1}^*
(Q(\theta))
\big)
&
\equiv
J_z^\lambda
\varphi 
\big(Q(\theta)\big)
\ \ \ \ \ \text{\rm and}
\\
J_c^\lambda h
\big(
{\sf L}_{x_1}
(Q(\theta))
\big)
&
\equiv
J_c^\lambda h
\big(
Q(\theta)
\big).
\endaligned
\end{equation}
As a consequence, for $2k$ even and for $2k+1$ odd, we have the
following four cancellation relations, useful below (we drop $\pi_z$
and $\pi_c$ after $J_z^\lambda \varphi$ and after $J_c^\lambda h$):
\def\theequation{11.30}\begin{equation}
\left\{
\aligned
J_z^\lambda\varphi\big(\Gamma_{2k}([xa]_{2k})\big)
&
\equiv
J_z^\lambda\varphi\big(\Gamma_{2k-1}([xa]_{2k-1})\big), 
\\
J_c^\lambda h\big(\Gamma_{2k}^*([a x]_{2k})\big)
&
\equiv
J_c^\lambda h\big(\Gamma_{2k-1}^*([a x]_{2k-1})\big), 
\\
J_z^\lambda\varphi\big(\Gamma_{2k+1}^*([a x]_{2k+1})\big)
&
\equiv
J_z^\lambda\varphi\big(\Gamma_{2k}^*([a x]_{2k})\big), 
\\
J_c^\lambda h\big(\Gamma_{2k+1}([xa]_{2k+1})\big)
&
\equiv
J_c^\lambda h\big(\Gamma_{2k}([xa]_{2k})\big).
\endaligned\right.
\end{equation}

We are now in position to state and to prove the second main
technical proposition.

\def\theproposition{11.31}\begin{proposition}
For every even chain-length $2k$ and for every jet-height $\lambda$, 
there exist two local $\K$-analytic maps, $\Phi_{ 2k}^\lambda$
valued in $\K^{ (n+ m) C_{ n+ m+ \lambda }^{ \lambda}}$, and $H_{
2k}^\lambda$ valued in $\K^{ (p+ m) C_{p+ m + \lambda }^\lambda }$
such that{\rm :}
\def\theequation{11.32}\begin{equation}
\left\{
\aligned
J_z^\lambda\varphi\left(
\Gamma_{2k}^*
\big([a x]_{2k}\big)
\right)
& 
\equiv
\Phi_{2k}^\lambda
\left(
[a x]_{2k} , \, 
J_z^{k(\kappa+
\kappa^*)+\lambda}\,
\varphi(0)
\right)
\ \ \ \ \ \text{\rm and}
\\
J_c^\lambda h 
\left(
\Gamma_{2k}
\big([xa]_{2k}\big) 
\right) 
&
\equiv
H_{2k}^\lambda
\left([xa]_{2k},\,J_c^{k(\kappa+\kappa^*)+\lambda}\,
\varphi(0)
\right).
\endaligned\right.
\end{equation}
Similarly, for every odd chain length $2k+1$ and for every jet eight
$\lambda$, there exist two local $\K$-analytic maps, $\Phi_{
2k+1}^\lambda$ valued in $\K^{ (n+ m) C_{ n+ m+ \lambda }^{ \lambda
}}$ and $H_{ 2k+1}^\lambda$ valued in $\K^{ (p+m) C_{ p+m + \lambda
}^\lambda }$, such that{\rm :}
\def\theequation{11.33}\begin{equation}
\left\{
\aligned
J_z^\lambda\varphi
\left(
\Gamma_{2k+1}\big([xa]_{2k+1}\big)
\right) 
&
\equiv 
\Phi_{2k+1}^\lambda
\left([xa]_{2k+1},\, 
J_c^{k\kappa+
(k+1)\kappa^*+\lambda}\,
h(0)
\right), 
\\
J_c^\lambda h\left(
\Gamma_{2k+1}^*\big([a x]_{2k+1}\big)
\right) 
&
\equiv 
H_{2k+1}^\lambda 
\left(
[a x]_{2k+1},\, 
J_z^{(k+1)\kappa+k\kappa^*+\lambda}\,
\varphi(0)
\right).
\endaligned\right.
\end{equation}
These maps depend only on $\Pi$, $\Pi^*$, $\Pi'$, ${\Pi'}^*$.
\end{proposition}

\proof
For $2 k+1 =1$, we replace $(z, c)$ by $\Gamma_1( [xa]_1)$ in the
first line of~\thetag{11.27} and by $\Gamma_1^*( [a x ]_1)$ in the
second line. Taking crucially account of the cancellation
properties~\thetag{11.29}, we get:
\def\theequation{11.34}\begin{equation}
\left\{
\aligned
J_z^\lambda\varphi
\big(\Gamma_1([xa]_1)\big)
&
\equiv
\Phi_0^\lambda\left(\Gamma_1([xa]_1),\, 
J_c^{\kappa^*+\lambda} 
h\big(\Gamma_1([xa]_1)\big)\right) 
\\
&
\equiv
\Phi_0^\lambda
\left(\Gamma_1([xa]_1), \, 
J_c^{\kappa^*+\lambda} h(0))\right) 
\\
&
=:
\Phi_1^\lambda
\left(
[xa]_1,\,J_c^{\kappa^*+\lambda}h(0)
\right), 
\\
J_c^\lambda h
\big(\Gamma_1^*([a x]_1)\big) 
& \
\equiv
H_0^\lambda\left(
\Gamma_1^*([a x]_1), \, J_z^{\kappa+\lambda} 
\varphi\big(\Gamma_1^*([a x]_1)\big)
\right) 
\\
& \
\equiv
H_0^\lambda\left(
\Gamma_1^*([ax]_1),\,J_z^{\kappa+\lambda} 
\varphi(0)\right) 
\\
&
=:
H_1^\lambda\left(
[ax]_1,\, 
J_z^{\kappa+\lambda} 
\varphi(0)
\right).
\endaligned\right.
\end{equation}
Here, the third line defines $\Phi_1^\lambda$ and the sixth line
defines $H_1^\lambda$. Thus, the proposition holds for $2k+1 =1$.

The rest of the proof proceeds by induction. We treat only the
induction step from an odd chain-length $2k+1$ to an even chain-length
$2k+2$, the other induction step being similar.

To this aim, we replace the variables $(z, c)$ in the first line
of~\thetag{11.27} by $\Gamma_{ 2k+2}^*( [a x]_{ 2k+2})$.
Taking account of the cancellation property
and of the induction assumption:
\begin{small}
\def\theequation{11.35}\begin{equation}
\aligned
J_z^\lambda
\varphi\left( 
\Gamma_{2k+2}^*
\big([ax]_{2k+2}\big)
\right)
&
\equiv 
\Phi_0^\lambda
\left(\Gamma_{2k+2}^*([ax]_{2k+2}),\, 
J_c^{\kappa^*+\lambda}h
\big(\Gamma_{2k+2}^*([ax]_{2k+2})\big)
\right) 
\\
&
\equiv 
\Phi_0^\lambda\left(\Gamma_{2k+2}^*([ax]_{2k+2}),\, 
J_c^{\kappa^*+\lambda}h
\big(\Gamma_{2k+1}^*([ax]_{2k+1})\big)
\right) 
\\
& \
\equiv 
\Phi_0^\lambda\left(\Gamma_{2k+2}^*([ax]_{2k+2}),\,
H_{2k+1}^{\kappa^*+\lambda}
\left([ax]_{2k+1},\,
J_c^{(k+1)(\kappa+\kappa^*)+\lambda}
\varphi(0)\right) \right)
\\
&
=:
\Phi_{2k+2}^\lambda\left(
[ax]_{2k+2},\,
J_c^{(k+1)(\kappa+\kappa^*)+\lambda}
\varphi(0)\right), 
\endaligned
\end{equation}
\end{small}
The last line defines $\Phi_{ 2k +2}^\lambda$. Similarly, we replace
$(z, c)$ in the second line of~\thetag{11.27} by $\Gamma_{ 2k+2 }([x a
]_{ 2k +2 })$. Taking account of the cancellation property and of the
induction assumption:
\begin{small}
\def\theequation{11.36}\begin{equation}
\aligned
J_c^\lambda h
\big(\Gamma_{2k+2}([xa]_{2k+2})\big) 
&
\equiv 
H_0^\lambda\left(
\Gamma_{2k+2}([xa]_{2k+2}),\, 
J_c^{\kappa+\lambda}\varphi
\big(\Gamma_{2k+2}([xa]_{2k+2})\big)
\right) 
\\
&
\equiv 
H_0^\lambda\left(
\Gamma_{2k+2}([xa]_{2k+2}),\, 
J_c^{\kappa+\lambda}\varphi\big(\Gamma_{2k+1}([xa]_{2k+1})\big)
\right) 
\\
&
\equiv 
H_0^\lambda\left(
\Gamma_{2k+2}([xa]_{2k+2}),\,
\Phi_{2k+1}^{\kappa+\lambda}\left(
[xa]_{2k+1},\, J_c^{(k+1)(\kappa+
\kappa^*)+\lambda}h(0)
\right)
\right) 
\\
&
=:
H_{2k+2}^\lambda\left(
[xa]_{2k+2},\,J_c^{(k+1)(\kappa+\kappa^*)+\lambda}
h(0)
\right).
\endaligned
\end{equation}
\end{small}
This completes the proof.
\endproof

\medskip\noindent
{\it End of the proof of Theorem~11.6.} With $(\mu, \mu^*)$ being the
type of $({\sf F}_{ \sf v}, \, {\sf F}_{ \sf p })$ and with $[a x]_{2
\mu^*}^0$ given by Corollary~10.31, the rank property~\thetag{ 10.32}
insures the existence of an affine $(n + m)$-dimensional space $H
\subset \K^{ \mu^* (p + n)}$ passing through $[a x]_{ 2 \mu^* }^0$ and
equipped with a local parametrization
\def\theequation{11.37}\begin{equation}
\K^{ n+m} 
\ni 
s 
\mapsto 
[ax]_{2\mu^*}(s)
\in 
H
\end{equation}
satisfying $[a x]_{ 2 \mu^*} (0)= [a x]_{2 \mu^*}^0$, such that the
map 
\def\theequation{11.38}\begin{equation}
\K^{n+m}\ni s \longmapsto 
\pi_z 
\big(
\Gamma_{2\mu^*}^*([a x]_{2\mu^*}(s))
\big)
=: 
z(s)\in\K^{n+m}
\end{equation}
is a local diffeomorphism fixing $0 \in \K^{ n+m}$. Replacing $z$ by
$z(s)$ in $\varphi (z)$ and applying the formula in the first line
of~\thetag{11.32} with $\lambda = 0$ and with $k = 2\mu^*$, we obtain
\def\theequation{11.39}\begin{equation}
\aligned
\varphi(z(s))
&
=
\varphi
\left(\pi_z\big(
\Gamma_{2\mu^*}^*([ax]_{2\mu^*}(s)\big)
\right)
\\
&
=
\varphi
\left(
\Gamma_{2\mu^*}^*\big([ax]_{2\mu^*}(s)\big)
\right)
\\
&
\equiv
\Phi_{2\mu^*}^0\left(
[a x]_{2\mu^*}(s), \, 
J_z^{\mu^*(\kappa+\kappa^*)}
\varphi(0)
\right).
\endaligned
\end{equation}
Inverting $s\mapsto z = z(s)$ as
$z \mapsto s = s(z)$, we finally get
\def\theequation{11.40}\begin{equation}
\aligned 
\varphi(z)
=
\varphi(z(s(z)))
&
\equiv 
\Phi_{2\mu^*}^0\left(
[a x]_{2\mu^*}(s(z)),\, 
J_z^{\mu^*(\kappa+\kappa^*)}\varphi(0)
\right) 
\\
&
=:
\Phi_\ell\left(
z,\,J_z^{\mu^*(\kappa+\kappa^*)}\varphi(0)
\right), 
\endaligned
\end{equation}
with $\ell := \mu^* (\kappa + \kappa^*)$, where the last line
defines $\Phi_\ell$. In conclusion, we have derived the first line
of~\thetag{ 11.7}. The second one is obtained similarly.

If $\Pi$, $\Pi^*$, $\Pi'$, ${\Pi'}^*$ are algebraic, so are
$\Gamma_k$, $\Gamma_k^*$, $\widehat{ H }$, $\widehat{ \Phi }$,
$\Phi_0^\lambda$, $H_0^\lambda$, $\Phi_k^\lambda$, $H_k^\lambda$ and
$\Phi_\ell$, $H_{ \ell^* }$.

The proof of Theorem~11.6 is complete.
\endproof

\newpage

$\:$\bigskip\bigskip

\begin{center}
{\large\bf
{\Large\bf II:}~Explicit prolongations of infinitesimal Lie symmetries
}
\end{center}

\begin{center}
\begin{minipage}[t]{12cm}
\baselineskip =0.35cm
{\tiny

\bigskip
\bigskip

\centerline{\bf Table of contents}

\medskip

{\bf 1.~Jet spaces 
and prolongations \dotfill 57.}

{\bf 2.~One independent variable and
one dependent variable \dotfill 62.}

{\bf 3.~Several independent variables and
one dependent variable \dotfill 68.}

{\bf 4.~One independent variable and
several dependent variables \dotfill 84.}

{\bf 5.~Several independent variables and
several dependent variables \dotfill 90.}

}\end{minipage}
\end{center}

\section*{\S1.~Jet spaces and prolongations}

\subsection*{1.1.~Choice of notations for the jet space variables}
Let $\K = \R$ or $\C$. Let $n\geqslant 1$ and $m\geqslant 1$ be two
positive integers and consider two sets of variables $x= (x^1, \dots,
x^n) \in \K^n$ and $y = (y^1, \dots, y^m)$. In the classical theory of
Lie symmetries of partial differential equations, one considers
certain differential systems whose (local) solutions should be
mappings of the form $y = y(x)$. We refer to~\cite{ ol1986} and
to~\cite{ bk1989} for an exposition of the fundamentals of the
theory. Accordingly, the variables $x$ are usually called {\sl
independent}, whereas the variables $y$ are called {\sl
dependent}. Not to enter in subtle regularity considerations (as
in~\cite{me2005b}), we shall assume $\mathcal{ C }^\infty$-smoothness
of all functions throughout this paper.

Let $\kappa \geqslant 1$ be a positive integer. For us, in a very
concrete way (without fiber bundles), the {\sl $\kappa$-th jet space}
$\mathcal{ J}_{ n, m}^\kappa$ consists of the space $\K^{ n+m+ m
\frac{ (n+m)!}{n! \ m!}}$ equipped with the affine coordinates
\def\theequation{1.2}\begin{equation}
\left(
x^i, y^j, y_{i_1}^j, y_{i_1,i_2}^j, 
\dots\dots, 
y_{i_1,i_2,\dots,i_\kappa}^j
\right),
\end{equation}
having the symmetries
\def\theequation{1.3}\begin{equation}
y_{i_1,i_2,\dots,i_\lambda}^j
=
y_{i_{\sigma(1)},i_{\sigma(2)},\dots,i_{\sigma(\lambda)}}^j,
\end{equation}
for every $\lambda$ with $1\leqslant \lambda \leqslant \kappa$ and for
every permutation $\sigma$ of the set $\{1, \dots, \lambda \}$. The
variable $y_{i_1, i_2, \dots, i_\lambda}^j$ is an independent
coordinate corresponding to the $\lambda$-th partial derivative
$\frac{ \partial^\lambda y^j}{ \partial x^{ i_1} \partial x^{ i_2}
\cdots \partial x^{ i_\lambda }}$. So the symmetries~\thetag{ 1.3} are
natural.

In the classical Lie theory (\cite{ ol1979},
\cite{ ol1986}, \cite{ bk1989}), all the
geometric objects: point transformations, vector fields, {\it etc.},
are local, defined in a neighborhood of some point lying in some
affine space $\K^N$. However, in this paper, the original geometric
motivations are rapidly forgotten in order to focus on combinatorial
considerations. Thus, to simplify the presentation, we shall not
introduce any special notation to speak of certain local open subsets
of $\K^{n+ m}$, or of the jet space $\mathcal{ J}_{ n, m}^\kappa =
\K^{n+m+ m \frac{ (n+m)! }{ n! \ m! }}$, {\it etc.}: we will always
work in global affine spaces $\K^N$.

\subsection*{ 1.4.~Prolongation $\varphi^{ (\kappa)}$ 
of a local diffeomorphism $\varphi$ to the $\kappa$-th jet space} In
this paragraph, we recall how the prolongation of a diffeomorphism to
the $\kappa$-th jet space is defined (\cite{ ol1979},
\cite{ ol1986}, \cite{ bk1989}).

Let $x_*\in \K^n$ be a central fixed point and let $\varphi : \K^{n+m}
\to \K^{n+m}$ be a diffeomorphism whose Jacobian matrix is close to
the identity matrix, at least in a small neighborhood of $x_*$. Let
\def\theequation{1.5}\begin{equation}
J_{x_*}^\kappa
:=
\left(
x_*^i,y_{*i_1}^j,
y_{*i_1,i_2}^j,
\dots\dots,
y_{*i_1,i_2,\dots,i\kappa}^j
\right)
\in
\left.
\mathcal{J}_{n,m}^\kappa
\right\vert_{x_*}
\end{equation}
be an arbitrary $\kappa$-jet based at $x_*$. The goal is to defined
its transformation $\varphi^{(\kappa)} ( J_{x_*}^\kappa)$ by
$\varphi$.

To this aim, choose an arbitrary mapping $\K^n\ni x \mapsto g(x) \in
\K^m$ defined at least in a neighborhood of $x_*$ and representing
this $\kappa$-th jet, {\it i.e.} satisfying
\def\theequation{1.6}\begin{equation}
y_{*i_1,\dots,i_\lambda}^j 
=
\frac{\partial^\lambda g^j}{
\partial x^{i_1}\cdots \partial x^{i_\lambda}}
(x_*),
\end{equation}
for every $\lambda\in \N$ with $0 \leqslant \lambda \leqslant \kappa$,
for all indices $i_1, \dots, i_\lambda$ with $1 \leqslant i_1, \dots,
i_\lambda \leqslant n$ and for every $j\in \N$ with $1\leqslant
j\leqslant m$. In accordance with the splitting $(x, y)\in \K^n\times
\K^m$ of coordinates, split the components of the diffeomorphism
$\varphi$ as $\varphi = (\phi, \psi) \in \K^n\times \K^m$. Write
$\left( \overline{x}, \overline{y} \right)$ the coordinates in the
target space, so that the diffeomorphism $\varphi$ is:
\def\theequation{1.7}\begin{equation}
\K^{n+m}
\ni
(x,y)
\longmapsto
\left(
\overline{x}, \overline{y}
\right) 
= 
\big(
\phi(x, y), \psi (x, y)
\big)
\in
\K^{n+m}.
\end{equation}
Restrict the variables $(x, y)$ to belong to the graph of $g$, namely
put $y:= g(x)$ above, which yields
\def\theequation{1.8}\begin{equation}
\left\{
\aligned
\overline{x} 
&
=
\phi(x, g(x)), \\
\overline{y}
&
=
\psi(x,g(x)).
\endaligned\right.
\end{equation}
As the differential of $\varphi$ at $x_*$ is close to the identity,
the first family of $n$ scalar equations may be solved with respect to
$x$, by means of the implicit function theorem. Denote $x
= \overline{ \chi}(\overline{ x})$ the resulting mapping, satisfying
by definition
\def\theequation{1.9}\begin{equation}
\overline{x}
\equiv 
\phi\left(
\overline{ \chi}(\overline{ x}), 
g(\overline{ \chi}(\overline{ x}))
\right).
\end{equation}
Replace $x$ by $\overline{ \chi}(\overline{ x})$ in the second family
of $m$ scalar equations~\thetag{1.8} above, which yields:
\def\theequation{1.10}\begin{equation}
\overline{y}
=
\psi\left(
\overline{ \chi}(\overline{ x}),
g(\overline{ \chi}(\overline{ x}))
\right).
\end{equation}
Denote simply by $\overline{ y} = \overline{ g} ( \overline{ x})$ this
last relation, where $\overline{ g} ( \cdot) := \psi \left( \overline{
\chi} (\cdot), g ( \overline{ \chi} (\cdot)) \right)$.

In summary, the graph $y=g(x)$ has been transformed to the graph
$\overline{y} = \overline{ g} (\overline{ x})$ by the diffeomorphism
$\varphi$.

Define then the {\sl transformed jet $\varphi^{(\kappa)} \left(
J_{x_*}^\kappa \right)$} to be the $\kappa$-th jet of $\overline{ g}$
at the point $\overline{ x}_* := \phi ( x_*)$, namely:
\def\theequation{1.11}\begin{equation}
\varphi^{(\kappa)}
\left(
J_{x_*}^\kappa
\right)
:= 
\left(
\frac{\partial^\lambda \overline{ g}^j}{
\partial \overline{ x}^{i_1} 
\cdots
\partial \overline{ x}^{i_\lambda}}
(\overline{ x}_*)
\right)_{
1\leqslant i_1, \dots, i_\lambda \leqslant n, \
0 \leqslant \lambda \leqslant \kappa
}^{
1\leqslant j\leqslant m}
\left.
\in\mathcal{J}_{n,m}^\kappa
\right\vert_{\overline{x}_*}.
\end{equation}
It may be shown that this jet does not depend on the choice of a local
graph $y = g(x)$ representing the $\kappa$-th jet $J_{x_*}^\kappa$ at
$x_*$. Furthermore, if $\pi_\kappa := \mathcal{ J}_{n, m}^\kappa \to
\K^m$ denotes the canonical projection onto the first factor, the
following diagram commutes:
$$
\diagram \mathcal{J}_{n,m}^\kappa \rto^{\varphi^{(\kappa)}} 
\dto_{\pi_\kappa} 
& \mathcal{J}_{n,m}^\kappa
\dto^{\pi_\kappa} \\
\K^{n+m} \rto^{\varphi} & 
\K^{n+m}
\enddiagram.
$$

\subsection*{1.12.~Inductive formulas for the 
$\kappa$-th prolongation $\varphi^{ (\kappa)}$} To present them, we
change our notations. Instead of $(\overline{ x}, \overline{ y})$, as
coordinates in the target space $\K^n \times \K^m$, we shall use
capital letters:
\def\theequation{1.13}\begin{equation}
\left(
X^1,\dots,X^n,Y^1,\dots,Y^m
\right).
\end{equation}
In the source space $\K^{n+m}$ equipped with the coordinates $(x, y)$,
we use the jet coordinates~\thetag{ 1.2} on the associated
$\kappa$-th jet space. In the target space $\K^{n+m}$ equipped with
the coordinates $(X, Y)$, we use the coordinates
\def\theequation{1.14}\begin{equation}
\left(
X^i, Y^j, Y_{X^{i_1}}^j, 
Y_{X^{i_1}X^{i_2}}^j,
\dots\dots,
Y_{X^{i_1}X^{i_2}\dots X^{i_\kappa}}^j
\right)
\end{equation}
on the associated $\kappa$-th jet space; to 
avoid confusion with $y_{i_1}, y_{ i_1, i_2}, \dots$ in subsequent
formulas,
we do not write 
$Y_{ i_1}, Y_{ i_1, i_2}, \dots$. In these notations, the
diffeomorphism $\varphi$ whose first order approximation is close to
the identity mapping in a neighborhood of $x_*$
may be written under the form:
\def\theequation{1.15}\begin{equation}
\varphi : 
\
\big(
x^{i'},y^{j'}
\big) 
\mapsto 
\left(
X^i,Y^j
\right) 
= 
\left(
X^i(x^{i'},y^{j'}), \ Y^j(x^{i'},y^{j'})
\right),
\end{equation}
for some $\mathcal{ C}^\infty$-smooth functions $X^i(x^{i'},y^{j'})$,
$i = 1, \dots, n$, and $Y^j(x^{i'},y^{j'})$, $j = 1, \dots, m$. The
first prolongation $\varphi^{(1)}$ of $\varphi$ may be written 
under the form:
\def\theequation{1.16}\begin{equation}
\varphi^{(1)}
: \
\left(
x^{i'},y^{j'}, y_{i_1'}^{j'}
\right) 
\longmapsto
\left(
X^i(x^{i'},y^{j'}), \
Y^j(x^{i'},y^{j'}), \
Y_{X^{i_1}}^j
\left(
x^{i'}, y^{j'}, y_{i_1'}^{j'}
\right)
\right),
\end{equation}
for some functions $Y_{X^{i_1}}^j \left( x^{i'}, y^{j'}, y_{i_1'}^{j'}
\right)$ which depend on the pure first jet variables
$y_{i_1'}^{j'}$. The way how these functions depend on the first order
partial derivatives functions $X_{x^{i'}}^i$, $X_{y^{j'}}^i$,
$Y_{x^{i'}}^j$, $Y_{y^{j'}}^j$ and on the pure first jet variables
$y_{i_1'}^{j'}$ is provided (in principle) by the following compact
formulas (\cite{ bk1989}):
\def\theequation{1.17}\begin{equation}
\left(
\begin{array}{c}
Y_{X^1}^j \\
\vdots \\
Y_{X^n}^j \\
\end{array}
\right)
=
\left(
\begin{array}{ccc}
D_1^1 X^1 & \cdots & D_1^1 X^n \\
\vdots & \cdots & \vdots \\
D_n^1 X^1 & \cdots & D_n^1 X^n \\
\end{array} 
\right)^{-1}
\left(
\begin{array}{c}
D_1^1 Y^j \\
\vdots \\
D_n^1 Y^j \\
\end{array}
\right), 
\end{equation}
where, for $i' = 1, \dots, n$, the symbol
$D_{ i'}^1$ denotes the {\sl $i'$-th
first order total differentiation operator}:
\def\theequation{1.18}\begin{equation}
D_{i'}^1 
:= 
\frac{\partial}{\partial x^{i'}}
+
\sum_{j'=1}^m\,y_{i'}^{j'}\,\frac{\partial}{\partial y^{j'}}.
\end{equation}
Striclty speaking, these formulas~\thetag{ 1.17} are not explicit,
because an inverse matrix is involved and because the terms $D_{i'}^1
X^i$, $D_{i'}^1 Y^j$ are not developed. However, it would be
feasible and
elementary to write down the corresponding totally explicit complete
formulas for the functions $Y_{X^{i_1}}^j = Y_{X^{i_1}}^j \left(
x^{i'}, y^{j'}, y_{i_1'}^{j'} \right)$.

Next, the second prolongation $\varphi^{ (2)}$ is of the form
\def\theequation{1.19}\begin{equation}
\varphi^{(2)}
: \
\left(
x^{i'},y^{j'}, y_{i_1'}^{j'},
y_{i_1',i_2'}^{j'}
\right) 
\longmapsto 
\left(
\varphi^{(1)}
\left(
x^{i'},y^{j'}, y_{i_1'}^{j'}
\right), 
\
Y_{X^{i_1}X^{i_2}}^j
\left(
x^{i'}, y^{j'}, y_{i_1'}^{j'}, y_{i_1',i_2'}^{j'}
\right)
\right),
\end{equation}
for some functions $Y_{X^{i_1} X^{i_2}}^j \left( x^{ i'}, y^{ j'}, y_{
i_1' }^{j'}, y_{i_1', i_2'}^{j'} \right)$ which depend on the pure
first and second jet variables. For $i = 1, \dots, n$, the
expressions of $Y_{X^{i_1}X^i}^j$ are given by the following compact
formulas (again \cite{ bk1989}):
\def\theequation{1.20}\begin{equation}
\left(
\begin{array}{c}
Y_{X^{i_1}X^1}^j \\
\vdots \\
Y_{X^{i_1}X^n}^j \\
\end{array}
\right)
=
\left(
\begin{array}{ccc}
D_1^1 X^1 & \cdots & D_1^1 X^n \\
\vdots & \cdots & \vdots \\
D_n^1 X^1 & \cdots & D_n^1 X^n \\
\end{array} 
\right)^{-1}
\left(
\begin{array}{c}
D_1^2 Y_{X^{i_1}}^j \\
\vdots \\
D_n^2 Y_{X^{i_1}}^j \\
\end{array}
\right), 
\end{equation}
where, for $i' = 1, \dots, n$, the
symbol $D_{i'}^2$ denotes the {\sl 
$i'$-th second order
total differentiation operator}:
\def\theequation{1.21}\begin{equation}
D_{i'}^2
:= 
\frac{\partial}{\partial x^{i'}}
+
\sum_{j'=1}^m\,y_{i'}^{j'}\,\frac{\partial}{\partial y^{j'}}
+
\sum_{j'=1}^m\,\sum_{i_1'=1}^n\,y_{i',i_1'}^{j'}\,
\frac{\partial}{\partial y_{i_1'}^{j'}}.
\end{equation}
Again, these formulas~\thetag{ 1.20} are not explicit in the sense
that an inverse matrix is involved and that the terms $D_{i'}^1 X^i$,
$D_{i'}^2 Y_{ X^{ i_1}}^j$ are not developed. It would already be a
nontrivial computational task to develope these expressions and to
find out some nice satisfying combinatorial formulas.

In order to present the general inductive non-explicit formulas for
the computation of the $\kappa$-th prolongation $\varphi^{(\kappa)}$,
we need some more notation. Let $\lambda \in \N$ be an arbitrary
integer. For $i' = 1, \dots, n$, let $D_{i'}^\lambda$ denotes the {\sl
$i'$-th $\lambda$-th order total differentiation operators}, defined
precisely by:
\def\theequation{1.22}\begin{equation}
\left\{
\aligned
D_{i'}^\lambda
&
:=
\frac{\partial}{\partial x^{i'}}
+
\sum_{j'=1}^m\,y_{i'}^{j'}\,\frac{\partial}{\partial y^{j'}}
+
\sum_{j'=1}^m\,\sum_{i_1'=1}^n\,y_{i',i_1'}^{j'}\,
\frac{\partial}{\partial y_{i_1'}^{j'}}
+
\sum_{j'=1}^m\,\sum_{i_1',i_2'=1}^n\,y_{i',i_1',i_2'}^{j'}\,
\frac{\partial}{\partial y_{i_1',i_2'}^{j'}}
+ \\
& \
\ \ \ \ \
+
\cdots
+
\sum_{j'=1}^m\,\sum_{i_1',i_2',\dots,i_{\lambda-1}'=1}^n\,
y_{i',i_1',i_2',\dots,i_{\lambda-1}'}^{j'}\,
\frac{\partial}{\partial y_{i_1',i_2',\dots,i_{\lambda-1}'}^{j'}}.
\endaligned\right.
\end{equation} 
Then, for $i = 1, \dots, n$, the expressions of $Y_{X^{i_1}\cdots
X^{i_{\lambda-1}} X^i }^j$ are given by the following compact formulas
(again \cite{ bk1989}):
\def\theequation{1.23}\begin{equation}
\left(
\begin{array}{c}
Y_{X^{i_1}\cdots X^{i_{\lambda-1}}X^1}^j \\
\vdots \\
Y_{X^{i_1}\cdots X^{i_{\lambda-1}}X^n}^j \\
\end{array}
\right)
=
\left(
\begin{array}{ccc}
D_1^1 X^1 & \cdots & D_1^1 X^n \\
\vdots & \cdots & \vdots \\
D_n^1 X^1 & \cdots & D_n^1 X^n \\
\end{array} 
\right)^{-1}
\left(
\begin{array}{c}
D_1^\lambda Y_{X^{i_1}\cdots X^{i_{\lambda-1}}}^j \\
\vdots \\
D_n^\lambda Y_{X^{i_1}\cdots X^{i_{\lambda-1}}}^j \\
\end{array}
\right).
\end{equation}
Again, these inductive formulas are incomplete and unsatisfactory.

\def\theproblem{1.24}\begin{problem}
Find totally explicit complete formulas
for the $\kappa$-th prolongation $\varphi^{(\kappa)}$.
\end{problem}

Except in the cases $\kappa = 1, 2$, we have not been able to solve
this problem. The case $\kappa = 1$ is elementary. Complete formulas
in the particular cases $\kappa = 2$, $n=1$, $m\geqslant 1$ and
$n\geqslant 1$, $m=1$ are implicitely provided in~\cite{ me2004} and in
Section~?(?), where one observes the appearance of some modifications
of the Jacobian determinant of the diffeomorphism $\varphi$, inserted
in a clearly understandable combinatorics. In fact, there is a nice
dictionary between the formulas for $\varphi^{ (2)}$ and the formulas
for the second prolongation $\mathcal{ L}^{ (2)}$ of a vector field
$\mathcal{ L}$ which were written in equation~\thetag{ 43} of~\cite{
gm2003a} ({\it see} also equations~\thetag{ 2.6}, \thetag{ 3.20},
\thetag{ 4.6} and~\thetag{ 5.3} in the next paragraphs). In the
passage from $\varphi^{ (2)}$ to $\mathcal{ L}^{ (2)}$, a sort of
formal first order linearization may be observed and the reverse
passage may be easily guessed. However, for $\kappa \geqslant 3$, the
formulas for $\varphi^{ (\kappa)}$ explode faster than the formulas
for the $\kappa$-th prolongation $\mathcal{ L}^{ (\kappa )}$ of a
vector field $\mathcal{ L}$. Also, the dictionary between $\varphi^{
(\kappa)}$ and $\mathcal{ L}^{ ( \kappa )}$ disappears. In fact, to
elaborate an appropriate dictionary, we believe that one should
introduce before a sort of formal $(\kappa-1)$-th order linearizations
of $\varphi^{ ( \kappa)}$, finer than the first order linearization
$\mathcal{ L}^{ (\kappa)}$. To be optimistic, we believe that the
final answer to Problem~1.24 is, nevertheless, accessible after hard
work.

The present article is devoted to present totally explicit complete
formulas for the $\kappa$-th prolongation $\mathcal{ L}^{ (\kappa )}$
of a vector field $\mathcal{ L}$ to $\mathcal{ J}_{n,m}^\kappa$, for
$n\geqslant 1$ arbitrary, for $m\geqslant 1$ arbitrary and for $\kappa
\geqslant 1$ arbitrary.

\subsection*{ 1.25.~Prolongation of a vector field to the
$\kappa$-th jet space} Consider a vector field
\def\theequation{1.26}\begin{equation}
\mathcal{ L} 
= 
\sum_{i=1}^n\mathcal{X}^i(x,y)\,
\frac{\partial}{\partial x^i}
+
\sum_{j=1}^m\,\mathcal{Y}^j(x,y)\,
\frac{\partial}{\partial y^j},
\end{equation}
defined in $\K^{ n+m}$. Its flow:
\def\theequation{1.27}\begin{equation}
\varphi_t ( x, y) 
:=
\exp 
\left(
t \mathcal{ L} 
\right) 
(x, y)
\end{equation}
constitutes a one-parameter family 
of diffeomorphisms of $\K^{ n+m} $ close
to the identity. The lift $(\varphi_t )^{ ( \kappa )}$ to the
$\kappa$-th jet space constitutes a one-parameter family of
diffeomorphisms of $\mathcal{ J}_{ n,m }^\kappa$. By definition, the
{\sl $\kappa$-th prolongation $\mathcal{ L }^{( \kappa)}$ of
$\mathcal{ L }$ to the jet space $\mathcal{ J}_{n, m }^\kappa$} is the
infinitesimal generator of $(\varphi_t )^{ (\kappa)}$, namely:
\def\theequation{1.28}\begin{equation}
\mathcal{ L}^{(\kappa)}
:= 
\left.
\frac{d}{dt}
\right\vert_{t=0}
\left[
(\varphi_t)^{(\kappa)}
\right].
\end{equation}

\subsection*{1.29.~Inductive formulas for the 
$\kappa$-th prolongation $\mathcal{ L}^{(\kappa)}$} As a vector field
defined in $\K^{n+m+ m \frac{ (n+m)! }{ n! \ m!}}$, the $\kappa$-th
prolongation $\mathcal{ L}^{ (\kappa) }$ may be written under the
general form:
\def\theequation{1.30}\begin{equation}
\left\{
\aligned
\mathcal{L}^{(\kappa)}
&
=
\sum_{i=1}^n\mathcal{X}^i\,\frac{\partial}{\partial x^i}
+
\sum_{j=1}^m\,\mathcal{Y}^j\,\frac{\partial}{\partial y^j}
+ \\
& \
\ \ \ \ \
+
\sum_{j=1}^m\,\sum_{i_1=1}^n\,{\bf Y}_{i_1}^j\,
\frac{\partial}{\partial y_{i_1}^j}
+
\sum_{j=1}^m\,\sum_{i_1,i_2=1}^n\,{\bf Y}_{i_1,i_2}^j\,
\frac{\partial}{\partial y_{i_1,i_2}^j}
+ 
\cdots
+ \\
& \
\ \ \ \ \ 
+
\sum_{j=1}^m\,\sum_{i_1,\dots,i_\kappa=1}^n\,
{\bf Y}_{i_1,\dots,i_\kappa}^j\,
\frac{\partial}{\partial y_{i_1,\dots,i_\kappa}^j}.
\endaligned\right.
\end{equation}
Here, the coefficients ${\bf Y}_{i_1}^j$, ${\bf Y}_{i_1, i_2}^j$,
$\dots$, ${\bf Y}_{i_1, i_2, \dots, i_\kappa}^j$ are uniquely
determined in terms of partial derivatives of the coefficients
$\mathcal{ X}^i$ and $\mathcal{ Y}^j$ of the original vector field
$\mathcal{ L}$, together with the pure jet variables $\left( y_{ i_1
}^j,\dots, y_{i_1, \dots, i_\kappa}^j \right)$, by means of the
following {\sl fundamental inductive formulas} (\cite{ ol1979}, 
\cite{ ol1986}, \cite{
bk1989}):
\def\theequation{1.31}\begin{equation}
\left\{
\aligned
{\bf Y}_{i_1}^j
&
:=
D_{i_1}^1
\left(
\mathcal{ Y}^j
\right)
-
\sum_{k=1}^n\,D_{i_1}^1
\left(
\mathcal{X}^k
\right)
\, y_k^j, 
\\
{\bf Y}_{i_1,i_2}^j
&
:=
D_{i_2}^2
\left(
{\bf Y}_{i_1}^j
\right)
-
\sum_{k=1}^n\,D_{i_2}^1
\left(
\mathcal{X}^k
\right)
\, y_{i_1,k}^j, 
\\
\cdots \cdots \cdots 
&
\cdots \cdots \cdots \cdots \cdots \cdots \cdots \cdots
\cdots \cdots \cdots \cdots \cdots \cdots \cdots \cdots
\\
{\bf Y}_{i_1,i_2,\dots,i_\kappa}^j
&
:=
D_{i_\kappa}^\kappa
\left(
{\bf Y}_{i_1,i_2,\dots,i_{\kappa-1}}^j
\right)
-
\sum_{k=1}^n\,D_{i_\kappa}^1
\left(
\mathcal{X}^k
\right)
\, y_{i_1,i_2,\dots,i_{\kappa-1},k}^j,
\endaligned\right.
\end{equation}
where, for every $\lambda \in \N$ with $0 \leqslant \lambda \leqslant
\kappa$, and for every $i\in \N$ with $1 \leqslant i' \leqslant n$,
the $i'$-th $\lambda$-th order total differentiation operator $D_{i'
}^\lambda$ was defined in~\thetag{ 1.22} above.

\def\theproblem{1.32}\begin{problem}
Applying these inductive formulas, find totally explicit complete
formulas for the $\kappa$-th prolongation $\mathcal{L}^{(\kappa)}$.
\end{problem}

The present article is devoted to provide all the desired formulas.

\subsection*{ 1.33.~Methodology of induction}
We have the intention of presenting our results in a purely inductive
style, based on several thorough visual comparisons between massive formulas
which will be written and commented in four different cases:

\smallskip

\begin{itemize}
\item[{\bf (i)}]
$n = 1$ and $m = 1$; $\kappa \geqslant 1$ arbitrary;
\item[{\bf (ii)}]
$n \geqslant 1$ and $m=1$; $\kappa \geqslant 1$ arbitrary;
\item[{\bf (iii)}]
$n=1$ and $m \geqslant 1$; $\kappa \geqslant 1$ arbitrary;
\item[{\bf (iv)}]
general case: $n \geqslant 1$ and $m \geqslant 1$; $\kappa
\geqslant 1$ arbitrary.
\end{itemize}

\smallskip

Accordingly, we shall particularize and slightly lighten our notations
in each of the three (preliminary) cases (i)
[Section~2], (ii) [Section~3] and (iii) [Section~4].

\section*{\S2.~One independent variable and one dependent variable}

\subsection*{2.1.~Simplified adapted notations}
Assume $n=1$ and $m=1$, let $\kappa\in \N$ with $\kappa \geqslant 1$
and simply denote the jet variables by:
\def\theequation{2.2}\begin{equation}
\left(
x,y,y_1,y_2,\dots,y_\kappa
\right) \in \mathcal{J}_{1,1}^\kappa.
\end{equation}
The $\kappa$-th prolongation of a vector field $\mathcal{ L} =
\mathcal{ X} \,\frac{ \partial }{ \partial x} + \mathcal{ Y}\,\frac{
\partial}{ \partial y}$ will be denoted by:
\def\theequation{2.3}\begin{equation}
\mathcal{L}^{(\kappa)}
=
\mathcal{X}\,\frac{\partial}{\partial x}
+
\mathcal{Y}\,\frac{\partial}{\partial y}
+
{\bf Y}_1\,\frac{\partial}{\partial y_1}
+
{\bf Y}_2\,\frac{\partial}{\partial y_2}
+ 
\cdots
+
{\bf Y}_\kappa\,\frac{\partial}{\partial y_\kappa}.
\end{equation}
The coefficients 
${\bf Y}_1$, ${\bf Y}_2$, $\dots$, ${\bf Y}_\kappa$
are computed by means of the inductive formulas:
\def\theequation{2.4}\begin{equation}
\left\{
\aligned
{\bf Y}_1
&
:= 
D^1(\mathcal{Y})
-
D^1(\mathcal{X})\,y_1, \\
{\bf Y}_2
&
:= 
D^2({\bf Y}_1)
-
D^1(\mathcal{X})\,y_2, \\
\cdots 
&
\cdots \cdots \cdots \cdots \cdots \cdots \cdots \cdots
\\
{\bf Y}_\kappa
&
:= 
D^\kappa({\bf Y}_{\kappa-1})
-
D^1(\mathcal{X})\,y_\kappa, \\
\endaligned\right.
\end{equation}
where, for $1 \leqslant \lambda \leqslant \kappa$:
\def\theequation{2.5}\begin{equation}
D^\lambda
:= 
\frac{\partial}{\partial x}
+
y_1\,\frac{\partial}{\partial y}
+ 
y_2\,\frac{\partial}{\partial y_1}
+
\cdots
+
y_\lambda\,\frac{\partial}{\partial y_{\lambda-1}}.
\end{equation}
By direct elementary computations, for $\kappa = 1$ and for $\kappa =
2$, we obtain the following two very classical formulas :
\def\theequation{2.6}\begin{equation}
\left\{
\aligned
{\bf Y}_1
&
=
\mathcal{Y}_x
+
\left[
\mathcal{Y}_y
-
\mathcal{X}_x
\right]
y_1
+
\left[
-\mathcal{X}_y
\right]
(y_1)^2, 
\\
{\bf Y}_2
&
=
\mathcal{Y}_{x^2}
+
\left[
2\,\mathcal{Y}_{xy}
-
\mathcal{X}_{x^2}
\right]
\,y_1
+
\left[
\mathcal{Y}_{y^2}
-
2\,\mathcal{X}_{xy}
\right](y_1)^2
+
\left[
-
\mathcal{X}_{y^2}
\right]
(y_1)^3
+ \\
& \
\ \ \ \ \
+
\left[
\mathcal{Y}_y
-
2\,\mathcal{X}_x
\right]\,y_2
+
\left[
-
3\,\mathcal{X}_y
\right]\,
y_1\,y_2.
\endaligned\right.
\end{equation}
Our main objective is to {\it devise the general combinatorics}.
Thus, to attain this aim, we have to achieve patiently formal
computations of the next coefficients ${\bf Y}_3$, ${\bf Y}_4$ and
${\bf Y}_5$. We systematically use parentheses $\left[ \cdot \right]$
to single out every coefficient of the polynomials ${\bf Y}_3$, ${\bf
Y}_4$ and ${\bf Y}_5$ in the pure jet variables $y_1, y_2, y_3, y_4$
and $y_5$, putting every sign inside these parentheses. We always put
the monomials in the pure jet variables $y_1, y_2, y_3, y_4$ and $y_5$
after the parentheses. For completeness, let us provide the
intermediate computation of the third coefficient ${\bf Y}_3$. In
detail:
$$
\aligned
{\bf Y}_3
&
=
D^3
\left(
{\bf Y}_2
\right)
-
D^1
\left(
\mathcal{ X}
\right)
y_3 
\\
&
=
\left(
\frac{\partial}{\partial x}
+
y_1\,\frac{\partial}{\partial y}
+
y_2\,\frac{\partial}{\partial y_1}
+
y_3\,\frac{\partial}{\partial y_2}
\right)
\Big(
\mathcal{Y}_{x^2}
+
\left[
2\,\mathcal{Y}_{xy}
-
\mathcal{X}_{x^2}
\right]
y_1
+ 
\\
& \
\ \ \ \ \
+
\left[
\mathcal{ Y}_{y^2}
-
2\,\mathcal{X}_{xy}
\right]
(y_1)^2
+
\left[
-
\mathcal{X}_{y^2}
\right]
(y_1)^3
+
\\
&
+
\left[
\mathcal{ Y}_y
-
2\,\mathcal{X}_x
\right]
y_2
+
\left[
-
3\,\mathcal{X}_y
\right]
y_1\,y_2
\Big)
\\
\endaligned
$$
\def\theequation{2.7}\begin{equation}
\aligned
&
=
\underline{
\mathcal{Y}_{x^3} 
}_{ \fbox{\tiny 1}}
+
\underline{
\left[
2\,\mathcal{Y}_{x^2y}
-
\mathcal{X}_{x^3}
\right]
y_1
}_{ \fbox{\tiny 2}}
+
\underline{
\left[
\mathcal{Y}_{xy^2}
-
2\,\mathcal{X}_{x^2y}
\right]
(y_1)^2
}_{ \fbox{\tiny 3}}
+
\underline{
\left[
-
\mathcal{X}_{xy^2}
\right]
(y_1)^3
}_{ \fbox{\tiny 4}}
+ \\
& \
\ \ \ \ \
+
\underline{ 
\left[
\mathcal{Y}_{xy}
-
2\,\mathcal{X}_{x^2}
\right]
y_2
}_{ \fbox{\tiny 6}}
+
\underline{ 
\left[
-
3\,\mathcal{X}_{xy}
\right]
y_1y_2
}_{ \fbox{\tiny 7}}
+
\underline{ 
\left[
\mathcal{Y}_{x^2y}
\right]
y_1
}_{ \fbox{\tiny 2}}
+ \\
& \
\ \ \ \ \
+
\underline{ 
\left[
2\,\mathcal{Y}_{xy^2}
-
\mathcal{X}_{x^2y}
\right]
(y_1)^2
}_{ \fbox{\tiny 3}}
+
\underline{ 
\left[
\mathcal{Y}_{y^3}
-
2\,\mathcal{X}_{xy^2}
\right]
(y_1)^3
}_{ \fbox{\tiny 4}}
+
\underline{ 
\left[
-
\mathcal{X}_{y^3}
\right]
(y_1)^4
}_{ \fbox{\tiny 5}}
+
\endaligned
\end{equation}
$$
\aligned
& \
\ \ \ \ \
+
\underline{ 
\left[
\mathcal{Y}_{y^2}
-
2\,\mathcal{X}_{xy}
\right]
y_1y_2
}_{ \fbox{\tiny 7}}
+
\underline{ 
\left[
-
3\,\mathcal{X}_{y^2}
\right]
(y_1)^2y_2
}_{ \fbox{\tiny 8}}
+
\underline{ 
\left[
2\,\mathcal{Y}_{xy}
-
\mathcal{X}_{x^2}
\right]
y_2
}_{ \fbox{\tiny 6}}
+ \\
& \
\ \ \ \ \
+
\underline{ 
\left[
\mathcal{Y}_{y^2}
-
2\,\mathcal{X}_{xy}
\right]
2\,y_1y_2
}_{ \fbox{\tiny 7}}
+
\underline{ 
\left[
-
\mathcal{X}_{y^2}
\right]
3(y_1)^2y_2
}_{ \fbox{\tiny 8}}
+
\underline{ 
\left[
-
3\,\mathcal{X}_y
\right]
(y_2)^2
}_{ \fbox{\tiny 9}}
+ \\
& \
\ \ \ \ \ 
+
\underline{ 
\left[
\mathcal{Y}_y
-
2\,\mathcal{X}_x
\right]
y_3
}_{ \fbox{\tiny 10}}
+
\underline{ 
\left[
-
3\,\mathcal{X}_y
\right]
y_1y_3
}_{ \fbox{\tiny 11}}
- \\
& \
\ \ \ \ \
-
\underline{ 
\left[
\mathcal{X}_x
\right]
y_3
}_{ \fbox{\tiny 10}}
-
\underline{ 
\left[
\mathcal{X}_y
\right]
y_1y_3
}_{ \fbox{\tiny 11}} \ .
\endaligned
$$
We have underlined all the terms with 
a number appended. Each number
refers to the order of appearance
of the terms in the final simplified
expression of ${\bf Y}_3$, also written in~\cite{ bk1989}
with different notations:
\def\theequation{2.8}\begin{equation}
\left\{
\aligned
{\bf Y}_3
& 
=
\mathcal{Y}_{x^3}
+
\left[
3\,\mathcal{Y}_{x^2y}
-
\mathcal{X}_{x^3}
\right]
y_1
+
\left[
3\,\mathcal{Y}_{xy^2}
-
3\,\mathcal{X}_{x^2y}
\right]
(y_1)^2
+ \\
& \
\ \ \ \ \
+
\left[
\mathcal{Y}_{y^3}
-
3\,\mathcal{X}_{xy^2}
\right]
(y_1)^3
+
\left[
-
\mathcal{X}_{y^3}
\right]
(y_1)^4
+
\left[
3\,\mathcal{Y}_{xy}
-
3\,\mathcal{X}_{x^2}
\right]
y_2
+ \\
& \
\ \ \ \ \
+
\left[
3\,\mathcal{Y}_{y^2}
-
9\,\mathcal{X}_{xy}
\right]
y_1y_2
+
\left[
-
6\,\mathcal{X}_{y^2}
\right]
(y_1)^2y_2
+
\left[
-
3\,\mathcal{X}_y
\right]
(y_2)^2
+ \\
& \
\ \ \ \ \
+
\left[
\mathcal{Y}_y
-
3\,\mathcal{X}_x
\right]
y_3
+
\left[
-
4\,\mathcal{X}_y
\right]
y_1y_3.
\endaligned\right.
\end{equation}
After similar manual computations, the intermediate details of which
we will not copy in this Latex file, we get the desired expressions of
${\bf Y}_4$ and of ${\bf Y}_5$.
Firstly:
\def\theequation{2.9}\begin{equation}
\small
\left\{
\aligned
{\bf Y}_4
&
=
\mathcal{Y}_{x^4}
+
\left[
4\,\mathcal{Y}_{x^3y}
-
\mathcal{X}_{x^4}
\right]
y_1
+
\left[
6\,\mathcal{Y}_{x^2y^2}
-
4\,\mathcal{X}_{x^3y}
\right]
(y_1)^2
+ \\
& \
\ \ \ \ \
+
\left[
4\,\mathcal{Y}_{xy^3}
-
6\,\mathcal{X}_{x^2y^2}
\right]
(y_1)^3
+
\left[
\mathcal{Y}_{y^4}
-
4\,\mathcal{X}_{xy^3}
\right]
(y_1)^4
+
\left[
-
\mathcal{X}_{y^4}
\right]
(y_1)^5
+ \\
& \
\ \ \ \ \
+
\left[
6\,\mathcal{Y}_{x^2y}
-
4\,\mathcal{X}_{x^3}
\right]
y_2
+
\left[
12\,\mathcal{Y}_{xy^2}
-
18\,\mathcal{X}_{x^2y}
\right]
y_1y_2
+ \\
& \
\ \ \ \ \
+
\left[
6\,\mathcal{Y}_{y^3}
-
24\,\mathcal{X}_{xy^2}
\right]
(y_1)^2y_2
+
\left[
-
10\,\mathcal{X}_{y^3}
\right]
(y_1)^3y_2
+ \\
& \
\ \ \ \ \
+
\left[
3\,\mathcal{Y}_{y^2}
-
12\,\mathcal{X}_{xy}
\right]
(y_2)^2
+
\left[
-
15\,\mathcal{X}_{y^2}
\right]
y_1(y_2)^2
+ \\
& \
\ \ \ \ \
+
\left[
4\,\mathcal{Y}_{xy}
-
6\,\mathcal{X}_{x^2}
\right]
y_3
+
\left[
4\,\mathcal{Y}_{y^2}
-
16\,\mathcal{X}_{xy}
\right]
y_1y_3
+
\left[
-
10\,\mathcal{X}_{y^2}
\right]
(y_1)^2y_3
+ \\
& \
\ \ \ \ \
+
\left[
-
10\,\mathcal{X}_y
\right]
y_2y_3
+
\left[
\mathcal{Y}_y
-
4\,\mathcal{X}_x
\right]
y_4
+
\left[
-
5\,\mathcal{X}_y
\right]
y_1y_4.
\endaligned\right.
\end{equation}
Secondly:
\def\theequation{2.10}\begin{equation}
\small
\left\{
\aligned
{\bf Y}_5
&
=
\mathcal{Y}_{x^5}
+
\left[
5\,\mathcal{Y}_{x^4y}
-
\mathcal{X}_{x^5}
\right]
y_1
+
\left[
10\,\mathcal{Y}_{x^3y^2}
-
5\,\mathcal{X}_{x^4y}
\right]
(y_1)^2
+ \\
& \
\ \ \ \ \
+
\left[
10\,\mathcal{Y}_{x^2y^3}
-
10\,\mathcal{X}_{x^3y^2}
\right]
(y_1)^3
+
\left[
5\,\mathcal{Y}_{xy^4}
-
10\,\mathcal{X}_{x^2y^3}
\right]
(y_1)^4
+ \\
& \
\ \ \ \ \
+
\left[
\mathcal{Y}_{y^5}
-
5\,\mathcal{X}_{xy^4}
\right]
(y_1)^5
+
\left[
-
\mathcal{X}_{y^5}
\right]
(y_1)^6
+
\left[
10\,\mathcal{Y}_{x^3y}
-
5\,\mathcal{X}_{x^4}
\right]
y_2
+ \\
& \
\ \ \ \ \
+
\left[
30\,\mathcal{Y}_{x^2y^2}
-
30\,\mathcal{X}_{x^3y}
\right]
y_1y_2
+
\left[
30\,\mathcal{Y}_{xy^3}
-
60\,\mathcal{X}_{x^2y^2}
\right]
(y_1)^2y_2
+ \\
& \
\ \ \ \ \
+
\left[
10\,\mathcal{Y}_{y^4}
-
50\,\mathcal{X}_{xy^3}
\right]
(y_1)^3y_2
+
\left[
-
15\,\mathcal{X}_{y^4}
\right]
(y_1)^4y_2
+ \\
& \
\ \ \ \ \
+
\left[
15\,\mathcal{Y}_{xy^2}
-
30\,\mathcal{X}_{x^2y}
\right]
(y_2)^2
+
\left[
15\,\mathcal{Y}_{y^3}
-
75\,\mathcal{X}_{xy^2}
\right]
y_1(y_2)^2
+ \\
& \
\ \ \ \ \
+
\left[
-
45\,\mathcal{X}_{y^3}
\right]
(y_1)^2(y_2)^2
+
\left[
-
15\,\mathcal{X}_{y^2}
\right]
(y_2)^3
+ \\
& \
\ \ \ \ \
+
\left[
10\,\mathcal{Y}_{x^2y}
-
10\,\mathcal{X}_{x^3}
\right]
y_3
+
\left[
20\,\mathcal{Y}_{xy^2}
-
40\,\mathcal{X}_{x^2y}
\right]
y_1y_3
+ \\
& \
\ \ \ \ \
+
\left[
10\,\mathcal{Y}_{y^3}
-
50\,\mathcal{X}_{xy^2}
\right]
(y_1)^2y_3
+
\left[
-
20\,\mathcal{X}_{y^3}
\right]
(y_1)^3y_3
+ \\
& \
\ \ \ \ \
+
\left[
10\,\mathcal{Y}_{y^2}
-
50\,\mathcal{X}_{xy}
\right]
y_2y_3
+
\left[
-
60\,\mathcal{X}_{y^2}
\right]
y_1y_2y_3
+
\left[
-
10\,\mathcal{X}_y
\right]
(y_3)^2
+ \\
& \
\ \ \ \ \
+
\left[
5\,\mathcal{Y}_{xy}
-
10\,\mathcal{X}_{x^2}
\right]
y_4
+
\left[
5\,\mathcal{Y}_{y^2}
-
25\,\mathcal{X}_{xy}
\right]
y_1y_4
+
\left[
-
15\,\mathcal{X}_{y^2}
\right]
(y_1)^2y_4
+ \\
& \
\ \ \ \ \
+
\left[
-
15\,\mathcal{X}_y
\right]
y_2y_4
+
\left[
\mathcal{Y}_y
-
5\,\mathcal{X}_y
\right]
y_5
+
\left[
-
6\,\mathcal{X}_y
\right]
y_1y_5.
\endaligned\right.
\end{equation}

\subsection*{2.11.~Formal inspection, formal intuition
and formal induction} Now, we have to comment these formulas. We have
written in length the five polynomials ${\bf Y}_1$, ${\bf Y}_2$, ${\bf
Y}_3$, ${\bf Y}_4$ and ${\bf Y}_5$ in the pure jet variables $y_1,
y_2, y_3, y_4$ and $y_5$. Except the first ``constant'' term
$\mathcal{ Y}_{x^\kappa}$, all the monomials in the expression of
${\bf Y}_\kappa$ are of the general form
\def\theequation{2.12}\begin{equation}
\left(
y_{\lambda_1}
\right)^{\mu_1}
\left(
y_{\lambda_2}
\right)^{\mu_2}
\cdots
\left(
y_{\lambda_d}
\right)^{\mu_d},
\end{equation}
for some positive integer $d\geqslant 1$, for some collection of
strictly increasing jet indices:
\def\theequation{2.13}\begin{equation}
1 \leqslant \lambda_1 < \lambda_2 < \cdots < 
\lambda_d \leqslant \kappa,
\end{equation}
and for some positive integers $\mu_1, \dots, \mu_d \geqslant 1$. This
and the next combinatorial facts may be confirmed by reading the
formulas giving ${\bf Y}_1$, ${\bf Y}_2$, ${\bf Y}_3$, ${\bf Y}_4$ and
${\bf Y}_5$. It follows that the integer $d$ satisfies the inequality
$d\leqslant \kappa+1$. To include the first ``constant'' term
$\mathcal{ Y}_{x^\kappa}$, we shall make the convention that putting
$d=0$ in the monomial~\thetag{ 2.12} yields the constant term $1$.
 
Furthermore, by inspecting the formulas giving ${\bf Y}_1$, ${\bf
Y}_2$, ${\bf Y}_3$, ${\bf Y}_4$ and ${\bf Y}_5$, we see that the
following inequality should be satisfied:
\def\theequation{2.14}\begin{equation}
\mu_1\lambda_1
+
\mu_2\lambda_2
+
\cdots
+
\mu_d\lambda_d 
\leqslant 
\kappa+1.
\end{equation}
For instance, in the expression of ${\bf Y}_4$, the two monomials
$(y_1)^3 y_2$ and $y_1 (y_2)^2$ do appear, but the two monomials
$(y_1)^4 y_2$ and $(y_1)^2 (y_2)^2$ cannot appear. All coefficients
of the pure jet monomials are of the general form:
\def\theequation{2.15}\begin{equation}
\left[
A\,
\mathcal{Y}_{x^\alpha y^{\beta+1}}
-
B\,
\mathcal{X}_{x^{\alpha+1}y^\beta}
\right],
\end{equation}
for some nonnegative integers $A, B, \alpha, \beta \in \N$. Sometimes
$A$ is zero, but $B$ is zero only for the (constant, with respect to
pure jet variables) term $\mathcal{ Y}_{x^\kappa}$. Importantly,
$\mathcal{ X}$ is differentiated once more with respect to $x$ and
$\mathcal{ Y}$ is differentiated once more with respect to $y$. Again,
this may be confirmed by reading all the terms in the formulas for
${\bf Y}_1$, ${\bf Y}_2$, ${\bf Y}_3$, ${\bf Y}_4$ and ${\bf Y}_5$.

In addition, we claim that there is a link between the couple
$(\alpha, \beta)$ and the collection $\{ \mu_1, \lambda_1, \dots,
\mu_d, \lambda_d \}$. To discover it, let us write some of the
monomials appearing in the expressions of ${\bf Y}_4$ (first column)
and of ${\bf Y}_5$ (second column), for instance:
\def\theequation{2.16}\begin{equation}
\left\{
\aligned
&
\left[6\,
\mathcal{Y}_{x^2y^2}
-
4\,\mathcal{X}_{x^3y}
\right]
(y_1)^2, 
\ \ \ \ \ \ \ \ \ \ \ \ \ \ \
&
\left[
5\,\mathcal{Y}_{xy^4}
-
10\,\mathcal{X}_{x^2y^3}
\right]
(y_1)^4,
&
\\
&
\left[
12\,\mathcal{Y}_{xy^2}
-
18\,\mathcal{X}_{x^2y}
\right]
y_1y_2,
\ \ \ \ \ \ \ \ \ \ \ \ \ \ \
&
\left[
30\,\mathcal{Y}_{xy^3}
-
60\,\mathcal{X}_{x^2y^2}
\right]
(y_1)^2y_2,
& 
\\
& 
\left[
-
10\,\mathcal{X}_{y^3}
\right]
(y_1)^3y_2, 
\ \ \ \ \ \ \ \ \ \ \ \ \ \ \
& 
\left[
-
15\,\mathcal{X}_{y^4}
\right]
(y_1)^4y_2,
&
\\
&
\left[
4\,\mathcal{Y}_{y^2}
-
16\,\mathcal{X}_{xy}
\right]
y_1y_3, 
\ \ \ \ \ \ \ \ \ \ \ \ \ \ \
& 
\left[
10\,\mathcal{Y}_{y^2}
-
50\,\mathcal{X}_{xy}
\right]
y_2y_3,
&
\\
&
\left[
-
10\,\mathcal{X}_{y^2}
\right]
(y_1)^2 y_3,
\ \ \ \ \ \ \ \ \ \ \ \ \ \ \
& 
\left[
-
60\,\mathcal{X}_{y^2}
\right]
y_1y_2y_3.
& \\
\endaligned\right.
\end{equation}
After some reflection, we discover the hidden intuitive rule: the
partial derivatives of $\mathcal{ Y}$ and of $\mathcal{ X}$ associated
with the monomial $(y_{\lambda_1})^{\mu_1} \cdots
(y_{\lambda_d})^{\mu_d}$ are, respectively:
\def\theequation{2.17}\begin{equation}
\left\{
\aligned
&
\mathcal{ Y}_{
x^{\kappa-\mu_1\lambda_1-\cdots-\mu_d\lambda_d}
\,
y^{\mu_1+\cdots+\mu_d}}, 
\\
&
\mathcal{ X}_{
x^{\kappa-\mu_1\lambda_1-\cdots-\mu_d\lambda_d+1}
\,
y^{\mu_1+\cdots+\mu_d-1}}. 
\\
\endaligned\right.
\end{equation}
This may be checked on each of the $10$ examples~\thetag{ 2.16}
above.

Now that we have explored and discovered the combinatorics of the pure
jet monomials, of the partial derivatives and of the complete sum
giving ${\bf Y}_\kappa$, we may express that it is of the following
general form:
\def\theequation{2.18}\begin{equation}
\left\{
\aligned
{\bf Y}_\kappa
&
=
\mathcal{ Y}_{x^\kappa}
+
\sum_{d=1}^{\kappa+1}
\ \
\sum_{1\leqslant\lambda_1<\cdots<\lambda_d\leqslant\kappa}
\ \
\sum_{\mu_1\geqslant 1,\dots,\mu_d\geqslant 1}
\
\sum_{
\mu_1\lambda_1
+
\cdots
+
\mu_d\lambda_d\leqslant \kappa+1} 
\\
& \
\ \ \ \ \
\left[
A_\kappa^{
(\mu_1, \lambda_1), \dots, (\mu_d, \lambda_d) }
\cdot
\mathcal{Y}_{
x^{\kappa-\mu_1\lambda_1-\cdots-\mu_d\lambda_d}
\,
y^{\mu_1+\cdots+\mu_d}
}
-
\right. \\
& \
\ \ \ \ \ \ \ \ \ \
\left.
-
B_\kappa^{
(\mu_1, \lambda_1), \dots, (\mu_d, \lambda_d) }
\cdot
\mathcal{X}_{
x^{\kappa-\mu_1\lambda_1-\cdots-\mu_d\lambda_d+1}
\,
y^{\mu_1+\cdots+\mu_d-1}
}
\right]
\cdot
\\
& \
\ \ \ \ \ \ \ \ \ \ \ \ \ \ \ \ \ \ \ \
\ \ \ \ \ \ \ \ \ \ \ \ \ \ \ \ \ \ \ \
\ \ \ \ \ \ \ \ \ \ \ \ \
\cdot
(y_{\lambda_1})^{\mu_1}
\cdots
(y_{\lambda_d})^{\mu_d}.
\endaligned\right.
\end{equation} 
Here, we separate the first term $\mathcal{ Y}_{x^\kappa}$ from the
general sum; it is the constant term in ${\bf Y}_\kappa$, which 
itself is a
polynomial with respect to the jet variables $y_\lambda$. In this
general formula, the only remaining unknowns are the nonnegative
integer coefficients $A_\kappa^{ (\mu_1, \lambda_1), \dots, (\mu_d,
\lambda_d) } \in \N$ and $B_\kappa^{ (\mu_1, \lambda_1), \dots,
(\mu_d, \lambda_d) } \in \N$. In Section~3 below, we shall explain how
we have discovered their exact value.

At present, even if we are unable to devise their explicit
expression, we may observe that the value of the special integer
coefficients $A^{(\mu_1, 1)}_{ \mu_1}$ and $B^{( \mu_1, 1)}_{ \mu_1}$
which are attached to the monomials ${\rm ct.}$, $y_1$, $(y_1)^2$,
$(y_1)^3$, $(y_1)^4$ and $(y_1)^5$ are simple. Indeed, by
inspecting the first terms in the expressions of ${\bf Y}_1$, ${\bf
Y}_2$, ${\bf Y}_3$, ${\bf Y}_4$ and ${\bf Y}_5$, we of course
recognize the binomial coefficients. In general:

\def\thelemma{2.19}\begin{lemma} For $\kappa \geqslant 1$, 
\def\theequation{2.20}\begin{equation}
\left\{
\aligned
{\bf Y}_\kappa
&
=
\mathcal{Y}_{x^\kappa}
+
\sum_{\lambda=1}^\kappa
\left[
\binom{\kappa}{\lambda}
\,\mathcal{Y}_{x^{\kappa-\lambda}y^\lambda}
-
\binom{\kappa}{\lambda-1}
\,\mathcal{X}_{x^{\kappa-\lambda+1}y^{\lambda-1}}
\right]
(y_1)^\lambda
+ \\
& \
\ \ \ \ \
+
\left[
-
\mathcal{X}_{y^\kappa}
\right]
(y_1)^\kappa
+
{\sf remainder},
\endaligned\right.
\end{equation}
where the term {\sf remainder} collects all remaining monomials in
the pure jet variables.
\end{lemma}

In addition, let us remind what we have observed and used in a
previous co-signed work.

\def\thelemma{2.21}\begin{lemma}
\text{\rm (\cite{ gm2003a}, p.~536)} 
For $\kappa \geqslant 4$, nine among
the monomials of ${\bf Y}_\kappa$ are of the following general
form{\rm :}
\def\theequation{2.22}\begin{equation}
\left\{
\aligned
{\bf Y}_\kappa
&
=
\mathcal{Y}_{x^\kappa}
+
\left[
C_\kappa^1 \, \mathcal{Y}_{x^{\kappa-1}y}
-
\mathcal{X}_{x^\kappa}
\right]
y_1
+
\left[
C_\kappa^2\,
\mathcal{Y}_{x^{\kappa-2}y}
-
C_\kappa^1 \,\mathcal{X}_{x^{\kappa-1}}
\right]
y_2
+ \\
& \
\ \ \ \ \
+
\left[
C_\kappa^2\,
\mathcal{Y}_{x^2y}
-
C_\kappa^3\,
\mathcal{X}_{x^3}
\right]
y_{\kappa-2}
+
\left[
C_\kappa^1 \,\mathcal{Y}_{xy}
-
C_\kappa^2\,
\mathcal{X}_{x^2}
\right]
y_{\kappa-1}
+ \\
& \
\ \ \ \ \
+
\left[
C_\kappa^1 \, \mathcal{Y}_{y^2}
-
\kappa^2 \, \mathcal{ X}_{xy}
\right]
y_1 y_{\kappa-1}
+
\left[
-C_\kappa^2 \, \mathcal{ X}_y
\right]
y_2y_{\kappa-1}
+ \\
& \
\ \ \ \ \
+
\left[
\mathcal{Y}_y
- 
C_\kappa^1\,\mathcal{X}_x
\right]
+
\left[
-C_{\kappa+1}^1 \, \mathcal{ X}_y
\right]
y_1 y_\kappa
+
{\sf remainder},
\endaligned\right.
\end{equation}
where the term {\sf remainder} denotes all the remaining monomials,
and where $C_\kappa^\lambda := \frac{ \kappa!}{(\kappa - \lambda)! \
\lambda !}$ is a notation for the binomial coefficient which occupies
less space in Latex ``equation mode'' than the classical notation
\def\theequation{2.23}\begin{equation}
\binom{\kappa}{\lambda}.
\end{equation}
\end{lemma}

Now, we state directly the final theorem, without further inductive or
intuitive information.

\def\thetheorem{2.24}\begin{theorem} 
For $\kappa \geqslant 1$, we have{\rm :}
\def\theequation{2.25}\begin{equation}
\boxed{
\aligned
{\bf Y}_\kappa
&
=
\mathcal{ Y}_{x^\kappa}
+
\sum_{d=1}^{\kappa+1}
\ \
\sum_{1\leqslant\lambda_1<\cdots<\lambda_d\leqslant\kappa}
\ \
\sum_{\mu_1\geqslant 1,\dots,\mu_d\geqslant 1} 
\
\sum_{
\mu_1\lambda_1
+
\cdots
+
\mu_d\lambda_d\leqslant \kappa+1} 
\\
& \
\ \ 
\left[
\frac{\kappa\cdots(\kappa-\mu_1\lambda_1-\cdots-\mu_d\lambda_d+1)}
{(\lambda_1!)^{\mu_1}\,\mu_1!
\cdots
(\lambda_d!)^{\mu_d}\,\mu_d!
}
\cdot
\mathcal{Y}_{
x^{\kappa-\mu_1\lambda_1-\cdots-\mu_d\lambda_d}
\,
y^{\mu_1+\cdots+\mu_d}
}
-
\right. \\
& \
\ \ \ \ \ \ \ \ \ \
\left.
-
\frac{\kappa\cdots(\kappa-\mu_1\lambda_1-\cdots-\mu_d\lambda_d+2)
(\mu_1\lambda_1+\cdots+\mu_d\lambda_d)}
{(\lambda_1!)^{\mu_1}\,\mu_1!
\cdots
(\lambda_d!)^{\mu_d}\,\mu_d!
}
\cdot
\right.
\\
& \
\ \ \ \ \ \ \ \ \ \ \ \ \ \ \ \
\cdot
\mathcal{X}_{
x^{\kappa-\mu_1\lambda_1-\cdots-\mu_d\lambda_d+1}
\,
y^{\mu_1+\cdots+\mu_d-1}
}
\Big]
(y_{\lambda_1})^{\mu_1}
\cdots
(y_{\lambda_d})^{\mu_d}.
\endaligned
}
\end{equation}
\end{theorem}

Once the correct theorem is formulated, its proof follows by
accessible induction arguments which will not be developed here. It is
better to continue through and to examine thorougly the case of
several variables, since it will help us considerably to explain how
we discovered the exact values of the integer coefficients $A_\kappa^{
(\mu_1, \lambda_1), \dots, (\mu_d, \lambda_d) }$ and $B_\kappa^{
(\mu_1, \lambda_1), \dots, (\mu_d, \lambda_d) }$.

\subsection*{2.26.~Verification and application} 
Before proceeding further, let us rapidly verify that the above
general formula~\thetag{ 2.25} is correct by inspecting two instances
extracted from ${\bf Y}_5$.

Firstly, the coefficient of $(y_1)^3 y_3$ in ${\bf Y}_5$ is obtained
by putting $\kappa = 5$, $d = 2$, $\lambda_1 = 1$, $\mu_1 = 3$,
$\lambda_2 = 3$ and $\mu_2 = 1$ in the general formula~\thetag{ 2.25},
which yields:
\def\theequation{2.27}\begin{equation} 
\left[ 
0
-
\frac{5\cdot 4\cdot 3\cdot 2\cdot 1\cdot 6}{ 
(1!)^3\ 3!\ (3!)^1\ 1!}
\,\mathcal{X}_{y^3} 
\right] = 
\left[ 
-
20\,\mathcal{X}_{y^3} 
\right].
\end{equation}
This value is the same as in the original formula~\thetag{ 2.10}:
confirmation.

Secondly, the coefficient of $y_1 (y_2)^2$ in ${\bf Y}_5$ is obtained
by $\kappa = 5$, $d = 2$, $\lambda_1 = 1$, $\mu_1 = 1$, $\lambda_2 =
2$ and $\mu_2 = 2$ in the general formula~\thetag{ 2.25}, which yields:
\def\theequation{2.28}\begin{equation}
\left[
\frac{5\cdot 4\cdot 3\cdot 2\cdot 1}{
(1!)^1\ 1!\ (2!)^2\ 2!}
\,\mathcal{Y}_{y^3}
-
\frac{5\cdot 4\cdot 3\cdot 2\cdot 5}{
(1!)^1\ 1!\ (2!)^2\ 2!}
\,\mathcal{X}_{xy^2}
\right]
= 
\left[
15\,\mathcal{Y}_{y^3}
-
75\,\mathcal{X}_{xy^2}
\right].
\end{equation}
This value is the same as in the original formula~\thetag{ 2.10};
again: confirmation.

Finally, applying our general formula~\thetag{ 2.25}, we deduce the
value of ${\bf Y}_6$ {\it without having to use ${\bf Y}_5$ and the
induction formulas~\thetag{ 2.4}}, which shortens substantially the
computations. 

For the pleasure, we obtain:
\def\theequation{2.29}\begin{equation}
\small
\left\{
\aligned
{\bf Y}_6
&
=
\mathcal{Y}_{x^6}
+
\left[
6\,\mathcal{Y}_{x^5y}
-
\mathcal{X}_{x^6}
\right]
y_1
+
\left[
15\,\mathcal{Y}_{x^4y^2}
-
6\,\mathcal{X}_{x^5y}
\right]
(y_1)^2
+ \\
& \
\ \ \ \ \
+
\left[
20\,\mathcal{Y}_{x^3y^3}
-
15\,\mathcal{X}_{x^4y^2}
\right]
(y_1)^3
+
\left[
15\,\mathcal{Y}_{x^2y^4}
-
20\,\mathcal{X}_{x^3y^3}
\right]
(y_1)^4
+ \\
& \
\ \ \ \ \
+
\left[
6\,\mathcal{Y}_{xy^5}
-
15\,\mathcal{X}_{x^2y^4}
\right]
(y_1)^5
+
\left[
\mathcal{Y}_{y^6}
-
6\,\mathcal{X}_{xy^5}
\right]
(y_1)^6
+
\left[
-
\mathcal{X}_{y^6}
\right]
(y_1)^7
+ \\
& \
\ \ \ \ \
+
\left[
15\,\mathcal{Y}_{x^4y}
-
6\,\mathcal{X}_{x^5}
\right]
y_2
+
\left[
60\,\mathcal{Y}_{x^3y^2}
-
45\,\mathcal{X}_{x^4y}
\right]
y_1y_2
+ \\
& \
\ \ \ \ \
+
\left[
90\,\mathcal{Y}_{x^2y^3}
-
120\,\mathcal{X}_{x^3y^2}
\right]
(y_1)^2y_2
+
\left[
60\,\mathcal{Y}_{xy^4}
-
150\,\mathcal{X}_{x^2y^3}
\right]
(y_1)^3y_2
+ \\
& \
\ \ \ \ \ 
+
\left[
15\,\mathcal{Y}_{y^5}
-
90\,\mathcal{X}_{xy^4}
\right]
(y_1)^4y_2
+
\left[
-
21\,\mathcal{X}_{y^5}
\right]
(y_1)^5y_2
+ \\
& \
\ \ \ \ \
+
\left[
45\,\mathcal{Y}_{x^2y^2}
-
60\,\mathcal{X}_{x^3y}
\right]
(y_2)^2
+
\left[
90\,\mathcal{Y}_{xy^3}
-
225\,\mathcal{X}_{x^2y^2}
\right]
y_1(y_2)^2
+ \\
& \
\ \ \ \ \
+
\left[
45\,\mathcal{Y}_{y^4}
-
270\,\mathcal{X}_{xy^3}
\right]
(y_1)^2(y_2)^2
+
\left[
-
210\,\mathcal{X}_{y^4}
\right]
(y_1)^3(y_2)^2
+ \\
& \
\ \ \ \ \
+
\left[
15\,\mathcal{Y}_{y^3}
-
90\,\mathcal{X}_{xy^2}
\right]
(y_2)^3
+
\left[
-
105\,\mathcal{X}_{y^3}
\right]
y_1(y_2)^3
+ \\
& \
\ \ \ \ \
+
\left[
20\,\mathcal{Y}_{x^3y}
-
15\,\mathcal{X}_{x^4}
\right]
y_3
+
\left[
60\,\mathcal{Y}_{x^2y^2}
-
80\,\mathcal{X}_{x^3y}
\right]
y_1y_3
+ \\
& \
\ \ \ \ \
+
\left[
60\,\mathcal{Y}_{xy^3}
-
150\,\mathcal{X}_{x^2y^2}
\right]
(y_1)^2y_3
+
\left[
20\,\mathcal{Y}_{y^4}
-
120\,\mathcal{X}_{xy^3}
\right]
(y_1)^3y_3
+ \\
& \
\ \ \ \ \
+
\left[
-
35\,\mathcal{X}_{y^4}
\right]
(y_1)^4y_3
+
\left[
60\,\mathcal{Y}_{xy^2}
-
150\,\mathcal{X}_{x^2y}
\right]
y_2y_3
+ \\
& \
\ \ \ \ \
+
\left[
60\,\mathcal{Y}_{y^3}
-
360\,\mathcal{X}_{xy^2}
\right]
y_1y_2y_3
+
\left[
-
210\,\mathcal{X}_{y^3}
\right]
(y_1)^2y_2y_3
+ \\
& \
\ \ \ \ \
+
\left[
-
105\,\mathcal{X}_{y^2}
\right]
(y_2)^2y_3
+
\left[
10\,\mathcal{Y}_{y^2}
-
60\,\mathcal{X}_{xy}
\right]
(y_3)^2
+ \\
& \
\ \ \ \ \
+
\left[
-
70\,\mathcal{X}_{y^2}
\right]
y_1(y_3)^2
+
\left[
15\,\mathcal{Y}_{x^2y}
-
20\,\mathcal{X}_{x^3}
\right]
y_4
+ \\
& \
\ \ \ \ \
+
\left[
30\,\mathcal{Y}_{xy^2}
-
75\,\mathcal{X}_{x^2y}
\right]
y_1y_4
+
\left[
15\,\mathcal{Y}_{y^3}
-
90\,\mathcal{X}_{xy^2}
\right]
(y_1)^2y_4
+ \\
& \
\ \ \ \ \
+
\left[
-
35\,\mathcal{X}_{y^3}
\right]
(y_1)^3y_4
+
\left[
15\,\mathcal{Y}_{y^2}
-
90\,\mathcal{X}_{xy}
\right]
y_2y_4
+ \\
& \
\ \ \ \ \
+ 
\left[
-
105\,\mathcal{X}_{y^2}
\right]
y_1y_2y_4
+ 
\left[
-
35\,\mathcal{X}_y
\right]
y_3y_4
+
\left[
6\,\mathcal{Y}_{xy}
-
15\,\mathcal{X}_{x^2}
\right]
y_5
+ \\
& \
\ \ \ \ \
+
\left[
6\,\mathcal{Y}_{y^2}
-
36\,\mathcal{X}_{xy}
\right]
y_1y_5+
\left[
-
21\,\mathcal{X}_{y^2}
\right]
(y_1)^2y_5
+
\left[
-
21\,\mathcal{X}_y
\right]
y_2y_5
+ \\
& \
\ \ \ \ \
+
\left[
\mathcal{Y}_y
-
6\,\mathcal{X}_y
\right]
y_6
+
\left[
-
7\,\mathcal{X}_y
\right]
y_1y_6.
\endaligned\right.
\end{equation}

\subsection*{ 2.30.~Deduction of the classical Fa\`a 
di Bruno formula} Let $x,y \in \K$ and let $g = g(x)$, $f = f ( y)$ be
two $\mathcal{ C }^\infty$-smooth functions $\K \to \K$. Consider the
composition $h := f \circ g$, namely $h(x) = f (g (x))$. For $\lambda
\in \N$ with $\lambda \geqslant 1$, simply denote by $g_\lambda$ the
$\lambda$-th derivative $\frac{ d^\lambda g}{d x^\lambda}$ and
similarly for $h_\lambda$. Also, abbreviate $f_\lambda := \frac{
d^\lambda f}{ d y^\lambda}$.

By the classical formula for the derivative of a composite function,
we have $h_1 = f_1 \, g_1$. Further computations provide the following
list of subsequent derivatives of $h$:
\def\theequation{2.31}\begin{equation}
\left\{
\aligned
h_1 
& 
=
f_1 \, g_1, 
\\
h_2
& 
=
f_2 \, (g_1)^2
+ 
f_1 \, g_2, 
\\
h_3
&
=
f_3\,(g_1)^3+3\,f_2\,g_1\,g_2
+
f_1\,g_3,
\\
h_4
&
=
f_4\,(g_1)^4
+
6\,f_3\,(g_1)^2\,g_2
+
3\,f_2\,(g_2)^2
+
4\,f_2\,g_1\,g_3
+
f_1\,g_4,
\\
h_5
&
=
f_5\,(g_1)^5\,
+
10\,f_4\,(g_1)^3\,g_2
+
15\,f_3\,(g_1)^2\,g_3
+
10\,f_3\,g_1\,(g_2)^2
+ \\
& \
\ \ \ \ \
+
10\,f_2\,g_2\,g_3
+
5\,f_2\,g_1\,g_4
+
f_1\,g_5, 
\\
h_6
&
=
f_6\,(g_1)^6\,
+
15\,f_5\,(g_1)^4\,g_2
+
45\,f_4\,(g_1)^2\,(g_2)^2
+
15\,f_3\,(g_2)^3
+ \\
& \
\ \ \ \ \
+
20\,f_4\,(g_1)^3\,g_3
+
60\,f_3\,g_1\,g_2\,g_3
+
10\,f_2\,(g_3)^2
+
15\,f_3\,(g_1)^2\,g_4
+ \\
& \
\ \ \ \ \
+
15\,f_2\,g_2\,g_4
+
6\,f_2\,g_1\,g_5
+
f_1\,g_6.
\endaligned\right.
\end{equation}

\def\thetheorem{2.32}\begin{theorem}
For every integer $\kappa \geqslant 1$, the $\kappa$-th derivative of the
composite function $h = f\circ g$ may be expressed as an explicit
polynomial in the partial derivatives of $f$ and of $g$ having integer
coefficients{\rm :}
\def\theequation{2.33}\begin{equation}
\boxed{
\aligned
\frac{ d^\kappa h}{dx^\kappa}
& 
=
\sum_{d=1}^\kappa
\
\sum_{1\leqslant\lambda_1<\cdots<\lambda_d\leqslant\kappa}
\
\sum_{\mu_1\geqslant 1,\dots,\mu_d\geqslant 1}
\
\sum_{\mu_1\lambda_1+\cdots+\mu_d\lambda_d=\kappa}
\\
& \
\ \ \ \ \ 
\frac{\kappa !}{(\lambda_1!)^{\mu_1}\ \mu_1! 
\cdots
(\lambda_d!)^{\mu_d}\ \mu_d!}
\
\frac{d^{\mu_1+\cdots+\mu_d} f}{
dy^{\mu_1+\cdots+\mu_d}}
\
\left(
\frac{d^{\lambda_1}g}{dx^{\lambda_1}}
\right)^{\mu_1}
\cdots\cdots
\left(
\frac{d^{\lambda_d}g}{dx^{\lambda_d}}
\right)^{\mu_d}.
\endaligned
}
\end{equation}
\end{theorem}

This is the classical {\it Fa\`a di Bruno formula}. Interestingly, we
observe that this formula is included as a subpart of the general
formula for ${\bf Y}_\kappa$, after a suitable translation. Indeed,
in the formulas for ${\bf Y}_1$, ${\bf Y}_2$, ${\bf Y}_3$, ${\bf
Y}_4$, ${\bf Y}_5$, ${\bf Y}_6$ and in the general sum for ${\bf
Y}_\kappa$, pick only the terms for which $\mu_1\lambda_1 + \cdots +
\mu_d \lambda_d = \kappa$ and drop $\mathcal{ X}$, which yields:
\def\theequation{2.34}\begin{equation}
\aligned
& \
\sum_{d=1}^\kappa
\
\sum_{1\leqslant\lambda_1<\cdots<\lambda_d\leqslant \kappa}
\
\sum_{\mu_1\geqslant 1,\dots,\mu_d\geqslant 1}
\
\sum_{\mu_1\lambda_1+\cdots+\mu_d\lambda_d=\kappa}
\\
& \
\left[
\frac{\kappa!}{
\mu_1!(\lambda_1!)^{\mu_1}
\cdots
\mu_d!(\lambda_d!)^{\mu_d}
}
\
\mathcal{Y}_{y^{\mu_1+\cdots+\mu_d}}
\right]
\left(
y_{\lambda_1}
\right)^{\mu_1}
\cdots
\left(
y_{\lambda_d}\right)^{\mu_d}.
\endaligned
\end{equation}
The similarity between the two formulas~\thetag{ 2.33} and~\thetag{
2.34} is now clearly visible.

The Fa\`a di Bruno formula may be established by means of
substitutions of power series (\cite{ f1969}, p.~222), by means of the
umbral calculus (\cite{ cs1996}), or by means of some induction
formulas, which we write for completeness. Define the differential
operators
\def\theequation{2.35}\begin{equation}
\small
\aligned
F_2
&
:=
g_2\,\frac{\partial}{\partial g_1}
+
g_1\left(
f_2\,\frac{\partial}{\partial f_1}
\right), 
\\
F_3
&
:=
g_2\,\frac{\partial}{\partial g_1}
+
g_3\,\frac{\partial}{\partial g_2}
+
g_1\left(
f_2\,\frac{\partial}{\partial f_1}
+
f_3\,\frac{\partial}{\partial f_2}
\right), 
\\
\cdots
&
\cdots\cdots
\cdots\cdots
\cdots\cdots
\cdots\cdots
\cdots\cdots
\cdots\cdots
\cdots\cdots
\\
F_\lambda
&
:=
g_2\,\frac{\partial}{\partial g_1}
+
g_3\,\frac{\partial}{\partial g_2}
+
\cdots
+
g_\lambda\,\frac{\partial}{\partial g_{\lambda-1}}
+
g_1\left(
f_2\,\frac{\partial}{\partial f_1}
+
f_3\,\frac{\partial}{\partial f_2}
+
\cdots
+
f_\lambda\,\frac{\partial}{\partial f_{\lambda-1}}
\right).
\endaligned
\end{equation}
Then we have 
\def\theequation{2.36}\begin{equation}
\aligned
h_2
&
=
F^2(h_1),
\\
h_3
&
=
F^3(h_2),
\\
\cdots
&
\cdots
\cdots
\cdots
\cdots
\\
h_\lambda
&
=
F^\lambda(h_{\lambda-1}).
\endaligned
\end{equation}

\section*{\S3.~Several independent variables and
one dependent variable}

\subsection*{3.1.~Simplified adapted notations}
As announced after the statement of Theorem~2.24, it is only after we
have treated the case of several independent variables that we will
understand perfectly the general formula~\thetag{ 2.25}, valid in the
case of one independent variable and one dependent variable. We will
discover massive formal computations, exciting our computational
intuition.

Thus, assume $n\geqslant 1$ and $m=1$, let $\kappa\in \N$ with $\kappa
\geqslant 1$ and simply denote (instead of~\thetag{ 1.2}) the jet
variables by:
\def\theequation{3.2}\begin{equation}
\left(
x^i,y,y_{i_1},y_{i_1,i_2},\dots,y_{i_1,i_2,\dots,i_\kappa}
\right).
\end{equation}
Also, instead of~\thetag{ 1.30}, denote
the $\kappa$-th prolongation of a vector field by:
\def\theequation{3.3}\begin{equation}
\left\{
\aligned
\mathcal{L}^{(\kappa)}
&
=
\sum_{i=1}^n\,\mathcal{X}^i\,\frac{\partial}{\partial x^i}
+
\mathcal{Y}\,\frac{\partial}{\partial y}
+
\sum_{i_1=1}^n\,{\bf Y}_{i_1}\,\frac{\partial}{\partial y_{i_1}}
+
\sum_{i_1,i_2=1}^n\,{\bf Y}_{i_1,i_2}\,
\frac{\partial}{\partial y_{i_1,i_2}}
+ \\
& \
\ \ \ \ \
+
\cdots
+
\sum_{i_1,i_2,\dots,i_\kappa=1}^n\,
{\bf Y}_{i_1,i_2,\dots,i_\kappa}\,
\frac{\partial}{\partial y_{i_1,i_2,\dots,i_\kappa}}.
\endaligned\right.
\end{equation}
The induction formulas are
\def\theequation{3.4}\begin{equation}
\left\{
\aligned
{\bf Y}_{i_1}
&
:=
D_{i_1}^1
\left(
\mathcal{ Y}
\right)
-
\sum_{k=1}^n\,D_{i_1}^1
\left(
\mathcal{X}^k
\right)
\, y_k, 
\\
{\bf Y}_{i_1,i_2}
&
:=
D_{i_2}^2
\left(
{\bf Y}_{i_1}
\right)
-
\sum_{k=1}^n\,D_{i_2}^1
\left(
\mathcal{X}^k
\right)
\, y_{i_1,k}, 
\\
\cdots \cdots \cdots 
&
\cdots \cdots \cdots \cdots \cdots \cdots \cdots \cdots
\cdots \cdots \cdots \cdots \cdots \cdots \cdots \cdots
\\
{\bf Y}_{i_1,i_2,\dots,i_\kappa}
&
:=
D_{i_\kappa}^\kappa
\left(
{\bf Y}_{i_1,i_2,\dots,i_{\kappa-1}}
\right)
-
\sum_{k=1}^n\,D_{i_\kappa}^1
\left(
\mathcal{X}^k
\right)
\, y_{i_1,i_2,\dots,i_{\kappa-1},k},
\endaligned\right.
\end{equation}
where the total differentiation operators $D_{ i' }^\lambda$ are
defined as in~\thetag{ 1.22}, dropping the sums $\sum_{j ' = 1}^m$ and
the indices $j'$.

\subsection*{3.5.~Two instructing explicit computations}
To begin with, let us compute ${\bf Y}_{i_1}$. With $D_{i_1}^1 =
\frac{ \partial }{\partial x^{i_1}} + y_{i_1}\, \frac{\partial
}{\partial y}$, we have:
\def\theequation{3.6}\begin{equation}
\aligned
{\bf Y}_{i_1}
&
=
D_{i_1}
\left(
\mathcal{Y}
\right)
-
\sum_{k_1=1}^n\,D_{i_1}^1
\left(
\mathcal{X}^{k_1}
\right)
y_{k_1} 
\\
& 
=
\mathcal{Y}_{x^{i_1}}
+
\mathcal{Y}_y\,y_{i_1}
-
\sum_{k_1=1}^n\,\mathcal{X}_{x^{i_1}}^{k_1}\,y_{k_1}
-
\sum_{k_1=1}^n\,\mathcal{X}_y^{k_1}\,y_{i_1}\,y_{k_1}.
\endaligned
\end{equation}
Searching for formal harmony and for coherence with the formula
$(2.6)_1$, we must include the term $\mathcal{ Y}_y\, y_{i_1}$ inside
the sum $\sum_{ k_1 =1 }^n\, \left[ \cdot \right] y_{ k_1}$. Using
the Kronecker symbol, we may write:
\def\theequation{3.7}\begin{equation}
\mathcal{Y}_y\,y_{i_1}
\equiv
\sum_{k_1=1}^n\,
\left[
\delta_{i_1}^{k_1}\,\mathcal{Y}_y
\right]
y_{k_1}.
\end{equation}
Also, we may rewrite the last term of~\thetag{ 3.6} with a double sum:
\def\theequation{3.8}\begin{equation}
-
\sum_{k_1=1}^n\,
\mathcal{X}_y^{k_1}\,y_{i_1}\,y_{k_1}
\equiv
\sum_{k_1,k_2=1}^n\,
\left[
-
\delta_{i_1}^{k_1}\,\mathcal{X}_y^{k_2}
\right]
y_{k_1}y_{k_2}.
\end{equation}
From now on and up to equation~\thetag{ 3.39}, we shall abbreviate any
sum $\sum_{k=1}^n$ from $1$ to $n$ as $\sum_k$. Putting everything
together, we get the final desired perfect expression of ${\bf
Y}_{i_1}$:
\def\theequation{3.9}\begin{equation}
{\bf Y}_{i_1}
=
\mathcal{Y}_{x^{i_1}}
+
\sum_{k_1}\,
\left[
\delta_{i_1}^{k_1}\,\mathcal{Y}_y
-
\mathcal{X}_{x^{i_1}}^{k_1}
\right]
y_{k_1}
+
\sum_{k_1,k_2}\,
\left[
-
\delta_{i_1}^{k_1}\,\mathcal{X}_y^{k_2}
\right]
y_{k_1}y_{k_2}.
\end{equation}
This completes the first explicit computation.

The second one is about ${\bf Y}_{ i_1, i_2}$. It becomes more
delicate, because several algebraic transformations must be achieved
until the final satisfying formula is obtained. Our goal is to
present each step very carefully, explaining every tiny
detail. Without such a care, it would be impossible to claim that some
of our subsequent computations, for which we will not provide the
intermediate steps, may be redone and verified. Consequently, we will
expose our rules of formal computation thoroughly.

Replacing the value of ${\bf Y}_1$ just obtained in the induction
formula $(3.4)_2$ and developing, we may conduct the very
first steps of the computation:
$$
\small
\aligned
{\bf Y}_{i_1,i_2}
&
=
D_{i_2}^2
\left(
{\bf Y}_{i_1}
\right)
-
\sum_{k_1}\,D_{i_2}^1
\left(
\mathcal{X}^{k_1}
\right)
y_{i_1,k_1} 
\\
& 
=
\left(
\frac{\partial}{\partial x^{i_2}}
+
y_{i_2}\,\frac{\partial}{\partial y}
+
\sum_{k_1}\,y_{i_2,k_1}\,
\frac{\partial}{\partial y_{k_1}}
\right)
\left(
\mathcal{Y}_{x^{i_1}}
+
\sum_{k_1}\,
\left[
\delta_{i_1}^{k_1}\,\mathcal{Y}_y
-
\mathcal{X}_{x^{i_1}}^{k_1}
\right]
y_{k_1}
+
\right.
\\
& \
\ \ \ \ \ \ \ \ \ \ \ \ \ \ \
\ \ \ \ \ \
\left.
+
\sum_{k_1,k_2}\,
\left[
-
\delta_{i_1}^{k_1}\,\mathcal{X}_y^{k_2}
\right]
y_{k_1}y_{k_2}
\right)
-
\sum_{k_1}\,
\left[
\mathcal{X}_{x^{i_2}}^{k_1}
+
y_{i_2}\,\mathcal{X}_y^{k_1}
\right]
y_{i_1,k_1}
\endaligned
$$
\def\theequation{3.10}\begin{equation}
\small
\aligned
&
=
\left(
\frac{\partial}{\partial x^{i_2}}
\right)
\left(
\mathcal{Y}_{x^{i_1}}
+
\sum_{k_1}\,
\left[
\delta_{i_1}^{k_1}\,\mathcal{Y}_y
-
\mathcal{X}_{x^{i_1}}^{k_1}
\right]
y_{k_1}
+
\sum_{k_1,k_2}\,
\left[
-
\delta_{i_1}^{k_1}\,\mathcal{X}_y^{k_2}
\right]
y_{k_1}y_{k_2}
\right)
+ \\
& \
\ \ \ \ \
+
\left(
y_{i_2}\,\frac{\partial}{\partial y}
\right)
\left(
\mathcal{Y}_{x^{i_1}}
+
\sum_{k_1}\,
\left[
\delta_{i_1}^{k_1}\,\mathcal{Y}_y
-
\mathcal{X}_{x^{i_1}}^{k_1}
\right]
y_{k_1}
+
\sum_{k_1,k_2}\,
\left[
-
\delta_{i_1}^{k_1}\,\mathcal{X}_y^{k_2}
\right]
y_{k_1}y_{k_2}
\right)
+ \\
& \
\ \ \ \ \
+
\left(
\sum_{k_1}\,y_{i_2,k_1}\,
\frac{\partial}{\partial y_{k_1}}
\right)
\left(
\mathcal{Y}_{x^{i_1}}
+
\sum_{k_1}\,
\left[
\delta_{i_1}^{k_1}\,\mathcal{Y}_y
-
\mathcal{X}_{x^{i_1}}^{k_1}
\right]
y_{k_1}
+
\sum_{k_1,k_2}\,
\left[
-
\delta_{i_1}^{k_1}\,\mathcal{X}_y^{k_2}
\right]
y_{k_1}y_{k_2}
\right) 
+ \\
& \
\ \ \ \ \
+
\sum_{k_1}\,
\left[
-
\mathcal{X}_{x^{i_2}}^{k_1}
\right]
y_{k_1,i_1}
+
\sum_{k_1}\,
\left[
-
\mathcal{X}_y^{k_1}
\right]
y_{i_2}y_{i_1,k_1}
\endaligned
\end{equation}
$$
\small
\aligned
& 
=
\mathcal{Y}_{x^{i_1}x^{i_2}}
+
\sum_{k_1}\,
\left[
\delta_{i_1}^{k_1}\,\mathcal{Y}_{x^{i_2}y}
-
\mathcal{X}_{x^{i_1}x^{i_2}}^{k_1}
\right]
y_{k_1}
+
\sum_{k_1,k_2}\,
\left[
-
\delta_{i_1}^{k_1}\,\mathcal{X}_{x^{i_2}y}^{k_2}
\right]
y_{k_1}y_{k_2}
+ \\
& \
\ \ \ \ \
+
\mathcal{Y}_{x^{i_1}y}\,y_{i_2}
+
\sum_{k_1}\,
\left[
\delta_{i_1}^{k_1}\,\mathcal{Y}_{yy}
-
\mathcal{X}_{x^{i_1}y}^{k_1}
\right]
y_{k_1}y_{i_2}
+
\sum_{k_1,k_2}\,
\left[
-
\delta_{i_1}^{k_1}\,\mathcal{X}_{yy}^{k_2}
\right]
y_{k_1}y_{k_2}y_{i_2}
+ \\
& \
\ \ \ \ \
+
\sum_{k_1}\,
\left[
\delta_{i_1}^{k_1}\,\mathcal{Y}_y
-
\mathcal{X}_{x^{i_1}}^{k_1}
\right]
y_{i_2,k_1}
+
\sum_{k_1,k_2}\,
\left[
-
\delta_{i_1}^{k_1}\,\mathcal{X}_y^{k_2}
\right]
y_{k_2}y_{i_2,k_1}
+
\sum_{k_1,k_2}\,
\left[
-
\delta_{i_1}^{k_1}\,\mathcal{X}_y^{k_2}
\right]
y_{k_1}y_{i_2,k_2}
+ \\
& \
\ \ \ \ \
+
\sum_{k_1}\,
\left[
-
\mathcal{X}_{x^{i_2}}^{k_1}
\right]
y_{k_1,i_1}
+
\sum_{k_1}
\left[
-
\mathcal{X}_y^{k_1}
\right]
y_{i_2}y_{i_1,k_1}.
\endaligned
$$
Some explanations are needed about the computation of the last two
terms of line 11, {\it i.e.} about the passage from line 7 of~\thetag{
3.10} just above to line 11. We have to compute:
\def\theequation{3.11}\begin{equation}
\footnotesize
\left(
\sum_{k_1}\,y_{i_2,k_1}\,
\frac{\partial}{\partial y_{k_1}}
\right)
\left(
\sum_{k_1,k_2}\,
\left[
-
\delta_{i_1}^{k_1}\,
\mathcal{X}_y^{k_2}
\right]
y_{k_1}y_{k_2}
\right).
\end{equation}
This term is of the form
\def\theequation{3.12}\begin{equation}
\footnotesize
\left(
\sum_{k_1}\,A_{k_1}\,
\frac{\partial}{\partial y_{k_1}}
\right)
\left(
\sum_{k_1,k_2}\,
\left[
B_{k_1,k_2}
\right]
y_{k_1}y_{k_2}
\right),
\end{equation}
where the terms $B_{k_1,k_2}$ are independent of the pure first jet
variables $y_{x^k}$.
By the rule of Leibniz for the differentiation 
of a product, we may write
\def\theequation{3.13}\begin{equation}
\footnotesize
\aligned
&
\left(
\sum_{k_1}\,A_{k_1}\,
\frac{\partial}{\partial y_{k_1}}
\right)
\left(
\sum_{k_1,k_2}\,
\left[
B_{k_1,k_2}
\right]
y_{k_1}y_{k_2}
\right)
= 
\\
&
=
\sum_{k_1,k_2}\,
\left[
B_{k_1,k_2}
\right]
y_{k_2}
\left(
\sum_{k_1'}\,A_{k_1'}\,
\frac{\partial}{\partial y_{k_1'}}
(y_{k_1})
\right)
+
\sum_{k_1,k_2}\,
\left[
B_{k_1,k_2}
\right]
y_{k_1}
\left(
\sum_{k_2'}\,A_{k_2'}\,
\frac{\partial}{\partial y_{k_2'}}
(y_{k_2})
\right)
\\
& 
=
\sum_{k_1,k_2}\,
\left[
B_{k_1,k_2}
\right]
y_{k_2}\,A_{k_1}
+
\sum_{k_1,k_2}\,
\left[
B_{k_1,k_2}
\right]
y_{k_1}\,A_{k_2}.
\endaligned
\end{equation}
This is how we have written line 11 of~\thetag{ 3.10}.

Next, the first term $\mathcal{ Y}_{ x^{i_1}y} \, y_{ i_2}$ in line 10
of~\thetag{ 3.10} is not in a suitable shape. For reasons of harmony
and coherence, we must insert it inside a sum of the form
$\sum_{k_1}\, \left[ \cdot \right] y_{k_1}$. Hence, using the
Kronecker symbol, we transform:
\def\theequation{3.14}\begin{equation}
\mathcal{ Y}_{x^{i_1}y} \, y_{i_2}
\equiv
\sum_{k_1}\,
\left[
\delta_{i_2}^{k_1}\,
\mathcal{Y}_{x^{i_1}y}
\right]
y_{k_1}.
\end{equation}
Also, we must ``summify'' the seven other terms, remaining in lines 10,
11 and 12 of~\thetag{ 3.10}. Sometimes, we use the symmetry $y_{i_2,
k_1} \equiv y_{ k_1, i_2}$ without mention. Similarly, we get:
$$
\aligned
\sum_{k_1}\,
\left[
\delta_{i_1}^{k_1}\,\mathcal{Y}_{yy}
-
\mathcal{X}_{x^{i_1y}}^{k_1}
\right]
y_{k_1}y_{i_2}
&
\equiv
\sum_{k_1,k_2}\,
\left[
\delta_{i_1}^{k_1}\,\delta_{i_2}^{k_2}\,
\mathcal{Y}_{yy}
-
\delta_{i_2}^{k_2}\,\mathcal{X}_{x^{i_1}y}^{k_1}
\right]
y_{k_1}y_{k_2}, \\
\sum_{k_1,k_2}\,
\left[
-
\delta_{i_1}^{k_1}\,
\mathcal{X}_{yy}^{k_2}
\right]
y_{k_1}y_{k_2}y_{i_2}
&
\equiv
\sum_{k_1,k_2,k_3}\,
\left[
-
\delta_{i_1}^{k_1}\,
\delta_{i_2}^{k_3}
\mathcal{X}_{yy}^{k_2}
\right]
y_{k_1}y_{k_2}y_{k_3}, \\
\sum_{k_1}
\left[
\delta_{i_1}^{k_1}\,
\mathcal{Y}_y
-
\mathcal{X}_{x^{i_1}}^{k_1}
\right]
y_{k_1,i_2}
&
\equiv
\sum_{k_1,k_2}
\left[
\delta_{i_1}^{k_1}\,
\delta_{i_2}^{k_2}\,
\mathcal{Y}_y
-
\delta_{i_2}^{k_2}\,
\mathcal{X}_{x^{i_1}}^{k_1}
\right]
y_{k_1,k_2},
\endaligned
$$
\def\theequation{3.15}\begin{equation}
\aligned
\sum_{k_1,k_2}\,
\left[
-
\delta_{i_1}^{k_1}\,\mathcal{X}_y^{k_2}
\right]
y_{k_2}y_{k_1,i_2}
&
=
\sum_{k_1,k_2}\,
\left[
-
\delta_{i_1}^{k_2}\,\mathcal{X}_y^{k_1}
\right]
y_{k_1}y_{k_2,i_2}
\\
& 
\equiv
\sum_{k_1,k_2,k_3}\,
\left[
-
\delta_{i_1}^{k_2}\,\delta_{i_2}^{k_3}\,
\mathcal{X}_y^{k_1}
\right]
y_{k_1}y_{k_2,k_3},
\endaligned
\end{equation}
$$
\aligned
\sum_{k_1,k_2}\,
\left[
-
\delta_{i_1}^{k_1}\,\mathcal{X}_y^{k_2}
\right]
y_{k_1}y_{k_2,i_2}
&
\equiv
\sum_{k_1,k_2,k_3}\,
\left[
-
\delta_{i_1}^{k_1}\,
\delta_{i_2}^{k_3}\,
\mathcal{X}_y^{k_2}
\right]
y_{k_1}y_{k_2,k_3}, 
\\
\sum_{k_1}\,
\left[
-
\mathcal{X}_{x^{i_2}}^{k_1}
\right]
y_{k_1,i_1}
&
\equiv
\sum_{k_1,k_2}\,
\left[
-
\delta_{i_1}^{k_2}
\mathcal{X}_{x^{i_2}}^{k_1}
\right]
y_{k_1,k_2}, 
\\
\sum_{k_1}\,
\left[
-
\mathcal{X}_y^{k_1}
\right]
y_{i_2}y_{k_1,i_1}
&
=
\sum_{k_2}\,
\left[
-
\mathcal{X}_y^{k_2}
\right]
y_{i_2}y_{k_2,i_1}
\\
& 
\equiv
\sum_{k_1,k_2,k_3}\,
\left[
-
\delta_{i_2}^{k_1}\,
\delta_{i_1}^{k_3}\,
\mathcal{X}_y^{k_2}
\right]
y_{k_1}y_{k_2,k_3}.
\endaligned
$$
In the sequel, for products of
Kronecker symbols, it will be convenient
to adopt the following self-evident contracted
notation:
\def\theequation{3.16}\begin{equation}
\delta_{i_1}^{k_1}\,
\delta_{i_2}^{k_2}
\equiv
\delta_{i_1,\,i_2}^{k_1,k_2}; 
\ \ \ \ \ \
{\rm generally:} 
\ \ \
\delta_{i_1}^{k_1}\,
\delta_{i_2}^{k_2}
\cdots
\delta_{i_\lambda}^{k_\lambda}
\equiv
\delta_{i_1,\,i_2,\,\cdots\,,i_\lambda}^{
k_1,k_2,\cdots,k_\lambda}.
\end{equation}
Re-inserting plainly these eight summified terms~\thetag{ 3.14},
\thetag{ 3.15} in the last expression~\thetag{ 3.10} of ${\bf Y}_{i_1,
i_2}$ (lines 10, 11 and 12), we get:
\def\theequation{3.17}\begin{equation}
\small
\aligned
{\bf Y}_{i_1,i_2}
&
=
\underline{
\mathcal{Y}_{x^{i_1}x^{i_2}} 
}_{ \fbox{\tiny 1}}
+
\underline{
\sum_{k_1}\,
\left[
\delta_{i_1}^{k_1}\,
\mathcal{Y}_{x^{i_2}y}
-
\mathcal{X}_{x^{i_1}x^{i_2}}^{k_1}
\right]
y_{k_1}
}_{ \fbox{\tiny 2}}
+
\underline{
\sum_{k_1,k_2}\,
\left[
-
\delta_{i_1}^{k_1}\,
\mathcal{X}_{x^{i_2}y}^{k_2}
\right]
y_{k_1}y_{k_2}
}_{ \fbox{\tiny 3}}
+ \\
& \
\ \ \ \ \
+
\underline{
\sum_{k_1}\,
\left[
\delta_{i_2}^{k_1}\,
\mathcal{Y}_{x^{i_1}y}
\right]
y_{k_1}
}_{ \fbox{\tiny 2}}
+ 
\underline{
\sum_{k_1,k_2}\,
\left[
\delta_{i_1,\,i_2}^{k_1,k_2}\,
\mathcal{Y}_{yy}
-
\delta_{i_2}^{k_2}\,
\mathcal{X}_{x^{i_1}y}^{k_1}
\right]
y_{k_1}y_{k_2}
}_{ \fbox{\tiny 3}}
+ \\
& \
\ \ \ \ \
+
\underline{
\sum_{k_1,k_2,k_3}\,
\left[
-
\delta_{i_1,\,i_2}^{k_1,k_3}\,
\mathcal{X}_{yy}^{k_2}
\right]
y_{k_1}y_{k_2}y_{k_3}
}_{ \fbox{\tiny 4}}
+
\underline{
\sum_{k_1,k_2}\,
\left[
\delta_{i_1,\,i_2}^{k_1,k_2}\,
\mathcal{Y}_y
-
\delta_{i_2}^{k_2}\,
\mathcal{X}_{x^{i_1}}^{k_1}
\right]
y_{k_1,k_2}
}_{ \fbox{\tiny 5}}
+ \\
& \
\ \ \ \ \
+
\underline{
\sum_{k_1,k_2,k_3}\,
\left[
-
\delta_{i_1,\,i_2}^{k_2,k_3}\,
\mathcal{X}_y^{k_1}
\right]
y_{k_1}y_{k_2,k_3}
}_{ \fbox{\tiny 6}}
+
\underline{
\sum_{k_1,k_2,k_3}\,
\left[
-
\delta_{i_1,\,i_2}^{k_1,k_3}\,
\mathcal{X}_y^{k_2}
\right]
y_{k_1}y_{k_2,k_3}
}_{ \fbox{\tiny 6}}
+ \\
& \
\ \ \ \ \
+
\underline{
\sum_{k_1,k_2}\,
\left[
-
\delta_{i_1}^{k_2}\,
\mathcal{X}_{x^{i_2}}^{k_1}
\right]
y_{k_1,k_2}
}_{ \fbox{\tiny 5}}
+
\underline{
\sum_{k_1,k_2,k_3}\,
\left[
-
\delta_{i_2,\,i_1}^{k_1,k_3}\,
\mathcal{X}_y^{k_2}
\right]
y_{k_1}y_{k_2,k_3}
}_{ \fbox{\tiny 6}}.
\endaligned
\end{equation}
Next, we gather the underlined terms, ordering them according to their
number. This yields 6 collections of
sums of monomials in the 
pure jet variables:
\def\theequation{3.18}\begin{equation}
\small
\aligned
{\bf Y}_{i_1,i_2}
&
=
\mathcal{Y}_{x^{i_1}x^{i_2}} 
+
\sum_{k_1}\,
\left[
\delta_{i_1}^{k_1}\,\mathcal{Y}_{x^{i_2}y}
+
\delta_{i_2}^{k_1}\,\mathcal{Y}_{x^{i_1}y}
-
\mathcal{X}_{x^{i_1}x^{i_2}}^{k_1}
\right]
y_{k_1}
+ \\
& \
\ \ \ \ \
+
\sum_{k_1,k_2}\,
\left[
\delta_{i_1,\,i_2}^{k_1,k_2}\,
\mathcal{Y}_{yy}
-
\delta_{i_1}^{k_1}\,
\mathcal{X}_{x^{i_2}y}^{k_2}
-
\delta_{i_2}^{k_2}\,
\mathcal{X}_{x^{i_1}y}^{k_1}
\right]
y_{k_1}y_{k_2}
+ \\
& \
\ \ \ \ \
+
\sum_{k_1,k_2,k_3}\,
\left[
-
\delta_{i_1,\,i_2}^{k_1,k_3}\,
\mathcal{X}_{yy}^{k_2}
\right]
y_{k_1}y_{k_2}y_{k_3}
+ \\
& \
\ \ \ \ \
+
\sum_{k_1,k_2}\,
\left[
\delta_{i_1,\,i_2}^{k_1,k_2}\,
\mathcal{Y}_y
-
\delta_{i_2}^{k_2}\,\mathcal{X}_{x^{i_1}}^{k_1}
-
\delta_{i_1}^{k_2}\,\mathcal{X}_{x^{i_2}}^{k_1}
\right]
y_{k_1,k_2}
+ \\
& \
\ \ \ \ \
+
\sum_{k_1,k_2,k_3}\,
\left[
-
\delta_{i_1,\,i_2}^{k_2,k_3}\,
\mathcal{X}_y^{k_1}
-
\delta_{i_1,\,i_2}^{k_1,k_3}\,
\mathcal{X}_y^{k_2}
-
\delta_{i_2,\,i_1}^{k_1,k_3}\,
\mathcal{X}_y^{k_2}
\right]
y_{k_1}y_{k_2,k_3}.
\endaligned
\end{equation}
To attain the real perfect harmony, this last expression has still to
be worked out a little bit.

\def\thelemma{3.19}\begin{lemma}
The final expression of ${\bf Y}_{i_1,i_2}$ is
as follows{\rm :}
\def\theequation{3.20}\begin{equation}
\left\{
\aligned
{\bf Y}_{i_1,i_2}
&
=
\mathcal{Y}_{x^{i_1}x^{i_2}} 
+
\sum_{k_1}\,
\left[
\delta_{i_1}^{k_1}\,\mathcal{Y}_{x^{i_2}y}
+
\delta_{i_2}^{k_1}\,\mathcal{Y}_{x^{i_1}y}
-
\mathcal{X}_{x^{i_1}x^{i_2}}^{k_1}
\right]
y_{k_1}
+ \\
& \
\ \ \ \ \
+
\sum_{k_1,k_2}\,
\left[
\delta_{i_1,\,i_2}^{k_1,k_2}\,
\mathcal{Y}_{yy}
-
\delta_{i_1}^{k_1}\,
\mathcal{X}_{x^{i_2}y}^{k_2}
-
\delta_{i_2}^{k_1}\,
\mathcal{X}_{x^{i_1}y}^{k_2}
\right]
y_{k_1}y_{k_2}
+ \\
& \
\ \ \ \ \
+
\sum_{k_1,k_2,k_3}\,
\left[
-
\delta_{i_1,\,i_2}^{k_1,k_2}\,
\mathcal{X}_{yy}^{k_3}
\right]
y_{k_1}y_{k_2}y_{k_3}
+ \\
& \
\ \ \ \ \
+
\sum_{k_1,k_2}\,
\left[
\delta_{i_1,\,i_2}^{k_1,k_2}\,
\mathcal{Y}_y
-
\delta_{i_1}^{k_1}\,\mathcal{X}_{x^{i_2}}^{k_2}
-
\delta_{i_2}^{k_1}\,\mathcal{X}_{x^{i_1}}^{k_2}
\right]
y_{k_1,k_2}
+ \\
& \
\ \ \ \ \
+
\sum_{k_1,k_2,k_3}\,
\left[
-
\delta_{i_1,\,i_2}^{k_1,k_2}\,
\mathcal{X}_y^{k_3}
-
\delta_{i_1,\,i_2}^{k_3,k_1}\,
\mathcal{X}_y^{k_2}
-
\delta_{i_1,\,i_2}^{k_2,k_3}\,
\mathcal{X}_y^{k_1}
\right]
y_{k_1}y_{k_2,k_3}.
\endaligned\right.
\end{equation}
\end{lemma}

\proof
As promised, we explain every tiny detail.

The first lines of~\thetag{ 3.18} and of~\thetag{ 3.20} are exactly
the same. For the transformations of terms in the second, in the third
and in the fourth lines, we use the following device. Let $\Upsilon_{
k_1, k_2}$ be an indexed quantity which is symmetric: $\Upsilon_{ k_1,
k_2} =\Upsilon_{ k_2, k_1}$. Let $A_{k_1, k_2}$ be an arbitrary
indexed quantity. Then obviously:
\def\theequation{3.21}\begin{equation}
\sum_{k_1,k_2}\,
A_{k_1, k_2}\,
\Upsilon_{k_1, k_2}
=
\sum_{k_1,k_2}\,
A_{k_2,k_1}\,
\Upsilon_{k_1,k_2}.
\end{equation}
Similar relations hold with a quantity $\Upsilon_{i_1, i_2, \dots,
i_\lambda}$ which is symmetric with respect to its $\lambda$ indices.
Consequently, in the second, in the third and in the fourth lines
of~\thetag{ 3.18}, we may permute freely certain indices in some of
the terms inside the brackets. This yields
the passage from lines 2, 3 and 4 of~\thetag{ 3.18}
to lines 2, 3 and 4 of~\thetag{ 3.20}.

It remains to explain how we pass from the fifth (last) line
of~\thetag{ 3.18} to the fifth (last) line of~\thetag{ 3.20}. The
bracket in the fifth line of~\thetag{ 3.18} contains three terms:
$\left[ -T_1 -T_2 -T_3 \right]$. The term $T_3$ involves the product
$\delta_{i_2,\,i_1}^{k_1,k_3}$, which we rewrite as
$\delta_{i_1,\,i_2}^{k_3,k_1}$, in order that $i_1$ appears before
$i_2$. Then, we rewrite the three terms in the new order $\left[ -T_2
-T_3 -T_1 \right]$, which yields:
\def\theequation{3.22}\begin{equation}
\sum_{k_1,k_2,k_3}\,
\left[
-
\delta_{i_1,\,i_2}^{k_1,k_3}\,
\mathcal{X}_y^{k_2}
-
\delta_{i_1,\,i_2}^{k_3,k_1}\,
\mathcal{X}_y^{k_2}
-
\delta_{i_1,\,i_2}^{k_2,k_3}\,
\mathcal{X}_y^{k_1}
\right]
y_{k_1}y_{k_2,k_3}.
\end{equation}
It remains to observe that we can permute $k_2$ and $k_3$ in the first
term $-T_2$, which yields the last
line of~\thetag{ 3.20}. The detailed proof is complete.
\endproof

\subsection*{3.23.~Final perfect expression of ${\bf Y}_{i_1, i_2,
i_3}$} Thanks to similar (longer) computations, we have obtained an
expression of ${\bf Y}_{i_1, i_2, i_3}$ which we consider to be in
final harmonious shape. Without copying the intermediate steps, let
us write down the result. The comments which 
are necessary to read it and to interpret it 
start just below.
$$
\small
\aligned
{\bf Y}_{i_1,i_2,i_3}
&
=
\mathcal{Y}_{x^{i_1}x^{i_2}x^{i_3}}
+
\sum_{k_1}\,
\left[
\delta_{i_1}^{k_1}\,\mathcal{Y}_{x^{i_2}x^{i_3}y}
+
\delta_{i_2}^{k_1}\,\mathcal{Y}_{x^{i_1}x^{i_3}y}
+
\delta_{i_3}^{k_1}\,\mathcal{Y}_{x^{i_1}x^{i_2}y}
-
\mathcal{X}_{x^{i_1}x^{i_2}x^{i_3}}^{k_1}
\right]
y_{k_1}
+ \\
& \
\ \ \ \ \
+
\sum_{k_1,k_2}\,
\left[
\delta_{i_1,\,i_2}^{k_1,k_2}\,\mathcal{Y}_{x^{i_3}y^2}
+
\delta_{i_3,\,i_1}^{k_1,k_2}\,\mathcal{Y}_{x^{i_2}y^2}
+
\delta_{i_2,\,i_3}^{k_1,k_2}\,\mathcal{Y}_{x^{i_1}y^2}
-
\right. \\
& \
\ \ \ \ \ \ \ \ \ \ \ \ \ \ \ \ \ \ \ \ 
\ \ \ \ \ \ \ \ \ \ \
\left.
-
\delta_{i_1}^{k_1}\,\mathcal{X}_{x^{i_2}x^{i_3}y}^{k_2}
-
\delta_{i_2}^{k_1}\,\mathcal{X}_{x^{i_1}x^{i_3}y}^{k_2}
-
\delta_{i_3}^{k_1}\,\mathcal{X}_{x^{i_1}x^{i_2}y}^{k_2}
\right]
y_{k_1}y_{k_2}
+ \\
& \
\ \ \ \ \
+
\sum_{k_1,k_2,k_3}\,
\left[
\delta_{i_1,\,i_2,\,i_3}^{k_1,k_2,k_3}\,\mathcal{Y}_{y^3}
-
\delta_{i_1,\,i_2}^{k_1,k_2}\,\mathcal{X}_{x^{i_3}y^2}^{k_3}
-
\delta_{i_1,\,i_3}^{k_1,k_2}\,\mathcal{X}_{x^{i_2}y^2}^{k_3}
-
\delta_{i_2,\,i_3}^{k_1,k_2}\,\mathcal{X}_{x^{i_1}y^2}^{k_3}
\right]
y_{k_1}y_{k_2}y_{k_3}
+ \\
& \
\ \ \ \ \
+
\sum_{k_1,k_2,k_3,k_4}\,
\left[
-\delta_{i_1,\,i_2,\,i_3}^{k_1,k_2,k_3}\,\mathcal{X}_{y^3}^{k_4}
\right]
y_{k_1}y_{k_2}y_{k_3}y_{k_4}
+ \\
\endaligned
$$
\def\theequation{3.24}\begin{equation}
\small
\aligned
& \
\ \ \ \ \
+
\sum_{k_1,k_2}\,
\left[
\delta_{i_1,\,i_2}^{k_1,k_2}\,\mathcal{Y}_{x^{i_3}y}
+
\delta_{i_3,\,i_1}^{k_1,k_2}\,\mathcal{Y}_{x^{i_2}y}
+
\delta_{i_2,\,i_3}^{k_1,k_2}\,\mathcal{Y}_{x^{i_1}y}
-
\right. 
\\
& \
\ \ \ \ \ \ \ \ \ \ \ \ \ \ \ \ \ \ \ \ \ \ 
\left.
-
\delta_{i_1}^{k_1}\,\mathcal{X}_{x^{i_2}x^{i_3}}^{k_2}
-
\delta_{i_2}^{k_1}\,\mathcal{X}_{x^{i_1}x^{i_3}}^{k_2}
-
\delta_{i_3}^{k_1}\,\mathcal{X}_{x^{i_1}x^{i_2}}^{k_2}
\right]
y_{k_1,k_2}
+ \\
& \
\ \ \ \ \
+
\sum_{k_1,k_2,k_3}\,
\left[
\delta_{i_1,\,i_2,\,i_3}^{k_1,k_2,k_3}\,\mathcal{Y}_{y^2}
+
\delta_{i_1,\,i_2,\,i_3}^{k_3,k_1,k_2}\,\mathcal{Y}_{y^2}
+
\delta_{i_1,\,i_2,\,i_3}^{k_2,k_3,k_1}\,\mathcal{Y}_{y^2}
-
\right.
\\
& \
\ \ \ \ \ \ \ \ \ \ \ \ \ \ \ \ \ \ \ \ \ \ \ \ \
\left.
-
\delta_{i_1,\,i_2}^{k_1,k_2}\,\mathcal{X}_{x^{i_3}y}^{k_3}
-
\delta_{i_1,\,i_2}^{k_3,k_1}\,\mathcal{X}_{x^{i_3}y}^{k_2}
-
\delta_{i_1,\,i_2}^{k_2,k_3}\,\mathcal{X}_{x^{i_3}y}^{k_1}
-
\right.
\\
& \
\ \ \ \ \ \ \ \ \ \ \ \ \ \ \ \ \ \ \ \ \ \ \ \ \
\left.
-
\delta_{i_1,\,i_3}^{k_1,k_2}\,\mathcal{X}_{x^{i_2}y}^{k_3}
-
\delta_{i_1,\,i_3}^{k_3,k_1}\,\mathcal{X}_{x^{i_2}y}^{k_2}
-
\delta_{i_1,\,i_3}^{k_2,k_3}\,\mathcal{X}_{x^{i_2}y}^{k_1}
-
\right.
\\
& \
\ \ \ \ \ \ \ \ \ \ \ \ \ \ \ \ \ \ \ \ \ \ \ \ \
\left.
-
\delta_{i_2,\,i_3}^{k_1,k_2}\,\mathcal{X}_{x^{i_1}y}^{k_3}
-
\delta_{i_2,\,i_3}^{k_3,k_1}\,\mathcal{X}_{x^{i_1}y}^{k_2}
-
\delta_{i_2,\,i_3}^{k_2,k_3}\,\mathcal{X}_{x^{i_1}y}^{k_1}
\right]
y_{k_1}y_{k_2,k_3}
+
\\
\endaligned
\end{equation}
$$
\small
\aligned
& \
\ \ \ \ \
+ 
\sum_{k_1,k_2,k_3,k_4}\,
\left[
-
\delta_{i_1,\,i_2,\,i_3}^{k_1,k_2,k_3}\,\mathcal{X}_{y^2}^{k_4}
-
\delta_{i_1,\,i_2,\,i_3}^{k_2,k_3,k_1}\,\mathcal{X}_{y^2}^{k_4}
-
\delta_{i_1,\,i_2,\,i_3}^{k_3,k_2,k_1}\,\mathcal{X}_{y^2}^{k_4}
- 
\right.
\\
& \
\ \ \ \ \ \ \ \ \ \ \ \ \ \ \ \ \ \ \ \ \ \ \ \ \
\left.
-
\delta_{i_1,\,i_2,\,i_3}^{k_3,k_4,k_1}\,\mathcal{X}_{y^2}^{k_2}
-
\delta_{i_1,\,i_2,\,i_3}^{k_3,k_1,k_4}\,\mathcal{X}_{y^2}^{k_2}
-
\delta_{i_1,\,i_2,\,i_3}^{k_1,k_3,k_4}\,\mathcal{X}_{y^2}^{k_2}
\right]
y_{k_1}y_{k_2}y_{k_3,k_4}
+ \\
& \
\ \ \ \ \
+
\sum_{k_1,k_2,k_3,k_4}\,
\left[
-
\delta_{i_1,\,i_2,\,i_3}^{k_1,k_2,k_3}\,\mathcal{X}_y^{k_4}
-
\delta_{i_1,\,i_2,\,i_3}^{k_2,k_3,k_1}\,\mathcal{X}_y^{k_4}
-
\delta_{i_1,\,i_2,\,i_3}^{k_3,k_1,k_2}\,\mathcal{X}_y^{k_4}
\right]
y_{k_1,k_2}y_{k_3,k_4}
+ \\
\endaligned
$$
$$
\small
\aligned
& \
\ \ \ \ \
+
\sum_{k_1,k_2,k_3}\,
\left[
\delta_{i_1,\,i_2,\,i_3}^{k_1,k_2,k_3}\,\mathcal{Y}_y
-
\delta_{i_1,\,i_2}^{k_1,k_2}\,\mathcal{X}_{x^{i_3}}^{k_3}
-
\delta_{i_1,\,i_3}^{k_1,k_2}\,\mathcal{X}_{x^{i_2}}^{k_3}
-
\delta_{i_2,\,i_3}^{k_1,k_2}\,\mathcal{X}_{x^{i_1}}^{k_3}
\right]
y_{k_1,k_2,k_3}
+ \\
& \
\ \ \ \ \
+
\sum_{k_1,k_2,k_3,k_4}\,
\left[
-
\delta_{i_1,\,i_2,\,i_3}^{k_1,k_2,k_3}\,\mathcal{X}_y^{k_4}
-
\delta_{i_1,\,i_2,\,i_3}^{k_4,k_1,k_2}\,\mathcal{X}_y^{k_3}
-
\delta_{i_1,\,i_2,\,i_3}^{k_3,k_4,k_1}\,\mathcal{X}_y^{k_2}
-
\delta_{i_1,\,i_2,\,i_3}^{k_2,k_3,k_4}\,\mathcal{X}_y^{k_1}
\right]
y_{k_1}y_{k_2,k_3,k_4}.
\endaligned
$$

\subsection*{3.25.~Comments, analysis and induction}
First of all, by comparing this expression of ${\bf Y}_{i_1, i_2,
i_3}$ with the expression~\thetag{ 2.8} of ${\bf Y}_3$, we easily
guess a part of the (inductional) dictionary beween the cases $n=1$
and the case $n \geqslant 1$. For instance, the three monomials $[\cdot ]
(y_1)^3$, $[\cdot]\, y_1 y_2$ and $[\cdot]\, (y_1 )^2\, y_2$ in ${\bf
Y }_3$ are replaced in ${\bf Y}_{ i_1, i_2, i_3}$ by the following
three sums:
\def\theequation{3.26}\begin{equation}
\small
\sum_{k_1,k_2,k_3}\,
\left[
\cdot
\right]
y_{k_1}y_{k_2}y_{k_3}, 
\ \ \ \ \ \ \ \ \ \ 
\sum_{k_1,k_2,k_3}\,
\left[
\cdot
\right]
y_{k_1}y_{k_2,k_3}, 
\ \ \ \ \ 
{\rm and}
\ \ \ \ \
\sum_{k_1,k_2,k_3,k_4}\,
\left[
\cdot
\right]
y_{k_1}y_{k_2}y_{k_3,k_4}. 
\end{equation}
Similar formal correspondences may be observed for all the monomials
of ${\bf Y}_1$, ${\bf Y}_{i_1}$, of ${\bf Y}_2$, ${\bf Y}_{i_1,i_2}$
and of ${\bf Y}_3$, ${\bf Y}_{i_1,i_2,i_3}$. Generally and
inductively speaking, the monomial
\def\theequation{3.27}\begin{equation}
\left[
\cdot
\right]
\left(
y_{\lambda_1}
\right)^{\mu_1}
\cdots
\left(
y_{\lambda_d}
\right)^{\mu_d}
\end{equation}
appearing in the expression~\thetag{ 2.25} of ${\bf Y }_\kappa$ should
be replaced by a certain multiple sum generalizing~\thetag{ 3.26}.
However, it is necessary to think, to pause and to search for an
appropriate formalism before writing down the desired multiple sum.

The jet variable $y_{ \lambda_1}$ should be replaced by a jet variable
corresponding to a $\lambda_1$-th partial derivative, say $y_{ k_1,
\dots, k_{ \lambda_1}}$, where $k_1, \dots,k_{ \lambda_1}= 1, \dots,
n$. For the moment, to simplify the discussion, we leave out the
presence of a sum of the form $\sum_{ k_1, \dots, k_{ \lambda_1}}$.
The $\mu_1$-th power $\left( y_{ \lambda_1 } \right)^{ \mu_1}$ should
be replaced {\it not}\, by $\left( y_{k_1, \dots, k_{ \lambda_1}}
\right)^{\mu_1}$, but by a product of $\mu_1$ different jet variables
$y_{k_1, \dots, k_{ \lambda_1}}$ of length $\lambda$, {\it with all
indices $k_\alpha = 1, \dots, n$ being distinct}. This rule may be
confirmed by inspecting the expressions of ${\bf Y}_{i_1}$, of ${\bf
Y}_{ i_1, i_2}$ and of ${\bf Y}_{i_1, i_2, i_3}$. So $y_{ k_1, \dots,
k_{ \lambda_1}}$ should be developed as a product of the form
\def\theequation{3.28}\begin{equation}
y_{k_1,\dots,k_{\lambda_1}}\,
y_{k_{\lambda_1+1},\dots,k_{2\lambda_1}}
\cdots \
y_{k_{(\mu_1-1)\lambda_1+1},\dots,k_{\mu_1\lambda_1}},
\end{equation}
where 
\def\theequation{3.29}\begin{equation}
k_1,\dots,k_{\lambda_1},\dots,k_{\mu_1\lambda_1}
=
1,\dots,n.
\end{equation}
Consider now the product $\left( y_{\lambda_1 } \right)^{ \mu_1}\left(
y_{\lambda_2 } \right)^{ \mu_2}$. How should it develope in the case
of several independent variables? For instance, in the expression of
${\bf Y}_{i_1,i_2,i_3}$, we have developed the product $(y_1)^2\,y_2$
as $y_{k_1} y_{k_2} y_{k_3,k_4}$. Thus, a reasonable proposal of
formalism would be that the product $\left( y_{\lambda_1 } \right)^{
\mu_1}\left( y_{\lambda_2 } \right)^{ \mu_2}$ should be developed as a
product of the form
\def\theequation{3.30}\begin{equation}
\aligned
&
y_{k_1,\dots,k_{\lambda_1}}\,
y_{k_{\lambda_1+1},\dots,k_{2\lambda_1}}
\cdots \
y_{k_{(\mu_1-1)\lambda_1+1},\dots,k_{\mu_1\lambda_1}} \\
&
y_{k_{\mu_1\lambda_1+1},\dots,k_{\mu_1\lambda_1+\lambda_2}}
\cdots \
y_{k_{\mu_1\lambda_1+(\mu_2-1)\lambda_2+1},\dots,
k_{\mu_1\lambda_1+\mu_2\lambda_2}},
\endaligned
\end{equation}
where
\def\theequation{3.31}\begin{equation}
k_1,\dots, k_{\lambda_1},\dots,k_{\mu_1\lambda_1},\dots,
k_{\mu_1\lambda_1+\mu_2\lambda_2}
=
1,\dots,n.
\end{equation}
However, when trying to write down the development of the general
monomial $\left( y_{\lambda_1 } \right)^{ \mu_1}\left( y_{\lambda_2 }
\right)^{ \mu_2} \cdots \left( y_{\lambda_d } \right)^{ \mu_d}$, we
would obtain the complicated product
\def\theequation{3.32}\begin{equation}
\small
\aligned
&
y_{k_1,\dots,k_{\lambda_1}}\,
y_{k_{\lambda_1+1},\dots,k_{2\lambda_1}}
\cdots \
y_{k_{(\mu_1-1)\lambda_1+1},\dots,k_{\mu_1\lambda_1}} 
\\
&
y_{k_{\mu_1\lambda_1+1},\dots,k_{\mu_1\lambda_1+\lambda_2}}
\dots
y_{k_{\mu_1\lambda_1+(\mu_2-1)\lambda_2+1},\dots,
k_{\mu_1\lambda_1+\mu_2\lambda_2}}
\\
&
\dots\dots\dots\dots\dots\dots
\dots\dots\dots\dots\dots\dots
\dots\dots\dots\dots\dots\dots
\\
&
y_{k_{\mu_1\lambda_1+\cdots+\mu_{d-1}\lambda_{d-1}+1},
\dots,k_{\mu_1\lambda_1+\cdots+\mu_{d-1}\lambda_{d-1}+\lambda_d}}
\cdots 
\\
& \
\ \ \ \ \ \ \ \ \ \ \ \ \ \ \ 
\cdots \
y_{k_{\mu_1\lambda_1+\cdots+\mu_{d-1}\lambda_{d-1}+
(\mu_d-1)\lambda_d+1},\dots,
k_{\mu_1\lambda_1+\cdots+\mu_d\lambda_d}}.
\endaligned
\end{equation} 
Essentially, this product is still readable. However, in it, some of
the integers $k_\alpha$ have a too long index $\alpha$, often
involving a sum. Such a length of $\alpha$ would be very inconvenient
in writing down and in reading the general Kronecker symbols $\delta_{
i_1, \dots \dots, i_\lambda }^{ k_{ \alpha_1}, \dots, k_{
\alpha_\lambda}}$ which should appear in the final expression of ${\bf
Y}_{ i_1, \dots, i_\kappa}$. One should read in advance Theorem~3.73
below to observe the presence of such multiple Kronecker symbols.
{\sf Consequently, for $\alpha = 1, \dots, \mu_1 \lambda_1, \dots,
\mu_1 \lambda_1 + \cdots + \mu_d \lambda_d$, we have to denote the
indices $k_\alpha$ differently}.

\def\thenotationalconvention{3.33}\begin{notationalconvention}
We denote $d$ collection of $\mu_d$ groups of $\lambda_d$ {\rm (a
priori distinct)} integers $k_\alpha = 1, \dots, n$ by
\def\theequation{3.34}\begin{equation}
\aligned
&
\underbrace{
\underbrace{
k_{1:1:1},\dots,k_{1:1:\lambda_1}}_{\lambda_1},
\dots,
\underbrace{
k_{1:\mu_1:1},\dots,k_{1:\mu_1:\lambda_1}}_{\lambda_1}}_{
\mu_1}, 
\\
&
\underbrace{
\underbrace{
k_{2:1:1},\dots,k_{2:1:\lambda_2}}_{\lambda_2},
\dots,
\underbrace{
k_{2:\mu_2:1},\dots,k_{2:\mu_2:\lambda_2}}_{\lambda_2}}_{
\mu_2},
\\
&
\text{\bf \ \
\dots\dots\dots\dots\dots\dots\dots\dots\dots\dots
\dots\dots\dots\dots\dots
}
\\
& 
\underbrace{
\underbrace{
k_{d:1:1},\dots,k_{d:1:\lambda_d}}_{\lambda_d},
\dots,
\underbrace{
k_{d:\mu_d:1},\dots,k_{d:\mu_d:\lambda_d}}_{\lambda_d}}_{
\mu_d}.
\endaligned
\end{equation}
Correspondingly, we identify the set 
\def\theequation{3.35}\begin{equation}
\small
\left\{
1,\dots,\lambda_1,\dots,\mu_1\lambda_1,
\dots\dots,
\mu_1\lambda_1+\mu_2\lambda_2,
\dots\dots,
\mu_1\lambda_1+\mu_2\lambda_2
+
\cdots
+
\mu_d\lambda_d
\right\}
\end{equation}
of all integers $\alpha$ from $1$ to $\mu_1 \lambda_1 + \mu_2
\lambda_2 + \cdots + \mu_d \lambda_d$ with the following specific set
\def\theequation{3.36}\begin{equation}
\small
\{
\underbrace{
\underbrace{
\underbrace{
\underbrace{
1\!\!:\!\!1\!\!:\!\!1,
\dots,
1\!\!:\!\!1\!\!:\!\!\lambda_1}_{\lambda_1},
\dots,
1\!\!:\!\!\mu_1\!\!:\!\!\lambda_1}_{\mu_1\lambda_1},
\dots,
2\!:\!\mu_2\!:\!\lambda_2}_{\mu_1\lambda_1+\mu_2\lambda_2},
\dots,
d\!:\!\mu_d\!:\!\lambda_d}_{\mu_1\lambda_1+\mu_2\lambda_2
+\cdots+\mu_d\lambda_d}
\},
\end{equation}
written in a lexicographic way which emphasizes clearly the
subdivision in $d$ collections of $\mu_d$ groups of $\lambda_d$
integers.
\end{notationalconvention}

With this notation at hand, we see that the development, in
several independent variables, of the general monomial $\left(
y_{\lambda_1 } \right)^{ \mu_1}
\cdots \left( y_{ \lambda_d } \right)^{ \mu_d }$, may be written as
follows:
\def\theequation{3.37}\begin{equation}
\aligned
y_{k_{1:1:1},\dots,k_{1:1:\lambda_1}}
\cdots \
y_{k_{1:\mu_1:1},\dots,k_{1:\mu_1:\lambda_1}}
\cdots \
y_{k_{d:1:1},\dots,k_{d:1:\lambda_d}}
\cdots\cdots\,
y_{k_{d:\mu_d:1},\dots,k_{d:\mu_d:\lambda_d}}.
\endaligned
\end{equation} 
Formally speaking, this expression is better than~\thetag{ 3.32}. Using
product symbols, we may even write it under the 
slightly more compact form
\def\theequation{3.38}\begin{equation}
\prod_{1\leqslant\nu_1\leqslant\mu_1}\,
y_{k_{1:\nu_1:1},\dots,k_{1:\nu_1:\lambda_1}}
\cdots \
\prod_{1\leqslant\nu_d\leqslant\mu_d}\,
y_{k_{d:\nu_d:1},\dots,k_{d:\nu_d:\lambda_d}}.
\end{equation}

Now that we have translated the monomial, we may add all the summation
symbols: the general expression of ${\bf Y}_\kappa$ (which generalizes
our three previous examples~\thetag{ 3.26}) will be of the form:
\def\theequation{3.39}\begin{equation}
\aligned
{\bf Y}_\kappa
&
=
\mathcal{Y}_{x^{i_1}\cdots x^{i_\kappa}}
+
\sum_{d=1}^{\kappa+1}
\ \
\sum_{1\leqslant\lambda_1<\cdots<\lambda_d\leqslant\kappa}
\ \
\sum_{\mu_1\geqslant 1,\dots,\mu_d\geqslant 1} 
\
\sum_{
\mu_1\lambda_1
+
\cdots
+
\mu_d\lambda_d\leqslant \kappa+1} 
\\
&
\sum_{k_{1:1:1},\dots,k_{1:1:\lambda_1}=1}^n
\cdots \
\sum_{k_{1:\mu_1:1},\dots,k_{1:\mu_1:\lambda_1}=1}^n
\cdots\cdots \
\sum_{k_{d:1:1},\dots,k_{d:1:\lambda_d}=1}^n
\cdots \
\sum_{k_{d:\mu_d:1},\dots,k_{d:\mu_d:\lambda_d}=1}^n
\\
&
\text{\bf[?]}
\prod_{1\leqslant\nu_1\leqslant\mu_1}\,
y_{k_{1:\nu_1:1},\dots,k_{1:\nu_1:\lambda_1}}
\ \cdots \
\prod_{1\leqslant\nu_d\leqslant\mu_d}\,
y_{k_{d:\nu_d:1},\dots,k_{d:\nu_d:\lambda_d}}.
\endaligned
\end{equation}
From now on, 
up to the end of the article, to be very precise, we will restitute
the bounds $\sum_{ k = 1 }^n$ of all the previously abbreviated sums
$\sum_k$. This is justified by the fact that, since we shall deal in
Section~5 below simultaneously with several independent variables
$(x^1, \dots, x^n)$ and with several dependent variables $(y^1, \dots,
y^m)$, we shall encounter sums $\sum_{ l = 1 }^m$, not to be confused
with sums $\sum_{ k = 1 }^n$.

\subsection*{3.40.~Combinatorics of the Kronecker symbols} 
Our next task is to determine what appears inside the brackets {\bf
[?]} of the above equation. We will treat this rather delicate
question very progressively. Inductively, we have to guess how we may
pass from the bracketed term of~\thetag{ 2.25}, namely from
\def\theequation{3.41}\begin{equation}
\aligned
&
\left[
\frac{\kappa\cdots(\kappa-\mu_1\lambda_1-\cdots-\mu_d\lambda_d+1)}
{(\lambda_1!)^{\mu_1}\,\mu_1!
\cdots
(\lambda_d!)^{\mu_d}\,\mu_d!
}
\cdot
\mathcal{Y}_{
x^{\kappa-\mu_1\lambda_1-\cdots-\mu_d\lambda_d}
\,
y^{\mu_1+\cdots+\mu_d}
}
-
\right. \\
& \
\ \ \ \ \ \ \ \ \ \
\left.
-
\frac{\kappa\cdots( 
\kappa-\mu_1\lambda_1-\cdots-\mu_d\lambda_d+2)
(\mu_1\lambda_1+\cdots+\mu_d\lambda_d)}
{(\lambda_1!)^{\mu_1}\,\mu_1!
\cdots
(\lambda_d!)^{\mu_d}\,\mu_d!
}
\cdot
\right.
\\
& \
\ \ \ \ \ \ \ \ \ \ \ \ \ \ \ \ \
\cdot
\mathcal{X}_{
x^{\kappa-\mu_1\lambda_1-\cdots-\mu_d\lambda_d+1}
\,
y^{\mu_1+\cdots+\mu_d-1}
}
\Big],
\endaligned
\end{equation}
to the corresponding (still unknown) bracketed term 
{\bf [?]}.

First of all, we examine the following term, extracted from the
complete expression of ${\bf Y}_{ i_1, i_2, i_3}$ (first line
of~\thetag{ 3.24}):
\def\theequation{3.42}\begin{equation}
\sum_{k_1=1}^n\,
\left[
\delta_{i_1}^{k_1}\,\mathcal{Y}_{x^{i_2}x^{i_3}y}
+
\delta_{i_2}^{k_1}\,\mathcal{Y}_{x^{i_1}x^{i_3}y}
+
\delta_{i_3}^{k_1}\,\mathcal{Y}_{x^{i_1}x^{i_2}y}
-
\mathcal{ X}_{x^{i_1}x^{i_2}x^{i_3}}^{k_1}
\right]
y_{k_1}.
\end{equation}
Here, the coefficient $\left[ 3\, \mathcal{ Y}_{ x^2 y} - \mathcal{
X}_{ x^3 } \right]$ of the monomial $y_1$ in ${\bf Y}_3$ is replaced
by the above bracketed terms. 

Let us precisely analyze the combinatorics. Here, $\mathcal{ X}_{x^3}$ is
replaced by $\mathcal{ X}_{x^{ i_1 }x^{ i_2}x^{i_3}}^{k_1}$, where the
lower indices $i_1, i_2, i_3$ come from ${\bf Y}_{i_1, i_2, i_3}$ and
where the upper index $k_1$ is the summation index. Also, the integer
$3$ in $3\, \mathcal{ Y}_{x^2 y}$ is replaced by a sum of exactly
three terms, each involving a single Kronecker symbol $\delta_i^k$, in
which the lower index is always an index $i= i_1, i_2, i_3$ and in
which the upper index is always equal to the summation index $k_1$.
By the way, more generally, we immediately observe that
all the successive positive integers 
\def\theequation{3.43}\begin{equation}
1, 3, 1, 3, 3, 1, 3, 1, 3, 3, 3, 9, 6, 3, 1, 3, 4
\end{equation}
appearing in the formula~\thetag{ 2.8} for ${\bf Y}_3$ are replaced,
in the formula~\thetag{ 3.24} for ${\bf Y}_{i_1, i_2, i_3}$, by sums
of exactly the same number of terms involving Kronecker
symbols. This observation will be a precious guide. Finally, in the
symbol $\delta_i^{k_1}$, if $i$ is chosen among the set $\{ i_1, i_2,
i_3\}$, for instance if $i = i_1$, it follows that the development of
$\mathcal{ Y}_{x^2y}$ necessarily involves the remaining indices, for
instance $\mathcal{ Y}_{x^{i_2}x^{i_3}y}$. Since there are three
choices for $i = i_1, i_2, i_3$, we recover the number $3$.

Next, comparing $\left[ \mathcal{ Y}_{yy}
-2\,\mathcal{ X}_{ xy} \right] (y_1)^2$ with 
the term
\def\theequation{3.44}\begin{equation}
\sum_{k_1,k_2=1}^n\,
\left[
\delta_{i_1,\,i_2}^{k_1,k_2}\,\mathcal{Y}_{yy}
-
\delta_{i_1}^{k_1}\,\mathcal{X}_{x^{i_2}y}^{k_1}
-
\delta_{i_2}^{k_1}\,\mathcal{X}_{x^{i_1}y}^{k_1}
\right]
y_{k_1}y_{k_2},
\end{equation}
extracted from the complete expression of ${\bf Y}_{i_1,i_2}$ (second
line of~\thetag{ 3.18}), we learn and we guess that the number of
Kronecker symbols before $\mathcal{ Y}_{x^\gamma y^\delta}$ must be
equal to the number of indices $k_\alpha$ minus $\gamma$. This rule
is confirmed by examining the term (second and third line of~\thetag{
3.24})
\def\theequation{3.45}\begin{equation}
\aligned
\sum_{k_1,k_2}
&
\left[
\delta_{i_1,\,i_2}^{k_1,k_2}\,\mathcal{Y}_{x^{i_3}y^2}
+
\delta_{i_3,\,i_1}^{k_1,k_2}\,\mathcal{Y}_{x^{i_2}y^2}
+
\delta_{i_2,\,i_3}^{k_1,k_2}\,\mathcal{Y}_{x^{i_1}y^2}
-
\right.
\\
& \
\ \ \ \ \
\left.
-
\delta_{i_1}^{k_1}\,\mathcal{X}_{x^{i_2}x^{i_3}y}^{k_2}
-
\delta_{i_2}^{k_1}\,\mathcal{X}_{x^{i_1}x^{i_3}y}^{k_2}
-
\delta_{i_3}^{k_1}\,\mathcal{X}_{x^{i_1}x^{i_2}y}^{k_2}
\right]
y_{k_1}y_{k_2},
\endaligned
\end{equation}
developing $\left[ 3\,\mathcal{ Y}_{xy^2} - 3\, \mathcal{
X}_{x^2y} \right] (y_1)^2$.

Also, we may examine the following term
\def\theequation{3.46}\begin{equation}
\aligned
\sum_{k_1,k_2=1}^n\,
&
\left[
\delta_{i_1,\,i_2}^{k_1,k_2}\,\mathcal{Y}_{x^{i_3}x^{i_4}y^2}
+
\delta_{i_1,\,i_3}^{k_1,k_2}\,\mathcal{Y}_{x^{i_2}x^{i_4}y^2}
+
\delta_{i_1,\,i_4}^{k_1,k_2}\,\mathcal{Y}_{x^{i_2}x^{i_3}y^2}
+ 
\right.
\\
& \
\ \ \ \ \
\left.
+
\delta_{i_2,\,i_3}^{k_1,k_2}\,\mathcal{Y}_{x^{i_1}x^{i_4}y^2}
+
\delta_{i_2,\,i_4}^{k_1,k_2}\,\mathcal{Y}_{x^{i_1}x^{i_3}y^2}
+
\delta_{i_3,\,i_4}^{k_1,k_2}\,\mathcal{Y}_{x^{i_1}x^{i_2}y^2}
-
\right.
\\
& \
\ \ \ \ \
\left.
-
\delta_{i_1}^{k_1}\,\mathcal{X}_{x^{i_2}x^{i_3}x^{i_4}y}^{k_1}
-
\delta_{i_2}^{k_1}\,\mathcal{X}_{x^{i_1}x^{i_2}x^{i_3}y}^{k_1}
-
\delta_{i_3}^{k_1}\,\mathcal{X}_{x^{i_1}x^{i_2}x^{i_4}y}^{k_1}
-
\right.
\\
& \
\ \ \ \ \
\left.
-
\delta_{i_4}^{k_1}\,\mathcal{X}_{x^{i_1}x^{i_2}x^{i_3}y}^{k_1}
\right]
y_{k_1}y_{k_2},
\endaligned
\end{equation} 
extracted from ${\bf Y}_{ i_1, i_2, i_3, i_4}$ and developing
$\left[ 6\, \mathcal{ Y}_{ x^2 y^2} - 4\, \mathcal{ X}_{ x^3 y}
\right] (y_1)^2$. We would like to mention that we have not written
the complete expression of ${\bf Y}_{ i_1, i_2, i_3, i_4}$, because it
would cover two and a half printed pages. 

By inspecting the way how the indices are permuted in the multiple
Kronecker symbols of the first two lines of this expression~\thetag{
3.46}, we observe that the six terms correspond exactly to the six
possible choices of two complementary ordered couples of integers in
the set $\{ 1, 2, 3, 4\}$, namely
\def\theequation{3.47}\begin{equation}
\aligned
&
\{1,2\}\cup\{3,4\}, 
\ \ \ \ \ \ \ \
\{1,3\}\cup\{2,4\}, 
\ \ \ \ \ \ \ \
\{1,4\}\cup\{2,3\}, 
\\
&
\{2,3\}\cup\{1,4\}, 
\ \ \ \ \ \ \ \
\{2,4\}\cup\{1,3\}, 
\ \ \ \ \ \ \ \
\{3,4\}\cup\{1,2\}.
\endaligned
\end{equation}
At this point, we start to devise the general combinatorics. Before
proceeding further, we need some notation.

\subsection*{ 3.48.~Permutation groups}
For every $p \in \N$ with $p \geqslant 1$, we denote by $\mathfrak{
S}_p$ the full permutation group of the set $\{ 1, 2, \dots, p-1,
p\}$. Its cardinal equals $p!$. The letters $\sigma$ and $\tau$ will
be used to denote an element of $\mathfrak{ S}_p$. If $p \geqslant 2$,
and if $q \in \N$ satisfies $1\leqslant q \leqslant p-1$, we denote by
$\mathfrak{ S}_p^q$ the subset of permutations $\sigma \in \mathfrak{
S}_p$ satisfying the two collections of inequalities
\def\theequation{3.49}\begin{equation}
\sigma(1)<\sigma(2)<\cdots<\sigma(q)
\ \ \ \ \ \ \ \ \ 
{\rm and}
\ \ \ \ \ \ \ \ \ 
\sigma(q+1)<\sigma(q+2)<\cdots<\sigma(p).
\end{equation}
The cardinal of $\mathfrak{ S}_p^q$ 
equals $C_p^q = \frac{ p!}{ q! \ (p-q) !}$.

\def\thelemma{3.50}\begin{lemma}
For $\kappa \geqslant 1$, the development of~\thetag{ 2.20} to several
independent variables $(x^1, \dots, x^n)$ is{\rm :}
\def\theequation{3.51}\begin{equation}
\small
\aligned
&
{\bf Y}_{i_1,i_2,\dots,i_\kappa}
=
\mathcal{Y}_{x^{i_1}x^{i_2}\cdots x^{i_\kappa}}
+
\sum_{k_1=1}^n\,
\left[
\sum_{\tau\in\mathfrak{S}_\kappa^1}\,
\delta_{i_{\tau(1)}}^{k_1}\,\mathcal{Y}_{x^{i_{\tau(2)}}
\cdots x^{i_{\tau(\kappa)}}y}
-
\mathcal{X}_{x^{i_1}x^{i_2}\cdots x^{i_{\kappa}}}^{k_1}
\right]
y_{k_1}
+ \\
& \
\ \ \ \ \
+
\sum_{k_1,k_2=1}^n\,
\left[
\sum_{\tau\in\mathfrak{S}_\kappa^2}\,
\delta_{i_{\tau(1)},i_{\tau(2)}}^{k_1,\ \ \ k_2}\,
\mathcal{Y}_{x^{i_{\tau(3)}}\cdots x^{i_{\tau(\kappa)}}y^2}
-
\sum_{\tau\in\mathfrak{S}_\kappa^1}\,
\delta_{i_{\tau(1)}}^{k_1}\,
\mathcal{X}_{x^{i_{\tau(2)}}\cdots x^{i_{\tau(\kappa)}}y}^{k_2}
\right]
y_{k_1}y_{k_2}
+
\\
& \
\ \ \ \ \
+
\sum_{k_1, k_2, k_3=1}^n
\left[
\sum_{\tau\in\mathfrak{S}_\kappa^3}\,
\delta_{i_{\tau(1)},i_{\tau(2)},i_{\tau(3)}}^{
k_1,\ \ \ k_2,\ \ \ k_3}\,
\mathcal{Y}_{x^{i_{\tau(4)}}\cdots x^{i_{\tau(\kappa)}}y^3}
-
\right.
\\
& \
\ \ \ \ \ \ \ \ \ \ \ \ \ \ \ 
\ \ \ \ \ \ \ \ \ \ \ \ \ \ \ 
\ \ \ \ \ \ 
\left.
-
\sum_{\tau\in\mathfrak{S}_\kappa^2}\,
\delta_{i_{\tau(1),i_{\tau(2)}}}^{k_1,\ \ \ k_2}
\mathcal{X}_{x^{i_{\tau(3)}}\cdots x^{i_{\tau(\kappa)}}y^2}^{k_3}
\right]
y_{k_1}y_{k_2}y_{k_3}
+ \\
& \
\ \ \ \ \
+
\cdots\cdots
+
\\
& \
\ \ \ \ \
+
\sum_{k_1,\dots,k_\kappa=1}^n
\left[
\delta_{i_1,\dots,\ i_\kappa}^{k_1,\dots,k_\kappa}\,
\mathcal{Y}_{y^\kappa}
-
\sum_{\tau\in\mathfrak{S}_\kappa^{\kappa-1}}\,
\delta_{i_{\tau(1)},\dots,i_{\tau(\kappa-1)}}^{
k_1,\dots\dots,k_{\kappa-1}}\,
\mathcal{X}_{x^{i_{\tau(\kappa)}}y^{\kappa-1}}^{k_\kappa}
\right]
y_{k_1}\cdots y_{k_\kappa}
+ \\
& \
\ \ \ \ \
+
\sum_{k_1,\dots,k_\kappa,k_{\kappa+1}=1}^n
\left[
-
\delta_{i_1,\dots,\ i_\kappa}^{k_1,\dots,k_\kappa}\,
\mathcal{X}_{y^\kappa}^{k_{\kappa+1}}
\right]
y_{k_1}\cdots y_{k_\kappa} y_{k_{\kappa+1}}
+
{\sf remainder}.
\endaligned
\end{equation}
Here, the term {\sf remainder} collects all
remaining monomials in the pure jet variables
$y_{ k_1, \dots, k_\lambda }$.
\end{lemma}

\subsection*{3.52.~Continuation}
Thus, we have devised how the part of ${\bf Y}_{i_1,\dots, i_\kappa}$
which involves only the jet variables $y_{k_\alpha}$ must be written.
To proceed further, we shall examine the following term, extracted
from ${\bf Y}_{i_1,i_2,i_3}$ (lines 12 and 13 of~\thetag{ 3.24})
\def\theequation{3.53}\begin{equation}
\aligned
& \
\ \ \ \ \
\sum_{k_1,k_2,k_3,k_4}\,
\left[
-
\delta_{i_1,\,i_2,\,i_3}^{k_1,k_2,k_3}\,\mathcal{X}_{y^2}^{k_4}
-
\delta_{i_1,\,i_2,\,i_3}^{k_2,k_3,k_1}\,\mathcal{X}_{y^2}^{k_4}
-
\delta_{i_1,\,i_2,\,i_3}^{k_3,k_2,k_1}\,\mathcal{X}_{y^2}^{k_4}
- 
\right.
\\
& \
\ \ \ \ \ \ \ \ \ \ \ \ \ \ \ \ \ \ \ \ \ \ \ \ \ \ \
\left.
-
\delta_{i_1,\,i_2,\,i_3}^{k_3,k_4,k_1}\,\mathcal{X}_{y^2}^{k_2}
-
\delta_{i_1,\,i_2,\,i_3}^{k_3,k_1,k_4}\,\mathcal{X}_{y^2}^{k_2}
-
\delta_{i_1,\,i_2,\,i_3}^{k_1,k_3,k_4}\,\mathcal{X}_{y^2}^{k_2}
\right]
y_{k_1}y_{k_2}y_{k_3,k_4},
\endaligned
\end{equation}
which developes the term $\left[ - 6\,\mathcal{ X}_{ y^2} \right]
(y_1)^2 y_2$ of ${\bf Y}_3$ (third line of~\thetag{ 2.8}). During the
computation which led us to the final expression~\thetag{ 3.24}, we
organized the formula in order that, in the six Kronecker symbols, the
lower indices $i_1,i_2,i_3$ are all written in the same order. But
then, {\it what is the rule for the appearance of the four upper
indices $k_1, k_2, k_3, k_4$}?

In April 2001, we discovered the rule by inspecting both~\thetag{
3.53} and the following complicated term, extracted from the complete
expression of ${\bf Y}_{i_1,i_2,i_3,i_4}$ written in one of our
manuscripts:
\def\theequation{3.54}\begin{equation}
\aligned
\sum_{k_1,k_2,k_3}\,
&
\left[
\delta_{i_1, \ i_2, \ i_3}^{k_1,k_2,k_3}\,\mathcal{Y}_{x^{i_4}y^2}
+
\delta_{i_1, \ i_2, \ i_3}^{k_2,k_1,k_3}\,\mathcal{Y}_{x^{i_4}y^2}
+
\delta_{i_1, \ i_2, \ i_3}^{k_2,k_3,k_1}\,\mathcal{Y}_{x^{i_4}y^2}
+
\right. 
\\
& \
\ \ 
+
\delta_{i_1, \ i_2, \ i_4}^{k_1,k_2,k_3}\,\mathcal{Y}_{x^{i_3}y^2}
+
\delta_{i_1, \ i_2, \ i_4}^{k_2,k_1,k_3}\,\mathcal{Y}_{x^{i_3}y^2}
+
\delta_{i_1, \ i_2, \ i_4}^{k_2,k_3,k_1}\,\mathcal{Y}_{x^{i_3}y^2}
+ \\
& \
\ \ 
+
\delta_{i_1, \ i_3, \ i_4}^{k_1,k_2,k_3}\,\mathcal{Y}_{x^{i_2}y^2}
+
\delta_{i_1, \ i_3, \ i_4}^{k_2,k_1,k_3}\,\mathcal{Y}_{x^{i_2}y^2}
+
\delta_{i_1, \ i_3, \ i_4}^{k_2,k_3,k_1}\,\mathcal{Y}_{x^{i_2}y^2}
+ \\
& \
\ \ 
+
\delta_{i_2, \ i_3, \ i_4}^{k_1,k_2,k_3}\,\mathcal{Y}_{x^{i_1}y^2}
+
\delta_{i_2, \ i_3, \ i_4}^{k_2,k_1,k_3}\,\mathcal{Y}_{x^{i_1}y^2}
+
\delta_{i_2, \ i_3, \ i_4}^{k_2,k_3,k_1}\,\mathcal{Y}_{x^{i_1}y^2}
- \\
& \
\ \ 
-
\delta_{i_1,\ i_2}^{k_1,k_2}\,\mathcal{X}_{x^{i_3}x^{i_4}y}^{k_3}
-
\delta_{i_1,\ i_2}^{k_2,k_1}\,\mathcal{X}_{x^{i_3}x^{i_4}y}^{k_3}
-
\delta_{i_1,\ i_2}^{k_2,k_3}\,\mathcal{X}_{x^{i_3}x^{i_4}y}^{k_1}
- \\
\endaligned
\end{equation}
\[
\aligned
& \
\ \ 
-
\delta_{i_1,\ i_3}^{k_1,k_2}\,\mathcal{X}_{x^{i_2}x^{i_4}y}^{k_3}
-
\delta_{i_1,\ i_3}^{k_2,k_1}\,\mathcal{X}_{x^{i_2}x^{i_4}y}^{k_3}
-
\delta_{i_1,\ i_3}^{k_2,k_3}\,\mathcal{X}_{x^{i_2}x^{i_4}y}^{k_1}
- \\
& \
\ \ 
-
\delta_{i_1,\ i_4}^{k_1,k_2}\,\mathcal{X}_{x^{i_2}x^{i_3}y}^{k_3}
-
\delta_{i_1,\ i_4}^{k_2,k_1}\,\mathcal{X}_{x^{i_2}x^{i_3}y}^{k_3}
-
\delta_{i_1,\ i_4}^{k_2,k_3}\,\mathcal{X}_{x^{i_2}x^{i_3}y}^{k_1}
- \\
& \
\ \ 
-
\delta_{i_2,\ i_3}^{k_1,k_2}\,\mathcal{X}_{x^{i_1}x^{i_4}y}^{k_3}
-
\delta_{i_2,\ i_3}^{k_2,k_1}\,\mathcal{X}_{x^{i_1}x^{i_4}y}^{k_3}
-
\delta_{i_2,\ i_3}^{k_2,k_3}\,\mathcal{X}_{x^{i_1}x^{i_4}y}^{k_1}
- \\
& \
\ \ 
-
\delta_{i_2,\ i_4}^{k_1,k_2}\,\mathcal{X}_{x^{i_1}x^{i_3}y}^{k_3}
-
\delta_{i_2,\ i_4}^{k_2,k_1}\,\mathcal{X}_{x^{i_1}x^{i_3}y}^{k_3}
-
\delta_{i_2,\ i_4}^{k_2,k_3}\,\mathcal{X}_{x^{i_1}x^{i_3}y}^{k_1}
- \\
& \
\ \
\left. 
-
\delta_{i_3,\ i_4}^{k_1,k_2}\,\mathcal{X}_{x^{i_1}x^{i_2}y}^{k_3}
-
\delta_{i_3,\ i_4}^{k_2,k_1}\,\mathcal{X}_{x^{i_1}x^{i_2}y}^{k_3}
-
\delta_{i_3,\ i_4}^{k_2,k_3}\,\mathcal{X}_{x^{i_1}x^{i_2}y}^{k_1}
\right]
y_{k_1}y_{k_2,k_3}.
\endaligned
\]
This sum developes the term $\left[ 12\, \mathcal{ Y}_{xy^2} -
18\,\mathcal{ X}_{ x^2 y} \right]y_1y_2$ of
${\bf Y}_3$ (third line of~\thetag{
2.9}). Let us explain what are the formal rules.

In the bracketed terms of~\thetag{ 3.53}, there are no permutation of
the indices $i_1,i_2,i_3$, but there is a certain unknown subset of
all the permutations of the four indices $k_1,k_2,k_3,k_4$. In the
bracketed terms of~\thetag{ 3.54}, two combinatorics are present:

\smallskip

\begin{itemize}
\item[$\bullet$]
there are some permutations of the indices $i_1,i_2,i_3,i_4$ and
\item[$\bullet$]
there are some permutations of the indices $k_1,k_2,k_3$. 
\end{itemize}

\smallskip

Here, the permutations of the indices $i_1,i_2,i_3,i_4$ are easily
guessed, since they are the same as the permutations which were
introduced in \S3.48 above. Indeed, in the first four
lines of~\thetag{ 3.54}, we
see the four decompositions 
\def\theequation{3.55}\begin{equation}
\{i_1,i_2,i_3\}\cup\{i_4\},
\ \ \ \ \ \ \ 
\{i_1,i_2,i_4\}\cup\{i_3\},
\ \ \ \ \ \ \ 
\{i_1,i_3,i_4\}\cup\{i_2\},
\ \ \ \ \ \ \ 
\{i_2,i_3,i_4\}\cup\{i_1\},
\end{equation}
of the set $\{i_1,i_2,i_3,i_4\}$, and in the 
last six lines of~\thetag{ 3.54}, we see the
six decompositions
\def\theequation{3.56}\begin{equation}
\aligned
&
\{i_1,i_2\}\cup\{i_3,i_4\},
\ \ \ \ \ \ \ 
\{i_1,i_3\}\cup\{i_2,i_4\},
\ \ \ \ \ \ \ 
\{i_1,i_4\}\cup\{i_2,i_3\},
\\
&
\{i_2,i_3\}\cup\{i_1,i_4\},
\ \ \ \ \ \ \ 
\{i_2,i_4\}\cup\{i_1,i_3\},
\ \ \ \ \ \ \ 
\{i_3,i_4\}\cup\{i_1,i_2\},
\endaligned
\end{equation}
so that~\thetag{ 3.54} may be written under the form
\def\theequation{3.57}\begin{equation}
\small
\aligned
\sum_{k_1,k_2,k_3}\,
\left[
\sum_{\tau\in\mathfrak{S}_4^3}\,
\sum_{\sigma\in\text{\bf ?}}\,
\delta_{i_{\tau(1)},i_{\tau(2)},i_{\tau(3)}}^{
k_{\tau(1)},k_{\tau(2)},k_{\tau(3)}}\,
\mathcal{Y}_{x^{i_{\tau(4)}}y^2}
-
\sum_{\tau\in\mathfrak{S}_4^2}\,
\sum_{\sigma\in\text{\bf ?}}\,
\delta_{i_{\tau(1)},i_{\tau(2)}}^{
k_{\tau(1)},k_{\tau(2)}}\,
\mathcal{X}_{x^{i_{\tau(3)}}x^{i_{\tau(4)}}y}^{k_{\tau(3)}}
\right]
y_{k_1}y_{k_2,k_3},
\endaligned
\end{equation}
where in the two above sums $\sum_{ \sigma \in \text{\bf ?}}$, the
letter $\sigma$ denotes a permutation of the set $\{1,2,3\}$ and where
the sign {\bf ?} refers to two (still unknown) subset of the full
permutation group $\mathfrak{ S}_3$. {\it The only remaining question
is to determine how the indices $k_\alpha$ are permuted in~\thetag{
3.53} and in~\thetag{ 3.54}}.

The answer may be guessed by looking at the permutations of the set
$\{k_1,k_2,k_3,k_4\}$ which stabilize the monomial
$y_{k_1}y_{k_2}y_{k_3,k_4}$ in~\thetag{ 3.53}: we clearly have the
following four symmetry relations between monomials:
\def\theequation{3.58}\begin{equation}
y_{k_1}y_{k_2}y_{k_3,k_4}
\equiv
y_{k_2}y_{k_1}y_{k_3,k_4}
\equiv
y_{k_1}y_{k_2}y_{k_4,k_3}
\equiv
y_{k_2}y_{k_1}y_{k_4,k_3},
\end{equation}
and nothing more. 
Then the number $6$ of bracketed terms in~\thetag{ 3.53}
is exactly equal to the cardinal $24 = 4!$ of the full permutation
group of the set $\{k_1,k_2,k_3,k_4\}$ divided by the number $4$ of
these symmetry relations. The set of permutations $\sigma$ of
$\{1,2,3,4\}$ satisfying these symmetry relations
\def\theequation{3.59}\begin{equation}
y_{k_{\sigma(1)}}y_{k_{\sigma(2)}}
y_{k_{\sigma(3)},k_{\sigma(4)}}
\equiv
y_{k_1}y_{k_2}y_{k_3,k_4}
\end{equation}
consitutes a subgroup of $\mathfrak{ S}_4$ which we will denote by
$\mathfrak{ H}_4^{(2,1),(1,2)}$. Furthermore, the coset
\def\theequation{3.60}\begin{equation}
\mathfrak{ F}_4^{(2,1),(1,2)} 
:=
\mathfrak{ S}_4 / \mathfrak{ H}_4^{(2,1),(1,2)}
\end{equation}
possesses the six representatives 
\def\theequation{3.61}\begin{equation}
\small
\aligned
\left(
\begin{array}{cccc}
1 & 2 & 3 & 4 \\
1 & 2 & 3 & 4 \\
\end{array}
\right), 
\ \ \ \ \ \ \ \ \ \ 
\left(
\begin{array}{cccc}
1 & 2 & 3 & 4 \\
2 & 3 & 1 & 4 \\
\end{array}
\right), 
\ \ \ \ \ \ \ \ \ \ 
\left(
\begin{array}{cccc}
1 & 2 & 3 & 4 \\
3 & 2 & 1 & 4 \\
\end{array}
\right), 
\\
\left(
\begin{array}{cccc}
1 & 2 & 3 & 4 \\
3 & 4 & 1 & 2 \\
\end{array}
\right), 
\ \ \ \ \ \ \ \ \ \ 
\left(
\begin{array}{cccc}
1 & 2 & 3 & 4 \\
3 & 1 & 4 & 2 \\
\end{array}
\right), 
\ \ \ \ \ \ \ \ \ \ 
\left(
\begin{array}{cccc}
1 & 2 & 3 & 4 \\
1 & 3 & 4 & 2 \\
\end{array}
\right), \\
\endaligned
\end{equation}
which exactly appear as the permutations of the upper indices of our
example~\thetag{ 3.53}. Of course, the question arises whether the
choice of such six representatives in the quotient $\mathfrak{ S}_4 /
\mathfrak{ H}_4^{ (2,1), (1,2)}$ is legitimate.

Fortunately, we observe that after conjugation by any permutation
$\sigma \in \mathfrak{ H }_4^{ (2,1), (1,2)}$, we do not perturb any of
the six terms of~\thetag{ 3.53}, for instance the third term
of~\thetag{ 3.53} is not perturbed, as shown by the following
computation
\def\theequation{3.62}\begin{equation}
\small
\aligned
&
\sum_{k_1,k_2,k_3,k_4}
\left[
-
\delta_{i_1, \ \ \ \ i_2, \ \ \ \ i_3}^{
k_{\sigma(3)}, k_{\sigma(2)}, k_{\sigma(1)}}\,
\mathcal{X}_{y^2}^{k_{\sigma(4)}}
\right]
y_{k_1}y_{k_2}y_{k_3,k_4}
= \\
& \
\ \ \ \ \
=
\sum_{k_1,k_2,k_3,k_4}
\left[
-
\delta_{i_1, \ i_2, \ i_3}^{
k_3, k_2, k_1}\,
\mathcal{X}_{y^2}^{k_{\sigma(4)}}
\right]
y_{k_{\sigma^{-1}(1)}}
y_{k_{\sigma^{-1}(2)}}
y_{k_{\sigma^{-1}(3)},k_{\sigma^{-1}(4)}}
\\
& \
\ \ \ \ \
=
\sum_{k_1,k_2,k_3,k_4}
\left[
-
\delta_{i_1, \ i_2, \ i_3}^{
k_3, k_2, k_1}\,
\mathcal{X}_{y^2}^{k_{\sigma(4)}}
\right]
y_{k_1}y_{k_2}y_{k_3,k_4}
\endaligned
\end{equation} 
thanks to the symmetry~\thetag{ 3.59}. Thus, as expected, the choice of
$6$ arbitrary representatives $\sigma \in \mathfrak{ F}_4^{(2,1), (1,
2)}$ in the bracketed terms of~\thetag{ 3.53} is free.
In conclusion, we have shown that~\thetag{ 3.53} may 
be written under the form:
\def\theequation{3.63}\begin{equation}
\aligned
\sum_{k_1,k_2,k_3,k_4}\,
\left[
-
\sum_{\sigma\in\mathfrak{F}_4^{(2,1),(1,2)}}\,
\delta_{i_1,\ \ \ \ i_2,\ \ \ \ i_3}^{
k_{\sigma(1)},k_{\sigma(2)},k_{\sigma(3)}}\,
\mathcal{X}_{y^2}^{k_{\sigma(4)}}
\right]
y_{k_1}y_{k_2}y_{k_3,k_4},
\endaligned
\end{equation}

This rule is confirmed by inspecting~\thetag{ 3.54} (as
well as all the other terms of ${\bf Y}_{i_1, i_2, i_3}$ and
of ${\bf Y}_{ i_1, i_2, i_3, i_4}$). Indeed, the
permutations $\sigma$ of the set $\{k_1, k_2, k_3\}$ which stabilize the
monomial $y_{k_1} y_{k_2, k_3}$ consist just of the identity
permutation and the transposition of $k_2$ and $k_3$. The coset
$\mathfrak{ S }_3 / \mathfrak{ H }_3^{ (1,1), (1,2)}$ has the three
representatives
\def\theequation{3.64}\begin{equation}
\small
\aligned
\left(
\begin{array}{ccc}
1 & 2 & 3 \\
1 & 2 & 3 \\
\end{array}
\right), 
\ \ \ \ \ \ \ \ \ \ 
\left(
\begin{array}{ccc}
1 & 2 & 3 \\
2 & 1 & 3 \\
\end{array}
\right), 
\ \ \ \ \ \ \ \ \ \ 
\left(
\begin{array}{cccc}
1 & 2 & 3 \\
2 & 3 & 1 \\
\end{array}
\right), 
\\
\endaligned
\end{equation}
which appear in the upper index position of each of the ten lines
of~\thetag{ 3.54}. It follows that~\thetag{ 3.54}
may be written under the form
\def\theequation{3.65}\begin{equation}
\aligned
\sum_{k_1,k_2,k_3}\,
&
\left[
\sum_{\tau\in\mathfrak{S}_4^3}\,
\sum_{\sigma\in\mathfrak{F}_3^{(1,1),(1,2)}}\,
\delta_{i_{\tau(1)},i_{\tau(2)},i_{\tau(3)}}^{
k_{\sigma(1)},k_{\sigma(2)},k_{\sigma(3)}}\,
\mathcal{Y}_{x^{i_{\tau(4)}}y^2}
- 
\right.
\\
& \
\left.
-
\sum_{\sigma\in\mathfrak{S}_4^2}\,
\sum_{\tau\in\mathfrak{F}_3^{(1,1),(1,2)}}\,
\delta_{i_{\tau(1)},i_{\tau(2)}}^{
k_{\sigma(1)},k_{\sigma(2)}}\,
\mathcal{X}_{x^{i_{\tau(3)}}
x^{i_{\tau(4)}}y}^{k_{\sigma(3)}}
\right]
y_{k_1}y_{k_2,k_3}.
\endaligned
\end{equation}

\subsection*{ 3.66.~General complete expression of ${\bf Y}_{i_1, \dots,
i_\kappa}$} As in the incomplete expression~\thetag{ 3.39} of ${\bf
Y}_{i_1, \dots, i_\kappa}$, consider integers $1\leqslant \lambda_1 <
\cdots < \lambda_d \leqslant \kappa$ and $\mu_1\geqslant 1, \dots, \mu_d\geqslant 1$
satisfying $\mu_1 \lambda_1 + \cdots + \mu_d \lambda_d \leqslant \kappa +
1$. By $\mathfrak{ H}_{ \mu_1 \lambda_1 + \cdots + \mathfrak{ H}_{
\mu_d \lambda_d}}$, we denote the subgroup of permutations $\tau \in
\mathfrak{ S}_{\mu_1 \lambda_1 + \cdots + \mathfrak{ H}_{ \mu_d
\lambda_d}}$ that leave unchanged the
general monomial~\thetag{ 3.38}, namely 
that satisfy
\def\theequation{3.67}\begin{equation}
\aligned
&
\prod_{1\leqslant\nu_1\leqslant\mu_1}\,
y_{k_{\sigma(1:\nu_1:1)},\dots,k_{\sigma(1:\nu_1:\lambda_1)}}
\cdots \
\prod_{1\leqslant\nu_d\leqslant\mu_d}\,
y_{k_{\sigma(d:\nu_d:1)},\dots,k_{\sigma(d:\nu_d:\lambda_d)}}
= \\
&
=
\prod_{1\leqslant\nu_1\leqslant\mu_1}\,
y_{k_{1:\nu_1:1},\dots,k_{1:\nu_1:\lambda_1}}
\cdots \
\prod_{1\leqslant\nu_d\leqslant\mu_d}\,
y_{k_{d:\nu_d:1},\dots,k_{d:\nu_d:\lambda_d}}.
\endaligned
\end{equation}
The structure of this group may be described as follows. For every $e
= 1, \dots, d$, an arbitrary permutation $\sigma$ of the set
\def\theequation{3.68}\begin{equation}
\{
\underbrace{
\underbrace{
e\!:1\!:\!1, \dots, e\!:1\!:\!\lambda_e}_{\lambda_e},
\underbrace{
e\!:2\!:\!1, \dots, e\!:2\!:\!\lambda_e}_{\lambda_e},
\cdots,
\underbrace{
e\!:\mu_e\!:\!1, \dots, e\!:\mu_e\!:\!\lambda_e}_{\lambda_e}}_{
\mu_e}
\}
\end{equation}
which leaves unchanged the monomial
\def\theequation{3.69}\begin{equation}
\prod_{1\leqslant\nu_e\leqslant\mu_e}\,
y_{k_{\sigma(e:\nu_e:1)},\dots,k_{\sigma(e:\nu_e:\lambda_e)}}=
\prod_{1\leqslant\nu_e\leqslant\mu_e}\,
y_{k_{e:\nu_e:1},\dots,k_{e:\nu_e:\lambda_e}}.
\end{equation}
uniquely decomposes as the composition of

\smallskip

\begin{itemize}
\item[$\bullet$]
$\mu_e$ arbitrary permutations of the $\mu_e$ groups
of $\lambda_e$ integers 
$\{e\!:\!\nu_e\!:\!1,
\dots,e\!:\!\nu_e\!:\!\lambda_e\}$, 
of total cardinal $(\lambda_e!)^{\mu_e}$;
\item[$\bullet$]
an arbitrary permutation
between these $\mu_e$ groups, of total
cardinal $\mu_e!$.
\end{itemize}

\smallskip

Consequently
\def\theequation{3.70}\begin{equation}
{\rm Card}
\left(
\mathfrak{H}_{\mu_1\lambda_1+\cdots+\mu_d\lambda_d}^{
(\mu_1,\lambda_1),\dots,(\mu_d,\lambda_d)}
\right)
=
\mu_1!(\lambda_1!)^{\mu_1}\cdots
\mu_d!(\lambda_d!)^{\mu_d}.
\end{equation}
Finally, define the coset
\def\theequation{3.71}\begin{equation}
\mathfrak{F}_{\mu_1\lambda_1+\cdots+\mu_d\lambda_d}^{
(\mu_1,\lambda_1),\dots,(\mu_d,\lambda_d)}
:=
\mathfrak{S}_{\mu_1\lambda_1+\cdots+\mu_d\lambda_d}
/
\mathfrak{H}_{\mu_1\lambda_1+\cdots+\mu_d\lambda_d}^{
(\mu_1,\lambda_1),\dots,(\mu_d,\lambda_d)}
\end{equation}
with
\def\theequation{3.72}\begin{equation}
\aligned
{\rm Card}
\left(
\mathfrak{F}_{\mu_1\lambda_1+\cdots+\mu_d\lambda_d}^{
(\mu_1,\lambda_1),\dots,(\mu_d,\lambda_d)}
\right)
&
=
\frac{
{\rm Card}
\left(
\mathfrak{S}_{\mu_1\lambda_1+\cdots+\mu_d\lambda_d}
\right)}{
{\rm Card}
\left(
\mathfrak{H}_{\mu_1\lambda_1+\cdots+\mu_d\lambda_d}^{
(\mu_1,\lambda_1),\dots,(\mu_d,\lambda_d)}
\right)}
\\
&
=
\frac{(\mu_1\lambda_1+\cdots+\mu_d\lambda_d)!}{
\mu_1!(\lambda_1!)^{\mu_1}\cdots
\mu_d!(\lambda_d!)^{\mu_d}}.
\endaligned
\end{equation}
In conclusion, by means of this formalism, we may write down the
complete generalization of Theorem~2.24 to several independent
variables.

\def\thetheorem{3.73}\begin{theorem}
For every $\kappa \geqslant 1$ and for every choice of $\kappa$ indices
$i_1,\dots, i_\kappa$ in the set $\{ 1, 2, \dots, n\}$, the general
expression of ${\bf Y}_{i_1, \dots, i_\kappa}$ is as follows{\rm :}
\def\theequation{3.74}\begin{equation}
\small
\boxed{
\aligned
{\bf Y}_{i_1, \dots, i_\kappa}
&
=
\mathcal{Y}_{x^{i_1}\cdots x^{i_\kappa}}
+
\sum_{d=1}^{\kappa+1}
\ \
\sum_{1\leqslant\lambda_1<\cdots<\lambda_d\leqslant\kappa}
\ \
\sum_{\mu_1\geqslant 1,\dots,\mu_d\geqslant 1} 
\
\sum_{
\mu_1\lambda_1
+
\cdots
+
\mu_d\lambda_d\leqslant \kappa+1} 
\\
&
\sum_{k_{1:1:1},\dots,k_{1:1:\lambda_1}=1}^n
\cdots \
\sum_{k_{1:\mu_1:1},\dots,k_{1:\mu_1:\lambda_1}=1}^n
\cdots\cdots \
\sum_{k_{d:1:1},\dots,k_{d:1:\lambda_d}=1}^n
\cdots \
\sum_{k_{d:\mu_d:1},\dots,k_{d:\mu_d:\lambda_d}=1}^n
\\
&
\left[
\aligned
& 
\sum_{\sigma\in\mathfrak{F}_{\mu_1\lambda_1+\cdots+\mu_d\lambda_d}^{
(\mu_1,\lambda_1),\dots,(\mu_d,\lambda_d)}}
\
\sum_{\tau\in\mathfrak{S}_\kappa^{
\mu_1\lambda_1+\cdots+\mu_d\lambda_d}}
\\
& \
\ \ \ \ \
\delta_{i_{\tau(1)},\dots,i_{\tau(\mu_1\lambda_1)},\dots,
i_{\tau(\mu_1\lambda_1+\cdots+\mu_d\lambda_d)}}^{
k_{\sigma(1:1:1)},\dots,k_{\sigma(1:\mu_1:\lambda_1)},
\dots,k_{\sigma(d:\mu_d:\lambda_d)}}\,
\frac{\partial^{\kappa-\mu_1\lambda_1-\cdots-\mu_d\lambda_d+
\mu_1+\cdots+\mu_d}
\mathcal{Y}}{
\partial x^{i_{\tau(\mu_1\lambda_1+\cdots+\mu_d\lambda_d+1)}}\cdots
\partial x^{i_{\tau(\kappa)}}
\left(\partial y\right)^{\mu_1+\cdots+\mu_d}}
- \\
& 
-
\sum_{\sigma\in\mathfrak{F}_{\mu_1\lambda_1+\cdots+\mu_d\lambda_d}^{
(\mu_1,\lambda_1),\dots,(\mu_d,\lambda_d)}}
\
\sum_{\tau\in\mathfrak{S}_\kappa^{
\mu_1\lambda_1+\cdots+\mu_d\lambda_d-1}}
\\
& \
\ \ \ \ \
\delta_{i_{\tau(1)},\dots,i_{\tau(\mu_1\lambda_1)},\dots,
i_{\tau(\mu_1\lambda_1+\cdots+\mu_d\lambda_d-1)}}^{
k_{\sigma(1:1:1)},\dots,k_{\sigma(1:\mu_1:\lambda_1)},
\dots,k_{\sigma(d:\mu_d:\lambda_d-1)}}\,
\frac{\partial^{\kappa-\mu_1\lambda_1-\cdots-\mu_d\lambda_d
+\mu_1+\cdots+\mu_d}\mathcal{X}^{k_{\sigma(d:\mu_d:\lambda_d)}}}{
\partial x^{i_{\tau(\mu_1\lambda_1+\cdots+\mu_d\lambda_d)}}\cdots
\partial x^{i_{\tau(\kappa)}}
\left(\partial y\right)^{\mu_1+\cdots+\mu_d-1}}
\endaligned
\right]
\cdot
\\
&
\ \ \ \ \ \ \ \ \ \ \ \ \ \ \ \ \ \ \ \ \ \ \ \ 
\cdot
\prod_{1\leqslant\nu_1\leqslant\mu_1}\,
y_{k_{1:\nu_1:1},\dots,k_{1:\nu_1:\lambda_1}}
\ \cdots \
\prod_{1\leqslant\nu_d\leqslant\mu_d}\,
y_{k_{d:\nu_d:1},\dots,k_{d:\nu_d:\lambda_d}}.
\endaligned
}
\end{equation}
\end{theorem}

Since the fundamental monomials appearing in the last line of~\thetag{
3.74} just above are not independent of each other, this formula has
still to be modified a little bit. We refer to Section~6 for details.

\subsection*{ 3.75.~Deduction of a multivariate Fa\`a
di Bruno formula} Let $n \in \N$ with $n\geqslant 1$, let $x = (x^1,\dots,
x^n) \in \K^n$, let $g = g( x^1, \dots, x^n)$ be a $\mathcal{
C}^\infty$-smooth function from $\K^n$ to $\K$, let $y \in \K$ and let
$f = f(y)$ be a $\mathcal{ C}^\infty$ function from $\K$ to $\K$. The
goal is to obtain an explicit formula for the partial derivatives of
the composition $h := f\circ g$, namely $h(x^1, \dots, x^n) := f (
g(x^1,\dots, x^n))$. For $\lambda \in \N$ with $\lambda \geqslant 1$ and
for arbitrary indices $i_1, \dots, i_\lambda = 1, \dots, n$, we
shall abbreviate the partial derivative $\frac{ \partial^\lambda g}{
\partial x^{i_1} \cdots \partial x^{i_\lambda}}$ by $g_{i_1,\dots,
i_\lambda}$ and similarly for $h_{i_1, \dots, i_\lambda}$. The
derivative $\frac{ d^\lambda f}{ d y^\lambda}$ will be
abbreviated by $f_\lambda$.

Appying the chain rule, we may compute:
\def\theequation{3.76}\begin{equation}
\small
\aligned
h_{i_1}
&
=
f_1 
\left[
g_{i_1}
\right], 
\\
h_{i_1,i_2}
&
=
f_2
\left[
g_{i_1}\,g_{i_2}
\right]
+
f_1
\left[
g_{i_1,i_2}
\right],
\\
h_{i_1,i_2,i_3}
&
=
f_3
\left[
g_{i_1}\,g_{i_2}\,g_{i_3}
\right]
+
f_2
\left[
g_{i_1}\,g_{i_2,i_3}
+
g_{i_2}\,g_{i_1,i_3}
+
g_{i_3}\,g_{i_1,i_2}
\right]
+
f_1
\left[
g_{i_1,i_2,i_3}
\right],
\\
h_{i_1,i_2,i_3,i_4}
&
=
f_4
\left[
g_{i_1}\,g_{i_2}\,g_{i_3}\,g_{i_4}
\right]
+
f_3
\left[
g_{i_2}\,g_{i_3}\,g_{i_1,i_4}
+
g_{i_3}\,g_{i_1}\,g_{i_2,i_4}
+
g_{i_1}\,g_{i_2}\,g_{i_3,i_4}
+
\right. 
\\
& \
\ \ \ \ \ \ \ \ \ \ \ \ \ \ \ \ \ \ \ \
\ \ \ \ \ \ \ \ \ \ \ \ \ \ \ \ \ \ \ \
\ \
\left.
+
g_{i_1}\,g_{i_4}\,g_{i_2,i_3}
+
g_{i_2}\,g_{i_4}\,g_{i_1,i_3}
+
g_{i_3}\,g_{i_4}\,g_{i_1,i_2}
\right]
+ \\
& \
\ \ \ \ \
+
f_2
\left[
g_{i_1,i_2}\,g_{i_3,i_4}
+
g_{i_1,i_3}\,g_{i_2,i_4}
+
g_{i_1,i_4}\,g_{i_2,i_3}
\right]
+ \\
& \
\ \ \ \ \
+
f_2
\left[
g_{i_1}\,g_{i_2,i_3,i_4}
+
g_{i_2}\,g_{i_1,i_3,i_4}
+
g_{i_3}\,g_{i_1,i_2,i_4}
+
g_{i_4}\,g_{i_1,i_2,i_3}
\right]
+ \\
& \
\ \ \ \ \
+
f_1
\left[
g_{i_1,i_2,i_3,i_4}
\right].
\endaligned
\end{equation}
Introducing the derivations
\def\theequation{3.77}\begin{equation}
\small
\aligned
F_i^2
&
:= 
\sum_{k_1=1}^n\,g_{k_1,i}\,
\frac{\partial}{\partial g_{k_1}}
+
g_i\left(
f_2\,\frac{\partial}{\partial f_1}
\right), \\
F_i^3
&
:=
\sum_{k_1=1}^n\,g_{k_1,i}\,
\frac{\partial}{\partial g_{k_1}}
+
\sum_{k_1,k_2=1}^n\,g_{k_1,k_2,i}\,
\frac{\partial}{\partial g_{k_1,k_2}}
+
g_i\left(
f_2\,\frac{\partial}{\partial f_1}
+
f_3\,\frac{\partial}{\partial f_2}
\right), \\
\endaligned
\end{equation}
\[
\small
\aligned
\text{\bf 
\dots\dots}
& \
\ \ \ \ \
\text{\bf 
\dots\dots\dots\dots\dots\dots\dots\dots\dots
\dots\dots\dots\dots\dots\dots\dots\dots\dots
\dots\dots\dots\dots\dots\dots\dots
} \\
F_i^\lambda
&
:=
\sum_{k_1=1}^n\,g_{k_1,i}\,
\frac{\partial}{\partial g_{k_1}}
+
\sum_{k_1,k_2=1}^n\,g_{k_1,k_2,i}\,
\frac{\partial}{\partial g_{k_1,k_2}}
+
\cdots
+ \\
& \
\ \ \ \ \ \ \ \ \ \ 
+
\sum_{k_1,\dots,k_{\lambda-1}=1}^n\,g_{k_1,\dots,k_{\lambda-1},i}\,
\frac{\partial}{\partial g_{k_1,\dots,k_{\lambda-1}}}
+ \\
& \
\ \ \ \ \ \ \ \ \ \ 
+
g_i\left(
f_2\,\frac{\partial}{\partial f_1}
+
f_3\,\frac{\partial}{\partial f_2}
+
\cdots
+
f_\lambda\,\frac{\partial}{\partial f_{\lambda-1}}
\right), \\
\endaligned
\]
we observe that the following
induction relations hold:
\def\theequation{3.78}\begin{equation}
\aligned
h_{i_1,i_2}
&
=
F_{i_2}^2
\left(
h_{i_1}
\right), 
\\
h_{i_1,i_2,i_3}
&
=
F_{i_3}^3
\left(
h_{i_1,i_2}
\right), 
\\
\text{\bf 
\dots\dots
}
& \
\ \ \ \ \
\text{\bf 
\dots\dots\dots\dots\dots
}
\\
h_{i_1,i_2,\dots,i_\lambda}
&
=
F_{i_\lambda}^\lambda
\left(
h_{i_1,i_2,\dots,i_{\lambda-1}}
\right).
\endaligned
\end{equation}

To obtain the explicit version of the Fa\`a di Bruno in the case of
several variables $(x^1, \dots, x^n)$ and one variable $y$, it
suffices to extract from the expression of ${\bf Y }_{ i_1,\dots,
i_\kappa}$ provided by Theorem~3.73 only the terms corresponding to
$\mu_1 \lambda_1 + \cdots + \mu_d\lambda_d = \kappa$, dropping all the
$\mathcal{ X}$ terms. After some simplifications and after a
translation by means of an elementary dictionary, we obtain a
statement.

\def\thetheorem{3.79}\begin{theorem} 
For every integer $\kappa \geqslant 1$ and for every choice of indices
$i_1, \dots, i_\kappa$ in the set $\{ 1, 2, \dots, n\}$, the
$\kappa$-th partial derivative of the composite function $h = h( x^1,
\dots, x^n) = f( g(x^1, \dots, x^n))$ with respect to the variables
$x^{i_1}, \dots, x^{i_\kappa}$ may be expressed as an explicit
polynomial depending on the derivatives of $f$, on the partial
derivatives of $g$ and having integer coefficients{\rm :}
\def\theequation{3.80}\begin{equation}
\boxed{
\aligned
\frac{\partial^\kappa h}{\partial x^{i_1}\cdots
\partial x^{i_\kappa}}
&
=
\sum_{d=1}^\kappa
\
\sum_{1\leqslant \lambda_1 < \cdots < \lambda_d \leqslant \kappa}
\
\sum_{\mu_1\geqslant 1,\dots,\mu_d\geqslant 1}
\
\sum_{\mu_1\lambda_1+\cdots+\mu_d\lambda_d=\kappa}
\
\frac{d^{\mu_1+\cdots+\mu_d} f}{
dy^{\mu_1+\cdots+\mu_d}}
\\
& 
\left[
\aligned
&
\sum_{\sigma\in\mathfrak{F}_\kappa^{
(\mu_1,\lambda_1),\dots,(\mu_d,\lambda_d)}}
\
\prod_{1\leqslant\nu_1\leqslant\mu_1}
\
\frac{\partial^{\lambda_1} g}{\partial
x^{i_{\sigma(1:\nu_1:1)}}\cdots
\partial x^{i_{\sigma(1:\nu_1:\lambda_1)}}}
\
\text{\bf \dots} 
\\
& \
\ \ \ \ \ \ \ \ \ \ \ \ \ \ \ \ \ \ \ \ \ \ \ \ \ \ \ \ \ \
\text{\bf \dots} 
\prod_{1\leqslant\nu_d\leqslant\mu_d}
\
\frac{\partial^{\lambda_d} g}{\partial
x^{i_{\sigma(d:\nu_d:1)}}\cdots
\partial x^{i_{\sigma(d:\nu_ds:\lambda_d)}}}
\endaligned
\right].
\endaligned
}
\end{equation}
\end{theorem}

In this formula, the coset $\mathfrak{ F }_\kappa^{
(\mu_1, \lambda_1 ),\dots, ( \mu_d, \lambda_d)}$
was defined in equation~\thetag{ 3.71}; we have made the identification{\rm :}
\def\theequation{3.81}\begin{equation}
\{1,\dots,\kappa\}
\equiv
\{
1\!:\!1\!:\!1,
\dots,
1\!:\!\mu_1\!:\!\lambda_1,
\text{\rm \dots\dots}, 
d\!:\!1\!:\!1,
\dots,
d\!:\!\mu_d\!:\!\lambda_d
\};
\end{equation}
and also, for the sake of clarity, we have restituted the complete
(not abbreviated) notation for the (partial) derivatives of $f$ and of
$g$. 

We refer to Section~6 for the final writing of the above
formula~\thetag{ 3.80}.

\section*{\S4.~One independent variable and
several dependent variables}


\subsection*{4.1.~Simplified adapted notations}
Assume $n = 1$ and $m \geqslant 1$, let $\kappa\in \N$ with $\kappa
\geqslant 1$ and simply denote the jet variables by (instead
of~\thetag{ 1.2}):
\def\theequation{4.2}\begin{equation}
\left(
x,y^j,y_1^j,y_2^j,\dots,y_\kappa^j
\right)
\in
\mathcal{J}_{1,m}^\kappa.
\end{equation}
Instead of~\thetag{ 1.30}, denote the $\kappa$-th prolongation
of a vector field by:
\def\theequation{4.3}\begin{equation}
\left\{
\aligned
\mathcal{L}^{(\kappa)}
&
=
\mathcal{X}\,\frac{\partial}{\partial x}
+
\sum_{j=1}^m\,\mathcal{Y}^j\,\frac{\partial}{\partial y^j}
+
\sum_{j=1}^m\,{\bf Y}_1^j\,\frac{\partial}{\partial y_1^j}
+
\sum_{j=1}^m\,{\bf Y}_2^j\,\frac{\partial}{\partial y_2^j}
+ \\
& \
\ \ \ \ \
+
\cdots
+
\sum_{j=1}^m\,{\bf Y}_\kappa^j\,\frac{\partial}{\partial y_\kappa^j}.
\endaligned\right.
\end{equation}
The induction formulas are:
\def\theequation{4.4}\begin{equation}
\left\{
\aligned
{\bf Y}_1^j
&
:=
D^1
\left(
\mathcal{Y}^j
\right)
-
D^1
\left(
\mathcal{ X}
\right)
y_1^j, 
\\
{\bf Y}_2^j
&
:=
D^2
\left(
{\bf Y}_1^j
\right)
-
D^1
\left(
\mathcal{ X}
\right)
y_2^j, 
\\
\cdots\cdots
&
\cdots\cdots
\\
{\bf Y}_\lambda^j
&
:=
D^\lambda
\left(
{\bf Y}_{\lambda-1}^j
\right)
-
D^1
\left(
\mathcal{ X}
\right)
y_\lambda^j, 
\endaligned\right.
\end{equation}
where the total
differentiation operators $D^\lambda$
are denoted by (instead of~\thetag{ 1.22}):
\def\theequation{4.5}\begin{equation}
\aligned
D^\lambda
:=
\frac{\partial}{\partial x}
+
\sum_{l=1}^m\,y_1^l\,
\frac{\partial}{\partial y^l}
+
\sum_{l=1}^m\,y_2^l\,
\frac{\partial}{\partial y_1^l}
+
\cdots
+ 
\sum_{l=1}^m\,y_\lambda^l\,
\frac{\partial}{\partial y_{\lambda-1}^l}.
\endaligned
\end{equation}
Applying the definitions in the first two lines of~\thetag{ 4.4}, we
compute, we simplify and we organize the results in a harmonious way,
using in an essential way the Kronecker symbol. Here, the
computations are more elementary than the computations of ${\bf
Y}_{i_1}$ and of ${\bf Y}_{i_1, i_2}$ achieved thoroughly 
in the previous
Section~3, so that we do not provide a Latex track of the
details. Firstly and secondly:
\def\theequation{4.6}\begin{equation}
\small
\left\{
\aligned
{\bf Y}_1^j
&
=
\mathcal{ Y}_x^j
+
\sum_{l_1=1}^m
\left[
\mathcal{Y}_{y^{l_1}}^j
-
\delta_{l_1}^j\,
\mathcal{X}_x
\right]
y_1^{l_1}
+
\sum_{l_1,l_2=1}^m
\left[
-
\delta_{l_1}^j\,
\mathcal{X}_{y^{l_2}}
\right]
y_1^{l_1}y_1^{l_2}, \\
{\bf Y}_2^j
&
=
\mathcal{ Y}_{x^2}^j
+
\sum_{l_1=1}^m
\left[
2\,\mathcal{Y}_{xy^{l_1}}^j
-
\delta_{l_1}^j\,
\mathcal{X}_{x^2}
\right]
y_1^{l_1}
+
\sum_{l_1,l_2=1}^m
\left[
\mathcal{Y}_{y^{l_1}y^{l_2}}^j
-
\delta_{l_1}^j\,2\,
\mathcal{X}_{xy^{l_2}}
\right]
y_1^{l_1}y_1^{l_2}
+ \\
& \
\ \ \ \ \
+
\sum_{l_1,l_2,l_3}
\left[
-
\delta_{l_1}^j\,
\mathcal{X}_{y^{l_2}y^{l_3}}
\right]
y_1^{l_1}y_1^{l_2}y_1^{l_3}
+
\sum_{l_1}
\left[
\mathcal{Y}_{y^{l_1}}^j
-
\delta_{l_1}^j\,2\,
\mathcal{X}_x
\right]
y_2^{l_1}
+ \\
& \
\ \ \ \ \ 
+
\sum_{l_1,l_2=1}^m
\left[
-
\delta_{l_1}^j\,
\mathcal{X}_{y^{l_2}}
-
\delta_{l_2}^j\,2\,
\mathcal{X}_{y^{l_1}}
\right]
y_1^{l_1}y_2^{l_2}.
\endaligned\right.
\end{equation}
Thirdly:
\def\theequation{4.7}\begin{equation}
\small
\aligned
{\bf Y}_3^j
&
=
\mathcal{ Y}_{x^3}^j
+
\sum_{l_1=1}^m
\left[
3\,\mathcal{Y}_{x^2y^{l_1}}^j
-
\delta_{l_1}^j\,
\mathcal{X}_{x^3}
\right]
y_1^{l_1}
+
\sum_{l_1,l_2=1}^m
\left[
3\,\mathcal{Y}_{xy^{l_1}y^{l_2}}^j
-
\delta_{l_1}^j\,
3\,
\mathcal{X}_{x^2y^{l_2}}
\right]
y_1^{l_1}y_1^{l_2}+ \\
& \
\ \ \ \ \
+
\sum_{l_1,l_2,l_3}
\left[
\mathcal{Y}_{y^{l_1}y^{l_2}y^{l_3}}^j
-
\delta_{l_1}^j\,
3\,\mathcal{X}_{xy^{l_2}y^{l_3}}
\right]
y_1^{l_1}y_1^{l_2}y_1^{l_3}
+ \\
& \ 
\ \ \ \ \
+
\sum_{l_1,l_2,l_3,l_4}
\left[
-
\delta_{l_1}^j\,
\mathcal{X}_{y^{l_2}y^{l_3}y^{l_4}}
\right]
y_1^{l_1}y_1^{l_2}y_1^{l_3}y_1^{l_4}
+
\sum_{l_1=1}^m
\left[
3\,\mathcal{Y}_{xy^{l_1}}^j
-
\delta_{l_1}^j\,
3\,\mathcal{X}_{x^2}
\right]
y_2^{l_1}
+ 
\endaligned
\end{equation}
$$
\small
\aligned
& \
\ \ \ \ \ 
+
\sum_{l_1,l_2=1}^m
\left[
3\,\mathcal{Y}_{y^{l_1}y^{l_2}}^j
-
\delta_{l_1}^j\,
3\,\mathcal{X}_{xy^{l_2}}
-
\delta_{l_2}^j\,
6\,\mathcal{X}_{xy^{l_1}}
\right]
y_1^{l_1}y_2^{l_2}
+ \\
& \
\ \ \ \ \
+
\sum_{l_1,l_2,l_3=1}^m
\left[
-
\delta_{l_1}^j\,
3\,\mathcal{X}_{y^{l_2}y^{l_3}}
-
\delta_{l_3}^j\,
3\,\mathcal{X}_{y^{l_1}y^{l_2}}
\right]
y_1^{l_1}y_1^{l_2}y_2^{l_3}
+ 
\sum_{l_1,l_2=1}^m
\left[
-
\delta_{l_3}^j\,
3\,\mathcal{X}_{y^{l_2}}
\right]
y_2^{l_1}y_2^{l_2}
+ \\
& \
\ \ \ \ \
+
\sum_{l_1=1}^m
\left[
\mathcal{Y}_{y^{l_1}}^j
-
\delta_{l_1}^j\,
3\,\mathcal{X}_x
\right]
y_3^{l_1}
+
\sum_{l_1,l_2=1}^m
\left[
-
\delta_{l_1}^j\,
\mathcal{X}_{y^{l_2}}
-
\delta_{l_2}^j\,
3\,\mathcal{X}_{y^{l_1}}
\right]
y_1^{l_1}y_3^{l_2}.
\endaligned
$$
Fourthly:
$$
\small
\aligned
{\bf Y}_4^j
&
=
\mathcal{ Y}_{x^4}^j
+
\sum_{l_1=1}^m
\left[
4\,\mathcal{Y}_{x^3y^{l_1}}^j
-
\delta_{l_1}^j
\mathcal{X}_{x^4}\,
\right]
y_1^{l_1}
+
\sum_{l_1,l_2=1}^m
\left[
6\,\mathcal{Y}_{x^2y^{l_1}y^{l_2}}^j
-
\delta_{l_1}^j\,
4\,
\mathcal{X}_{x^3y^{l_2}}
\right]
y_1^{l_1}y_1^{l_2}
+ \\
& \
\ \ \ \ \
+
\sum_{l_1,l_2,l_3=1}^m
\left[
4\,\mathcal{Y}_{xy^{l_1}y^{l_2}y^{l_3}}^j
-
\delta_{l_1}^j\,
6\,\mathcal{X}_{x^2y^{l_2}y^{l_3}}
\right]
y_1^{l_1}y_1^{l_2}y_1^{l_3}
+ \\
& \
\ \ \ \ \
+
\sum_{l_1,l_2,l_3,l_4=1}^m
\left[
\mathcal{Y}_{xy^{l_1}y^{l_2}y^{l_3}y^{l_4}}^j
-
\delta_{l_1}^j\,
4\,\mathcal{X}_{xy^{l_2}y^{l_3}y^{l_4}}
\right]
y_1^{l_1}y_1^{l_2}y_1^{l_3}y_1^{l_4}
+ \\
& \ 
\ \ \ \ \
+
\sum_{l_1,l_2,l_3,l_4,l_5=1}^m
\left[
-
\delta_{l_1}^j\,
\mathcal{X}_{y^{l_2}y^{l_3}y^{l_4}y^{l_5}}
\right]
y_1^{l_1}y_1^{l_2}y_1^{l_3}y_1^{l_4}y_1^{l_5}
+
\sum_{l_1=1}^m
\left[
6\,\mathcal{Y}_{x^2y^{l_1}}^j
-
\delta_{l_1}^j\,
4\,\mathcal{X}_{x^3}
\right]
y_2^{l_1}
+
\endaligned
$$
\def\theequation{4.8}\begin{equation}
\small
\aligned
& \
\ \ \ \ \ 
+
\sum_{l_1,l_2=1}^m
\left[
12\,\mathcal{Y}_{xy^{l_1}y^{l_2}}^j
-
\delta_{l_1}^j\,
6\,\mathcal{X}_{x^2y^{l_2}}
-
\delta_{l_2}^j\,
12\,\mathcal{X}_{x^2y^{l_1}}
\right]
y_1^{l_1}y_2^{l_2}
+ \\
& \
\ \ \ \ \ 
+
\sum_{l_1,l_2,l_3=1}^m
\left[
6\,\mathcal{Y}_{y^{l_1}y^{l_2}y^{l_3}}^j
-
\delta_{l_1}^j\,
12\,\mathcal{X}_{xy^{l_2}y^{l_3}}
-
\delta_{l_3}^j\,
12\,\mathcal{X}_{xy^{l_1}y^{l_2}}
\right]
y_1^{l_1}y_1^{l_2}y_2^{l_3}
+ \\
& \
\ \ \ \ \ 
+
\sum_{l_1,l_2,l_3,l_4=1}^m
\left[
-
\delta_{l_1}^j\,
6\,\mathcal{X}_{y^{l_2}y^{l_3}y^{l_4}}
-
\delta_{l_4}^j\,
4\,\mathcal{X}_{y^{l_1}y^{l_2}y^{l_3}}
\right]
y_1^{l_1}y_1^{l_2}y_1^{l_3}y_2^{l_4}
+ \\
\endaligned
\end{equation}
$$
\small
\aligned
& \
\ \ \ \ \
+
\sum_{l_1,l_2=1}^m
\left[
3\,\mathcal{Y}_{y^{l_1}y^{l_2}}^j
-
\delta_{l_1}^j\,
12\,\mathcal{X}_{xy^{l_2}}
\right]
y_2^{l_1}y_2^{l_2}
+ \\
& \
\ \ \ \ \
+
\sum_{l_1,l_2,l_3=1}^m
\left[
-
\delta_{l_1}^j\,
3\,\mathcal{X}_{y^{l_2}y^{l_3}}
-
\delta_{l_2}^j\,
12\,\mathcal{X}_{y^{l_1}y^{l_3}}
\right]
y_1^{l_1}y_2^{l_2}y_2^{l_3}
+
\sum_{l_1=1}^m
\left[
4\,\mathcal{Y}_{xy^{l_1}}^j
-
\delta_{l_1}^j\,
6\,\mathcal{X}_{x^2}
\right]
y_3^{l_1}
+ \\
& \
\ \ \ \ \
+
\sum_{l_1,l_2=1}^m
\left[
4\,\mathcal{Y}_{y^{l_1}y^{l_2}}^j
-
\delta_{l_1}^j\,
4\,\mathcal{X}_{xy^{l_2}}
-
\delta_{l_2}^j\,
12\,\mathcal{X}_{xy^{l_1}}
\right]
y_1^{l_1}y_3^{l_2}
+ \\
& \
\ \ \ \ \
+
\sum_{l_1,l_2,l_3=1}^m
\left[
-
\delta_{l_1}^j\,
4\,\mathcal{X}_{y^{l_2}y^{l_3}}
-
\delta_{l_3}^j\,
6\,\mathcal{X}_{y^{l_1}y^{l_2}}
\right]
y_1^{l_1}y_1^{l_2}y_3^{l_3}
+ \\
& \
\ \ \ \ \
+
\sum_{l_1,l_2=1}^m
\left[
-
\delta_{l_1}^j\,
4\,\mathcal{X}_{y^{l_2}}
-
\delta_{l_2}^j\,
6\,\mathcal{X}_{y^{l_1}}
\right]
y_2^{l_1}y_3^{l_2}
+ \\
& \
\ \ \ \ \
+
\sum_{l_1=1}^m
\left[
\mathcal{Y}_{y^{l_1}}^j
-
\delta_{l_1}^j\,
4\,\mathcal{X}_x
\right]
y_4^{l_1}
+
\sum_{l_1,l_2=1}^m
\left[
-
\delta_{l_1}^j\,
\mathcal{X}_{y^{l_2}}
-
\delta_{l_2}^j\,
4\,\mathcal{X}_{y^{l_1}}
\right]
y_1^{l_1}y_4^{l_2}.
\endaligned
$$

\subsection*{4.9.~Inductive elaboration of the
general formula} Now we compare the formula~\thetag{ 2.9} for ${\bf
Y}_4$ with the above formula~\thetag{ 4.8} for ${\bf Y}_4^j$. The goal
is to find the rules of transformation and of development by
inspecting several instances, in order to devise how to transform and
to develope the formula~\thetag{ 2.25} to several dependent variables.

First of all, we have to develope the general monomial $(y_{ \lambda_1
})^{ \mu_1} \cdots (y_{ \lambda_d })^{ \mu_d}$. In every monomial
present in the expressions of ${\bf Y}_1^j$, of ${\bf Y}_2^j$, of
${\bf Y}_3^j$ and of ${\bf Y}_4^j$ above, we see that the number
$\alpha$ of indices $l_\beta$ appearing in all the sums $\sum_{l_1,
\dots, l_\alpha = 1}^m$ is exactly equal to $\mu_1 + \dots+ \mu_d$. To
denote these $\mu_1 + \cdots + \mu_d$ indices $l_\beta$, we shall use
the notation:
\def\theequation{4.10}\begin{equation}
\underbrace{
\underbrace{
l_{1:1},\dots,l_{1:\mu_1}}_{\mu_1},
\dots,
\underbrace{
l_{d:1},\dots,l_{d:\mu_d}}_{\mu_d}}_{
\mu_1+\cdots+\mu_d},
\end{equation} 
inspired by Convention~3.33. With such a choice of notation, we may
avoid long subscripts in the indices $l_\beta$, like $l_{\mu_1+\cdots+
\mu_d}$. It follows that the development of the general monomial $(y_{
\lambda_1 })^{ \mu_1} \cdots (y_{ \lambda_d })^{ \mu_d}$ to several
dependent variables yields $m^{\mu_1+ \cdots + \mu_d}$ possible
choices:
\def\theequation{4.11}\begin{equation}
\prod_{1\leqslant\nu_1\leqslant\mu_1}
y_{\lambda_1}^{l_{1:\nu_1}}
\ \cdots\cdots
\prod_{1\leqslant\nu_d\leqslant\mu_d}
y_{\lambda_d}^{l_{d:\nu_d}},
\end{equation}
where the indices $l_{1 :1}, \dots, l_{1: \mu_1 }, \dots, l_{ d:1},
\dots, l_{d: \mu_d }$ take their values in the set $\{ 1, 2, \dots, m
\}$. Consequently, the general expression of ${\bf Y }_\kappa^j$ must
be of the form:
\def\theequation{4.12}\begin{equation}
\small
\aligned
{\bf Y}_\kappa^j
&
=
\mathcal{Y}_{x^\kappa}^j
+
\sum_{d=1}^{\kappa+1}
\
\sum_{1\leqslant\lambda_1<\cdots<\lambda_d\leqslant\kappa}
\
\sum_{\mu_1\geqslant 1,\dots,\mu_d\geqslant 1}
\
\sum_{\mu_1\lambda_1+\cdots+\mu_d\lambda_d\leqslant\kappa+1} 
\
\\
& \
\ \ \ \ \ \ \ \ \ \ \ \ \ \ \
\sum_{l_{1:1}=1}^m
\cdots
\sum_{l_{1:\mu_1}=1}^m
\cdots\cdots
\sum_{l_{d:1}=1}^m
\cdots
\sum_{l_{d:\mu_d}=1}^m
\
\text{\bf [?]}
\\
& \
\ \ \ \ \ \ \ \ \ \ \ \ \ \ \ \ \ \ \
\prod_{1\leqslant\nu_1\leqslant\mu_1}
y_{\lambda_1}^{l_{1:\nu_1}}
\ \cdots\cdots
\prod_{1\leqslant\nu_d\leqslant\mu_d}
y_{\lambda_d}^{l_{d:\nu_d}},
\endaligned
\end{equation}
where the term in brackets {\bf [?]} is still unknown. To determine
it, let us examine a few instances.

According to~\thetag{ 4.8} (fourth line), the term $\left[ 6\,
\mathcal{ Y}_{ x^2y} - 4\, \mathcal{ X}_{ x^3} \right] y_2$ of ${\bf
Y}_4$ developes as $\sum_{ l_1 =1}^m \, \left[ 6\, \mathcal{ Y}_{ x^2
y^{l_1}}^j - \delta_{ l_1}^j \, 4 \, \mathcal{ X}_{ x^3} \right] y_2^{
l_1}$ in ${\bf Y}_4^j$. Here, 
$6\, \mathcal{ Y}_{ x^2y}$ just becomes
$6\, \mathcal{ Y}_{ x^2y^{l_1}}^j$. Thus, we suspect that the term
$\frac{\kappa\cdots(\kappa-\mu_1\lambda_1-\cdots-\mu_d\lambda_d+1)}
{(\lambda_1!)^{\mu_1}\,\mu_1! \cdots (\lambda_d! )^{\mu_d}\,\mu_d! }
\cdot \mathcal{ Y}_{ x^{ \kappa-\mu_1\lambda_1-\cdots-\mu_d\lambda_d}
\, y^{\mu_1+\cdots+\mu_d} } $ of the second line of~\thetag{ 2.25}
should simply be developed as
\def\theequation{4.13}\begin{equation}
\small
\aligned
& 
\frac{\kappa(\kappa-1)\cdots(\kappa-\mu_1\lambda_1-\cdots
-\mu_d\lambda_d+1)}{
(\lambda_1!)^{\mu_1}\ \mu_1!\cdots
(\lambda_d!)^{\mu_d}\ \mu_d!}
\cdot
\\
& \
\ \ \ \ \ \ \ \
\cdot
\frac{\partial^{\kappa-\mu_1\lambda_1-\cdots
-\mu_d\lambda_d+\mu_1+\cdots+\mu_d}\mathcal{Y}^j}{
(\partial x)^{\kappa-\mu_1\lambda_1-\cdots
-\mu_d\lambda_d}
\partial y^{l_{1:1}}
\cdots
\partial y^{l_{1:\mu_1}}
\cdots
\partial y^{l_{d:1}}
\cdots
\partial y^{l_{d:\mu_d}}
}.
\endaligned
\end{equation}
This rule is confirmed by inspecting all the other monomials of ${\bf
Y}_1^j$, of ${\bf Y}_2^j$, of ${\bf Y}_3^j$ and of ${\bf Y}_4^j$.

It remains to determine how we must develope the term in $\mathcal{
X}$ appearing in the last two lines of~\thetag{ 2.25}. To begin with,
let us rewrite in advance this term in the slightly different
shape, emphasizing a factorization:
\def\theequation{4.14}\begin{equation}
\small
\aligned
\frac{\kappa\cdots(
\kappa-\mu_1\lambda_1-\cdots-\mu_d\lambda_d+2)}
{(\lambda_1!)^{\mu_1}\,\mu_1!
\cdots
(\lambda_d!)^{\mu_d}\,\mu_d!
}
\left[
(\mu_1\lambda_1+\cdots+\mu_d\lambda_d)
\mathcal{X}_{
x^{\kappa-\mu_1\lambda_1-\cdots-\mu_d\lambda_d+1}
\,
y^{\mu_1+\cdots+\mu_d-1}
}
\right].
\endaligned
\end{equation}
Then we examine four instances extracted from the complete expression
of ${\bf Y}_4^j$:
\def\theequation{4.15}\begin{equation}
\small
\left\{
\aligned
&
\sum_{l_1,l_2,l_3=1}^m
\left[
4\,\mathcal{Y}_{xy^{l_1}y^{l_2}y^{l_3}}^j
-
\delta_{l_1}^j\,
6\,\mathcal{X}_{x^2y^{l_2}y^{l_3}}
\right]
y_1^{l_1}y_1^{l_2}y_1^{l_3}, 
\\
&
\sum_{l_1,l_2=1}^m
\left[
12\,\mathcal{Y}_{xy^{l_1}y^{l_2}}^j
-
\delta_{l_1}^j\,
6\,\mathcal{X}_{x^2y^{l_2}}
-
\delta_{l_2}^j\,
12\,\mathcal{X}_{x^2y^{l_1}}
\right]
y_1^{l_1}y_2^{l_2},
\\
&
\sum_{l_1,l_2,l_3,l_4=1}^m
\left[
-
\delta_{l_1}^j\,
6\,\mathcal{X}_{y^{l_2}y^{l_3}y^{l_4}}
-
\delta_{l_4}^j\,
4\,\mathcal{X}_{y^{l_1}y^{l_2}y^{l_3}}
\right]
y_1^{l_1}y_1^{l_2}y_1^{l_3}y_2^{l_4},
\\
&
\sum_{l_1,l_2,l_3=1}^m
\left[
-
\delta_{l_1}^j\,
4\,\mathcal{X}_{y^{l_2}y^{l_3}}
-
\delta_{l_3}^j\,
6\,\mathcal{X}_{y^{l_1}y^{l_2}}
\right]
y_1^{l_1}y_1^{l_2}y_3^{l_3},
\endaligned\right.
\end{equation}
and we compare them to the corresponding terms of ${\bf Y}_4$:
\def\theequation{4.16}\begin{equation}
\small
\left\{
\aligned
&
\left[
4\,\mathcal{Y}_{xy^3}
-
6\,\mathcal{X}_{x^2y^2}
\right]
(y_1)^3, 
\\
&
\left[
12\,\mathcal{Y}_{xy^2}
-
18\,\mathcal{X}_{x^2y}
\right]
y_1y_2,
\\
&
\left[
-
10\,\mathcal{X}_{y^3}
\right]
(y_1)^3y_2,
\\
&
\left[
-
10\,\mathcal{X}_{y^2}
\right]
(y_1)^2y_3.
\endaligned\right.
\end{equation}
In the development from~\thetag{ 4.16} to~\thetag{ 4.15}, we see that
the four integers just before $\mathcal{ X}$, namely $6 = 6$, $18 = 6
+ 12$, $10 = 6 + 4$ and $10 = 4 + 6$, are split in a certain
manner. Also, a single Kronecker symbol $\delta_{l_\alpha}^j$ is added
as a factor. {\it What are the rules}?

In the second splitting $18 = 6 + 12$, we see that the relative weight
of $6$ and of $12$ is the same as the relative weight of $1$ and $2$
in the splitting $3 = 1 + 2$ issued from the lower indices of the
corresponding monomial $y_1^{l_1} y_2^{l_2}$. Similarly, in the third
splitting $10 = 6 + 4$, the relative weight of $6$ and of $4$ is the
same as the relative weight of $1+1+1$ and of $2$ issued from the
lower indices of the corresponding monomial $y_1^{l_1} y_1^{l_2}
y_1^{l_3} y_2^{l_4}$. This rule may be confirmed by inspecting all
the other monomials of ${\bf Y}_2$, ${\bf Y}_2^j$, of ${\bf Y}_3$,
${\bf Y}_3^j$ and of ${\bf Y}_4$, ${\bf Y}_4^j$. For a general $\kappa
\geqslant 1$, the splitting of integers just amounts to decompose the sum
appearing inside the brackets of~\thetag{ 4.14} as $\mu_1\lambda_1,
\mu_2\lambda_2, \dots, \mu_d\lambda_d$. In fact, when we 
wrote~\thetag{ 4.14}, we emphasized in advance the decomposable
factor $(\mu_1 \lambda_1 + \cdots + \mu_d \lambda_d)$.

Next, we have to determine what is the subscript $\alpha$ in the
Kronecker symbol $\delta_{l_\alpha}^j$. We claim that in the four
instances~\thetag{ 4.15}, the subscript $\alpha$ is intrinsically
related to weight splitting. Indeed, recall that in the second line
of~\thetag{ 4.15}, the number $6$ of the splitting $18 = 6 + 12$ is
related to the number $1$ in the splitting $3 = 1 + 2$ of the lower
indices of the monomial $y_1^{l_1} y_2^{l_2}$. It follows that the
index $l_\alpha$ {\it must be}\, the index $l_1$ of the monomial
$y_1^{l_1}$. Similarly, also in the second line of~\thetag{ 4.15}, the
number $12$ of the splitting $18 = 6 + 12$ being related to the number
$2$ in the splitting $3 = 1 + 2$ of the lower indices of the monomial
$y_1^{l_1} y_2^{l_2}$, it follows that the index $l_\alpha$ attached
to the second $\mathcal{ X}$ term must be the index $l_2$ of the
monomial $y_2^{l_2}$.

This rule is still ambiguous. Indeed, let us examine the third line
of~\thetag{ 4.15}. We have the splitting $10 = 6 + 4$, homologous to
the splitting of relative weights $5 = (1+1+1) + 2$ in the monomial
$y_1^{ l_1} y_1^{ l_2} y_1^{ l_3} y_2^{ l_4}$. Of course, it is clear
that we must choose the index $l_4$ for the Kronecker symbol
associated to the second $\mathcal{ X}$ term $-4\, \mathcal{ X}_{
y^3}$, thus obtaining $-\delta_{ l_4}^j \, 4 \, \mathcal{ X}_{ y^{
l_1} y^{ l_2} y^{ l_3}}$. However, since the monomial $y_1^{ l_1}
y_1^{ l_2} y_1^{ l_3}$ has three indices $l_1$, $l_2$ and $l_3$, there
arises a question: {\it what index $l_\alpha$ must we choose for the
Kronecker symbol $\delta_{ l_\alpha}^j$ attached to the first
$\mathcal{ X}$ term $6\,\mathcal{ X}_{y^3}${\rm :} the index $l_1$,
the index $l_2$ or the index $l_3$}?

The answer is simple: {\it any of the three indices $l_1$, $l_2$ or
$l_3$ works}. Indeed, since the monomial $y_1^{ l_1} y_1^{ l_2}
y_1^{ l_3}$ is symmetric with respect to all permutations of the set
of three indices $\{ l_1, l_2, l_3\}$, we have
\def\theequation{4.17}\begin{equation}
\small
\aligned
\sum_{ l_1, l_2, l_3, l_4 = 1}^m\,
\left[
-
\delta_{l_1}^j\,6\,\mathcal{X}_{y^{l_2}y^{l_3}y^{l_4}}
\right]
y_1^{l_1}y_1^{l_2}y_1^{l_3}y_2^{l_4}
&
=
\sum_{ l_1, l_2, l_3, l_4 = 1}^m\,
\left[
-
\delta_{l_2}^j\,6\,\mathcal{X}_{y^{l_1}y^{l_3}y^{l_4}}
\right]
y_1^{l_1}y_1^{l_2}y_1^{l_3}y_2^{l_4}
= 
\\
&
=
\sum_{ l_1, l_2, l_3, l_4 = 1}^m\,
\left[
-
\delta_{l_3}^j\,6\,\mathcal{X}_{y^{l_1}y^{l_2}y^{l_4}}
\right]
y_1^{l_1}y_1^{l_2}y_1^{l_3}y_2^{l_4}.
\endaligned
\end{equation}
In fact, we have systematically used such symmetries during the
intermediate computations (not exposed here) which we achieved
manually to obtain the final expressions of ${\bf Y}_1^j$, of ${\bf
Y}_2^j$, of ${\bf Y}_3^j$ and of ${\bf Y}_4^j$. To fix ideas, we have
always choosen the first index. Here, the first index is $l_1$; in the
first sum of line~9 of~\thetag{ 4.8}, the first index $l_\alpha$ for
the second weight $12$ is $l_2$.

This rule may be confirmed by inspecting all the monomials of ${\bf
Y}_2^j$, of ${\bf Y}_3^j$, of ${\bf Y}_4^j$ (and also of ${\bf
Y}_5^j$, which we have computed in a 
manuscript, but not copied in this Latex file).

From these considerations, we deduce that for the general formula, the
weight decomposition is simply $\mu_1\lambda_1, \dots, \mu_d\lambda_d$
and that the Kronecker symbol $\delta_\alpha^j$ is intrinsically
associated to the weights. In conclusion,
building on inductive reasonings, we have formulated the following
statement.

\def\thetheorem{4.18}\begin{theorem}
For one independent variable $x$, for several dependent variables
$(y^1, \dots, y^m)$ and for $\kappa \geqslant 1$, the general expression of
the coefficient ${\bf Y }_\kappa^j$ of the prolongation~\thetag{ 4.3} 
of a vector field is{\rm :}
\def\theequation{4.19}\begin{equation}
\boxed{
\small
\aligned
{\bf Y}_\kappa^j
&
=
\mathcal{Y}_{x^\kappa}^j
+
\sum_{d=1}^{\kappa+1}
\
\sum_{1\leqslant\lambda_1<\cdots<\lambda_d\leqslant\kappa}
\
\sum_{\mu_1\geqslant 1,\dots,\mu_d\geqslant 1}
\
\sum_{\mu_1\lambda_1+\cdots+\mu_d\lambda_d\leqslant\kappa+1} 
\\
& \
\ \ \ \ \
\sum_{l_{1:1}=1}^m
\cdots
\sum_{l_{1:\mu_1}=1}^m
\cdots\cdots
\sum_{l_{d:1}=1}^m
\cdots
\sum_{l_{d:\mu_d}=1}^m
\
\frac{\kappa(\kappa-1)\cdots
(\kappa-\mu_1\lambda_1+\cdots+\mu_d\lambda_d+2)}{
(\lambda_1!)^{\mu_1}\ \mu_1!\cdots(\lambda_d!)^{\mu_d}\ \mu_d!
}
\\
& \
\left[
\aligned
&
(\kappa-\mu_1\lambda_1-\cdots-\mu_d\lambda_d+1)
\frac{\partial^{\kappa-\mu_1\lambda_1-\cdots-\mu_d\lambda_d+
\mu_1+\cdots+\mu_d}\mathcal{Y}^j}{
(\partial x)^{\kappa-\mu_1\lambda_1-\cdots-\mu_d\lambda_d}
\partial y^{l_{1:1}}
\cdots
\partial y^{l_{1:\mu_1}}
\cdots
\partial y^{l_{d:1}}
\cdots
\partial y^{l_{d:\mu_d}}}
- \\
& \
-
\delta_{l_{1:1}}^j\,\mu_1\lambda_1\,
\frac{\partial^{\kappa-\mu_1\lambda_1-\cdots-\mu_d\lambda_d+
\mu_1+\cdots+\mu_d}\mathcal{X}}{
(\partial x)^{\kappa-\mu_1\lambda_1-\cdots-\mu_d\lambda_d+1}
\widehat{\partial y^{l_{1:1}}}
\cdots
\partial y^{l_{1:\mu_1}}
\cdots
\partial y^{l_{d:1}}
\cdots
\partial y^{l_{d:\mu_d}}}
-
\\
& 
\ \ \ \ \ \ \ \ \
-
\cdots
-
\\
& \
-
\delta_{l_{d:1}}^j\,\mu_d\lambda_d\,
\frac{\partial^{\kappa-\mu_1\lambda_1-\cdots-\mu_d\lambda_d+
\mu_1+\cdots+\mu_d}\mathcal{X}}{
(\partial x)^{\kappa-\mu_1\lambda_1-\cdots-\mu_d\lambda_d+1}
\partial y^{l_{1:1}}
\cdots
\partial y^{l_{1:\mu_1}}
\cdots
\widehat{\partial y^{l_{d:1}}}
\cdots
\partial y^{l_{d:\mu_d}}}
\endaligned
\right]
\cdot
\\
& \
\ \ \ \ \ \ \ \ \ \ \ \ \ \ \ \ \ \ \ 
\ \ \ \ \ \ \ \ \ \ \ \ \ \ \ \ \ \ \ 
\cdot
\prod_{1\leqslant\nu_1\leqslant\mu_1}
y_{\lambda_1}^{l_{1:\nu_1}}
\ \cdots\cdots
\prod_{1\leqslant\nu_d\leqslant\mu_d}
y_{\lambda_d}^{l_{d:\nu_d}}.
\endaligned
}
\end{equation}
Here, the notation $\widehat{ \partial y^l}$ means 
that the partial derivative is dropped.
\end{theorem}

Since the fundamental monomials appearing in the last line of~\thetag{
4.19} just above are not independent of each other, this formula has
still to be modified a little bit. We refer to Section~6 for details.

\subsection*{ 4.20.~Deduction of a multivariate Fa\`a
di Bruno formula} Let $m \in \N$ with $m\geqslant 1$, let $y = (y^1,\dots,
y^m) \in \K^m$, let $f = f( y^1, \dots, y^m)$ be a $\mathcal{
C}^\infty$-smooth function from $\K^m$ to $\K$, let $x \in \K$ and let
$g^1 = g^1(x), \dots, g^m = g^m( x)$ be $\mathcal{ C}^\infty$
functions from $\K$ to $\K$. The goal is to obtain an explicit formula
for the derivatives, with respect to $x$, of the composition $h :=
f\circ g$, namely $h(x) := f \left( g^1(x), \dots, g^m(x)
\right)$. For $\lambda \in \N$ with $\lambda \geqslant 1$, and for $j= 1,
\dots, m$, we shall abbreviate the derivative $\frac{ d^\lambda g^j}{
dx^\lambda}$ by $g_\lambda^j$ and similarly for $h_\lambda$. The
partial derivatives $\frac{ \partial^\lambda f}{ \partial
y^{l_1}\cdots \partial y^{l_\lambda}}$ will be abbreviated by $f_{l_1,
\dots, l_\lambda }$.

Appying the chain rule, we may compute:
\def\theequation{4.21}\begin{equation}
\small
\aligned
h_1
& 
=
\sum_{l_1=1}^m\,f_{l_1}\,g_1^{l_1}, 
\\
h_2
& 
=
\sum_{l_1,l_2=1}^m\,f_{l_1,l_2}\,g_1^{l_1}\,g_1^{l_2}
+
\sum_{l_1=1}^m\,f_{l_1}\,g_2^{l_1}, 
\\
h_3
& 
=
\sum_{l_1,l_2,l_3=1}^m\,f_{l_1,l_2,l_3}\,
g_1^{l_1}\,g_1^{l_2}\,g_1^{l_3}
+
3\,\sum_{l_1,l_2=1}^m\,f_{l_1,l_2}\,
g_1^{l_1}\,g_2^{l_2}
+
\sum_{l_1=1}^m\,f_{l_1}\,g_3^{l_1},
\\
h_4
& 
=
\sum_{l_1,l_2,l_3,l_4=1}^m\,f_{l_1,l_2,l_3,l_4}\,
g_1^{l_1}\,g_1^{l_2}\,g_1^{l_3}\,g_1^{l_4}
+
6\,\sum_{l_1,l_2,l_3=1}^m\,f_{l_1,l_2,l_3}\,
g_1^{l_1}\,g_1^{l_2}\,g_2^{l_3}
+ \\
& \
\ \ \ \ \
+
3\,\sum_{l_1,l_2=1}^m\,f_{l_1,l_2}\,
g_2^{l_1}\,g_2^{l_2}
+
4\,\sum_{l_1,l_2=1}^m\,f_{l_1,l_2}\,
g_1^{l_1}\,g_3^{l_2}
+
\sum_{l_1=1}^m\,
f_{l_1}\,g_4^{l_1}, 
\\
h_5
& 
=
\sum_{l_1,l_2,l_3,l_4,l_5=1}^m\,f_{l_1,l_2,l_3,l_4,l_5}\,
g_1^{l_1}\,g_1^{l_2}\,g_1^{l_3}\,g_1^{l_4}\,g_1^{l_5}
+
10\,\sum_{l_1,l_2,l_3,l_4=1}^m\,f_{l_1,l_2,l_3,l_4}\,
g_1^{l_1}\,g_1^{l_2}\,g_1^{l_3}\,g_2^{l_4}
+ \\
& \
\ \ \ \ \
+
15\,\sum_{l_1,l_2,l_3=1}^m\,f_{l_1,l_2,l_3}\,
g_1^{l_1}\,g_2^{l_2}\,g_2^{l_3}
+
10\,\sum_{l_1,l_2,l_3=1}^m\,f_{l_1,l_2,l_3}\,
g_1^{l_1}\,g_1^{l_2}\,g_3^{l_3}
+ \\
& \
\ \ \ \ \
+
10\,\sum_{l_1,l_2=1}^m\,f_{l_1,l_2}\,
g_2^{l_1}\,g_3^{l_2}
+
5\,\sum_{l_1,l_2=1}^m\,f_{l_1,l_2}\,
g_1^{l_1}\,g_4^{l_2}
+
\sum_{l_1=1}^m\,
f_{l_1}\,g_5^{l_1}.
\endaligned
\end{equation}
Introducing the derivations
\def\theequation{4.22}\begin{equation}
\small
\aligned
F^2
&
:= 
\sum_{l_1=1}^m\,g_2^{l_1}\,
\frac{\partial}{\partial g_1^{l_1}}
+
\sum_{l_1=1}^m\,g_1^{l_1}
\left(
\sum_{l_2=1}^m\,
f_{l_1,l_2}\,\frac{\partial}{\partial f_{l_2}}
\right), \\
F^3
&
:=
\sum_{l_1=1}^m\,g_2^{l_1}\,
\frac{\partial}{\partial g_1^{l_1}}
+
\sum_{l_1=1}^m\,g_3^{l_1}\,
\frac{\partial}{\partial g_2^{l_1}}
+
\sum_{l_1=1}^m\,g_1^{l_1}
\left(
\sum_{l_2=1}^m\,
f_{l_1,l_2}\,\frac{\partial}{\partial f_{l_2}}
+
\sum_{l_2,l_3=1}^m\,
f_{l_1,l_2,l_3}\,\frac{\partial}{\partial f_{l_2,l_3}}
\right), 
\\
& \
\ \ \ \ \
\text{\bf 
\dots\dots\dots\dots\dots\dots\dots\dots\dots
\dots\dots\dots\dots\dots\dots\dots\dots\dots
\dots\dots\dots\dots\dots\dots\dots\dots\dots
} \\
F^\lambda
&
:=
\sum_{l_1=1}^m\,g_2^{l_1}\,
\frac{\partial}{\partial g_1^{l_1}}
+
\sum_{l_1=1}^m\,g_3^{l_1}\,
\frac{\partial}{\partial g_2^{l_1}}
+
\cdots
+
\sum_{l_1=1}^m\,g_\lambda^{l_1}\,
\frac{\partial}{\partial g_{\lambda-1}^{l_1}}
+ \\
& \
\ \ \ \ \
+
\sum_{l_1=1}^m\,g_1^{l_1}
\left(
\sum_{l_2=1}^m\,
f_{l_1,l_2}\,\frac{\partial}{\partial f_{l_2}}
+
\sum_{l_2,l_3=1}^m\,
f_{l_1,l_2,l_3}\,\frac{\partial}{\partial f_{l_2,l_3}}
+
\cdots
+
\sum_{l_2,\dots,l_\lambda=1}^m\,
f_{l_1,l_2,\dots,l_\lambda}\,
\frac{\partial}{\partial f_{l_2,\dots,l_\lambda}}
\right),
\endaligned
\end{equation}
we observe that the following
induction relations hold:
\def\theequation{4.23}\begin{equation}
\aligned
h_2
&
=
F^2
\left(
h_1
\right), 
\\
h_3
&
=
F^3
\left(
h_2
\right), 
\\
\text{\bf 
\dots\dots
}
& \
\ \ \ \ \
\text{\bf 
\dots\dots\dots
}
\\
h_\lambda
&
=
F^\lambda
\left(
h_{\lambda-1}
\right).
\endaligned
\end{equation}
To obtain the explicit version of the Fa\`a di Bruno in the case of
one variable $x$ and several variables $(y^1, \dots, y^m)$, it
suffices to extract from the expression of ${\bf Y}_\kappa^j$ provided
by Theorem~4.18 only the terms corresponding to $\mu_1 \lambda_1 +
\cdots + \mu_d\lambda_d = \kappa$, dropping all the $\mathcal{ X}$
terms. After some simplifications and after a translation by means of
an elementary dictionary, we may formulate a statement.

\def\thetheorem{4.24}\begin{theorem} 
For every integer $\kappa \geqslant 1$, the $\kappa$-th partial derivative
of the composite function $h = h( x) = f \left( 
g^1(x), \dots, g^m(x) \right)$ with
respect to $x$ may be expressed as an explicit polynomial depending on
the partial derivatives of $f$, on the derivatives of $g$ and having
integer coefficients{\rm:}
\def\theequation{4.25}\begin{equation}
\boxed{
\aligned
\frac{d^\kappa h}{dx^\kappa}
&
=
\sum_{d=1}^\kappa
\
\sum_{1\leqslant \lambda_1 < \cdots < \lambda_d \leqslant \kappa}
\
\sum_{\mu_1\geqslant 1,\dots,\mu_d\geqslant 1}
\
\sum_{\mu_1\lambda_1+\cdots+\mu_d\lambda_d=\kappa}
\
\frac{\kappa!}{
(\lambda_1!)^{\mu_1}\ \mu_1! 
\cdots
(\lambda_d!)^{\mu_d}\ \mu_d!
}
\\
& \
\ \ \ \ \ \ \ \ \ \ \ \ \ \ \
\sum_{l_{1:1},\dots,l_{1:\mu_1}=1}^m
\ \cdots \
\sum_{l_{d:1},\dots,l_{d:\mu_d}=1}^m
\\
\\
& \
\ \ \ \ \
\frac{\partial^{\mu_1+\cdots+\mu_d}f}{
\partial y^{l_{1:1}} 
\cdots
\partial y^{l_{1:\mu_1}}
\cdots
\partial y^{l_{d:1}} 
\cdots
\partial y^{l_{d:\mu_d}}
}
\
\prod_{1\leqslant\nu_1\leqslant\mu_1}
\frac{d^{\lambda_1} g^{l_{1:\nu_1}}}{d x^{\lambda_1}}
\ \cdots
\prod_{1\leqslant\nu_d\leqslant\mu_d}
\frac{d^{\lambda_d} g^{l_{d:\nu_d}}}{d x^{\lambda_d}}.
\endaligned
}
\end{equation}
\end{theorem} 

We refer to Section~6 for the final writing of the above
formula~\thetag{ 4.25}.

\section*{\S5.~Several independent variables and
several dependent variables}

\subsection*{5.1.~Expression of ${\bf Y}_{i_1}^j$, 
of ${\bf Y}_{i_1,i_2}^j$ and of ${\bf Y}_{i_1,i_2,i_3}^j$} Applying
the induction~\thetag{1.31} and working out the obtained formulas
until they take a perfect shape, we obtain firstly:
\def\theequation{5.2}\begin{equation}
\small
{\bf Y}_{i_1}^j
=
\mathcal{Y}_{x^{i_1}}^j
+
\sum_{l_1=1}^m\ \sum_{k_1=1}^n
\left[
\delta_{i_1}^{k_1}\,\mathcal{Y}_{y^{l_1}}^j
-
\delta_{l_1}^j\,
\mathcal{X}_{x^{i_1}}^{k_1}
\right]
y_{k_1}^{l_1}
+ 
\sum_{l_1,l_2=1}^m\ \sum_{k_1,k_2=1}^n
\left[
-
\delta_{l_2}^j\,
\delta_{i_1}^{k_1}\,\mathcal{X}_{y^{l_1}}^{k_2}
\right]
y_{k_1}^{l_1}y_{k_2}^{l_2}.
\end{equation}
Secondly:
\def\theequation{5.3}\begin{equation}
\small
\aligned{\bf Y}_{i_1,i_2}^j
&
=
\mathcal{Y}_{x^{i_1}x^{i_2}}^j
+
\sum_{l_1=1}^m\ \sum_{k_1=1}^n
\left[
\delta_{i_1}^{k_1}\,\mathcal{Y}_{x^{i_2}y^{l_1}}^j
+
\delta_{i_2}^{k_1}\,\mathcal{Y}_{x^{i_1}y^{l_1}}^j
-
\delta_{l_1}^j\,
\mathcal{X}_{x^{i_1}x^{i_2}}^{k_1}
\right]
y_{k_1}^{l_1}
+ \\
& \
\ \ \ \ \
+
\sum_{l_1,l_2=1}^m\ \sum_{k_1,k_2=1}^n
\left[
\delta_{i_1, \ i_2}^{k_1,k_2}\,
\mathcal{Y}_{y^{l_1}y^{l_2}}^j
-
\delta_{l_2}^j\,\delta_{i_1}^{k_1}\,
\mathcal{X}_{x^{i_2}y^{l_1}}^{k_2}
-
\delta_{l_2}^j\,\delta_{i_2}^{k_1}\,
\mathcal{X}_{x^{i_1}y^{l_1}}^{k_2}
\right]
y_{k_1}^{l_1}y_{k_2}^{l_2}
+ \\
& \
\ \ \ \ \
+
\sum_{l_1,l_2,l_3=1}^m\ \sum_{k_1,k_2,k_3=1}^n
\left[
-
\delta_{l_3}^j\,\delta_{i_1,\ i_2}^{k_1,k_2}\,
\mathcal{X}_{y^{l_1}y^{l_2}}^{k_3}
\right]
y_{k_1}^{l_1}y_{k_2}^{l_2}y_{k_3}^{l_3}
+ \\
& \
\ \ \ \ \
+
\sum_{l_1=1}^m\ \sum_{k_1,k_2=1}^n
\left[
\delta_{i_1,\ i_2}^{k_1,k_2}\,
\mathcal{Y}_{y^{l_1}}^j
-
\delta_{l_1}^j\,\delta_{i_1}^{k_1}\,
\mathcal{X}_{x^{i_2}}^{k_2}
-
\delta_{l_1}^j\,\delta_{i_2}^{k_1}\,
\mathcal{X}_{x^{i_1}}^{k_2}
\right]
y_{k_1,k_2}^{l_1}
+ \\
& \
\ \ \ \ \
+
\sum_{l_1,l_2=1}^m\ \sum_{k_1,k_2,k_3=1}^n
\left[
-
\delta_{l_1}^j\,\delta_{i_1,\ i_2}^{k_2,k_3}\,
\mathcal{X}_{y^{l_2}}^{k_1}
-
\delta_{l_2}^j\,\delta_{i_1,\ i_2}^{k_3,k_1}\,
\mathcal{X}_{y^{l_1}}^{k_2}
-
\delta_{l_2}^j\,\delta_{i_1,\ i_2}^{k_1,k_2}\,
\mathcal{X}_{y^{l_1}}^{k_3}
\right]
y_{k_1}^{l_1}y_{k_2}^{l_2}y_{k_3}^{l_3}.
\endaligned
\end{equation}
Thirdly:
$$
\small
\aligned
{\bf Y}_{i_1,i_2,i_3}^j
&
=
\mathcal{Y}_{x^{i_1}x^{i_2}x^{i_3}}^j
+
\sum_{l_1=1}^m\ \sum_{k_1=1}^n
\left[
\delta_{i_1}^{k_1}\,\mathcal{Y}_{x^{i_2}x^{i_3}y^{l_1}}^j
+
\delta_{i_2}^{k_1}\,\mathcal{Y}_{x^{i_1}x^{i_3}y^{l_1}}^j
+
\delta_{i_3}^{k_1}\,\mathcal{Y}_{x^{i_1}x^{i_2}y^{l_1}}^j
-
\delta_{l_1}^j\,
\mathcal{X}_{x^{i_1}x^{i_2}x^{i_3}}^{k_1}
\right]
y_{k_1}^{l_1}
+ \\
& \
\ \ \ \ \
+
\sum_{l_1,l_2=1}^m\ \sum_{k_1,k_2=1}^n
\left[
\delta_{i_1, \ i_2}^{k_1,k_2}\,
\mathcal{Y}_{x^{i_3}y^{l_1}y^{l_2}}^j
+
\delta_{i_3, \ i_1}^{k_1,k_2}\,
\mathcal{Y}_{x^{i_2}y^{l_1}y^{l_2}}^j
+
\delta_{i_2, \ i_3}^{k_1,k_2}\,
\mathcal{Y}_{x^{i_1}y^{l_1}y^{l_2}}^j
- 
\right. 
\\
& \
\ \ \ \ \ \ \ \ \ \ \ \ \ \ \ \ \ \ \ \ 
\ \ \ \ \ \ \ \ \ \ \ \ \ \ \ \ \
\left.
-
\delta_{l_2}^j\,\delta_{i_1}^{k_1}\,
\mathcal{X}_{x^{i_2}x^{i_3}y^{l_1}}^{k_2}
-
\delta_{l_2}^j\,\delta_{i_2}^{k_1}\,
\mathcal{X}_{x^{i_1}x^{i_3}y^{l_1}}^{k_2}
-
\delta_{l_2}^j\,\delta_{i_3}^{k_1}\,
\mathcal{X}_{x^{i_1}x^{i_2}y^{l_1}}^{k_2}
\right]
y_{k_1}^{l_1}y_{k_2}^{l_2}
+ \\
& \
\ \ \ \ \
+
\sum_{l_1,l_2,l_3=1}^m\ \sum_{k_1,k_2,k_3=1}^n
\left[
\delta_{i_1, \ i_2, \ i_3}^{k_1,k_2,k_3}\,
\mathcal{Y}_{y^{l_1}y^{l_2}y^{l_3}}^j
-
\delta_{l_3}^j\,\delta_{i_1,\ i_2}^{k_1,k_2}\,
\mathcal{X}_{x^{i_3}y^{l_1}y^{l_2}}^{k_3}
-
\right.
\\
& \
\ \ \ \ \ \ \ \ \ \ \ \ \ \ \ \ \ \ \ \ \
\ \ \ \ \ \ \ \ \ \ \ \ \ \ \ \ \ \ \ \ \
\ \ \ \ 
\left.
-
\delta_{l_3}^j\,\delta_{i_1,\ i_3}^{k_1,k_2}\,
\mathcal{X}_{x^{i_2}y^{l_1}y^{l_2}}^{k_3}
-
\delta_{l_3}^j\,\delta_{i_2,\ i_3}^{k_1,k_2}\,
\mathcal{X}_{x^{i_1}y^{l_1}y^{l_2}}^{k_3}
\right]
y_{k_1}^{l_1}y_{k_2}^{l_2}y_{k_3}^{l_3}
+
\endaligned
$$
$$
\small
\aligned
& \
\ \ \ \ \ 
+
\sum_{l_1,l_2,l_3,l_4=1}^m\ \sum_{k_1,k_2,k_3,k_4=1}^n
\left[
-
\delta_{l_4}^j\,\delta_{i_1,\ i_2,\ i_3}^{k_1,k_2,k_3}\,
\mathcal{X}_{y^{l_1}y^{l_2}y^{l_3}}^{k_4}
\right]
y_{k_1}^{l_1}y_{k_2}^{l_2}y_{k_3}^{l_3}y_{k_4}^{l_4}
+ \\
& \
\ \ \ \ \
+
\sum_{l_1=1}^m\ \sum_{k_1,k_2=1}^n
\left[
\delta_{i_1,\ i_2}^{k_1,k_2}\,
\mathcal{Y}_{x^{i_3}y^{l_1}}^j
+
\delta_{i_3,\ i_1}^{k_1,k_2}\,
\mathcal{Y}_{x^{i_2}y^{l_1}}^j
+
\delta_{i_2,\ i_3}^{k_1,k_2}\,
\mathcal{Y}_{x^{i_1}y^{l_1}}^j
- 
\right. \\
& \
\ \ \ \ \ \ \ \ \ \ \ \ \ \ \ \ \ \ \ \
\ \ \ \ \ \ \ \ \ \ \ \ \ \
\left.
-
\delta_{l_1}^j\,\delta_{i_1}^{k_1}\,
\mathcal{X}_{x^{i_2}x^{i_3}}^{k_2}
-
\delta_{l_1}^j\,\delta_{i_2}^{k_1}\,
\mathcal{X}_{x^{i_1}x^{i_3}}^{k_2}
-
\delta_{l_1}^j\,\delta_{i_3}^{k_1}\,
\mathcal{X}_{x^{i_1}x^{i_2}}^{k_2}
\right]
y_{k_1,k_2}^{l_1}
+ 
\endaligned
$$
\def\theequation{5.4}\begin{equation}
\small
\aligned
& \
\ \ \ \ \
+
\sum_{l_1,l_2=1}^m\ \sum_{k_1,k_2,k_3=1}^n
\left[
\delta_{i_1,\ i_2,\ i_3}^{k_1,k_2,k_3}\,
\mathcal{Y}_{y^{l_1}y^{l_2}}^j
+
\delta_{i_1,\ i_2,\ i_3}^{k_3,k_1,k_2}\,
\mathcal{Y}_{y^{l_1}y^{l_2}}^j
+
\delta_{i_1,\ i_2,\ i_3}^{k_2,k_3,k_1}\,
\mathcal{Y}_{y^{l_1}y^{l_2}}^j
-
\right. 
\\
& \
\ \ \ \ \ \ \ \ \ \ \ \ \ \ \ \ \ \ \ \
\ \ \ \ \ \ \ \ \ \ \ \ \ \ \ \ \ \ \ \
\ \ 
\left.
-
\delta_{l_1}^j\,\delta_{i_1,\ i_2}^{k_2,k_3}\,
\mathcal{X}_{x^{i_3}y^{l_2}}^{k_1}
-
\delta_{l_1}^j\,\delta_{i_1,\ i_3}^{k_2,k_3}\,
\mathcal{X}_{x^{i_2}y^{l_2}}^{k_1}
-
\delta_{l_1}^j\,\delta_{i_2,\ i_3}^{k_2,k_3}\,
\mathcal{X}_{x^{i_1}y^{l_2}}^{k_1}
-
\right. 
\\
& \
\ \ \ \ \ \ \ \ \ \ \ \ \ \ \ \ \ \ \ \
\ \ \ \ \ \ \ \ \ \ \ \ \ \ \ \ \ \ \ \
\ \ 
\left.
-
\delta_{l_2}^j\,\delta_{i_1,\ i_2}^{k_3,k_1}\,
\mathcal{X}_{x^{i_3}y^{l_1}}^{k_2}
-
\delta_{l_2}^j\,\delta_{i_1,\ i_3}^{k_3,k_1}\,
\mathcal{X}_{x^{i_2}y^{l_1}}^{k_2}
-
\delta_{l_2}^j\,\delta_{i_2,\ i_3}^{k_3,k_1}\,
\mathcal{X}_{x^{i_1}y^{l_1}}^{k_2}
-
\right. 
\\
& \
\ \ \ \ \ \ \ \ \ \ \ \ \ \ \ \ \ \ \ \
\ \ \ \ \ \ \ \ \ \ \ \ \ \ \ \ \ \ \ \
\ \ 
\left.
-
\delta_{l_2}^j\,\delta_{i_1,\ i_2}^{k_1,k_2}\,
\mathcal{X}_{x^{i_3}y^{l_1}}^{k_3}
-
\delta_{l_2}^j\,\delta_{i_1,\ i_3}^{k_1,k_2}\,
\mathcal{X}_{x^{i_2}y^{l_1}}^{k_3}
-
\delta_{l_2}^j\,\delta_{i_2,\ i_3}^{k_1,k_2}\,
\mathcal{X}_{x^{i_1}y^{l_1}}^{k_3}
\right]
y_{k_1}^{l_1}y_{k_2,k_3}^{l_2}
+ \\
\endaligned
\end{equation}
$$
\small
\aligned
& \
\ \ \ \ \
+
\sum_{l_1,l_2,l_3=1}^m\ \sum_{k_1,k_2,k_3,k_4=1}^n
\left[
-
\delta_{l_3}^j\,\delta_{i_1,\ i_2,\ i_3}^{k_1,k_2,k_3}\,
\mathcal{X}_{y^{l_1}y^{l_2}}^{k_4}
-
\delta_{l_3}^j\,\delta_{i_1,\ i_2,\ i_3}^{k_2,k_3,k_1}\,
\mathcal{X}_{y^{l_1}y^{l_2}}^{k_4}
-
\delta_{l_3}^j\,\delta_{i_1,\ i_2,\ i_3}^{k_3,k_2,k_1}\,
\mathcal{X}_{y^{l_1}y^{l_2}}^{k_4}
- 
\right. 
\\
& \
\ \ \ \ \ \ \ \ \ \ \ \ \ \ \ \ \ \ \ \
\ \ \ \ \ \ \ \ \ \ \ \ \ \ \ \ \ \ \ \ 
\ \ \ \ \ \ \ \ \ \
\left.
-
\delta_{l_2}^j\,\delta_{i_1,\ i_2,\ i_3}^{k_3,k_4,k_1}\,
\mathcal{X}_{y^{l_1}y^{l_3}}^{k_2}
-
\delta_{l_2}^j\,\delta_{i_1,\ i_2,\ i_3}^{k_3,k_1,k_4}\,
\mathcal{X}_{y^{l_1}y^{l_3}}^{k_2}
-
\right. 
\\
& \
\ \ \ \ \ \ \ \ \ \ \ \ \ \ \ \ \ \ \ \
\ \ \ \ \ \ \ \ \ \ \ \ \ \ \ \ \ \ \ \ 
\ \ \ \ \ \ \ \ \ \ \ \ \ \ \ \ \ \ \ \
\ \ \ \ \ \ \ \ \ \ \ \ \ \ \ \ \ \ \ \
\ \ \ \ \ \ \ \ \ \ \ \ \ \ \ \ \ \ \ \ 
\left.
-
\delta_{l_2}^j\,\delta_{i_1,\ i_2,\ i_3}^{k_1,k_3,k_4}\,
\mathcal{X}_{y^{l_1}y^{l_3}}^{k_2}
\right]
y_{k_1}^{l_1}y_{k_2}^{l_2}y_{k_3,k_4}^{l_3}
+ \\
& \
\ \ \ \ \
+
\sum_{l_1,l_2=1}^m\ \sum_{k_1,k_2,k_3,k_4=1}^n
\left[
-
\delta_{l_2}^j\,\delta_{i_1,\ i_2,\ i_3}^{k_1,k_2,k_3}\,
\mathcal{X}_{y^{l_1}}^{k_3}
-
\delta_{l_2}^j\,\delta_{i_1,\ i_2,\ i_3}^{k_2,k_4,k_1}\,
\mathcal{X}_{y^{l_1}}^{k_3}
-
\delta_{l_2}^j\,\delta_{i_1,\ i_2,\ i_3}^{k_4,k_1,k_2}\,
\mathcal{X}_{y^{l_1}}^{k_3}
\right]
y_{k_1,k_2}^{l_1}y_{k_3,k_4}^{l_2}
+ \\
\endaligned
$$
$$
\small
\aligned
& \
\ \ \ \ \
+
\sum_{l_1=1}^m\ \sum_{k_1,k_2,k_3=1}^n
\left[
\delta_{i_1,\ i_2, \i_3}^{k_1,k_2,k_3}\,
\mathcal{Y}_{y^{l_1}}^j
-
\delta_{l_1}^j\,\delta_{i_1,\ i_2}^{k_1,k_2}\,
\mathcal{X}_{x^{i_3}}^{k_3}
-
\delta_{l_1}^j\,\delta_{i_1,\ i_3}^{k_1,k_2}\,
\mathcal{X}_{x^{i_2}}^{k_3}
-
\delta_{l_1}^j\,\delta_{i_2,\ i_3}^{k_1,k_2}\,
\mathcal{X}_{x^{i_1}}^{k_3}
\right]
y_{k_1,k_2,k_3}^{l_1}
+ \\
& \
\ \ \ \ \
+
\sum_{l_1,l_2=1}^m\ \sum_{k_1,k_2,k_3,k_4=1}^n
\left[
-
\delta_{l_2}^j\,\delta_{i_1,\ i_2,\ i_3}^{k_1,k_2,k_3}\,
\mathcal{X}_{y^{l_1}}^{k_4}
-
\delta_{l_2}^j\,\delta_{i_1,\ i_2,\ i_3}^{k_4,k_1,k_2}\,
\mathcal{X}_{y^{l_1}}^{k_3}
-
\delta_{l_2}^j\,\delta_{i_1,\ i_2,\ i_3}^{k_3,k_4,k_1}\,
\mathcal{X}_{y^{l_1}}^{k_2}
-
\right.
\\
& \
\ \ \ \ \ \ \ \ \ \ \ \ \ \ \ \ \ \ \ \
\ \ \ \ \ \ \ \ \ \ \ \ \ \ \ \ \ \ \ \
\ \ \ \ \ \ \ \
\left.
-
\delta_{l_1}^j\,\delta_{i_1,\ i_2,\ i_3}^{k_2,k_3,k_4}\,
\mathcal{X}_{y^{l_2}}^{k_1}
\right]
y_{k_1}^{l_1}y_{k_2,k_3,k_4}^{l_2}.
\endaligned
$$

\subsection*{5.5.~Final synthesis}
To obtain the general formula for ${\bf Y}_{ i_1, \dots, i_\kappa}^j$,
we have to achieve the synthesis between the two formulas~\thetag{
3.74} and~\thetag{ 4.19}. We start with~\thetag{ 3.74} and we modify
it until we reach the final formula for ${\bf Y}_{ i_1, \dots,
i_\kappa }^j$.

We have to add the $\mu_1+\cdots+\mu_d$ 
sums $\sum_{ l_{ 1:1} =1 }^m \cdots \sum_{
l_{ 1: \mu_1 } =1 }^m \cdots \cdots \sum_{ l_{ d:1} =1}^m \cdots
\sum_{l_{ d: \mu_d }= 1}^m$, together with various indices
$l_\alpha$. About these indices, the only point which is not obvious
may be analyzed as follows.

A permutation $\sigma \in \mathfrak{ F}_{
\mu_1\lambda_1 + \cdots + \mu_d \lambda_d }^{ (\mu_1, \lambda_1),
\dots, (\mu_d, \lambda_d)}$ yields the list:
\def\theequation{5.6}\begin{equation}
\aligned
& \
\sigma(1\!:\!1\!:\!1),\dots,\sigma(1\!:\!1\!:\!\lambda_1),
\dots
\sigma(1\!:\!\mu_1\!:\!1),\dots,\sigma(1\!:\!\mu_1\!:\!\lambda_1),
\dots
\\
& \
\ \ \ \ \ \ \ \ \ 
\dots,
\sigma(d\!:\!1\!:\!1),\dots,\sigma(1\!:\!1\!:\!\lambda_d),
\dots
\sigma(d\!:\!\mu_d\!:\!1),\dots,\sigma(d\!:\!\mu_d\!:\!\lambda_d),
\endaligned
\end{equation} 
In the sixth line of~\thetag{ 3.74}, the last term $\sigma ( d\! : \!
\mu_d \! : \! \lambda_d )$ of the above list appears as the subscript
of the upper index $k_{ \sigma (d: \mu_d: \lambda_d )}$ of the term
$\mathcal{ X }^{ k_{ \sigma( d:\mu_d: \lambda_d) }}$. According to the
formal rules of Theorem~4.19, we have to multiply the partial
derivative of $\mathcal{ X }^{ k_{ \sigma( d:\mu_d: \lambda_d ) }}$ by
a certain Kronecker symbol $\delta_{ l_\alpha}^j$. The question is:
{\it what is the subscript $\alpha$ and how
to denote it}?

As explained before the statement of Theorem~4.19, the subscript
$\alpha$ is obtained as follows. The term $\sigma ( d\! : \! \mu_d \!
: \! \lambda_d )$ is of the form $(e\! : \! \nu_d \! : \! \gamma_e
)$, for some $e$ with $1\leqslant e \leqslant d$, for some $\nu_e$ with $1\leqslant
\nu_e \leqslant \mu_e$ and for some $\gamma_e$ with $1\leqslant \gamma_e \leqslant
\lambda_e$. The single pure jet variable
\def\theequation{5.7}\begin{equation}
\aligned
y_{k_{e:\nu_e:1},\dots,
k_{e:\nu_e:\gamma_e},\dots,
k_{e:\nu_e:\lambda_e}}^{l_{e:\nu_e}}
\endaligned
\end{equation}
appears inside the total monomial
\def\theequation{5.8}\begin{equation}
\aligned
\prod_{1\leqslant\nu_1\leqslant\mu_1}\,
y_{k_{1:\nu_1:1},\dots,k_{1:\nu_1:\lambda_1}}^{l_{1:\nu_1}}
\ \cdots \
\prod_{1\leqslant\nu_d\leqslant\mu_d}\,
y_{k_{d:\nu_d:1},\dots,k_{d:\nu_d:\lambda_d}}^{l_{d:\nu_d}},
\endaligned
\end{equation}
placed at the end of the formula for ${\bf Y }_{ i_1, \dots,
i_\kappa}^j$ ({\it see} in advance formula~\thetag{ 5.13} below; this
total monomial generalizes to several dependent variables the total
monomial appearing in the last line of~\thetag{ 3.74}). According to
the rule explained before the statement of Theorem~4.18, the index
$l_\alpha$ must be equal to $l_{e : \nu_e }$, since
$l_{e : \nu_e }$ is attached to the monomial~\thetag{ 5.7}.
Coming back to the term
$\sigma ( d\! : \! \mu_d \! : \! \lambda_d )$, we shall denote this
index by
\def\theequation{5.9}\begin{equation}
\aligned
l_{e:\nu_e}
=:
l_{\pi(e:\nu_e:\gamma_e)}
=:
l_{\pi\sigma(d:\mu_d:\lambda_d)},
\endaligned
\end{equation}
where the symbol $\pi$ denotes the projection
from the set 
\def\theequation{5.10}\begin{equation}
\aligned
\{
1\!:\!1\!:\!1,\dots,1\!:\!\mu_1\!:\!\lambda_1,
\dots\dots,
d\!:\!1\!:\!1,\dots,d\!:\!\mu_d\!:\!\lambda_d
\}
\endaligned
\end{equation}
to the set
\def\theequation{5.11}\begin{equation}
\aligned
\{
1\!:\!1,\dots,1\!:\!\mu_1,
\dots,
d\!:\!1,\dots,d\!:\!\mu_d
\}
\endaligned
\end{equation}
simply defined by $\pi(e \! : \! \nu_e \! : \! \gamma_e) := (e \! : \!
\nu_e)$.

In conclusion, by means of this formalism, we may write down the
complete generalization of Theorems~2.24, 3.73 and~4.18 to several
independent variables and several dependent variables

\def\thetheorem{5.12}\begin{theorem}
For $j = 1, \dots, m$, for every $\kappa \geqslant 1$ and for every choice
of $\kappa$ indices $i_1,\dots, i_\kappa$ in the set $\{ 1, 2, \dots,
n\}$, the general expression of ${\bf Y}_{i_1, \dots, i_\kappa }^j$ is
as follows{\rm :}
\def\theequation{5.13}\begin{equation}
\small
\boxed{
\aligned
{\bf Y}_{i_1, \dots, i_\kappa}^j
&
=
\mathcal{Y}_{x^{i_1}\cdots x^{i_\kappa}}^j
+
\sum_{d=1}^{\kappa+1}
\ \
\sum_{1\leqslant\lambda_1<\cdots<\lambda_d\leqslant\kappa}
\ \
\sum_{\mu_1\geqslant 1,\dots,\mu_d\geqslant 1} 
\
\sum_{
\mu_1\lambda_1
+
\cdots
+
\mu_d\lambda_d\leqslant \kappa+1} 
\\
& \
\ \ \ \ \ \ \ \ \ \ \ \ \ \ \ \ \ \ \ \ 
\ \ \ \ \ \ \ \ \ \ \ \ \ 
\sum_{l_{1:1}=1}^m
\cdots
\sum_{l_{1:\mu_1}=1}^m
\cdots\cdots
\sum_{l_{d:1}=1}^m
\cdots
\sum_{l_{d:\mu_d}=1}^m
\\
&
\sum_{k_{1:1:1},\dots,k_{1:1:\lambda_1}=1}^n
\cdots \
\sum_{k_{1:\mu_1:1},\dots,k_{1:\mu_1:\lambda_1}=1}^n
\cdots\cdots \
\sum_{k_{d:1:1},\dots,k_{d:1:\lambda_d}=1}^n
\cdots \
\sum_{k_{d:\mu_d:1},\dots,k_{d:\mu_d:\lambda_d}=1}^n
\\
&
\left[
\aligned
& 
\sum_{\sigma\in\mathfrak{F}_{\mu_1\lambda_1+\cdots+\mu_d\lambda_d}^{
(\mu_1,\lambda_1),\dots,(\mu_d,\lambda_d)}}
\
\sum_{\tau\in\mathfrak{S}_\kappa^{
\mu_1\lambda_1+\cdots+\mu_d\lambda_d}}\,
\delta_{i_{\tau(1)},\dots,i_{\tau(\mu_1\lambda_1)},\dots,
i_{\tau(\mu_1\lambda_1+\cdots+\mu_d\lambda_d)}}^{
k_{\sigma(1:1:1)},\dots,k_{\sigma(1:\mu_1:\lambda_1)},
\dots,k_{\sigma(d:\mu_d:\lambda_d)}}
\cdot
\\
& \
\ \ \ \ \
\cdot
\frac{\partial^{\kappa-\mu_1\lambda_1-\cdots-\mu_d\lambda_d+
\mu_1+\cdots+\mu_d}
\mathcal{Y}^j}{
\partial x^{i_{\tau(\mu_1\lambda_1+\cdots+\mu_d\lambda_d+1)}}\cdots
\partial x^{i_{\tau(\kappa)}}
\partial y^{l_{1:1}}\cdots\partial y^{l_{d:\mu_d}}}\
- \\
& 
-
\sum_{\sigma\in\mathfrak{F}_{\mu_1\lambda_1+\cdots+\mu_d\lambda_d}^{
(\mu_1,\lambda_1),\dots,(\mu_d,\lambda_d)}}
\
\sum_{\tau\in\mathfrak{S}_\kappa^{
\mu_1\lambda_1+\cdots+\mu_d\lambda_d-1}}\,
\delta_{i_{\tau(1)},\dots,i_{\tau(\mu_1\lambda_1)},\dots,
i_{\tau(\mu_1\lambda_1+\cdots+\mu_d\lambda_d-1)}}^{
k_{\sigma(1:1:1)},\dots,k_{\sigma(1:\mu_1:\lambda_1)},
\dots,k_{\sigma(d:\mu_d:\lambda_d-1)}}
\cdot
\\
& \
\ \ \ \ \
\cdot
\delta_{l_{\pi\sigma(d:\mu_d:\lambda_d)}}^j
\cdot
\frac{\partial^{\kappa-\mu_1\lambda_1-\cdots-\mu_d\lambda_d
+\mu_1+\cdots+\mu_d}\mathcal{X}^{k_{\sigma(d:\mu_d:\lambda_d)}}}{
\partial x^{i_{\tau(\mu_1\lambda_1+\cdots+\mu_d\lambda_d)}}\cdots
\partial x^{i_{\tau(\kappa)}}
\partial y^{l_{1:1}}\cdots
\widehat{\partial y^{l_{\pi\sigma(d:\mu_d:\lambda_d)}}}
\cdots\partial y^{l_{d:\mu_d}}}
\endaligned
\right]
\cdot
\\
& \
\ \ \ \ \ \ \ \ \ \ \ \ \ \ \ \ \ \ \
\cdot
\prod_{1\leqslant\nu_1\leqslant\mu_1}\,
y_{k_{1:\nu_1:1},\dots,k_{1:\nu_1:\lambda_1}}^{l_{1:\nu_1}}
\ \cdots \
\prod_{1\leqslant\nu_d\leqslant\mu_d}\,
y_{k_{d:\nu_d:1},\dots,k_{d:\nu_d:\lambda_d}}^{l_{d:\nu_d}}.
\endaligned
}
\end{equation}
\end{theorem}

In this formula, the coset $\mathfrak{ F }_{ \mu_1 \lambda_1 + \cdots
+ \mu_d \lambda_d }^{ ( \mu_1, \lambda_1 ),\dots, ( \mu_d,
\lambda_d)}$ was defined in equation~\thetag{ 3.71};
as in Theorem~3.73, we have made the
identification{\rm :}
\def\theequation{5.14}\begin{equation}
\{1,\dots,\kappa\}
\equiv
\{
1\!:\!1\!:\!1,
\dots,
1\!:\!\mu_1\!:\!\lambda_1,
\text{\rm \dots\dots}, 
d\!:\!1\!:\!1,
\dots,
d\!:\!\mu_d\!:\!\lambda_d
\}.
\end{equation}

Since the fundamental monomials appearing in the last line of~\thetag{
4.19} just above are not independent of each other, this formula has
still to be modified a little bit. We refer to Section~6 for details.

\subsection*{ 5.15.~Deduction of the 
most general multivariate Fa\`a di Bruno formula} Let $n \in \N$ with
$n\geqslant 1$, let $x = (x^1,\dots, x^n) \in \K^n$, let $m\in \N$ with
$m\geqslant 1$, let $g^j = g^j ( x^1, \dots, x^n)$, $j=1, \dots, m$, be
$\mathcal{ C}^\infty$-smooth functions from $\K^n$ to $\K^m$, let $y =
(y^1, \dots, y^m) \in \K^m$ and let $f = f(y^1, \dots, y^m)$ be a
$\mathcal{ C}^\infty$ function from $\K^m$ to $\K$. The goal is to
obtain an explicit formula for the partial derivatives of the
composition $h := f\circ g$, namely
\def\theequation{5.16}\begin{equation}
h(x^1, \dots, x^n) := f
\left(
g^1(x^1,\dots, x^n), \dots, g^m(x^1,\dots, x^n) 
\right).
\end{equation}
For $j= 1,\dots, m$, for $\lambda \in \N$ with $\lambda \geqslant 1$ and
for arbitrary indices $i_1, \dots, i_\lambda = 1, \dots, n$, we shall
abbreviate the partial derivative $\frac{ \partial^\lambda g^j}{
\partial x^{i_1} \cdots \partial x^{i_\lambda}}$ by $g_{i_1,\dots,
i_\lambda }^j$ and similarly for $h_{i_1, \dots, i_\lambda}$. For
arbitrary indices $l_1, \dots, l_\lambda = 1, \dots, m$, the partial
derivative $\frac{ \partial^\lambda f}{ \partial y^{l_1} \cdots
\partial y^{l_\lambda }}$ will be abbreviated by $f_{ l_1, \dots,
l_\lambda }$.

Appying the chain rule, we may compute:
\def\theequation{5.17}\begin{equation}
\small
\aligned
h_{i_1}
&
=
\sum_{l_1=1}^m\,f_{l_1}
\left[
g_{i_1}^{l_1}
\right],
\\
h_{i_1,i_2}
&
=
\sum_{l_1,l_2=1}^m\,f_{l_1,l_2}
\left[
g_{i_1}^{l_1}\,g_{i_2}^{l_2}
\right]
+
\sum_{l_1=1}^m\,f_{l_1}
\left[
g_{i_1,i_2}^{l_1}
\right],
\\
h_{i_1,i_2,i_3}
&
=
\sum_{l_1,l_2,l_3=1}^m\,f_{l_1,l_2,l_3}
\left[
g_{i_1}^{l_1}\,g_{i_2}^{l_2}\,g_{i_3}^{l_3}
\right]
+
\sum_{l_1,l_2=1}^m\,f_{l_1,l_2}
\left[
g_{i_1}^{l_1}\,g_{i_2,i_3}^{l_2}
+
g_{i_2}^{l_1}\,g_{i_1,i_3}^{l_2}
+
g_{i_3}^{l_1}\,g_{i_1,i_2}^{l_2}
\right]
+ \\
& \
\ \ \ \ \
+
\sum_{l_1=1}^m\,f_{l_1}
\left[
g_{i_1,i_2,i_3}^{l_1}
\right],
\endaligned
\end{equation}
$$
\small
\aligned
h_{i_1,i_2,i_3,i_4}
&
=
\sum_{l_1,l_2,l_3,l_4=1}^m\,f_{l_1,l_2,l_3,l_4}
\left[
g_{i_1}^{l_1}\,g_{i_2}^{l_2}\,g_{i_3}^{l_3}\,g_{i_4}^{l_4}
\right]
+ \\
& \
\ \ \ \ \ 
\sum_{l_1,l_2,l_3=1}^m\,f_{l_1,l_2,l_3}
\left[
g_{i_2}^{l_1}\,g_{i_3}^{l_2}\,g_{i_1,i_4}^{l_3}
+
g_{i_3}^{l_1}\,g_{i_1}^{l_2}\,g_{i_2,i_4}^{l_3}
+
g_{i_1}^{l_1}\,g_{i_2}^{l_2}\,g_{i_3,i_4}^{l_3}
+
\right.
\\
& \
\ \ \ \ \ \ \ \ \ \ \ \ \ \ \ \ \ \ \ 
\ \ \ \ \ \ \ \ \ \ \ \ \ \ \ \ \ \
\left.
+
g_{i_1}^{l_1}\,g_{i_4}^{l_2}\,g_{i_2,i_3}^{l_3}
+
g_{i_2}^{l_1}\,g_{i_4}^{l_2}\,g_{i_3,i_1}^{l_3}
+
g_{i_3}^{l_1}\,g_{i_4}^{l_2}\,g_{i_1,i_2}^{l_3}
\right]
+ \\
& \
\ \ \ \ \
+
\sum_{l_1,l_2=1}^m\,f_{l_1,l_2}
\left[
g_{i_1,i_2}^{l_1}\,g_{i_3,i_4}^{l_2}
+
g_{i_1,i_3}^{l_1}\,g_{i_2,i_4}^{l_2}
+
g_{i_1,i_4}^{l_1}\,g_{i_2,i_3}^{l_2}
\right]
+ \\
& \
\ \ \ \ \
+
\sum_{l_1,l_2=1}^m\,f_{l_1,l_2}
\left[
g_{i_1}^{l_1}\,g_{i_2,i_3,i_4}^{l_2}
+
g_{i_2}^{l_1}\,g_{i_1,i_3,i_4}^{l_2}
+
g_{i_3}^{l_1}\,g_{i_1,i_2,i_4}^{l_2}
+
g_{i_4}^{l_1}\,g_{i_1,i_2,i_3}^{l_2}
\right]
+ \\
& \
\ \ \ \ \
+
\sum_{l_1=1}^m\,f_{l_1}
\left[
g_{i_1,i_2,i_3,i_4}^{l_1}
\right].
\endaligned
$$
Introducing the derivations
\begin{small}
\[
\aligned
F_i^2
&
:= 
\sum_{k_1=1}^n\,\sum_{l_1=1}^m\,g_{k_1,i}^{l_1}\,
\frac{\partial}{\partial g_{k_1}^{l_1}}
+
\sum_{l_1=1}^m\,g_i^{l_1}\left(
\sum_{l_2=1}^m\,f_{l_1,l_2}\,\frac{\partial}{\partial f_{l_2}}
\right), \\
F_i^3
&
:=
\sum_{k_1=1}^n\,\sum_{l_1=1}^m\,g_{k_1,i}^{l_1}\,
\frac{\partial}{\partial g_{k_1}^{l_1}}
+
\sum_{k_1,k_2=1}^n\,\sum_{l_1=1}^m\,g_{k_1,k_2,i}^{l_1}\,
\frac{\partial}{\partial g_{k_1,k_2}^{l_1}}
+ \\
& \
\ \ \ \ \ 
+
\sum_{l_1=1}^m\,g_i^{l_1}
\left(
\sum_{l_2=1}^m\,f_{l_1,l_2}\,\frac{\partial}{\partial f_{l_2}}
+
\sum_{l_2,l_3=1}^m\,f_{l_1,l_2,l_3}\,
\frac{\partial}{\partial f_{l_2,l_3}}
\right),
\endaligned
\]
\def\theequation{5.18}\begin{equation}
\\
\text{\bf 
\dots\dots}
\ \ \ \ \
\text{\bf 
\dots\dots\dots\dots\dots\dots\dots\dots\dots
\dots\dots\dots\dots\dots\dots\dots\dots\dots
\dots
}
\end{equation}
\[
\aligned
F_i^\lambda
&
:=
\sum_{k_1=1}^n\,\sum_{l_1=1}^m\,g_{k_1,i}^{l_1}\,
\frac{\partial}{\partial g_{k_1}^{l_1}}
+
\sum_{k_1,k_2=1}^n\,\sum_{l_1=1}^m\,g_{k_1,k_2,i}^{l_1}\,
\frac{\partial}{\partial g_{k_1,k_2}^{l_1}}
+ 
\cdots
+ \\
& \
\ \ \ \ \ 
+
\sum_{k_1,k_2,\dots,k_{\lambda-1}=1}^n\,\sum_{l_1=1}^m\,
g_{k_1,k_2,\dots,k_{\lambda-1},i}
\
\frac{\partial}{\partial g_{k_1,\dots,k_{\lambda-1}}^{l_1}}
+ \\
& \
\ \ \ \ \ 
+
\sum_{l_1=1}^m\,g_i^{l_1}
\left(
\sum_{l_2=1}^m\,f_{l_1,l_2}\,\frac{\partial}{\partial f_{l_2}}
+
\sum_{l_2,l_3=1}^m\,f_{l_1,l_2,l_3}\,
\frac{\partial}{\partial f_{l_2,l_3}}
+ 
\right. 
\\
& \
\ \ \ \ \ \ \ \ \ \ \ \ \ \ \ \ \ \ \ \ 
\ \ \ \ \ \ \ \ \ \
\left.
+
\cdots
+
\sum_{l_2,l_3,\dots,l_\lambda}\,
f_{l_1,l_2,l_3,\dots,l_\lambda}\,
\frac{\partial}{\partial f_{l_2,l_3,\dots,l_\lambda}}
\right),
\endaligned
\]
\end{small}

\noindent
we observe that the following
induction relations hold:
\def\theequation{5.19}\begin{equation}
\aligned
h_{i_1,i_2}
&
=
F_{i_2}^2
\left(
h_{i_1}
\right), 
\\
h_{i_1,i_2,i_3}
&
=
F_{i_3}^3
\left(
h_{i_1,i_2}
\right), 
\\
\text{\bf 
\dots\dots
}
& \
\ \ \ \ \
\text{\bf 
\dots\dots\dots\dots\dots
}
\\
h_{i_1,i_2,\dots,i_\lambda}
&
=
F_{i_\lambda}^\lambda
\left(
h_{i_1,i_2,\dots,i_{\lambda-1}}
\right).
\endaligned
\end{equation}

To obtain the explicit version of the Fa\`a di Bruno in the case of
several variables $(x^1, \dots, x^n)$ and several variables $(y^1,
\dots, y^m)$, it suffices to extract from the expression of ${\bf Y
}_{ i_1,\dots, i_\kappa }^j$ provided by Theorem~5.12 only the terms
corresponding to $\mu_1 \lambda_1 + \cdots + \mu_d\lambda_d = \kappa$,
dropping all the $\mathcal{ X}$ terms. After some simplifications and
after a translation by means of an elementary dictionary, we obtain
the fourth and the most general multivariate Fa\`a di Bruno formula.

\def\thetheorem{5.20}\begin{theorem} 
For every integer $\kappa \geqslant 1$ and for every choice of indices
$i_1, \dots, i_\kappa$ in the set $\{ 1, 2, \dots, n\}$, the
$\kappa$-th partial derivative of the composite function 
\def\theequation{5.21}\begin{equation}
h = h( x^1,
\dots, x^n) = 
f \left( g^1 (x^1, \dots, x^n),\dots, g^m (x^1, \dots,
x^n) \right) 
\end{equation}
with respect to the variables $x^{i_1}, \dots,
x^{i_\kappa}$ may be expressed as an explicit polynomial depending on
the partial derivatives of $f$, on the partial derivatives of the
$g^j$ and having integer coefficients{\rm :}
\def\theequation{5.22}\begin{equation}
\boxed{
\aligned
\frac{\partial^\kappa h}{\partial x^{i_1}\cdots
\partial x^{i_\kappa}}
&
=
\sum_{d=1}^\kappa
\
\sum_{1\leqslant \lambda_1 < \cdots < \lambda_d \leqslant \kappa}
\
\sum_{\mu_1\geqslant 1,\dots,\mu_d\geqslant 1}
\
\sum_{\mu_1\lambda_1+\cdots+\mu_d\lambda_d=\kappa}
\\
& \
\ \ \ \ \
\sum_{l_{1:1},\dots,l_{1:\mu_1}=1}^m
\cdots
\sum_{l_{d:1},\dots,l_{d:\mu_d}=1}^m
\
\frac{\partial^{\mu_1+\cdots+\mu_d}f}{
\partial y^{l_{1:1}}\cdots
\partial y^{l_{1:\mu_1}}\cdots 
\partial y^{l_{d:1}}\cdots
\partial y^{l_{d:\mu_d}}
}
\\
& 
\ \ \ \ \ \ \ \ \ \ 
\left[
\aligned
&
\sum_{\sigma\in\mathfrak{F}_\kappa^{
(\mu_1,\lambda_1),\dots,(\mu_d,\lambda_d)}}
\
\prod_{1\leqslant\nu_1\leqslant\mu_1}
\
\frac{\partial^{\lambda_1} g^{l_{1:\nu_1}}}{\partial
x^{i_{\sigma(1:\nu_1:1)}}\cdots
\partial x^{i_{\sigma(1:\nu_1:\lambda_1)}}}
\
\text{\bf \dots} 
\\
& \
\ \ \ \ \ \ \ \ \ \ \ \ \ \ 
\ \ \ \ \ \ \ \ \ \ \ \ \ \
\text{\bf \dots} 
\prod_{1\leqslant\nu_d\leqslant\mu_d}
\
\frac{\partial^{\lambda_d} g^{l_{d:\nu_d}}}{\partial
x^{i_{\sigma(d:\nu_d:1)}}\cdots
\partial x^{i_{\sigma(d:\nu_d:\lambda_d)}}}
\endaligned
\right].
\endaligned
}
\end{equation}
\end{theorem}

\newpage

$\:$\bigskip\bigskip

\begin{center}
{\large\bf
{\Large\bf III:}~Systems of second order
}
\end{center}

\begin{center}
\begin{minipage}[t]{12cm}
\baselineskip =0.35cm
{\scriptsize

\centerline{\bf Table of contents}

\smallskip

{\bf 1.~Explicit characterizations of flatness \dotfill 96.}

{\bf 2.~Completely integrable systems of second order ordinary
differential equations \dotfill 100.}

{\bf 3.~First and second auxiliary systems \dotfill 107.} 

}\end{minipage}
\end{center}

\bigskip
\bigskip

\section*{\S1.~Explicit characterizations of flatness}

In 1883, S.~Lie obtained the following explicit characterization of
the local equivalence of a second order ordinary differential equation
\thetag{ $\mathcal{ E}_1$}: 
$y_{xx } = F( x, y, y_x)$ to the Newtonian free particle equation with
one degree of freedom $Y_{X X } = 0$. All the functions are assumed to
be analytic.

\def\thetheorem{1.1}\begin{theorem}
{\rm (\cite{ lie1883}, pp.~362--365)} 
Let $\K=\R$ of $\C$. Let $x\in \K$ and $y\in \K$. A local second
order ordinary differential equation $y_{xx } = F( x, y, y_x)$ is
equivalent under an invertible point transformation $(x, y)\mapsto (X
(x,y), Y (x,y))$ to the free particle equation $Y_{ XX } =0$ if and
only if the following two conditions are satisfied{\rm :}

\begin{itemize}

\smallskip\item[{\bf (i)}]
$F_{y_x y_x y_x y_x }=0$, or equivalently $F$ is a degree three
polynomial in $y_x$, namely there exist four functions $G$, $H$, $L$
and $M$ of $(x,y)$ such that $F$ can be written as
\def\theequation{1.2}\begin{equation}
F(x,y,y_x)
=
G(x,y)+y_x\cdot H(x,y)+
(y_x)^2\cdot L(x,y)+(y_x)^3\cdot M(x,y);
\end{equation}

\smallskip\item[{\bf (ii)}]
the four functions $G$, $H$, $L$ and $M$ satisfy the following system
of two second order quasi-linear partial differential
equations{\rm :}
\def\theequation{1.3}\begin{equation}
\left\{
\aligned
0 = 
-2\, G_{yy} & \
+\frac{ 4}{3}\, H_{xy}
-
\frac{ 2}{3} \, L_{xx}+ \\
& \
+ 2\, (G\, L)y-2\, G_x\, M-4\, G\, M_x+\frac{ 2}{3}
\, H\, L_x-\frac{ 4}{3}\, H\, H_y, \\
0 = 
-\frac{ 2}{3}\, H_{yy} & \
+\frac{ 4}{3}\, L_{xy} 
-
2\, M_{xx}+\\
& \
+2\, G\, M_y+4\, G_y\, M-
2\, (H\, M)_x-\frac{ 2}{3}\, 
H_y\, L+\frac{ 4}{3}\, L\, L_x.
\endaligned\right.
\end{equation}
\end{itemize}
\end{theorem}

\def\theopenquestion{1.4}\begin{openquestion}
Deduce an explicit necessary and sufficient condition for the
associated submanifold of solutions $y = \Pi (x, a, b)$ to be
locally equivalent to $Y = B + XA$.
\end{openquestion}

Assuming $F = F (x, y_x)$ to be independent of $y$, or equivalently
assuming $\mathcal{ M}_{ (\mathcal{ E}_1)}$ to be:
\def\theequation{1.5}\begin{equation}
y 
= 
b+\Pi (x, a),
\end{equation}
the author has checked that equivalence to $Y = B + XA$ holds if and
only if two differential rational expressions annihilate:
\def\theequation{1.6}\begin{equation}
\aligned
0
&
=
\frac{\Pi_{x^2a^4}}{
\big(\Pi_{xa}\big)^4}
-
6\,
\frac{\Pi_{x^2a^3}\,\Pi_{xa^2}}{
\big(\Pi_{xa}\big)^5}
+
15\,
\frac{\Pi_{x^2a^2}\,\big(\Pi_{xa^2}\big)^2}{
\big(\Pi_{xa}\big)^6}
-
4\,
\frac{\Pi_{x^2a^2}\,\Pi_{xa^3}}{
\big(\Pi_{xa}\big)^5}
\\
&\ \ \ \ \
-
\frac{\Pi_{x^2a}\,\Pi_{xa^4}}{
\big(\Pi_{xa}\big)^5}
+
10\,
\frac{\Pi_{xa^3}\,\Pi_{x^2a}\,\Pi_{xa^2}}{
\big(\Pi_{xa}\big)^6}
-
15\,
\frac{\Pi_{x^2a}\,\big(\Pi_{xa^2}\big)^3}{
\big(\Pi_{xa}\big)^7}
\ \ \ \ \ \ \ \text{\rm and}
\\
0
&
=
\frac{\Pi_{x^4a^2}}{
\big(\Pi_{xa}\big)^2}
-
6\,
\frac{\Pi_{x^3a^2}\,\Pi_{x^2a}}{
\big(\Pi_{xa}\big)^3}
-
4\,
\frac{\Pi_{x^3a}\,\Pi_{x^2a^2}}{
\big(\Pi_{xa}\big)^3}
-
\frac{\Pi_{x^4a}\,\Pi_{xa^2}}{
\big(\Pi_{xa}\big)^3}
+
\\
&\ \ \ \ \
+
15\,
\frac{\Pi_{x^2a^2}\,\big(\Pi_{x^2a}\big)^2}{
\big(\Pi_{xa}\big)^4}
+
10\,
\frac{\Pi_{x^3a}\,\Pi_{x^2a}\,\Pi_{xa^2}}{
\big(\Pi_{xa}\big)^4}
-
15\,
\frac{\big(\Pi_{x^2a}\big)^3\,\Pi_{xa^2}}{
\big(\Pi_{xa}\big)^5}.
\endaligned
\end{equation}
As an application, this characterizes local sphericity of a rigid
hypersurface $w = \bar w + i\, \Theta (z, \bar z)$ of $\C^2$. The
answer for a general $y = \Pi (x, a, b)$, together with a proof, will
appear elsewhere.

\smallskip

A modern restitution of Lie's original proof of Theorem~1.1 may be
found in~\cite{ me2004}. In this reference, we generalize Theorem~1.1
to several dependent variables $y = (y^1, y^2, \dots, y^m)$. In the
present Part~III, we will instead pass to several independent variables
$x = (x^1, x^2, \dots, x^n)$.

\def\thetheorem{1.7}\begin{theorem}
Let $\K = \R$ or $\C$, let $n\in\N$,
\underline{\rm suppose $n \geqslant 2$} and consider a system of
completely integrable partial differential equations in $n$
independent variables $x= (x^1, \dots, x^n) \in \K^n$ and in one
dependent variable $y \in \K$ of the form{\rm :}
\def\theequation{1.8}\begin{equation}
y_{x^{j_1}x^{j_2}}(x) 
=
F^{j_1, j_2} 
\big(
x, y(x),y_{x^1}(x),\dots,y_{x^n}(x)
\big), 
\ \ \ \ \ \ \ \
1\leqslant j_1, j_2 \leqslant n,
\end{equation}
where $F^{j_1, j_2} = F^{ j_2, j_1}$. Under a local change of
coordinates $(x, y) \mapsto (X, Y) = (X(x, y), Y(x, y))$, this
system~\thetag{1.8} is equivalent to the {\rm simplest}
``flat'' system
\def\theequation{1.9}\begin{equation}
Y_{X^{j_1} X^{j_2}}
=
0, 
\ \ \ \ \ \ \ \
1\leqslant j_1, j_2 \leqslant n,
\end{equation}
{\rm if and only if} there exist {\rm arbitrary} functions $G_{j_1,
j_2}$, $H_{ j_1, j_2}^{ k_1}$, $L_{ j_1}^{ k_1}$ and $M^{k_1}$ of the
variables $(x, y)$, for $1\leqslant j_1,j_2,k_1 \leqslant n$,
satisfying the two symmetry conditions $G_{j_1, j_2}= G_{j_2, j_1}$
and $H_{j_1, j_2}^{ k_1}= H_{j_2,j_1}^{k_1}$, such that the
equation~\thetag{1.8} is of the specific cubic polynomial form{\rm :}
\def\theequation{1.10}\begin{equation}
y_{x^{j_1}x^{j_2}}
=
G_{j_1,j_2}+\sum_{k_1=1}^n\, 
y_{x^{k_1}}\left(
H_{j_1,j_2}^{k_1} +
\frac{ 1}{2}\, y_{x^{j_1}} \, L_{j_2}^{k_1}+
\frac{ 1}{2}\, y_{x^{j_2}}\, L_{j_1}^{k_1} +
y_{x^{j_1}}
y_{x^{j_2}}\, M^{k_1}
\right), 
\end{equation}
for $j_1, j_2 = 1, \dots, n$.
\end{theorem}

It may seem quite paradoxical and counter-intuitive (or even false?)
that {\it every} system~\thetag{ 1.10}, for {\it arbitrary} choices of
functions $G_{j_1, j_2}$, $H_{j_1, j_2}^{k_1}$, $L_{j_1}^{k_1}$ and
$M^{k_1}$, is {\it automatically} equivalent to $Y_{X^{j_1} X^{j_2}}=
0$. However, a strong hidden assumption holds: that of {\sl complete
integrability}. Shortly, this crucial condition amounts to say that
\def\theequation{1.11}\begin{equation}
D_{j_3}
\left(
F^{j_1,j_2}
\right)=D_{j_2}\left(
F^{j_1,j_3}
\right),
\end{equation}
for all $j_1, j_2, j_3 = 1, \dots, n$, where, for $j =1, \dots, n$, the
$D_j$ are the {\sl total differentiation operators} defined by
\def\theequation{1.12}\begin{equation}
D_j:=
\frac{ \partial }{\partial x^j}+
y_{x^j} \, 
\frac{ \partial }{\partial y}+
\sum_{l=1}^n\,
F^{j,l}\, 
\frac{ \partial }{\partial y_{x^l}}.
\end{equation}
These conditions are non-void precisely when $n\geqslant 2$. More
concretely, developing out~\thetag{ 1.11} when the
$F^{j_1, j_2}$ are of the specific cubic polynomial form~\thetag{
1.10}, after some nontrivial manual computation, we obtain the
complicated cubic differential polynomial in the
variables $y_{x^k}$. Equating to zero all the
coefficients of this cubic polynomial, we obtain four familes (I'),
(II'), (III') and (IV') of \underline{first order} partial
differential equations satisfied by $G_{j_1, j_2}$, $H_{j_1,
j_2}^{k_1}$, $L_{j_1}^{k_1}$ and $M^{k_1}$:
\def\theequation{I'}\begin{equation}
\left\{
0 =
G_{j_1,j_2, x^{j_3}}- G_{j_1,j_3, x^{j_2}}
+
\sum_{k_1=1}^n\, H_{j_1,j_2}^{k_1}\, G_{k_1, j_3}
-
\sum_{k_1=1}^n\, H_{j_1, j_3}^{k_1}\, G_{k_1, j_2}.
\right.
\end{equation}
\def\theequation{II'}\begin{equation}
\left\{
\aligned
0
&
=
\delta_{j_3}^{k_1}\, G_{j_1,j_2, y}
-
\delta_{j_2}^{k_1}\, G_{j_1,j_3, y}
+
H_{j_1, j_2, x^{j_3}}^{k_1}
- 
H_{j_1, j_3, x^{j_2}}^{k_1}
+ 
\\
& \
\ \ \ \ \
+
\frac{ 1}{2}\,G_{j_1, j_3}\, L_{j_2}^{k_1}
- 
\frac{ 1}{2}\,G_{j_1, j_2}\, L_{j_3}^{k_1}
+ \\
& \
\ \ \ \ \
+
\frac{1}{2}\,\delta_{j_1}^{k_1}\,\sum_{k_2=1}^n\,
G_{k_2,j_3}\,L_{j_2}^{k_2}
-
\frac{ 1}{2}\,\delta_{j_1}^{k_1}\,\sum_{k_2=1}^n\, 
G_{k_2,j_2}\, L_{j_3}^{k_2}
+ \\
& \
\ \ \ \ \
+
\frac{ 1}{2}\,\delta_{j_2}^{k_1}\,\sum_{k_2=1}^n\, 
G_{k_2, j_3}\, L_{j_1}^{k_2}
-
\frac{ 1}{2}\,\delta_{j_3}^{k_1}\,\sum_{k_2=1}^n\, 
G_{k_2,j_2}\,L_{j_1}^{k_2}
+ \\
& \
\ \ \ \ \
+
\sum_{k_2=1}^n\,H_{k_2,j_3}^{k_1}\,H_{j_1, j_2}^{k_2}
- 
\sum_{k_2=1}^n\,H_{k_2,j_2}^{k_1}\,H_{j_1,j_3}^{k_2}.
\endaligned\right.
\end{equation}
\def\theequation{III'}\begin{equation}
\left\{
\aligned
0
&
=
\sum_{\sigma\in\mathfrak{S}_2}\left(
\delta_{j_3}^{k_{2)}}\,H_{j_1,j_2,y}^{k_{\sigma(1)}}
-
\delta_{j_2}^{k_{\sigma(2)}}\,H_{j_1,j_3,y}^{k_{\sigma(1)}}
+
\right. \\
& \
\ \ \ \ \ \ \ \ \ \ \ \ \ \ \
\left.
+
\frac{1}{2}\,\delta_{j_2}^{k_{\sigma(2)}}\,L_{j_1,x^{j_3}}^{k_{\sigma(1)}}
-
\frac{1}{2}\,\delta_{j_3}^{k_{\sigma(2)}}\,L_{j_1,x^{j_2}}^{k_{\sigma(1)}}
+
\right. \\
& \
\ \ \ \ \ \ \ \ \ \ \ \ \ \ \
\left.
+
\frac{1}{2}\,\delta_{j_1}^{k_{\sigma(2)}}\,L_{j_2,x^{j_3}}^{k_{\sigma(1)}}
-
\frac{1}{2}\,\delta_{j_1}^{k_{\sigma(2)}}\,L_{j_3,x^{j_2}}^{k_{\sigma(1)}}
+
\right. \\
& \
\ \ \ \ \ \ \ \ \ \ \ \ \ \ \
\left.
+
\delta_{j_2}^{k_{\sigma(2)}}\,G_{j_1,j_3}\,M^{k_{\sigma(1)}}
-
\delta_{j_3}^{k_{\sigma(2)}}\,G_{j_1,j_2}\,M^{k_{\sigma(1)}}
+
\right.
\\
& \ 
\ \ \ \ \ \ \ \ \ \ \ \ \ \ \
\left.
+
\delta_{j_1,\ \ \ \ j_2}^{k_{\sigma(1)},k_{\sigma(2)}}\,
\sum_{k_3=1}^n\,G_{k_3,j_3}\,M^{k_3}
-
\delta_{j_1,\ \ \ \ j_3}^{k_{\sigma(1)},k_{\sigma(2)}}\,
\sum_{k_3=1}^n\,G_{k_3,j_2}\,M^{k_3}
+ 
\right. \\
& \
\ \ \ \ \ \ \ \ \ \ \ \ \ \ \
+
\left.
\frac{1}{2}\,\delta_{j_1}^{k_{\sigma(1)}}\,\sum_{k_3=1}^n\, 
H_{k_3,j_3}^{k_{\sigma(2)}}\,L_{j_2}^{k_3} 
-
\frac{1}{2}\,\delta_{j_1}^{k_{\sigma(1)}}\,\sum_{k_3=1}^n\,
H_{k_3,j_2}^{k_{\sigma(2)}}\,L_{j_3}^{k_3}
+ 
\right. \\
& \
\ \ \ \ \ \ \ \ \ \ \ \ \ \ \
\left.
+
\frac{1}{2}\,\delta_{j_2}^{k_{\sigma(1)}}\,\sum_{k_3=1}^n\,
H_{k_3,j_3}^{k_{\sigma(2)}}\,L_{j_1}^{k_3}
-
\frac{1}{2}\,\delta_{j_3}^{k_{\sigma(1)}}\,\sum_{k_3=1}^n\, 
H_{k_3,j_2}^{k_{\sigma(2)}}\,L_{j_1}^{k_3}
+ 
\right. \\
& \
\ \ \ \ \ \ \ \ \ \ \ \ \ \ \
\left.
+
\frac{1}{2}\,\delta_{j_3}^{k_{\sigma(1)}}\,\sum_{k_3=1}^n\,
H_{j_1,j_2}^{k_3}\,L_{k_3}^{k_{\sigma(2)}}
-
\frac{1}{2}\,\delta_{j_2}^{k_{\sigma(1)}}\,\sum_{k_3=1}^n\, 
H_{j_1,j_3}^{k_3}\ L_{k_3}^{k_{\sigma(2)}}
\right).
\endaligned\right.
\end{equation}
\def\theequation{IV'}\begin{equation}
\left\{
\aligned
0
&
=
\sum_{\sigma\in\mathfrak{S}_3}\left( 
\frac{1}{2}\,\delta_{j_3,\ \ \ \ j_1}^{k_{\sigma(3)},k_{\sigma(2)}}\,
L_{j_2,y}^{k_{\sigma(1)}}
- 
\frac{1}{2}\,\delta_{j_2,\ \ \ \ j_1}^{k_{\sigma(3)},k_{\sigma(2)}}\,
L_{j_3,y}^{k_{\sigma(1)}}
+ 
\right. \\
& \ 
\ \ \ \ \ \ \ \ \ \ \ \
\left.
+
\delta_{j_2,\ \ \ \ j_1}^{k_{\sigma(3)},k_{\sigma(2)}}\,
M_{x^{j_3}}^{k_{\sigma(1)}}
- 
\delta_{j_3,\ \ \ \ j_1}^{k_{\sigma(3)},k_{\sigma(2)}}\,
M_{x^{j_2}}^{k_{\sigma(1)}}
+ 
\right. \\
& \
\ \ \ \ \ \ \ \ \ \ \ \
\left.
+
\delta_{j_2,\ \ \ \ j_1}^{k_{\sigma(3)},k_{\sigma(1)}}\, 
\sum_{k_4=1}^n\,H_{k_4,j_3}^{k_{\sigma(2)}}\,M^{k_4}
-
\delta_{j_3,\ \ \ \ j_1}^{k_{\sigma(3)},k_{\sigma(1)}}\, 
\sum_{k_4=1}^n\,H_{k_4,j_2}^{k_{\sigma(2)}}\,M^{k_4}
+ 
\right. \\
& \
\ \ \ \ \ \ \ \ \ \ \ \
\left.
+
\frac{1}{4}\,\delta_{j_1,\ \ \ \ j_3}^{k_{\sigma(1)},k_{\sigma(3)}}\,
\sum_{k_4=1}^n\,L_{k_4}^{k_{\sigma(2)}}\,L_{j_2}^{k_4}
-
\frac{1}{4}\,\delta_{j_1,\ \ \ \ j_2}^{k_{\sigma(1)},k_{\sigma(3)}}\, 
\sum_{k_4=1}^n\,L_{k_4}^{k_{\sigma(2)}}\, L_{j_3}^{k_4}
\right).
\endaligned\right.
\end{equation}
(These systems (I'), (II'), (III') and (IV') should be distinguished
from the systems (I), (II), (III) and (IV) of Theorem~1.7 in~\cite{
me2004}, although they are quite similar.) Here, the indices $j_1, j_2,
j_3, k_1, k_2, k_3$ vary in $\{1, 2, \dots, n\}$. By $\mathfrak{ S}_2$
and by $\mathfrak{ S}_3$, we denote the permutation group of $\{1,
2\}$ and of $\{1, 2, 3\}$. To facilitate hand- and Latex-writing,
partial derivatives are denoted as indices after a comma; for
instance, $G_{j_1, j_2, x^{j _3}}$ is an abreviation for $\partial
G_{j_1, j_2}/ \partial x^{j _3}$. To deduce (I'), (II'), (III') and
(IV') from equation~\thetag{ 1.11}, we use the fact that every cubic
polynomial equation of the form
\def\theequation{1.13}\begin{equation}
\aligned
0 
&
\equiv
A 
+
\sum_{k_1=1}^n\,B_{k_1}\cdot\,
y_{x^{k_1}}
+
\sum_{k_1=1}^n\,\sum_{k_2=1}^n\,C_{k_1,k_2}\cdot\, 
y_{x^{k_1}}\,y_{x^{k_2}}
+ \\
& \ 
\ \ \ \ \ \ \ 
+
\sum_{k_1=1}^n\,\sum_{k_2=1}^n\,\sum_{k_2=1}^n\,D_{k_1,k_2,k_3}\cdot\,
y_{x^{k_1}}\,y_{x^{k_2}}\,y_{x^{k_3}}
\endaligned
\end{equation}
is equivalent to the
annihilation of the following symmetric sums of its coefficients:
\def\theequation{1.14}\begin{equation}
\left\{
\aligned
0 
&
=
A, 
\\
0 
&
=
B_{k_1}, 
\\
0 
&
=
C_{k_1,k_2}
+
C_{k_2,k_1}, 
\\ 
0 
&
=
D_{k_1,k_2,k_3}
+
D_{k_3,k_1,k_2}
+
D_{k_2,k_3,k_1}
+
D_{k_2,k_1,k_3}
+
D_{k_3,k_2,k_1}
+
D_{k_1,k_3,k_2}.
\endaligned\right.
\end{equation}
for all $k_1, k_2, k_3 = 1, \dots, n$.

In conclusion, the functions $G_{j_1, j_2 }$, $H_{ j_1, j_2 }^{ k_1}$,
$L_{ j_1}^{ k_1}$ and $M^{ k_1}$ in the statement of Theorem~1.7 are
far from being arbitrary: they satisfy the complicated system of first
order partial differential equations (I'), (II'), (III') and (IV')
above.

Our proof of Theorem~1.7 is similar to the one provided in~\cite{
me2004}, in the case of systems of second order ordinary differential
equations, so that most steps of the proof will be summarized.

In the end of this paper, we will delineate a complicated system of
{\it second}\, order partial differential equations satisfied by
$G_{j_1, j_2 }$, $H_{j_1, j_2 }^{ k_1}$, $L_{ j_1 }^{ k_1}$ and $M^{
k_1 }$ which is the exact analog of the system described in the
abstract. The main technical part of the proof of Theorem~1.7 will be
to establish that this second order system is a consequence, by linear
combinations and by differentiations, of the first order system (I'),
(II'), (III') and (IV').

\def\theopenquestion{1.15}\begin{openquestion}
Are Theorems~1.1 and 1.7 true under weaker smoothness assumptions,
namely with a $\mathcal{ C }^2$ or a $W_{\rm loc }^{ 1, \infty }$
right-hand side~?
\end{openquestion}

We refer to~\cite{ ma2003} for inspiration and appropriate tools.

\def\theopenquestion{1.16}\begin{openquestion}
Deduce from Theorem~1.7 an explicit necessary and sufficient condition
for the associated submanifold of solutions $y = b + \Pi (x^i, 
a^k, b)$ to
be locally equivalent to 
$Y = B + X^1A^1 + \cdots + X^n A^n$.
\end{openquestion}

As an application, this would characterize local sphericity of a Levi
nondegenerate hypersurface $M \subset \C^{ n+ 1}$ with $n \geqslant
2$.

\medskip

Generalizing the Lie-Tresse classification would be a great
achievement.

\def\theopenproblem{1.17}\begin{openproblem}
For $n = 2$ establish a complete list of normal forms of all possible
systems~\thetag{ 1.?} according to their Lie symmetry group.
In case of success, classify Levi nondegenerate
real analytic hypersurfaces of $\C^3$ up
to biholomorphisms.
\end{openproblem}

\section*{ \S2.~Completely integrable systems of \\
second order ordinary differential equations}

\subsection*{ 2.1.~Prolongation of a point transformation to 
the second order jet space} Let $\K = \R$ or $\C$, let
$n\in \N$, 
\underline{suppose
$n\geqslant 2$}, let $x = (x^1, \dots, x^n) \in \K^n$ and let $y
\in \K$. According to the main assumption of Theorem~1.7, we have to
consider a local $\K$-analytic diffeomorphism of the form
\def\theequation{2.2}\begin{equation}
\left(
x^{j_1}, y 
\right)
\longmapsto 
\left(
 X^j (x^{j_1}, y), 
\
Y( x^{j_1}, y)
\right),
\end{equation}
which transforms the system~\thetag{ 1.8} to the system $Y_{X^{i_1}
X^{i_2}}= 0$, $1\leqslant j_1, j_2 \leqslant n$. Without loss of
generality, we shall assume that this transformation is close to the
identity. To obtain the precise expression~\thetag{ 2.35} of the
transformed system~\thetag{ 1.8}, we have to prolong the above
diffeomorphism to the second order jet space. We introduce the
coordinates $\left( x^j, y, y_{x^{j_1}}, y_{x^{j_1}x^{j_2}} \right)$
on the second order jet space. Let
\def\theequation{2.3}\begin{equation}
D_k
:= 
\frac{\partial}{\partial x^k}
+
y_{x^k}\,\frac{\partial}{\partial y}
+
\sum_{l=1}^n\,y_{x^kx^l}\,\frac{\partial}{\partial y_{x^l}},
\end{equation}
be the $k$-th total differentiation operator. According to~\cite{
ol1986, bk1989, ol1995}, for the first order partial derivatives, one
has the (implicit, compact) expression:
\def\theequation{2.4}\begin{equation}
\left(
\begin{array}{c}
Y_{X^1} \\
\vdots \\
Y_{X^n}
\end{array}
\right)
=
\left(
\begin{array}{ccc}
D_1X^1 & \cdots & D_1X^n \\
\vdots & \cdots & \vdots \\
D_nX^1 & \cdots & D_nX^n
\end{array}
\right)^{-1} 
\left(
\begin{array}{c}
D_1Y \\
\vdots \\
D_nY
\end{array}
\right),
\end{equation}
where $(\cdot)^{ -1}$ denotes the inverse matrix, which exists, since
the transformation~\thetag{ 2.2} is close to the identity. For the
second order partial derivatives, again according to~\cite{ 
ol1986, bk1989, ol1995},
one has the (implicit, compact) expressions:
\def\theequation{2.5}\begin{equation}
\left(
\begin{array}{c}
Y_{X^jX^1} \\
\vdots \\
Y_{X^jX^n}
\end{array}
\right)
=
\left(
\begin{array}{ccc}
D_1X^1 & \cdots & D_1X^n \\
\vdots & \cdots & \vdots \\
D_nX^1 & \cdots & D_nX^n
\end{array}
\right)^{-1} 
\left(
\begin{array}{c}
D_1Y_{X^j} \\
\vdots \\
D_nY_{X^j}
\end{array}
\right),
\end{equation}
for $j = 1, \dots, n$. Let $DX$ denote the matrix $\left( D_iX^j
\right)_{ 1\leqslant i\leqslant n}^{ 1\leqslant j\leqslant n}$, where
$i$ is the index of lines and $j$ the index of columns, let $Y_X$
denote the column matrix $\left( Y_{X^i} \right)_{ 1\leqslant i
\leqslant n}$ and let $DY$ be the column matrix $\left( D_i Y
\right)_{ 1\leqslant i \leqslant n}$.

By inspecting~\thetag{ 2.5} above, we see that the
equivalence between {\bf (i)}, {\bf (ii)}
and {\bf (iii)} just below is obvious:

\def\thelemma{2.6}\begin{lemma}
The following conditions are equivalent{\rm :}
\begin{itemize}
\item[{\bf (i)}]
the differential equations $Y_{X^jX^k} = 0$ hold for $1\leqslant j,k \leqslant
n${\rm ;}
\item[{\bf (ii)}]
the matrix equations $D_k( Y_X)= 0$ hold for $1\leqslant k\leqslant n${\rm ;}
\item[{\bf (iii)}]
the matrix equations $DX \cdot D_k( Y_X)= 0$ hold for $1\leqslant 
k\leqslant n${\rm ;}
\item[{\bf (iv)}]
the matrix equations
$0 = D_k( DX) \cdot Y_X - D_k (DY)$ hold for $1\leqslant k\leqslant n$.
\end{itemize}
\end{lemma}

Formally, in the sequel, it will be more convenient to achieve the
explicit computations starting from condition {\bf (iv)}, since no
matrix inversion at all is involved in it.

\proof
Indeed, applying the
total differentiation operator $D_k$ to 
the matrix equation~\thetag{ 2.4} written under the equivalent
form $0 = DX\cdot Y_X - DY$, we get:
\def\theequation{2.7}\begin{equation}
0 
=
D_k(DX)\cdot Y_X 
+
DX\cdot D_k(Y_X)
-
D_k(DY),
\end{equation}
so that the equivalence between {\bf (iii)} and
{\bf (iv)} is now clear.
\endproof

\subsection*{ 2.8.~An explicit formula in the case $n = 2$}
Thus, we can start to develope explicitely the 
matrix equations
\def\theequation{2.9}\begin{equation}
0 
=
D_k (DX) \cdot Y_X 
- 
D_k(DY).
\end{equation}
In it, some huge formal expressions are hidden behind the symbol
$D_k$. Proceeding inductively, we start by examinating the case $n=2$
thoroughly. By direct computations which require to be clever, we
reconstitute some $3\times 3$ determinants in the four (in 
fact three) developed
equations~\thetag{ 2.9}. After some work, the first equation is:
\def\theequation{2.10}\begin{equation}
\aligned
0
& 
= 
y_{x^1x^1}
\cdot
\left\vert
\begin{array}{ccc}
X_{x^1}^1 & X_{x^2}^1 & X_y^1 \\
X_{x^1}^2 & X_{x^2}^2 & X_y^2 \\
Y_{x^1} & Y_{x^2} & Y_y \\
\end{array}
\right\vert
+
\left\vert
\begin{array}{ccc}
X_{x^1}^1 & X_{x^2}^1 & X_{x^1x^1}^1 \\
X_{x^1}^2 & X_{x^2}^2 & X_{x^1x^1}^2 \\
Y_{x^1} & Y_{x^2} & Y_{x^1x^1} \\
\end{array}
\right\vert
+ \\
& \
\ \ \ \ \
+
y_{x^1}
\cdot
\left\{
2\,
\left\vert
\begin{array}{ccc}
X_{x^1}^1 & X_{x^2}^1 & X_{x^1y}^1 \\
X_{x^1}^2 & X_{x^2}^2 & X_{x^1y}^2 \\
Y_{x^1} & Y_{x^2} & Y_{x^1y} \\
\end{array}
\right\vert
-
\left\vert
\begin{array}{ccc}
X_{x^1x^1}^1 & X_{x^2}^1 & X_y^1 \\
X_{x^1x^1}^2 & X_{x^2}^2 & X_y^2 \\
Y_{x^1x^1} & Y_{x^2} & Y_y \\
\end{array}
\right\vert
\right\}
+ \\
& \
\ \ \ \ \
+
y_{x^2}
\cdot
\left\{
-
\left\vert
\begin{array}{ccc}
X_{x^1}^1 & X_{x^1x^1}^1 & X_y^1 \\
X_{x^1}^2 & X_{x^1x^1}^2 & X_y^2 \\
Y_{x^1} & Y_{x^1x^1} & Y_y \\
\end{array}
\right\vert
\right\}
+ \\
\endaligned
\end{equation}
$$
\aligned
& \
\ \ \ \ \
+
y_{x^1}\,y_{x^1}
\cdot
\left\{
\left\vert
\begin{array}{ccc}
X_{x^1}^1 & X_{x^2}^1 & X_{yy}^1 \\
X_{x^1}^2 & X_{x^2}^2 & X_{yy}^2 \\
Y_{x^1} & Y_{x^2} & Y_{yy} \\
\end{array}
\right\vert
-
2\,
\left\vert
\begin{array}{ccc}
X_{x^1y}^1 & X_{x^2}^1 & X_y^1 \\
X_{x^1y}^2 & X_{x^2}^2 & X_y^2 \\
Y_{x^1y} & Y_{x^2} & Y_y \\
\end{array}
\right\vert
\right\}
+ \\
& \
\ \ \ \ \ 
+
y_{x^1}\,y_{x^2}
\cdot
\left\{
-2\,
\left\vert
\begin{array}{ccc}
X_{x^1}^1 & X_{x^1y}^1 & X_y^1 \\
X_{x^1}^2 & X_{x^1y}^2 & X_y^2 \\
Y_{x^1} & Y_{x^1y} & Y_y \\
\end{array}
\right\vert
\right\}
+ \\
& \
\ \ \ \ \
+
y_{x^1}\,y_{x^1}\,y_{x^1}
\cdot
\left\{
-
\left\vert
\begin{array}{ccc}
X_{yy}^1 & X_{x^2}^1 & X_y^1 \\
X_{yy}^2 & X_{x^2}^2 & X_y^2 \\
Y_{yy} & Y_{x^2} & Y_y \\
\end{array}
\right\vert
\right\}
+ \\
& \
\ \ \ \ \
+
y_{x^1}\,y_{x^1}\,y_{x^2}
\cdot
\left\{
-
\left\vert
\begin{array}{ccc}
X_{x^1}^1 & X_{yy}^1 & X_y^1 \\
X_{x^1}^2 & X_{yy}^2 & X_y^2 \\
Y_{x^1} & Y_{yy} & Y_y \\
\end{array}
\right\vert
\right\}.
\endaligned
$$
This formula and the two next~\thetag{ 2.22}, 
\thetag{ 2.23} have been checked by Sylvain Neut
and Michel Petitot with the help of Maple.

\subsection*{2.11.~Comparison with the 
coefficients of the second prolongation of a vector field} At present,
it is useful to make an illuminating digression which will help us to
devise what is the general form of the development of the
equations~\thetag{ 2.9}. Consider an arbitrary vector field of the
form
\def\theequation{2.12}\begin{equation}
\mathcal{ L}
:= 
\sum_{k=1}^n\,\mathcal{X}^k\,\frac{\partial}{\partial x^k}
+
\mathcal{Y}\,\frac{\partial}{\partial y},
\end{equation}
where the coefficients $\mathcal{ X}^k$ and $\mathcal{ Y}$ are
functions of $(x^i, y)$. According to~\cite{ ol1986, bk1989, ol1995},
there exists a unique prolongation $\mathcal{ L}^{(2)}$ of this vector
field to the second order jet space, of the form
\def\theequation{2.13}\begin{equation}
\mathcal{ L}^{(2)} 
:=
\mathcal{L}
+
\sum_{j_1=1}^n\,
{\bf Y}_{j_1}\,\frac{\partial}{\partial y_{x^{j_1}}}
+
\sum_{j_1=1}^n\,\sum_{j_2=1}^n\,
{\bf Y}_{j_1,j_2}\,\frac{\partial}{\partial y_{x^{j_1}x^{j_2}}},
\end{equation}
where the coefficients ${\bf Y}_{ j_1}$, ${\bf Y}_{ j_1, j_2}$ may be
computed by means of formulas~\thetag{ 3.4} of Section~3(II). In
Part~II, we obtained the following perfect formulas:
\def\theequation{2.14}\begin{equation}
\left\{
\aligned
{\bf Y}_{j_1,j_2}
&
=
\mathcal{Y}_{x^{j_1}x^{j_2}}
+
\sum_{k_1=1}^n\,y_{x^{k_1}}
\cdot
\left\{
\delta_{j_1}^{k_1}\,\mathcal{Y}_{x^{j_2}y}
+
\delta_{j_2}^{k_1}\,\mathcal{Y}_{x^{j_1}y}
-
\mathcal{X}_{x^{j_1}x^{j_2}}^{k_1}
\right\}
+ \\
& \
\ \ \ \ \
+
\sum_{k_1=1}^n\,\sum_{k_2=1}^n\,y_{x^{k_1}}\,y_{x^{k_2}}
\cdot
\left\{
\delta_{j_1,\,j_2}^{k_1,k_2}\,\mathcal{Y}_{yy}
-
\delta_{j_1}^{k_1}\,\mathcal{X}_{x^{j_2}y}^{k_2}
-
\delta_{j_2}^{k_1}\,\mathcal{X}_{x^{j_1}y}^{k_2}
\right\}
+ \\
& \
\ \ \ \ \
+
\sum_{k_1=1}^n\,\sum_{k_2=1}^n\,\sum_{k_3=1}^n\,
y_{x^{k_1}}\,y_{x^{k_2}}\,y_{x^{k_3}}
\cdot
\left\{
-\delta_{j_1,\,j_2}^{k_1,k_2}\,\mathcal{X}_{yy}^{k_3}
\right\},
\endaligned\right.
\end{equation} 
for $j_1, j_2 = 1,\dots, n$.
The expression of ${\bf Y}_{j_1}$ does not
matter for us here.
Specifying this formula to the the case $n=2$ and taking account of
the symmetry ${\bf Y}_{ 1, 2} = {\bf Y}_{ 2,1}$ we get the following
three second order coefficients:
\def\theequation{2.15}\begin{equation}
\left\{
\aligned
{\bf Y}_{1,1}
&
=
\mathcal{Y}_{x^1x^1}
+
y_{x^1}
\cdot
\left\{
2\,\mathcal{Y}_{x^1y}
-
\mathcal{X}_{x^1x^1}^1
\right\}
+
y_{x^2}
\cdot
\left\{
-
\mathcal{X}_{x^1x^1}^2
\right\}
+ \\
& \
\ \ \ \ \
+
y_{x^1}\,y_{x^1}
\cdot
\left\{
\mathcal{Y}_{yy}
-
2\,\mathcal{X}_{x^1y}^1
\right\}
+
y_{x^1}\,y_{x^2}
\cdot
\left\{
-
2\,\mathcal{X}_{x^1y}^2
\right\}
+ \\
& \
\ \ \ \ \ 
+
y_{x^1}\,y_{x^1}\,y_{x^1}
\cdot
\left\{
-
\mathcal{X}_{yy}^1
\right\}
+
y_{x^1}\,y_{x^1}\,y_{x^2}
\cdot
\left\{
-
\mathcal{X}_{yy}^2
\right\}, \\
{\bf Y}_{1,2}
&
=
\mathcal{Y}_{x^1x^2}
+
y_{x^1}
\cdot
\left\{
\mathcal{Y}_{x^2y}
-
\mathcal{X}_{x^1x^2}^1
\right\}
+
y_{x^2}
\cdot
\left\{
\mathcal{Y}_{x^1y}
-
\mathcal{X}_{x^1x^2}^2
\right\}
+ \\
& \
\ \ \ \ \
+
y_{x^1}\,y_{x^1}
\cdot
\left\{
-
\mathcal{X}_{x^2y}^1
\right\}
+
y_{x^1}\,y_{x^2}
\cdot
\left\{
\mathcal{Y}_{yy}
-
\mathcal{X}_{x^1y}^1
-
\mathcal{X}_{x^2y}^2
\right\}
+ \\
& \
\ \ \ \ \
+
y_{x^2}\,y_{x^2}
\cdot
\left\{
-
\mathcal{X}_{x^1y}^2
\right\}
+ \\
& \
\ \ \ \ \ 
+
y_{x^1}\,y_{x^1}\,y_{x^2}
\cdot
\left\{
-
\mathcal{X}_{yy}^1
\right\}
+
y_{x^1}\,y_{x^2}\,y_{x^2}
\cdot
\left\{
-
\mathcal{X}_{yy}^2
\right\}, \\
{\bf Y}_{2,2}
&
=
\mathcal{Y}_{x^2x^2}
+
y_{x^1}
\cdot
\left\{
-
\mathcal{X}_{x^2x^2}^1
\right\}
+
y_{x^2}
\cdot
\left\{
2\,\mathcal{Y}_{x^2y}
-
\mathcal{X}_{x^2x^2}^2
\right\}
+ \\
& \
\ \ \ \ \
+
y_{x^1}\,y_{x^2}
\cdot
\left\{
-
2\,\mathcal{X}_{x^2y}^1
\right\}
+
y_{x^2}\,y_{x^2}
\cdot
\left\{
\mathcal{Y}_{yy}
-
2\,\mathcal{X}_{x^2y}^2
\right\}
+ \\
& \
\ \ \ \ \ 
+
y_{x^1}\,y_{x^2}\,y_{x^2}
\cdot
\left\{
-
\mathcal{X}_{yy}^1
\right\}
+
y_{x^2}\,y_{x^2}\,y_{x^2}
\cdot
\left\{
-
\mathcal{X}_{yy}^2
\right\}.
\endaligned\right.
\end{equation}
We would like to mention that the computation of ${\bf Y}_{j_1, j_2}$,
$1\leqslant j_1, j_2 \leqslant 2$, above is easier than the verification
of~\thetag{ 2.10}. Based on the three formulas~\thetag{ 2.15}, we claim
that we can guess the second and the third equations, 
which would be obtained by developing and by simplifying~\thetag{ 2.9},
namely with
$y_{x^1x^2}$ and with $y_{ x^2 x^2}$ instead of $y_{x^1 x^2}$
in~\thetag{ 2.10}. Our dictionary to translate from the first
formula~\thetag{ 2.15} to~\thetag{ 2.10} may be described as
follows. Begin with the {\sl Jacobian determinant}
\def\theequation{2.16}\begin{equation}
\left\vert
\begin{array}{ccc}
X_{x^1}^1 & X_{x^2}^1 & X_y^1 \\
X_{x^1}^2 & X_{x^2}^2 & X_y^2 \\
Y_{x^1} & Y_{x^2} & Y_y \\
\end{array}
\right\vert
\end{equation}
of the change of coordinates~\thetag{ 2.2}. Since this change of
coordinates is close to the identity, we may consider that the following
Jacobian matrix approximation holds:
\def\theequation{2.17}\begin{equation}
\left(
\begin{array}{ccc}
X_{x^1}^1 & X_{x^2}^1 & X_y^1 \\
X_{x^1}^2 & X_{x^2}^2 & X_y^2 \\
Y_{x^1} & Y_{x^2} & Y_y \\
\end{array}
\right)
\cong
\left(
\begin{array}{ccc}
1 & 0 & 0 \\
0 & 1 & 0 \\
0 & 0 & 1\\
\end{array}
\right).
\end{equation}
The jacobian matrix has three columns. There are six possible second
order derivatives with respect to the variables $(x^1, x^2, y)$,
namely
\def\theequation{2.18}\begin{equation}
(\cdot)_{x^1x^1}, 
\ \ \ \ \
(\cdot)_{x^1x^2}, 
\ \ \ \ \
(\cdot)_{x^2x^2}, 
\ \ \ \ \
(\cdot)_{x^1y}, 
\ \ \ \ \
(\cdot)_{x^2y}, 
\ \ \ \ \
(\cdot)_{yy}.
\end{equation}
In the Jacobian determinant~\thetag{ 2.16}, by replacing any one of
the three columns of first order derivatives with a column of second
order derivatives, we obtain exactly $3\times 6 =18$ possible
determinants. For instance, by replacing the third column by the
second order derivative $(\cdot )_{ x^1y}$ or the first column by the
second order derivative $(\cdot )_{ x^1x^1}$, we get:
\def\theequation{2.19}\begin{equation}
\left\vert
\begin{array}{ccc}
X_{x^1}^1 & X_{x^2}^1 & X_{x^1y}^1 \\
X_{x^1}^2 & X_{x^2}^2 & X_{x^1y}^2 \\
Y_{x^1} & Y_{x^2} & Y_{x^1y} \\
\end{array}
\right\vert 
\ \ \ \ \ \ \ \ \ \ 
{\rm or}
\ \ \ \ \ \ \ \ \ \
\left\vert
\begin{array}{ccc}
X_{x^1x^1}^1 & X_{x^2}^1 & X_y^1 \\
X_{x^1x^1}^2 & X_{x^2}^2 & X_y^2 \\
Y_{x^1x^1} & Y_{x^2} & Y_y \\
\end{array}
\right\vert.
\end{equation}
We recover the two determinants appearing in the second line
of~\thetag{ 2.10}. On the other hand, according to the 
approximation~\thetag{ 2.17}, these two determinants are
essentially equal to
\def\theequation{2.20}\begin{equation}
\left\vert
\begin{array}{ccc}
1 & 0 & X_{x^1y}^1 \\
0 & 1 & X_{x^1y}^2 \\
0 & 0 & Y_{x^1y} \\
\end{array}
\right\vert 
=
Y_{x^1y}
\ \ \ \ \ \ \ \ \ \ 
\text{\rm or to}
\ \ \ \ \ \ \ \ \ \
\left\vert
\begin{array}{ccc}
X_{x^1x^1}^1 & 0 & 0 \\
X_{x^1x^1}^2 & 1 & 0 \\
Y_{x^1x^1} & 0 & 1 \\
\end{array}
\right\vert
=
X_{x^1x^1}^1.
\end{equation}
Consequently, in the second line of~\thetag{ 2.10}, up to a change to
calligraphic letters, we recover the coefficient
\def\theequation{2.21}\begin{equation}
2\,\mathcal{Y}_{x^1y}
-
\mathcal{X}_{x^1x^1}^1
\end{equation}
of $y_{x_1}$ in the expression of ${\bf Y}_{1,1}$ in~\thetag{ 2.15}.
In conclusion, we have discovered how to pass symbolically from the
first equation~\thetag{ 2.15} to the equation~\thetag{ 2.10}
and conversely.

Translating the second equation~\thetag{ 2.15}, we deduce, {\it
without any further computation}, that the second equation which would
be obtained by developing~\thetag{ 2.9} in length, is:
\def\theequation{2.22}\begin{equation}
\aligned
0
& 
= 
y_{x^1x^2}
\cdot
\left\vert
\begin{array}{ccc}
X_{x^1}^1 & X_{x^2}^1 & X_y^1 \\
X_{x^1}^2 & X_{x^2}^2 & X_y^2 \\
Y_{x^1} & Y_{x^2} & Y_y \\
\end{array}
\right\vert
+
\left\vert
\begin{array}{ccc}
X_{x^1}^1 & X_{x^2}^1 & X_{x^1x^2}^1 \\
X_{x^1}^2 & X_{x^2}^2 & X_{x^1x^2}^2 \\
Y_{x^1} & Y_{x^2} & Y_{x^1x^2} \\
\end{array}
\right\vert
+ \\
& \
\ \ \ \ \
+
y_{x^1}
\cdot
\left\{
\left\vert
\begin{array}{ccc}
X_{x^1}^1 & X_{x^2}^1 & X_{x^2y}^1 \\
X_{x^1}^2 & X_{x^2}^2 & X_{x^2y}^2 \\
Y_{x^1} & Y_{x^2} & Y_{x^2y} \\
\end{array}
\right\vert
-
\left\vert
\begin{array}{ccc}
X_{x^1x^2}^1 & X_{x^2}^1 & X_y^1 \\
X_{x^1x^2}^2 & X_{x^2}^2 & X_y^2 \\
Y_{x^1x^2} & Y_{x^2} & Y_y \\
\end{array}
\right\vert
\right\}
+ \\
& \
\ \ \ \ \
+
y_{x^2}
\cdot
\left\{
\left\vert
\begin{array}{ccc}
X_{x^1}^1 & X_{x^2}^1 & X_{x^1y}^1 \\
X_{x^1}^2 & X_{x^2}^2 & X_{x^1y}^2 \\
Y_{x^1} & Y_{x^2} & Y_{x^1y} \\
\end{array}
\right\vert
-
\left\vert
\begin{array}{ccc}
X_{x^1}^1 & X_{x^1x^2}^1 & X_y^1 \\
X_{x^1}^2 & X_{x^1x^2}^2 & X_y^2 \\
Y_{x^1} & Y_{x^1x^2} & Y_y \\
\end{array}
\right\vert
\right\}
+ \\
\endaligned
\end{equation}
$$
\aligned
& \
\ \ \ \ \
+
y_{x^1}\,y_{x^1}
\cdot
\left\{
-
\left\vert
\begin{array}{ccc}
X_{x^2y}^1 & X_{x^2}^1 & X_y^1 \\
X_{x^2y}^2 & X_{x^2}^2 & X_y^2 \\
Y_{x^2y} & Y_{x^2} & Y_y \\
\end{array}
\right\vert
\right\}
+ \\
& \
\ \ \ \ \ 
+
y_{x^1}\,y_{x^2}
\cdot
\left\{
\left\vert
\begin{array}{ccc}
X_{x^1}^1 & X_{x^2}^1 & X_{yy}^1 \\
X_{x^1}^2 & X_{x^2}^2 & X_{yy}^2 \\
Y_{x^1} & Y_{x^2} & Y_{yy} \\
\end{array}
\right\vert
-
\left\vert
\begin{array}{ccc}
X_{x^1y}^1 & X_{x^2}^1 & X_y^1 \\
X_{x^1y}^2 & X_{x^2}^2 & X_y^2 \\
Y_{x^1y} & Y_{x^2} & Y_y \\
\end{array}
\right\vert
-
\right. \\
& \
\ \ \ \ \
\left.
-
\left\vert
\begin{array}{ccc}
X_{x^1}^1 & X_{x^2y}^1 & X_y^1 \\
X_{x^1}^2 & X_{x^2y}^2 & X_y^2 \\
Y_{x^1} & Y_{x^2y} & Y_y \\
\end{array}
\right\vert
\right\}
+
y_{x^2}\,y_{x^2}
\left\{
-
\left\vert
\begin{array}{ccc}
X_{x^1}^1 & X_{x^1y}^1 & X_y^1 \\
X_{x^1}^2 & X_{x^1y}^2 & X_y^2 \\
Y_{x^1} & Y_{x^1y} & Y_y \\
\end{array}
\right\vert
\right\}
+ \\
& \
\ \ \ \ \
+
y_{x^1}\,y_{x^1}\,y_{x^2}
\cdot
\left\{
-
\left\vert
\begin{array}{ccc}
X_{yy}^1 & X_{x^2}^1 & X_y^1 \\
X_{yy}^2 & X_{x^2}^2 & X_y^2 \\
Y_{yy} & Y_{x^2} & Y_y \\
\end{array}
\right\vert
\right\}
+ 
\\
& \
\ \ \ \ \
+
y_{x^1}\,y_{x^2}\,y_{x^2}
\cdot
\left\{
-
\left\vert
\begin{array}{ccc}
X_{x^1}^1 & X_{yy}^1 & X_y^1 \\
X_{x^1}^2 & X_{yy}^2 & X_y^2 \\
Y_{x^1} & Y_{yy} & Y_y \\
\end{array}
\right\vert
\right\}.
\endaligned
$$ 
Using the third equation~\thetag{ 2.15}, we also deduce, {\it without
any further computation}, that the third equation which would be
obtained by developing~\thetag{ 2.9} in length, is:
\def\theequation{2.23}\begin{equation}
\aligned
0
& 
= 
y_{x^2x^2}
\cdot
\left\vert
\begin{array}{ccc}
X_{x^1}^1 & X_{x^2}^1 & X_y^1 \\
X_{x^1}^2 & X_{x^2}^2 & X_y^2 \\
Y_{x^1} & Y_{x^2} & Y_y \\
\end{array}
\right\vert
+
\left\vert
\begin{array}{ccc}
X_{x^1}^1 & X_{x^2}^1 & X_{x^2x^2}^1 \\
X_{x^1}^2 & X_{x^2}^2 & X_{x^2x^2}^2 \\
Y_{x^1} & Y_{x^2} & Y_{x^2x^2} \\
\end{array}
\right\vert
+ \\
& \
\ \ \ \ \
+
y_{x^1}
\cdot
\left\{
-
\left\vert
\begin{array}{ccc}
X_{x^2x^2}^1 & X_{x^2}^1 & X_y^1 \\
X_{x^1x^2}^2 & X_{x^2}^2 & X_y^2 \\
Y_{x^2x^2} & Y_{x^2} & Y_y \\
\end{array}
\right\vert
\right\}
+ \\
& \
\ \ \ \ \
+
y_{x^2}
\cdot
\left\{
2\,
\left\vert
\begin{array}{ccc}
X_{x^1}^1 & X_{x^2}^1 & X_{x^2y}^1 \\
X_{x^1}^2 & X_{x^2}^2 & X_{x^2y}^2 \\
Y_{x^1} & Y_{x^2} & Y_{x^2y} \\
\end{array}
\right\vert
-
\left\vert
\begin{array}{ccc}
X_{x^1}^1 & X_{x^2x^2}^1 & X_y^1 \\
X_{x^1}^2 & X_{x^2x^2}^2 & X_y^2 \\
Y_{x^1} & Y_{x^2x^2} & Y_y \\
\end{array}
\right\vert
\right\}
+ \\
\endaligned
\end{equation}
$$
\aligned
& \
\ \ \ \ \ 
+
y_{x^1}\,y_{x^2}
\cdot
\left\{
-2\,
\left\vert
\begin{array}{ccc}
X_{x^2y}^1 & X_{x^2}^1 & X_y^1 \\
X_{x^2y}^2 & X_{x^2}^2 & X_y^2 \\
Y_{x^2y} & Y_{x^2} & Y_y \\
\end{array}
\right\vert
\right\}
+ \\
& \
\ \ \ \ \
+
y_{x^2}\,y_{x^2}
\cdot
\left\{
\left\vert
\begin{array}{ccc}
X_{x^1}^1 & X_{x^2}^1 & X_{yy}^1 \\
X_{x^1}^2 & X_{x^2}^2 & X_{yy}^2 \\
Y_{x^1} & Y_{x^2} & Y_{yy} \\
\end{array}
\right\vert
-
2\,
\left\vert
\begin{array}{ccc}
X_{x^1}^1 & X_{x^2y}^1 & X_y^1 \\
X_{x^1}^2 & X_{x^2y}^2 & X_y^2 \\
Y_{x^1} & Y_{x^2y} & Y_y \\
\end{array}
\right\vert
\right\}
+ \\
& \
\ \ \ \ \
+
y_{x^1}\,y_{x^2}\,y_{x^2}
\cdot
\left\{
-
\left\vert
\begin{array}{ccc}
X_{yy}^1 & X_{x^2}^1 & X_y^1 \\
X_{yy}^2 & X_{x^2}^2 & X_y^2 \\
Y_{yy} & Y_{x^2} & Y_y \\
\end{array}
\right\vert
\right\}
+ \\
& \
\ \ \ \ \
+
y_{x^2}\,y_{x^2}\,y_{x^2}
\cdot
\left\{
-
\left\vert
\begin{array}{ccc}
X_{x^1}^1 & X_{yy}^1 & X_y^1 \\
X_{x^1}^2 & X_{yy}^2 & X_y^2 \\
Y_{x^1} & Y_{yy} & Y_y \\
\end{array}
\right\vert
\right\}.
\endaligned
$$

\subsection*{ 2.24.~Appropriate formalism}
To describe the combinatorics underlying formulas~\thetag{ 2.10},
\thetag{ 2.22} and~\thetag{ 2.23}, as in~\cite{ me2004}, let us
introduce the following notation for the Jacobian determinant:
\def\theequation{2.25}\begin{equation}
\Delta(x^1\vert x^2\vert y)
:=
\left\vert
\begin{array}{ccc}
X_{x^1}^1 & X_{x^2}^1 & X_y^1 \\
X_{x^1}^2 & X_{x^2}^2 & X_y^2 \\
Y_{x^1} & Y_{x^2} & Y_y
\end{array}
\right\vert.
\end{equation}
Here, in the notation $\Delta(x^1\vert x^2\vert y)$, the three spaces
between the two vertical lines $\vert$ refer to the three columns of
the Jacobian determinant, and the terms $x^1$, $x^2$, $y$ in $(x^1
\vert x^2 \vert y)$ designate the partial derivatives appearing in
each column. Accordingly, in the following two examples of {\sl modified
Jacobian determinants}:
\def\theequation{2.26}\begin{equation}
\left\{
\aligned
&
\Delta( \underline{ x^1 x^2} \vert x^2 \vert y)
:=
\left\vert
\begin{array}{ccc}
X_{\underline{x^1x^2}}^1 & X_{x^2}^1 & X_y^1 \\
X_{\underline{x^1x^2}}^2 & X_{x^2}^2 & X_y^2 \\
Y_{\underline{x^1x^2}} & Y_{x^2} & Y_y
\end{array}
\right\vert 
\ \ \ \ \
\ \ \ \ \
{\rm and}
\\
&
\Delta(x^1\vert x^2\vert \underline{x^1y})
:=
\left\vert
\begin{array}{ccc}
X_{x^1}^1 & X_{x^2}^1 & X_{\underline{x^1y}}^1 \\
X_{x^1}^2 & X_{x^2}^2 & X_{\underline{x^1y}}^2 \\
Y_{x^1} & Y_{x^2} & Y_{\underline{x^1y}}
\end{array}
\right\vert,
\endaligned\right.
\end{equation}
we simply mean which column of first order derivatives
is replaced by a column of second order derivatives
in the original Jacobian determinant.

As there are $6$ possible second order derivatives $(\cdot )_{ x^1
x^1}$, $(\cdot )_{ x^1 x^2}$, $(\cdot )_{ x^1 x^y}$, $(\cdot )_{ x^2
x^2}$, $(\cdot )_{ x^2 y}$ and $(\cdot )_{ y y}$ together with $3$
columns, we obtain $3\times 6 = 18$ possible modified Jacobian
determinants:
\def\theequation{2.27}\begin{equation}
\left\{
\aligned
& \
\Delta(x^1x^1 \vert x^2 \vert y) \ \ \ \ \
&
\Delta(x^1 \vert x^1x^1 \vert y) \ \ \ \ \ 
&
\Delta(x^1 \vert x^2 \vert x^1x^1) \ \ \ \ \ 
\\
& \
\Delta(x^1x^2 \vert x^2 \vert y) \ \ \ \ \
&
\Delta(x^1 \vert x^1x^2 \vert y) \ \ \ \ \ 
&
\Delta(x^1 \vert x^2 \vert x^1x^2) \ \ \ \ \ 
\\
& \
\Delta(x^1y \vert x^2 \vert y) \ \ \ \ \
&
\Delta(x^1 \vert x^1y \vert y) \ \ \ \ \ 
&
\Delta(x^1 \vert x^2 \vert x^1y) \ \ \ \ \ 
\\
& \
\Delta(x^2x^2 \vert x^2 \vert y) \ \ \ \ \
&
\Delta(x^1 \vert x^2x^2 \vert y) \ \ \ \ \ 
&
\Delta(x^1 \vert x^2 \vert x^2x^2) \ \ \ \ \ 
\\
& \
\Delta(x^2y \vert x^2 \vert y) \ \ \ \ \
&
\Delta(x^1 \vert x^2y \vert y) \ \ \ \ \ 
&
\Delta(x^1 \vert x^2 \vert x^2y) \ \ \ \ \ 
\\
& \
\Delta(yy \vert x^2 \vert y) \ \ \ \ \
&
\Delta(x^1 \vert yy \vert y) \ \ \ \ \ 
&
\Delta(x^1 \vert x^2 \vert yy). \ \ \ \ \ 
\\
\endaligned\right.
\end{equation}

Next, we observe that if we want to solve with
respect to $y_{ x^1x^1}$ in~\thetag{ 2.10}, 
with respect to $y_{ x^1x^2}$ in~\thetag{ 2.22}
and with respect to $y_{ x^2x^2}$ in~\thetag{ 2.23}, 
we have to divide by the Jacobian determinant 
$\Delta( x^1 \vert x^2\vert y)$. 
Consequently, we introduce 
$18$ new {\sl square functions} as follows:
\def\theequation{2.28}\begin{equation}
\small
\left\{
\aligned
\square_{x^1x^1}^1:=&\ \frac{ \Delta(x^1x^1\vert x^2 \vert y)}{
\Delta(x^1\vert x^2 \vert y)} 
\ \ \ \ \ & 
\square_{x^1x^2}^1:=&\ \frac{ \Delta(x^1x^2\vert x^2 \vert y)}{
\Delta(x^1\vert x^2 \vert y)} 
\ \ \ \ \ & 
\square_{x^1y}^1:=&\ \frac{ \Delta(x^1y\vert x^2 \vert y)}{
\Delta(x^1\vert x^2 \vert y)} 
\ \ \ \ \ 
\\
\square_{x^2x^2}^1:=&\ \frac{ \Delta(x^2x^2\vert x^2 \vert y)}{
\Delta(x^1\vert x^2 \vert y)} 
\ \ \ \ \ & 
\square_{x^2y}^1:=&\ \frac{ \Delta(x^2y\vert x^2 \vert y)}{
\Delta(x^1\vert x^2 \vert y)} 
\ \ \ \ \ & 
\square_{yy}^1:=&\ \frac{ \Delta(yy\vert x^2 \vert y)}{
\Delta(x^1\vert x^2 \vert y)} 
\ \ \ \ \ 
\\
\square_{x^1x^1}^2:=&\ \frac{ \Delta(x^1\vert x^1x^1 \vert y)}{
\Delta(x^1\vert x^2 \vert y)} 
\ \ \ \ \ & 
\square_{x^1x^2}^2:=&\ \frac{ \Delta(x^1\vert x^1x^2 \vert y)}{
\Delta(x^1\vert x^2 \vert y)} 
\ \ \ \ \ & 
\square_{x^1y}^2:=&\ \frac{ \Delta(x^1\vert x^1y \vert y)}{
\Delta(x^1\vert x^2 \vert y)} 
\ \ \ \ \
\\
\square_{x^2x^2}^2:=&\ \frac{ \Delta(x^1\vert x^2x^2 \vert y)}{
\Delta(x^1\vert x^2 \vert y)} 
\ \ \ \ \ & 
\square_{x^2y}^2:=&\ \frac{ \Delta(x^1\vert x^2y \vert y)}{
\Delta(x^1\vert x^2 \vert y)} 
\ \ \ \ \ & 
\square_{yy}^2:=&\ \frac{ \Delta(x^1\vert yy \vert y)}{
\Delta(x^1\vert x^2 \vert y)} 
\ \ \ \ \ 
\\
\square_{x^1x^1}^3:=&\ \frac{ \Delta(x^1\vert x^2 \vert x^1x^1)}{
\Delta(x^1\vert x^2 \vert y)} 
\ \ \ \ \ & 
\square_{x^1x^2}^3:=&\ \frac{ \Delta(x^1\vert x^2 \vert x^1x^2)}{
\Delta(x^1\vert x^2 \vert y)} 
\ \ \ \ \ & 
\square_{x^1y}^3:=&\ \frac{ \Delta(x^1\vert x^2 \vert x^1y)}{
\Delta(x^1\vert x^2 \vert y)} 
\ \ \ \ \ 
\\
\square_{x^2x^2}^3:=&\ \frac{ \Delta(x^1\vert x^2 \vert x^2x^2)}{
\Delta(x^1\vert x^2 \vert y)} 
\ \ \ \ \ & 
\square_{x^2y}^3:=&\ \frac{ \Delta(x^1\vert x^2 \vert x^2y)}{
\Delta(x^1\vert x^2 \vert y)} 
\ \ \ \ \ & 
\square_{yy}^3:=&\ \frac{ \Delta(x^1\vert x^2 \vert yy)}{
\Delta(x^1\vert x^2 \vert y)}.
\ \ \ \ \ 
\\
\endaligned\right.
\end{equation} 

Thanks to these notations, we can rewrite the three equations~\thetag{
2.10}, \thetag{ 2.22} and~\thetag{ 2.23} in a more compact style.

\def\thelemma{2.29}\begin{lemma}
A completely integrable system of
{\rm three} second order partial 
differential equations 
\def\theequation{2.30}\begin{equation}
\left\{
\aligned
y_{x^1x^1} (x)
&
= 
F^{1,1}
\left(
x^1, x^2, y(x), y_{x^1}(x), y_{x^2}(x) 
\right), 
\\
y_{x^1x^2} (x)
& 
= 
F^{1,2}
\left(
x^1, x^2, y(x), y_{x^1}(x), y_{x^2}(x) 
\right), 
\\
y_{x^2x^2} (x)
&
= 
F^{2,2}
\left(
x^1, x^2, y(x), y_{x^1}(x), y_{x^2}(x) 
\right), 
\\
\endaligned\right.
\end{equation}
is equivalent to the simplest
system $Y_{X^1 X^1} = 0$, 
$Y_{ X^1 X^2} = 0$, $Y_{X^2, X^2} = 0$,
{\rm if and only if}
there exist local 
$\K$-analytic functions 
$X^1$, $X^2$, $Y$ such that it
may be written under the
specific form{\rm :}
\def\theequation{2.31}\begin{equation}
\small
\left\{
\aligned
y_{x^1x^1}
& 
=
- 
\square_{x^1x^1}^3
+
y_{x^1}
\cdot
\left(
-
2\,\square_{x^1y}^3
+
\square_{x^1x^1}^1
\right)
+
y_{x^2}
\cdot
\left(
\square_{x^1x^1}^2
\right)
+ \\
& \
\ \ \ \ \ 
+
y_{x^1}\,y_{x^1}
\cdot
\left(
-
\square_{yy}^3
+
2\,\square_{x^1y}^1
\right)
+
y_{x^1}\,y_{x^2}
\cdot
\left(
2\,\square_{x^1y}^2
\right)
+ \\
& \
\ \ \ \ \ 
+
y_{x^1}\,y_{x^1}\,y_{x^1}
\cdot
\left(
\square_{yy}^1
\right)
+
y_{x^1}\,y_{x^1}\,y_{x^2}
\cdot
\left(
\square_{yy}^2
\right),
\\
y_{x^1x^2}
& 
=
- 
\square_{x^1x^2}^3
+
y_{x^1}
\cdot
\left(
-
\square_{x^2y}^3
+
\square_{x^1x^2}^1
\right)
+
y_{x^2}
\cdot
\left(
-
\square_{x^1y}^3
+
\square_{x^1x^2}^2
\right)
+ \\
& \
\ \ \ \ \
+
y_{x^1}\,y_{x^1}
\cdot
\left(
\square_{x^2y}^1
\right)
+
y_{x^1}\,y_{x^2}
\cdot
\left(
-
\square_{yy}^3
+
\square_{x^1y}^1
+
\square_{x^2y}^2
\right)
+ \\
& \
\ \ \ \ \ 
+
y_{x^2}\,y_{x^2}
\cdot
\left(
\square_{x^1y}^2
\right)
+
y_{x^1}\,y_{x^1}\,y_{x^2}
\cdot
\left(
\square_{yy}^1
\right)
+
y_{x^1}\,y_{x^2}\,y_{x^2}
\cdot
\left(
\square_{yy}^2
\right),
\\
y_{x^2x^2}
& 
=
- 
\square_{x^2x^2}^3
+
y_{x^1}
\cdot
\left(
\square_{x^2x^2}^1
\right)
+
y_{x^2}
\cdot
\left(
-
2\,\square_{x^2y}^3
+
\square_{x^2x^2}^2
\right)
+ \\
& \
\ \ \ \ \ 
+
y_{x^1}\,y_{x^2}
\cdot
\left(
2\,\square_{x^2y}^1
\right)
+
y_{x^2}\,y_{x^2}
\cdot
\left(
-
\square_{yy}^3
+
2\,\square_{x^2y}^2
\right)
+ \\
& \
\ \ \ \ \ 
+
y_{x^1}\,y_{x^2}\,y_{x^2}
\cdot
\left(
\square_{yy}^1
\right)
+
y_{x^2}\,y_{x^2}\,y_{x^2}
\cdot
\left(
\square_{yy}^2
\right).
\endaligned\right.
\end{equation}
\end{lemma}

\subsection*{ 2.32.~General formulas}
The formal dictionary between the original determinantial
formulas~\thetag{ 2.10}, \thetag{ 2.22}, \thetag{ 2.23}, between the
coefficients~\thetag{ 2.15} of the second order prolongation of a
vector field and between the new square formulas~\thetag{ 2.31} above
is evident. Consequently, {\it without any computation}, just by
translating the family of formulas~\thetag{ 2.14}, we may deduce the
exact formulation of the desired generalization of Lemma~2.29 above.

\def\thelemma{2.33}\begin{lemma}
A completely integrable system of second order partial 
differential equations of the form
\def\theequation{2.34}\begin{equation}
y_{x^{j_1}x^{j_2}}(x) 
=
F^{j_1, j_2} 
\left(
x, y(x), y_{x^1}(x),\dots,y_{x^n}(x)
\right), 
\ \ \ \ \ \ \ \ 
j_1,j_2=1,\dots n, 
\end{equation}
is equivalent to the simplest system $Y_{X^{j_1} X^{j_2}} = 0$, $j_1,
j_2 = 1, \dots, n$, {\rm if and only if} there exist local
$\K$-analytic functions $X^{l}$, $Y$ such that it may be written under
the specific form{\rm :}
\def\theequation{2.35}\begin{equation}
\left\{
\aligned
y_{x^{j_1}x^{j_2}}
&
=
-
\square_{x^{j_1}x^{j_2}}^{n+1}
+
\sum_{k_1=1}^n\,y_{x^{k_1}}
\cdot
\left\{
\left(
\square_{x^{j_1}x^{j_2}}^{k_1}
-
\delta_{j_1}^{k_1}\,\square_{x^{j_2}y}^{n+1}
-
\delta_{j_2}^{k_1}\,\square_{x^{j_1}y}^{n+1}
\right)
+ 
\right.
\\
& \
\ \ \ \ \
\left.
+
y_{x^{j_1}}
\cdot
\left(
\square_{x^{j_2}y}^{k_1}
-
\frac{1}{2}\,\delta_{j_2}^{k_1}\,\square_{yy}^{n+1}
\right)
+
y_{x^{j_2}}
\cdot
\left(
\square_{x^{j_1}y}^{k_1}
-
\frac{1}{2}\,\delta_{j_1}^{k_1}\,\square_{yy}^{n+1}
\right)
+ 
\right.
\\
& \
\ \ \ \ \
\left.
+
y_{x^{j_1}}\,y_{x^{j_2}}
\cdot
\left(
\square_{yy}^{k_1}
\right)
\right\}.
\endaligned\right.
\end{equation}
\end{lemma}

Of course, to define the square functions in the context of $n
\geqslant 2$
independent variables $(x^1, x^2, \dots, x^n)$, we introduce the
Jacobian determinant
\def\theequation{2.36}\begin{equation}
\Delta(x^1 \vert x^2 \vert \cdots \vert x^n \vert y)
:= 
\left\vert
\begin{array}{cccc}
X_{x^1}^1 & \cdots & X_{x^n}^1 & X_y^1 \\
\vdots & \cdots & \vdots & \vdots \\
X_{x^1}^n & \cdots & X_{x^n}^n & X_y^n \\
Y_{x^1} & \cdots & Y_{x^n} & Y_y \\
\end{array}
\right\vert,
\end{equation}
together with its modifications
\def\theequation{2.37}\begin{equation}
\Delta
\left(
x^1\vert\cdots\vert^{k_1}\,x^{j_1}\,x^{j_2} \vert \cdots \vert y
\right),
\end{equation}
in which the $k_1$-th column of partial first order derivatives
$\vert^{k_1}\, x^{k_1} \vert$ is replaced by the column $\vert^{k_1}
\, x^{j_1} x^{j_2} \vert$ of partial derivatives. Here, the indices
$k_1$, $j_1$, $j_2$ satisfy $1\leqslant k_1, j_1, j_2 \leqslant n+1$,
with the convention that we adopt the notational equivalence
\def\theequation{2.38}\begin{equation}
\fbox{
$x^{n+1}\equiv y$
}.
\end{equation}
This convention will be convenient to write some of our general
formulas in the sequel.

As we promised to only summarize the proof of Theorem~1.7 in this
paper, we will not develope the proof of Lemma~2.33: 
it is similar to the proof of
Lemma~3.32 in~\cite{ me2004}.

\section*{ \S3.~First and second auxiliary system}

\subsection*{3.1.~Functions $G_{j_1, j_2}$, 
$H_{j_1, j_2}^{k_1}$, $L_{j_1}^{k_1}$ and $M^{k_1}$} To discover the
four families of functions appearing in the statement of Theorem~1.7,
by comparing~\thetag{ 2.35} and~\thetag{ 1.10}, it suffices (of course)
to set:
\def\theequation{3.2}\begin{equation}
\left\{
\aligned
G_{j_1,j_2}
&
:= 
-
\square_{x^{j_1}x^{j_2}}^{n+1}, 
\\
H_{j_1,j_2}^{k_1}
&
:=
\square_{x^{j_1}x^{j_2}}^{k_1}
-
\delta_{j_1}^{k_1}\,\square_{x^{j_2}y}^{n+1}
-
\delta_{j_2}^{k_1}\,\square_{x^{j_1}y}^{n+1}, 
\\
L_{j_1}^{k_1}
&
:=
2\,\square_{x^{j_1}y}^{k_1}
-
\delta_{j_1}^{k_1}\,\square_{yy}^{n+1},
\\
M^{k_1}
&
:=
\square_{yy}^{k_1}.
\endaligned\right.
\end{equation}

Consequently, we have shown the ``only if'' part of Theorem~1.7, which
is the easiest implication.

To establish the ``if'' part, by far the most difficult implication,
the very main lemma can be stated as follows.

\def\thelemma{3.3}\begin{lemma}
The partial differerential relations {\rm (I')}, {\rm (II')}, {\rm
(III')} and {\rm (IV')} which express in length the compatibility
conditions~\thetag{ 1.11} are necessary and sufficient
for the existence of functions $X^l$, $Y$ of $(x^{l_1}, y)$ satisfying
the second order nonlinear system of partial differential
equations~\thetag{ 3.2} above.
\end{lemma}

Indeed, the collection of equations~\thetag{ 3.2} is a system of
partial differential equations with unknowns $X^l$, $Y$, by virtue of
the definition of the square functions.

\subsection*{ 3.4.~First auxiliary system}
To proceed further, we observe that there are $(m+1)$ more square
functions than functions $G_{j_1, j_2}$, $H_{j_1, j_2}^{k_1}$,
$L_{j_1}^{k_1}$ and $M^{k_1}$. Indeed, a simple counting yields:
\def\theequation{3.5}\begin{equation}
\left\{
\aligned
{}
& \
\# \{\square_{x^{j_1}x^{j_2}}^{k_1}\}
=
\frac{n^2(n+1)}{2},
& \
\# \{\square_{x^{j_1}y}^{k_1} \}
=
n^2, \\
& \
\# \{\square_{yy}^{k_1}\}
=
n, 
& \
\# \{\square_{x^{j_1}x^{j_2}}^{n+1} \}
= 
\frac{n(n+1)}{2}, \\
& \
\# \{\square_{x^{j_1}y}^{n+1}\}
= 
n, 
\ \ \ \ \ \ \ \ \ \ \ 
\ \ \ \ \ \ \ \ \ \ \ 
& \
\# \{\square_{yy}^{n+1} \}
= 
1, \\
\endaligned\right.
\end{equation}
whereas
\def\theequation{3.6}\begin{equation}
\left\{
\aligned
{}
& \
\# \{G^{j_1,j_2}\}
= 
\frac{n(n+1)}{2}, 
& \ 
\#
\{H_{j_1,j_2}^{k_1} \}
=
\frac{n^2(n+1)}{2}, \\
& \
\# \{L_{j_1}^{k_1} \}
=
n^2, 
& \ 
\ \ \ \ \ \ \ \ \ \ \ 
\ \ \ \ \ \ \ \ \ \ \
\# \{M^{k_1}\}
=
n.
\endaligned\right.
\end{equation}
Here, the indices $j_1$, $j_2$, $k_1$ satisfy $1\leqslant j_1, j_2,
k_1 \leqslant n$. Similarly as in~\cite{ me2004}, to transform the
system~\thetag{ 3.2} in a true complete system, let us introduce
functions $\Pi_{j_1,j_2}^{k_1}$ of $(x^{l_1}, y)$, where $1\leqslant
j_1, j_2, k_1 \leqslant n+1$, which satisfy the symmetry
$\Pi_{j_1,j_2}^{k_1} = \Pi_{j_1, j_1}^{k_1}$, and let us introduce the
following {\sl first auxiliary system}:
\def\theequation{3.7}\begin{equation}
\left\{
\aligned
\square_{x^{j_1}x^{j_2}}^{k_1}
=
\Pi_{j_1,j_2}^{k_1}, 
\ \ \ \ \ \ \ &
\square_{x^{j_1}y}^{k_1}
=
\Pi_{j_1, n+1}^{k_1}, 
\ \ \ \ \ \ \ &
\square_{yy}^{k_1}
=
\Pi_{n+1,n+1}^{k_1}, 
\\
\square_{x^{j_1}x^{j_2}}^{n+1}
=
\Pi_{j_1,j_2}^{n+1}, 
\ \ \ \ \ \ \ &
\square_{x^{j_1}y}^{n+1}
=
\Pi_{j_1, n+1}^{n+1}, 
\ \ \ \ \ \ \ &
\square_{yy}^{n+1}
=
\Pi_{n+1,n+1}^{n+1}. 
\endaligned\right.
\end{equation}
It is complete. The necessary and sufficient
conditions for the existence of solutions
$X^l$, $Y$ follow by cross
differentiations.

\def\thelemma{3.8}\begin{lemma}
For all $j_1, j_2, j_3, k_1 = 1, 2, \dots, n+1$, 
we have the cross differentiation relations
\def\theequation{3.9}\begin{equation}
\left(
\square_{x^{j_1}x^{j_2}}^{k_1}
\right)_{x^{j_3}}
-
\left(
\square_{x^{j_1}x^{j_3}}^{k_1}
\right)_{x^{j_2}}
=
-
\sum_{k_2=1}^{n+1}\,
\square_{x^{j_1}x^{j_2}}^{k_2}\,\square_{x^{j_3}x^{k_2}}^{k_1}
+
\sum_{k_2=1}^{n+1}\,
\square_{x^{j_1}x^{j_3}}^{k_2}\,\square_{x^{j_2}x^{k_2}}^{k_1}.
\end{equation}
\end{lemma}

The proof of this lemma is exactly the same as the proof of Lemma~3.40
in~\cite{ me2004}.

As a direct consequence, we deduce that a necessary and sufficient
condition for the existence of solutions $\Pi_{j_1, j_2}^{ k_1}$ to
the first auxiliary system is that they satisfy the following
compatibility partial differential relations:
\def\theequation{3.10}\begin{equation}
\frac{\partial \Pi_{j_1,j_2}^{k_1}}{\partial x^{j_3}}
-
\frac{\partial \Pi_{j_1,j_3}^{k_1}}{\partial x^{j_2}}
=
-
\sum_{k_2=1}^{n=1}\,\Pi_{j_1,j_2}^{k_2}
\cdot
\Pi_{j_3,k_2}^{k_1}
+
\sum_{k_2=1}^{n=1}\,\Pi_{j_1,j_3}^{k_2}
\cdot
\Pi_{j_2,k_2}^{k_1},
\end{equation}
for all $j_1, j_2, j_3, k_1 = 1, \dots, n+1$.

We shall have to specify this system in length according to the
splitting $\{1, 2, \dots, n\}$ and $\{n+1\}$ of the indices of
coordinates. We obtain six families of equations equivalent
to~\thetag{ 3.10} just above:
\def\theequation{3.11}\begin{equation}
\left\{
\aligned
\left(
\Pi_{j_1,j_2}^{n+1}
\right)_{x^{j_3}}
-
\left(
\Pi_{j_1,j_3}^{n+1}
\right)_{x^{j_2}}
&
=
-
\sum_{k_2=1}^n\,\Pi_{j_1,j_2}^{k_2}\,\Pi_{j_3,k_2}^{n+1}
-
\Pi_{j_1,j_2}^{n+1}\,\Pi_{j_3,n+1}^{n+1}
+ \\
& \
\ \ \ \ \ \ \ \ \ \
+
\sum_{k_2=1}^n\,\Pi_{j_1,j_3}^{k_2}\,\Pi_{j_2,k_2}^{n+1}
+
\Pi_{j_1,j_3}^{n+1}\,\Pi_{j_2,n+1}^{n+1}, 
\\
\left(
\Pi_{j_1,j_2}^{n+1}
\right)_y
-
\left(
\Pi_{j_1,n+1}^{n+1}
\right)_{x^{j_2}}
&
=
-
\sum_{k_2=1}^n\,\Pi_{j_1,j_2}^{k_2}\,\Pi_{n+1,k_2}^{n+1}
-
\Pi_{j_1,j_2}^{n+1}\,\Pi_{n+1,n+1}^{n+1}
+ \\
& \
\ \ \ \ \ \ \ \ \ \
+
\sum_{k_2=1}^n\,\Pi_{j_1,n+1}^{k_2}\,\Pi_{j_2,k_2}^{n+1}
+
\Pi_{j_1,n+1}^{n+1}\,\Pi_{j_2,n+1}^{n+1}, 
\\
\left(
\Pi_{j_1,n+1}^{n+1}
\right)_y
-
\left(
\Pi_{n+1,n+1}^{n+1}
\right)_{x^{j_1}}
&
=
-
\sum_{k_2=1}^n\,\Pi_{j_1,n+1}^{k_2}\,\Pi_{n+1,k_2}^{n+1}
-
\underline{
\Pi_{j_1,n+1}^{n+1}\,\Pi_{n+1,n+1}^{n+1} 
}_{\circ \! \! \! \! \! \tiny{\sf a}}
+ 
\\
& \
\ \ \ \ \ \ \ \ \ \
+
\sum_{k_2=1}^n\,\Pi_{n+1,n+1}^{k_2}\,\Pi_{j_1,k_2}^{n+1}
+
\underline{
\Pi_{n+1,n+1}^{n+1}\,\Pi_{j_1,n+1}^{n+1}
}_{\circ \! \! \! \! \! \tiny{\sf a}} \ , 
\\
\left(
\Pi_{j_1,j_2}^{k_1}
\right)_{x^{j_3}}
-
\left(
\Pi_{j_1,j_3}^{k_1}
\right)_{x^{j_2}}
&
=
-
\sum_{k_2=1}^n\,\Pi_{j_1,j_2}^{k_2}\,\Pi_{j_3,k_2}^{k_1}
-
\Pi_{j_1,j_2}^{n+1}\,\Pi_{j_3,n+1}^{k_1}
+ \\
& \
\ \ \ \ \ \ \ \ \ \
+
\sum_{k_2=1}^n\,\Pi_{j_1,j_3}^{k_2}\,\Pi_{j_2,k_2}^{k_1}
+
\Pi_{j_1,j_3}^{n+1}\,\Pi_{j_2,n+1}^{k_1}, 
\\
\left(
\Pi_{j_1,j_2}^{k_1}
\right)_y
-
\left(
\Pi_{j_1,n+1}^{k_1}
\right)_{x^{j_2}}
&
=
-
\sum_{k_2=1}^n\,\Pi_{j_1,j_2}^{k_2}\,\Pi_{n+1,k_2}^{k_1}
-
\Pi_{j_1,j_2}^{n+1}\,\Pi_{n+1,n+1}^{k_1}
+ \\
& \
\ \ \ \ \ \ \ \ \ \
+
\sum_{k_2=1}^n\,\Pi_{j_1,n+1}^{k_2}\,\Pi_{j_2,k_2}^{k_1}
+
\Pi_{j_1,n+1}^{n+1}\,\Pi_{j_2,n+1}^{k_1}, 
\\
\left(
\Pi_{j_1,n+1}^{k_1}
\right)_y
-
\left(
\Pi_{n+1,n+1}^{k_1}
\right)_{x^{j_1}}
&
=
-
\sum_{k_2=1}^n\,\Pi_{j_1,n+1}^{k_2}\,\Pi_{n+1,k_2}^{k_1}
-
\Pi_{j_1,n+1}^{n+1}\,\Pi_{n+1,n+1}^{k_1}
+ \\
& \
\ \ \ \ \ \ \ \ \ \
+
\sum_{k_2=1}^n\,\Pi_{n+1,n+1}^{k_2}\,\Pi_{j_1,k_2}^{k_1}
+
\Pi_{n+1,n+1}^{n+1}\,\Pi_{j_1,n+1}^{k_1}. 
\\
\endaligned\right.
\end{equation}
where the indices $j_1, j_2, j_3, k_1$ vary in the set $\{ 1, 2, 1,
\dots, n\}$.

\subsection*{3.12.~Principal unknowns}
As there are $(m+1)$ more square (or Pi) functions than the functions
$G_{j_1, j_2}$, $H_{j_1, j_2}^{k_1}$, $L_{j_1}^{k_1}$ and $M^{ k_1}$,
we cannot invert directly the linear system~\thetag{ 3.2}. To
quasi-inverse it, we choose the $(m+1)$ specific square functions
\def\theequation{3.13}\begin{equation}
\Theta^1
:= 
\square_{x^1x^1}^1, 
\ \ \ \ \
\Theta^2
:= 
\square_{x^2x^2}^2,
\cdots\cdots,
\Theta^{n+1}
:= 
\square_{x^{n+1}x^{n+1}}^{n+1},
\end{equation}
calling them {\sl principal unknowns}, 
and we get the quasi-inversion:
\def\theequation{3.14}\begin{equation}
\left\{
\aligned
\Pi_{j_1,j_2}^{k_1}
&
=
\square_{x^{j_1}x^{j_2}}^{k_1}
=
H_{j_1,j_2}^{k_1}
-
\frac{1}{2}\,\delta_{j_1}^{k_1}\,H_{j_2,j_2}^{j_2}
-
\frac{1}{2}\,\delta_{j_2}^{k_1}\,H_{j_1,j_1}^{j_1}
+
\frac{1}{2}\,\delta_{j_1}^{k_1}\,\Theta^{j_2}
+
\frac{1}{2}\,\delta_{j_2}^{k_1}\,\Theta^{j_1}, 
\\
\Pi_{j_1,n+1}^{k_1}
&
=
\square_{x^{j_1}y}^{k_1}
=
\frac{1}{2}\,L_{j_1}^{k_1}
+
\frac{1}{2}\,\delta_{j_1}^{k_1}\,\Theta^{n+1}, 
\\
\Pi_{n+1,n+1}^{k_1}
&
=
\square_{yy}^{k_1}
=
M^{k_1},
\\
\Pi_{j_1,j_2}^{n+1}
&
=
\square_{x^{j_1}x^{j_2}}^{n+1}
=
-
G_{j_1,j_2},
\\
\Pi_{j_1,n+1}^{n+1}
&
=
\square_{x^{j_1}y}^{n+1}
=
-
\frac{1}{2}\,H_{j_1, j_1}^{j_1}
+
\frac{1}{2}\,\Theta^{j_1}.
\endaligned\right.
\end{equation}

\subsection*{ 3.15.~Second auxiliary system}
Replacing the five families of functions $\Pi_{j_1, j_2}^{k_1}$,
$\Pi_{j_1, n+1}^{k_1}$, $\Pi_{ n+1, n+1}^{k_1}$, $\Pi_{j_1,
j_2}^{n+1}$, $\Pi_{j_1, n+1}^{n+1}$ by their values obtained
in~\thetag{ 3.14} just above together with the
principal unknowns
\def\theequation{3.16}\begin{equation}
\Pi_{j_1, j_1}^{j_1} 
=
\Theta^{j_1}, 
\ \ \ \ \ \ \ \ \ \
\Pi_{n+1,n+1}^{n+1}
= 
\Theta^{n+1},
\end{equation}
in the six equations $(3.11)_1$, $(3.11)_2$, $(3.11)_3$, $(3.11)_4$,
$(3.11)_5$ and $(3.11)_6$, after hard computations that we will not
reproduce here, we obtain six families of equations. From now on, we
abbreviate every sum $\sum_{k=1}^n$ as $\sum_{k_1}$.

Firstly:
\def\theequation{3.17}\begin{equation}
0
=
\underline{
G_{j_1,j_2,x^{j_3}}
-
G_{j_1,j_3,x^{j_2}}
}
+
\sum_{k_1}\,G_{j_3,k_1}\,H_{j_1,j_2}^{k_1}
-
\sum_{k_1}\,G_{j_2,k_1}\,H_{j_1,j_3}^{k_1}.
\end{equation}
This is (I') of Theorem~1.7.
Just above and below, we plainly underline the
monomials involving a first order derivative. Secondly:
\def\theequation{3.18}\begin{equation}
\left\{
\aligned
\Theta_{x^{j_2}}^{j_1}
&
=
-
\underline{
2\,G_{j_1,j_2,y}
+
H_{j_1,j_1,x^{j_2}}^{j_1}
}
+ \\
& \
\ \ \ \ \
+
\sum_{k_1}\,G_{j_2,k_1}\,L_{j_1}^{k_1}
+
\frac{1}{2}\,H_{j_1,j_1}^{j_1}\,H_{j_2,j_2}^{j_2}
-
\sum_{k_1}\,H_{j_1,j_2}^{k_1}\,H_{k_1,k_1}^{k_1}
- \\
& \
\ \ \ \ \
-
G_{j_1,j_2}\,\Theta^{n+1}
-
\frac{1}{2}\,H_{j_1,j_1}^{j_1}\,\Theta^{j_2}
-
\frac{1}{2}\,H_{j_2,j_2}^{j_2}\,\Theta^{j_1}
+
\sum_{k_1}\,H_{j_1,j_2}^{k_1}\,\Theta^{k_1}
+ \\
& \
\ \ \ \ \
+ 
\frac{1}{2}\,\Theta^{j_1}\,\Theta^{j_2}.
\endaligned\right.
\end{equation}
Thirdly:
\def\theequation{3.19}\begin{equation}
\left\{
\aligned
-
\Theta_{x^{j_1}}^{n+1}
+
\frac{1}{2}\,\Theta_y^{j_1}
& 
=
\underline{
\frac{1}{2}\,H_{j_1,j_1,y}^{j_1}
}
- \\
& \
\ \ \ \ \ \
-
\sum_{k_1}\,G_{j_1,k_1}\,M^{k_1}
+
\frac{1}{4}\,\sum_{k_1}\,H_{k_1,k_1}^{k_1}\,L_{j_1}^{k_1}
+ \\
& \
\ \ \ \ \ 
+
\frac{1}{4}\,H_{j_1,j_1}^{j_1}\,\Theta^{n+1}
-
\frac{1}{4}\,\sum_{k_1}\,L_{j_1}^{k_1}\,\Theta^{k_1}
-
\frac{1}{4}\,\Theta^{j_1}\,\Theta^{n+1}.
\endaligned\right.
\end{equation}
Fourtly:
\def\theequation{3.20}\begin{equation}
\left\{
\aligned
&
\frac{1}{2}\,\delta_{j_1}^{k_1}\,\Theta_{x^{j_3}}^{j_2}
-
\frac{1}{2}\,\delta_{j_1}^{k_1}\,\Theta_{x^{j_2}}^{j_3}
+
\frac{1}{2}\,\delta_{j_2}^{k_1}\,\Theta_{x^{j_3}}^{j_1}
-
\frac{1}{2}\,\delta_{j_3}^{k_1}\,\Theta_{x^{j_2}}^{j_1}
= \\
& \
\ \ \ \ \
=
-
\underline{
H_{j_1,j_2,x^{j_3}}^{k_1}
+
H_{j_1,j_3,x^{j_2}}^{k_1}
-
\frac{1}{2}\,\delta_{j_1}^{k_1}\,H_{j_3,j_3,x^{j_2}}^{j_3}
+
\frac{1}{2}\,\delta_{j_1}^{k_1}\,H_{j_2,j_2,x^{j_3}}^{j_2}
}
- \\
& \
\ \ \ \ \ \ \ \ \ \ 
-
\underline{
\frac{1}{2}\,\delta_{j_3}^{k_1}\,H_{j_1,j_1,x^{j_2}}^{j_1}
+
\frac{1}{2}\,\delta_{j_2}^{k_1}\,H_{j_1,j_1,x^{j_3}}^{j_1}
}
- \\
& \
\ \ \ \ \ \ \ \ \ \ 
-
\frac{1}{2}\,G_{j_1,j_2}\,L_{j_3}^{k_1}
+
\frac{1}{2}\,G_{j_1,j_3}\,L_{j_2}^{k_1}
-
\frac{1}{4}\,\delta_{j_3}^{k_1}\,H_{j_1,j_1}^{j_1}\,H_{j_2,j_2}^{j_2}
+
\frac{1}{4}\,\delta_{j_2}^{k_1}\,H_{j_1,j_1}^{j_1}\,H_{j_3,j_3}^{j_3}
- \\
& \
\ \ \ \ \ \ \ \ \ \ 
-
\sum_{k_2}\,H_{j_1,j_2}^{k_2}\,H_{j_3,k_2}^{k_1}
+
\sum_{k_2}\,H_{j_1,j_3}^{k_2}\,H_{j_2,k_2}^{k_1}
-
\frac{1}{2}\,\delta_{j_2}^{k_1}\,H_{j_1,j_3}^{k_2}\,H_{k_2,k_2}^{k_2}
+
\frac{1}{2}\,\delta_{j_3}^{k_1}\,H_{j_1,j_2}^{k_2}\,H_{k_2,k_2}^{k_2}
- \\
& \
\ \ \ \ \ \ \ \ \ \ 
-
\frac{1}{2}\,\delta_{j_2}^{k_1}\,G_{j_1,j_3}\,\Theta^{n+1}
+
\frac{1}{2}\,\delta_{j_3}^{k_1}\,G_{j_1,j_2}\,\Theta^{n+1}
- \\
& \
\ \ \ \ \ \ \ \ \ \ 
-
\frac{1}{4}\,\delta_{j_2}^{k_1}\,H_{j_1,j_1}^{j_1}\,\Theta^{j_3}
+
\frac{1}{4}\,\delta_{j_3}^{k_1}\,H_{j_1,j_1}^{j_1}\,\Theta^{j_2}
-
\frac{1}{4}\,\delta_{j_2}^{k_1}\,H_{j_3,j_3}^{j_3}\,\Theta^{j_1}
+
\frac{1}{4}\,\delta_{j_3}^{k_1}\,H_{j_2,j_2}^{j_2}\,\Theta^{j_1}
- \\
& \
\ \ \ \ \ \ \ \ \ \ 
-
\frac{1}{2}\,\delta_{j_3}^{k_1}\,
\sum_{k_1}\,H_{j_1,j_2}^{k_2}\,\Theta^{k_2}
+
\frac{1}{2}\,\delta_{j_2}^{k_1}\,
\sum_{k_1}\,H_{j_1,j_3}^{k_2}\,\Theta^{k_2}- \\
& \
\ \ \ \ \ \ \ \ \ \ 
-
\frac{1}{4}\,\delta_{j_3}^{k_1}\,\Theta^{j_1}\,\Theta^{j_2}
+
\frac{1}{4}\,\delta_{j_2}^{k_1}\,\Theta^{j_1}\,\Theta^{j_3}.
\endaligned\right.
\end{equation}
Fifthly:
\def\theequation{3.21}\begin{equation}
\left\{
\aligned
&
\frac{1}{2}\,\delta_{j_1}^{k_1}\,\Theta_y^{j_2}
+
\frac{1}{2}\,\delta_{j_2}^{k_1}\,\Theta_y^{j_1}
-
\frac{1}{2}\,\delta_{j_1}^{k_1}\,\Theta_{x^{j_2}}^{n+1}
= \\
& \
\ \ \ \ \
=
G_{j_1,j_2}\,M^{k_1}
+
\frac{1}{2}\,\sum_{k_2}\,H_{j_1,k_2}^{k_1}\,L_{j_1}^{k_2}
-
\frac{1}{2}\,\sum_{k_2}\,H_{j_1,j_2}^{k_2}\,L_{k_2}^{k_1}
-
\frac{1}{4}\,\delta_{j_2}^{k_1}\,\sum_{k_2}\,
H_{k_2,k_2}^{k_2}\,L_{j_1}^{k_2}
- \\
& \
\ \ \ \ \ \ \ \ \ \
-
\frac{1}{4}\,\delta_{j_2}^{k_1}\,H_{j_1,j_1}^{j_1}\,\Theta^{n+1}
+
\frac{1}{4}\,\delta_{j_2}^{k_1}\,\sum_{k_2}\,
L_{j_1}^{k_2}\,\Theta^{k_2}
+
\frac{1}{4}\,\delta_{j_2}^{k_1}\,\Theta^{j_1}\,\Theta^{n+1}.
\endaligned\right.
\end{equation}
Sixthly:
\def\theequation{3.22}\begin{equation}
\left\{
\aligned
\delta_{j_1}^{k_1}\,\Theta_y^{n+1}
&
=
-
\underline{
L_{j_1, y}^{k_1}
+
2\,M_{x^{j_1}}^{k_1}
}
+ \\
& \
\ \ \ \ \ 
+
2\,\sum_{k_2}\,H_{j_1,k_2}^{k_1}\,M^{k_2}
-
\delta_{j_1}^{k_1}\,\sum_{k_2}\,H_{k_2,k_2}^{k_2}\,M^{k_2}
-
\frac{1}{2}\,\sum_{k_2}\,L_{j_1}^{k_2}\,L_{k_2}^{k_1}
+ \\
& \
\ \ \ \ \
+
\delta_{j_1}^{k_1}\,\sum_{k_2}\,M^{k_2}\,\Theta^{k_2}
+
\frac{1}{2}\,\delta_{j_1}^{k_1}\,\Theta^{n+1}\,\Theta^{n+1}.
\endaligned\right.
\end{equation}

\subsection*{3.23.~Solving $\Theta_{ x^{j_2 } }^{ j_1}$, 
$\Theta_y^{j_1}$, $\Theta_{x^{ j_1}}^{ n+1}$ and $\Theta_y^{ n+1 }$}
From the six families of equations~\thetag{ 3.17}, \thetag{ 3.18},
\thetag{ 3.19}, \thetag{ 3.20}, \thetag{ 3.21} and~\thetag{ 3.22}, we
can solve $\Theta_{ x^{j_2 }}^{j_1}$, $\Theta_y^{ j_1}$, $\Theta_{ x^{
j_1}}^{ n+1}$ and $\Theta_y^{ n+1}$. Not mentioning the (hard)
intermediate computations, we obtain firstly:
\def\theequation{3.24}\begin{equation}
\left\{
\aligned
\Theta_{x^{j_2}}^{j_1}
&
=
-
\underline{
2\,G_{j_1,j_2,y}
+
H_{j_1,j_1,x^{j_2}}^{j_1}
}
+
\sum_l\,G_{j_2,l}\,L_{j_1}^l
+
\frac{1}{2}\,H_{j_1,j_1}^{j_1}\,H_{j_2,j_2}^{j_2}
-
\sum_l\,H_{j_1,j_2}^l\,H_{l,l}^l
- \\
& \
\ \ \ \ \
-
G_{j_1,j_2}\,\Theta^{n+1}
-
\frac{1}{2}\,H_{j_1,j_1}\,\Theta^{j_1}
-
\frac{1}{2}\,H_{j_2,j_2}^{j_2}\,\Theta^{j_1}
+
\sum_l\,H_{j_1,j_2}^l\,\Theta^l
+
\frac{1}{2}\,\Theta^{j_1}\,\Theta^{j_2}.
\endaligned\right.
\end{equation}
Secondly:
\def\theequation{3.25}\begin{equation}
\left\{
\aligned
\Theta_y^{j_1}
&
=
-
\underline{
\frac{1}{3}\,H_{j_1,j_1,y}^{j_1}
+
\frac{2}{3}\,L_{j_1,x^{j_1}}^{j_1}
}
+
\frac{4}{3}\,G_{j_1,j_1}\,M^{j_1}
+
\frac{2}{3}\,\sum_l\,G_{j_1,l}\,M^l
-
\frac{1}{2}\,\sum_l\,H_{l,l}^l\,L_{j_1}^l
+ \\
& \
\ \ \ \ \
+
\frac{2}{3}\,\sum_l\,H_{j_1,l}^{j_1}\,L_{j_1}^l
-
\frac{2}{3}\,\sum_l\,H_{j_1,j_1}^l\,L_l^{j_1}
-
\frac{1}{2}\,H_{j_1,j_1}^{j_1}\,\Theta^{n+1}
+
\frac{1}{2}\,\sum_l\,L_{j_1}^l\,\Theta^l
+ \\
& \
\ \ \ \ \
+
\frac{1}{2}\,\Theta^{j_1}\,\Theta^{n+1}.
\endaligned\right.
\end{equation}
Thirdly:
\def\theequation{3.26}\begin{equation}
\left\{
\aligned
\Theta_{x^{j_1}}^{n+1}
&
=
-
\underline{
\frac{2}{3}\,H_{j_1,j_1,y}^{j_1}
+
\frac{1}{3}\,L_{j_1,x^{j_1}}^{j_1}
}
+
\frac{2}{3}\,G_{j_1,j_1}\,M^{j_1}
+
\frac{4}{3}\,\sum_l\,G_{j_1,l}\,M^l
-
\frac{1}{2}\,\sum_l\,H_{l,l}^l\,L_{j_1}^l
+ \\
& \
\ \ \ \ \ 
+
\frac{1}{3}\,\sum_l\,H_{j_1,l}^{j_1}\,L_{j_1}^l
-
\frac{1}{3}\,\sum_l\,H_{j_1,j_1}^l\,L_l^{j_1}
-
\frac{1}{2}\,H_{j_1,j_1}^{j_1}\,\Theta^{n+1}
+
\frac{1}{2}\,\sum_l\,L_{j_1}^l\,\Theta^l
+ \\
& \
\ \ \ \ \
+
\frac{1}{2}\,\Theta^{j_1}\,\Theta^{n+1}.
\endaligned\right.
\end{equation}
Fourtly:
\def\theequation{3.27}\begin{equation}
\left\{
\aligned
\Theta_y^{n+1}
&
=
-
\underline{
L_{j_1,y}^{j_1}
+
2\,M_{x^{j_1}}^{j_1}
}
+
2\,\sum_l\,H_{j_1,l}^{j_1}\,M^l
-
\sum_l\,H_{l,l}^l\,M^l
-
\frac{1}{2}\,\sum_l\,L_{j_1}^l\,L_l^{j_1}
+ \\
& \
\ \ \ \ \ 
+
\sum_l\,M^l\,\Theta^l
+
\frac{1}{2}\,\Theta^{n+1}\,\Theta^{n+1}.
\endaligned\right.
\end{equation}
These four families of partial differential equations constitute the
{\sl second auxiliary system}. By replacing these solutions in the
three remaining families of equations~\thetag{ 3.20}, \thetag{ 3.21}
and~\thetag{ 3.22}, we obtain supplementary equations
(which we do not copy) that are direct consequences of 
(I'), (II'), (III'), (IV').

To complete the proof of the main Lemma~3.3 above, it suffices now to
establish the first implication of the following list, since the other
three have been already established.

\begin{itemize}
\item[$\bullet$]
Some given functions $G_{j_1, j_2}$, $H_{j_1, j_2}^{k_1}$,
$L_{j_1}^{k_1}$ and $M^{k_1}$ of $(x^{l_1},y)$ satisfy the four
families of partial differential equations (I'), (II'), (III') and
(IV') of Theorem~1.7.
\item[$\Downarrow$]
\item[$\bullet$]
There exist functions $\Theta^{j_1}$, $\Theta^{n+1}$ satisfying the
second auxiliary system~\thetag{ 3.24}, \thetag{ 3.25}, \thetag{ 3.26}
and~\thetag{ 3.27}.
\item[$\Downarrow$]
\item[$\bullet$]
These solution functions $\Theta^{j_1}$, $\Theta^{n+1}$ satisfy the
six families of partial differential equations~\thetag{ 3.17},
\thetag{ 3.18}, \thetag{ 3.19}, \thetag{ 3.20}, \thetag{ 3.21}
and~\thetag{ 3.22}.
\item[$\Downarrow$]
\item[$\bullet$]
There exist functions $\Pi_{ j_1, j_2}^{k_1}$ of $(x^{l_1}, y)$, $1
\leqslant j_1, j_2, k_1 \leqslant m+1$, satisfying the first auxiliary
system~\thetag{ 3.7} of partial differential equations.
\item[$\Downarrow$]
\item[$\bullet$]
There exist functions $X^{l_2}$, $Y$ of $(x^{l_1},y)$ transforming the
system $y_{x^{j_1}x^{j_2}} = F^{j_1,j_2} (x^{l_1}, y, y_{x^{l_2}})$,
$j_1, j_2 = 1, \dots, n$, to the simplest system $Y_{X^{j_1}X^{j_2}} =
0$, $j_1, j_1 = 1, \dots, n$.
\end{itemize}

\smallskip

\subsection*{ 3.28.~Compatibility conditions for the second
auxiliary system} We notice that the second auxiliary system is also a
complete system. Thus, to establish the first above implication, it
suffices to show that the four families of compatibility conditions:
\def\theequation{3.29}\begin{equation}
\left\{
\aligned
0
& 
=
\left(
\Theta_{x^{j_2}}^{j_1}
\right)_{x^{j_3}}
-
\left(
\Theta_{x^{j_3}}^{j_1}
\right)_{x^{j_2}}, 
\\
0
& 
=
\left(
\Theta_{x^{j_2}}^{j_1}
\right)_y
-
\left(
\Theta_y^{j_1}
\right)_{x^{j_2}}, 
\\
0
& 
=
\left(
\Theta_{x^{j_1}}^{n+1}
\right)_{x^{j_2}}
-
\left(
\Theta_{x^{j_2}}^{n+1}
\right)_{x^{j_1}}, 
\\
0
& 
=
\left(
\Theta_{x^{j_2}}^{n+1}
\right)_y
-
\left(
\Theta_y^{n+1}
\right)_{x^{j_2}}, 
\endaligned\right.
\end{equation}
are a consequence of (I'), (I''), (III'), (IV').

For instance, in $(3.29)_1$, replacing $\Theta_{x^{j_2}}^{j_1}$ by its
expression~\thetag{ 3.24}, differentiating it with respect to
$x^{j_3}$, replacing $\Theta_{x^{j_3}}^{j_1}$ by its
expression~\thetag{ 3.24}, differentiating it with respect to
$x^{j_2}$ and substracting, we get: 
\def\theequation{3.30}\begin{equation}
\left\{
\aligned
0
& 
=
-
2\,G_{j_1,j_2,yx^{j_3}}
+
2\,G_{j_1,j_3,yx^{j_2}}
+
\underline{ 
H_{j_1,j_1,x^{j_2}x^{j_3}}^{j_1}
}_{\, \tiny{\sf a}}
-
\underline{
H_{j_1,j_1,x^{j_3}x^{j_2}}^{j_1}
}_{\, \tiny{\sf a}}
+ \\
& \
\ \ \ \ \
+
\frac{1}{2}\,
\underline{\Theta_{x^{j_3}}^{j_1}}\,\Theta^{j_2}
+
\frac{1}{2}\,\Theta^{j_1}\,
\underline{\Theta_{x^{j_3}}^{j_2}}
-
\frac{1}{2}\,
\underline{\Theta_{x^{j_2}}^{j_1}}\,\Theta^{j_3}
-
\frac{1}{2}\,\Theta^{j_1}\,
\underline{\Theta_{x^{j_2}}^{j_3}}
- \\
& \
\ \ \ \ \
-
\frac{1}{2}\,H_{j_1,j_1,x^{j_3}}^{j_1}\,\Theta^{j_2}
-
\frac{1}{2}\,H_{j_1,j_1}^{j_1}\,
\underline{\Theta_{x^{j_3}}^{j_2}}
+
\frac{1}{2}\,H_{j_1,j_1,x^{j_2}}^{j_1}\,\Theta^{j_3}
+
\frac{1}{2}\,H_{j_1,j_1}^{j_1}\,
\underline{\Theta_{x^{j_2}}^{j_3}}
- \\
& \
\ \ \ \ \ 
-
\frac{1}{2}\,H_{j_2,j_2,x^{j_3}}^{j_2}\,\Theta^{j_1}
-
\frac{1}{2}\,H_{j_2,j_2}^{j_2}\,
\underline{\Theta_{x^{j_3}}^{j_1}}
+
\frac{1}{2}\,H_{j_3,j_3,x^{j_2}}^{j_3}\,\Theta^{j_1}
+
\frac{1}{2}\,H_{j_3,j_3}^{j_3}\,
\underline{\Theta_{x^{j_2}}^{j_1}}
- \\
& \
\ \ \ \ \
-
G_{j_1,j_2,x^{j_3}}\,\Theta^{n+1}
-
G_{j_1,j_2}\,\underline{\Theta_{x^{j_3}}^{n+1}}
+
G_{j_1,j_3,x^{j_2}}\,\Theta^{n+1}
+
G_{j_1,j_3}\,\underline{\Theta_{x^{j_2}}^{n+1}}
+ \\
& \
\ \ \ \ \
+
\sum_l\,H_{j_1,j_2,x^{j_3}}^l\,\Theta^l
+
\sum_l\,H_{j_1,j_2}^l\,
\underline{\Theta_{x^{j_3}}^l}
-
\sum_l\,H_{j_1,j_3,x^{j_2}}^l\,\Theta^l
-
\sum_l\,H_{j_1,j_3}^l\,
\underline{\Theta_{x^{j_2}}^l}
+ \\
& \
\ \ \ \ \
+
\frac{1}{2}\,H_{j_1,j_1,x^{j_3}}^{j_1}\,H_{j_2,j_2}^{j_2}
+
\frac{1}{2}\,H_{j_1,j_1}^{j_1}\,H_{j_2,j_2,x^{j_3}}^{j_2}
-
\frac{1}{2}\,H_{j_1,j_1,x^{j_2}}^{j_1}\,H_{j_3,j_3}^{j_3}
-
\frac{1}{2}\,H_{j_1,j_1}^{j_1}\,H_{j_3,j_3,x^{j_2}}^{j_3}
- \\
& \
\ \ \ \ \ 
-
\sum_l\,H_{j_1,j_2,x^{j_3}}^l\,H_{l,l}^l
-
\sum_l\,H_{j_1,j_2}^l\,H_{l,l,x^{j_3}}^l
+
\sum_l\,H_{j_1,j_3,x^{j_2}}^l\,H_{l,l}^l
+
\sum_l\,H_{j_1,j_3}^l\,H_{l,l,x^{j_2}}^l
+ \\
& \
\ \ \ \ \ 
+
\sum_l\,G_{j_2,l,x^{j_3}}\,L_{j_1}^l
+
\sum_l\,G_{j_2,l}\,L_{j_1,x^{j_3}}^l
-
\sum_l\,G_{j_3,l,x^{j_2}}\,L_{j_1}^l
-
\sum_l\,G_{j_3,l}\,L_{j_1,x^{j_2}}^l.
\endaligned\right.
\end{equation}
Next, replacing the twelve first order partial derivatives underlined
just above:
\def\theequation{3.31}\begin{equation}
\left\{
\aligned
&
\underline{\Theta_{x^{j_3}}^{j_1}}, 
\ \ \ \ \
\underline{\Theta_{x^{j_3}}^{j_2}}, 
\ \ \ \ \
\underline{\Theta_{x^{j_2}}^{j_1}}, 
\ \ \ \ \
\underline{\Theta_{x^{j_2}}^{j_3}},
\ \ \ \ \
\underline{\Theta_{x^{j_3}}^{j_2}}, 
\ \ \ \ \
\underline{\Theta_{x^{j_2}}^{j_3}},
\\
&
\underline{\Theta_{x^{j_3}}^{j_1}},
\ \ \ \ \
\underline{\Theta_{x^{j_2}}^{j_1}},
\ \ \ \ \
\underline{\Theta_{x^{j_3}}^{n+1}},
\ \ \ \ \
\underline{\Theta_{x^{j_2}}^{n+1}},
\ \ \ \ \
\underline{\Theta_{x^{j_3}}^l},
\ \ \ \ \
\underline{\Theta_{x^{j_2}}^l}.
\endaligned\right.
\end{equation}
by their values issued from~\thetag{ 3.24}, \thetag{ 3.26} and
adapting the summation indices, we get the explicit developed form of
the first family of compatibility conditions $(3.29)_1$:
\def\theequation{3.32}\begin{equation}
\left\{
\aligned
0 
= 
?
&
=
-
\underline{\underline{
2\,G_{j_1,j_2,x^{j_3}y}
+
2\,G_{j_1,j_3,x^{j_2}y}
}}
- \\
& \
\ \ \ \ \
-
\underline{
\sum_l\,G_{j_3,l,x^{j_2}}\,L_{j_1}^l
+
\sum_l\,G_{j_2,l,x^{j_3}}\,L_{j_1}^l
-
G_{j_1,j_2,y}\,H_{j_3,j_3}^{j_3}
+
G_{j_1,j_3,y}\,H_{j_2,j_2}^{j_2}
}
- \\
& \
\ \ \ \ \
-
\underline{
2\,\sum_l\,G_{l,j_3}\,H_{j_1,j_2}^l
+
2\,\sum_l\,G_{l,j_2}\,H_{j_1,j_3}^l
-
\sum_l\,H_{j_1,j_2,x^{j_3}}^l\,H_{l,l}^l
+
\sum_l\,H_{j_1,j_3,x^{j_2}}^l\,H_{l,l}^l
}
- \\
& \
\ \ \ \ \
-
\underline{
\frac{2}{3}\,H_{j_2,j_2,y}^{j_2}\,G_{j_1,j_3}
+
\frac{2}{3}\,H_{j_3,j_3,y}^{j_3}\,G_{j_1,j_2}
-
\frac{2}{3}\,L_{j_3,x^{j_3}}^{j_3}\,G_{j_1,j_2}
+
\frac{2}{3}\,L_{j_2,x^{j_2}}^{j_2}\,G_{j_1,j_3}
}
- \\
& \
\ \ \ \ \
-
\underline{
\sum_l\,L_{j_1,x^{j_2}}^l\,G_{j_3,l}
+
\sum_l\,L_{j_1,x^{j_3}}^l\,G_{j_2,l}
}
- \\
& \
\ \ \ \ \
-
\frac{2}{3}\,G_{j_1,j_2}\,G_{j_3,j_3}\,M^{j_3}
+
\frac{2}{3}\,G_{j_1,j_3}\,G_{j_2,j_2}\,M^{j_2}
-
\frac{4}{3}\,\sum_l\,G_{j_1,j_2}\,G_{j_3,l}\,M^l
+ \\
& \
\ \ \ \ \
+
\frac{4}{3}\,\sum_l\,G_{j_1,j_3}\,G_{j_2,l}\,M^l
-
\frac{1}{2}\,\sum_l\,G_{j_3,l}\,H_{j_1,j_1}^{j_1}\,L_{j_2}^l
+
\frac{1}{2}\,\sum_l\,G_{j_2,l}\,H_{j_1,j_1}^{j_1}\,L_{j_3}^l
- \\
& \
\ \ \ \ \
-
\frac{1}{2}\,\sum_l\,G_{j_3,l}\,H_{j_2,j_2}^{j_2}\,L_{j_1}^l
+
\frac{1}{2}\,\sum_l\,G_{j_2,l}\,H_{j_3,j_3}^{j_3}\,L_{j_1}^l
-
\frac{1}{2}\,\sum_l\,G_{j_1,j_3}\,H_{l,l}^l\,L_{j_2}^l
+ \\
& \
\ \ \ \ \
+
\frac{1}{2}\,\sum_l\,G_{j_1,j_2}\,H_{l,l}^l\,L_{j_3}^l
-
\frac{1}{3}\,\sum_l\,G_{j_1,j_2}\,H_{j_3,l}^{j_3}\,L_{j_3}^l
+
\frac{1}{3}\,\sum_l\,G_{j_1,j_3}\,H_{j_2,l}^{j_2}\,L_{j_2}^l
- \\
& \
\ \ \ \ \
-
\frac{1}{3}\,G_{j_1,j_3}\,H_{j_2,j_2}^l\,L_l^{j_2}
+
\frac{1}{3}\,G_{j_1,j_2}\,H_{j_3,j_3}^l\,L_l^{j_3}
- \\
& \
\ \ \ \ \
-
\sum_l\,\sum_p\,G_{j_2,p}\,H_{j_1,j_3}^l\,L_l^p
+
\sum_l\,\sum_p\,G_{j_3,p}\,H_{j_1,j_2}^l\,L_l^p
- \\
& \
\ \ \ \ \
-
\sum_l\,\sum_p\,H_{j_1,j_2}^l\,H_{l,j_3}^p\,H_{p,p}^p
+
\sum_l\,\sum_p\,H_{j_1,j_3}^l\,H_{l,j_2}^p\,H_{p,p}^p.
\endaligned\right.
\end{equation}

\def\thelemma{3.33}\begin{lemma}
{\rm (\cite{ me2003, me2004})}
This first family of compatibility conditions for the second auxiliary
system obtained by developing $(3.29)_1$ in length, together with the
three remaining families obtained by developing $(3.29)_2$,
$(3.29)_3$, $(3.29)_4$ in length, are consequences, by linear
combinations and by differentiations, of {\rm (I')}, {\rm (II')}, 
{\rm (III')}, {\rm (IV')}, of Theorem~1.7.
\end{lemma}

The summarized proof of Theorem~1.7 is complete.
\qed

\newpage

\begin{center}
{\large\bf
{\Large\bf IV:}~Bibliography
}
\end{center}

\vfill\end{document}